\documentclass[11pt]{amsart}
\usepackage{amsmath,amscd,latexsym,verbatim,amssymb}
\usepackage[frenchb]{babel}
\usepackage{times}
\DeclareMathOperator{\End}{End}
\DeclareMathOperator{\GR}{Gr}
\DeclareMathOperator{\Hom}{Hom}
\DeclareMathOperator{\Img}{Im}
\DeclareMathOperator{\Ind}{Ind}
\DeclareMathOperator{\Red}{Red}

\newcommand{\num}{{\operatorname{num}}}

\newcommand{\kim}{{\operatorname{kim}}}

\newcommand{\rat}{{\operatorname{rat}}}
\newcommand{\snum}{{\operatorname{snum}}}
\newcommand{\rad}{{\operatorname{rad}}}
\newcommand{\red}{{\operatorname{red}}}
\newcommand{\comp}{{\operatorname{comp}}}

\newcommand{\Spec}{\operatorname{Spec}}
\newcommand{\diag}{\operatorname{diag}}
\newcommand{\uHom}{{\underline{\operatorname{Hom}}}}
\newcommand{\resp}{{\it resp.} }
\newcommand{\cf}{{\it cf.} }
\newcommand{\ie}{{\it i.e.} }
\newcommand{\eg}{{\it e.g.} }

\newcommand{\loccit}{{\it loc. cit.} }
\newcommand{\If}{\Longrightarrow}
\newcommand{\bA}{\mathbf{A}}

\newcommand{\F}{\mathbf{F}}
\newcommand{\N}{\mathbf{N}}
\newcommand{\bP}{\mathbf{P}}
\newcommand{\Q}{\mathbf{Q}}
\newcommand{\R}{\mathbf{R}}
\newcommand{\boP}{\mathbf{P}}

\newcommand{\U}{\mathbf{U}}
\newcommand{\Z}{\mathbf{Z}}
\newcommand{\bG}{\mathbb{G}}
\newcommand{\bI}{\mathbb{I}}
\newcommand{\bT}{\mathbb{T}}
\newcommand{\Sec}{\mathbf{Sec}}
\newcommand{\frg}{\mathfrak{g}}
\newcommand{\fsl}{\mathfrak{sl}}
\newcommand{\sA}{{\mathcal{A}}}
\newcommand{\sB}{{\mathcal{B}}}
\newcommand{\sC}{{\mathcal{C}}}
\newcommand{\sE}{{\mathcal{E}}}
\newcommand{\sG}{{\mathcal{G}}}
\newcommand{\sH}{{\mathcal{H}}}
\newcommand{\sI}{{\mathcal{I}}}
\newcommand{\sJ}{{\mathcal{J}}}

\newcommand{\sN}{{\mathcal{N}}}
\newcommand{\sO}{{\mathcal{O}}}

\newcommand{\sR}{{\mathcal{R}}}
\newcommand{\sS}{{\mathcal{S}}}
\newcommand{\sU}{{\mathcal{U}}}
\newcommand{\sV}{{\mathcal{V}}}
\newcommand{\sX}{{\mathcal{X}}}

\newcommand{\To}{\longrightarrow}
\newcommand{\inj}{\hookrightarrow}

\newcommand{\surj}{\rightarrow\!\!\!\!\!\rightarrow}
\newcommand{\Surj}{\relbar\joinrel\surj}
\newcommand{\iso}{\overset{\sim}{\longrightarrow}}
\newcommand{\osi}{\overset{\sim}{\longleftarrow}}
\newcommand{\rr}{\rightrightarrows}

\newcommand{\Cat}{\operatorname{\mathbf{Cat}}}
\newcommand{\Gr}{\operatorname{\mathbf{Gr}}}
\newcommand{\Acg}{\operatorname{\mathbf{Acg}}}
\newcommand{\Hmg}{\operatorname{\mathbf{Hmg}}}
\newcommand{\un}{{\bf 1}}
\newcommand{\simtimes}{\mathrel{\overset{\otimes}{\simeq}}}
\newcommand{\point}[1]{\overset{\bullet}{#1}}
\newcommand{\pointplus}[1]{\overset{\oplus}{#1}}

\renewcommand{\epsilon}{\varepsilon}
\renewcommand{\phi}{\varphi}

\renewcommand{\lim}{\varprojlim}
\newcommand{\colim}{\varinjlim}

\newcommand{\s}{\scriptstyle}
\renewcommand{\d}{\displaystyle}

\newcommand{\Ker}{\operatorname{Ker}}
\newcommand{\IM}{\operatorname{Im}}
\newcommand{\Coker}{\operatorname{Coker}}
\newcommand{\Ext}{\operatorname{Ext}}
\newcommand{\Tor}{\operatorname{Tor}}

\newcommand{\car}{\operatorname{car}}
\newcommand{\Aff}{\operatorname{Aff}}
\newcommand{\Gaff}{\operatorname{Gaff}}
\newcommand{\Gred}{\operatorname{Gred}}
\newcommand{\Pred}{\operatorname{{}^pRed}}
\newcommand{\PU}{\operatorname{{}^pU}}
\newcommand{\bGaff}{\overline{\operatorname{Gaff}}}
\newcommand{\bGred}{\overline{\operatorname{Gred}}}
\newcounter{spec}
\newenvironment{thlist}{\begin{list}{\rm{(\roman{spec})}}%
{\usecounter{spec}\labelwidth=20pt\itemindent=0pt\labelsep=10pt}}%
{\end{list}}

\newtheorem{Thm}{Th\'eor\`eme}
\newtheorem{lem}{Lemma}[section]
\newtheorem{Th}{Theorem}[section]

\swapnumbers
\newtheorem{thm}{Th\'eor\`eme}[subsection]
\newtheorem{lemme}[thm]{Lemme}
\newtheorem{prop}[thm]{Proposition}
\newtheorem{sco}[thm]{Scolie}
\newtheorem{cor}[thm]{Corollaire}

\newtheorem{sorite}[thm]{Sorite}
\theoremstyle{definition}
\newtheorem{defn}[thm]{D\'efinition}
\newtheorem{nota}[thm]{Notation}
\newtheorem{ex}[thm]{Exemple}
\newtheorem{meg}[thm]{Mise en garde}
\newtheorem{qn}[thm]{Question}

\newtheorem{contrex}[thm]{Contre-exemple}
\newtheorem{rem}[thm]{Remarque}
\newtheorem{rems}[thm]{Remarques}

\newtheorem{hyps}[thm]{Hypoth\`eses}

\newtheorem{para}[thm]{}

\numberwithin{equation}{section}

 at10pt

\newcommand{\prf}{\noindent {\bf D\'emonstration. }}
\renewcommand{\qed}{\hfill $\square$\medskip}

\begin{document}

\title{Nilpotence, radicaux et structures mo\-no\-\"{\i}\-dales}
\author[Yves
Andr\'e et Bruno Kahn]{Yves
Andr\'e et Bruno Kahn\\ avec un
appendice de Peter O'Sullivan}
  \address{D\'epartement de Math\'ematiques et Applications, \'Ecole Normale Sup\'erieure, \\ 
45 rue d'Ulm,  75230
  Paris cedex 05\\France.}
\email{andre@dmi.ens.fr}
 \address{Institut de Math\'ematiques de
 Jussieu\\175--179 rue du
 Chevaleret\\ \break 75013
 Paris\\France.}
\email{kahn@math.jussieu.fr}
 \address{39 Edgecliffe Ave\\
         Coogee, NSW 2034\\
         Australia}
\email{pjosullivan@optusnet.com.au}
\date{10 octobre 2002}
\subjclass{16N, 16D, 18D10, 18E, 14L,
16G60, 13E10, 17C}
\maketitle

\tableofcontents

\renewcommand{\abstractname}{Summary}

\
\bigskip\bigskip\bigskip\bigskip\bigskip\bigskip\bigskip 
\bigskip\bigskip

\begin{abstract}
For $K$ a field, a \emph{Wedderburn $K$-linear category} is a
$K$-linear category $\sA$ whose radical $\sR$ is locally nilpotent and
such that
$\bar \sA:=\sA/\sR$ is semi-simple and remains so after any extension of
scalars. We prove existence and uniqueness results for sections of the
projection $\sA\to \bar\sA$, in the vein of the theorems of Wedderburn.
There are two such results: one in the general case and one when $\sA$
has a monoidal structure for which $\sR$ is a monoidal ideal. The latter
applies notably to Tannakian categories over a field of characteristic
zero, and we get a generalisation of the Jacobson-Morozov theorem: the
existence of a \emph{pro-reductive envelope} $\Pred(G)$ associated to any
affine group scheme $G$ over $K$ ($\Pred(\bG_a)=SL_2$, and $\Pred(G)$ is
infinite-dimensional for any bigger unipotent group). Other applications
are given in this paper as well as in the note \cite{ak(note)} on motives.
\end{abstract}

\newpage

\section*{Introduction}

\bigskip\bigskip

En alg\`ebre non commutative, on rencontre d'abord, par ordre de
complexit\'e croissant, les anneaux semi-simples, puis les anneaux
semi-pri\-mai\-res, c'est-\`a-dire les extensions d'un anneau semi-simple
par un id\'eal nilpotent (le radical). Parmi ces derniers, le cas o\`u
l'anneau semi-simple est en fait une alg\`ebre s\'eparable est
particuli\`erement agr\'eable puisque, d'apr\`es un th\'eor\`eme
classique de Wedderburn, l'extension par le radical se scinde.

Par ailleurs, en renversant la tautologie qui identifie tout anneau \`a
une ca\-t\'e\-go\-rie additive
\`a un seul objet, on peut consid\'erer les ca\-t\'e\-go\-ries additives comme des
anneaux \`a plusieurs objets. Ce point de vue, popu\-laris\'e par B.
Mitchell et R. Street, permet de s'inspirer largement de l'alg\`ebre non
commutative dans les questions cat\'egoriques. Comme tous premiers
exemples, on obtient la notion d'id\'eal, et celle de radical, d'une
ca\-t\'e\-go\-rie additive.

Suivant ce point de vue, le premier th\`eme de cet article est
l'\'etude des ca\-t\'e\-go\-ries semi-primaires, extensions d'une ca\-t\'e\-go\-rie
semi-simple par un id\'eal v\'erifiant une condition convenable de
nilpotence (\cf d\'efini\-tion \ref{D4}), et de l'analogue cat\'egorique du
th\'eor\`eme de Wedderburn. Le second th\`eme est l'\'etude du
radical en pr\'esence d'une structure mo\-no\-\"{\i}\-dale. Ces th\`emes se
rejoignent dans une version mo\-no\-\"{\i}\-dale du th\'eor\`eme de Wedderburn
(th\'e\-o\-r\`e\-mes \ref{T2} et \ref{t4}).

\medskip Nous donnons deux applications principales de ces r\'esultats.
La premi\`ere est la construction de l'\emph{enveloppe pro-r\'eductive}
$\Pred(G)$ d'un groupe al\-g\'e\-bri\-que lin\'eaire
quelconque $G$ en ca\-rac\-t\'e\-ris\-ti\-que nulle. Les re\-pr\'e\-sen\-ta\-tions
ind\'e\-com\-po\-sa\-bles de $G$ sont en bijection avec les re\-pr\'e\-sen\-ta\-tions
ir\-r\'e\-duc\-ti\-bles de $\Pred(G)$, et ceci caract\'erise $\Pred(G)$ si
$K$ est alg\'ebriquement clos (proposition \ref{p14.1}). Le
pro\-to\-ty\-pe est
$\Pred(\bG_a)=SL_2$ (th\'eor\`eme de Jacobson-Morozov). Toutefois,
$\Pred(G)$ n'est en g\'en\'eral pas de
dimension finie: cela arrive d\'ej\`a pour $G=\bG_a\times \bG_a$.
L'existence m\^eme de $\Pred(G)$ n'en implique pas moins une s\'erie de
r\'esultats concrets sur les re\-pr\'e\-sen\-ta\-tions ind\'e\-com\-po\-sa\-bles des
groupes alg\'ebriques.

La seconde application concerne la
ca\-t\'e\-go\-rie des motifs purs  
construits en termes d'une
\'equiva\-lence ad\'equate quelconque (fix\'ee) pour 
les cycles
alg\'ebriques. Nous d\'e\-tail\-le\-rons ailleurs (voir d\'ej\`a
  
\cite{ak(note)}), en nous contentant dans ce texte de
br\`eves allusions, 
notamment dans la deuxi\`eme partie. Indiquons ici
seulement que l'on s'attend
\`a ce qu'une telle ca\-t\'e\-go\-rie de motifs soit semi-primaire, de radical
compatible
\`a la structure mo\-no\-\"{\i}\-dale (cela d\'ecou\-le\-rait de la 
conjecture
de Beilinson-Murre et des conjectures standard de Grothendieck); par
ailleurs, nous nous appuyons dans \cite{ak(note)} sur notre version
mo\-no\-\"{\i}\-dale du th\'eor\`eme de Wedderburn pour cons\-truire
inconditionnellement les groupes de Galois motiviques.

En fait, l'un des 
projets directeurs de ce travail a \'et\'e de
s\'eparer nettement ce qui dans la 
th\'eorie
encore largement conjecturale des motifs purs est de nature purement
cat\'egorique, et ce qui est de nature g\'eom\'etrique.

\medskip
D\'ecrivons plus en d\'etail les quatre parties de ce travail.

Les deux premiers paragraphes (et l'appendice) contiennent une discussion des 
notions cat\'egoriques fondamentales d'id\'eal, de radical (not\'e $\rad\sA$),
semi-simplicit\'e, semi-primarit\'e, etc...

On y introduit la notion de $K$-{\it ca\-t\'e\-go\-rie de Wedderburn}, qui joue un 
r\^ole important dans la suite ($K$ \'etant un corps): en supposant pour 
simplifier l'existence de biproduits finis, il s'agit d'une ca\-t\'e\-go\-rie
$\sA$ dont les ensembles de mor\-phis\-mes $\sA(A,B)$ sont des $K$-espaces
vectoriels (la composition \'etant $K$-lin\'eaire), telle que les
$\rad\sA(A,A)$ soient des id\'eaux nilpotents des alg\`ebres
d'endomor\-phis\-mes $\sA(A,A)$, et telle que les $K$-alg\`ebres
$(\sA/\rad\sA)(A,A)$ soient s\'eparables, \ie absolument semi-simples.  

Cette condition de nilpotence est assez faible, et on montre au \S
\ref{radrenf} qu'on ne peut gu\`ere la renforcer sans restreindre
consid\'erablement le champ d'application de la th\'eorie. Dans le cas de
ca\-t\'e\-go\-ries de modules, l'analyse de ces conditions de nilpotence
s'av\`ere intimement li\'ee \`a celle des carquois d'Auslander-Reiten. 

Bien que le radical se comporte ``mal" en g\'en\'eral par extension des
sca\-lai\-res, on montre au \S \ref{s3} que la situation est
meilleure dans le cas d'une $K$-ca\-t\'e\-go\-rie de Wedderburn.
Cette premi\`ere partie se termine par l'\'etude et la comparaison (\S
\ref{ext}) de deux notions d'extensions des scalaires pour les
$K$-ca\-t\'e\-go\-ries que l'on rencontre dans la litt\'erature.

\medskip
L'objectif de la seconde partie est d'\'etudier en d\'etail le radical
d'une ca\-t\'e\-go\-rie mo\-no\-\"{\i}\-dale sym\'e\-trique rigide (\ie dont les
objets ont des duaux). Il se trouve qu'il n'est pas compatible \`a la
structure mo\-no\-\"{\i}\-dale en g\'en\'eral (un exemple de ce ph\'enom\`ene
est donn\'e par la ca\-t\'e\-go\-rie des re\-pr\'e\-sen\-ta\-tions de $GL_p$ en
caract\'eristique $p>0$). Nous comparons chemin faisant le radical au
plus grand id\'eal propre compatible \`a la structure mo\-no\-\"{\i}\-dale, et
analysons le quotient de la ca\-t\'e\-go\-rie par son radical. Actions du
groupe sy\-m\'e\-tri\-que, puissances ext\'erieures et sy\-m\'e\-tri\-ques sont mises
en \oe uvre dans cette optique. Nous nous sommes largement inspir\'es du
r\'ecent travail de Kimura \cite{ki} sur les motifs de Chow ``de
dimension finie".

On montre que le radical d'une ca\-t\'e\-go\-rie tannakienne sur un corps de 
ca\-rac\-t\'e\-ris\-ti\-que nulle est toujours mo\-no\-\"{\i}\-dal, et
qu'une telle ca\-t\'e\-go\-rie  est toujours de Wedderburn. Les
th\'eor\`emes
\ref{absJannsen} et \ref{nouveau}  de ce paragraphe, beaucoup plus
g\'en\'eraux, contiennent une version abstraite  des r\'esultats de
Jannsen sur les motifs num\'eriques et de Kimura sur les  motifs de Chow
\cite{jannsen,ki}. Une attention particuli\`ere a
\'et\'e port\'ee aux variantes $\Z/2$-gradu\'ees, en vue des applications
aux motifs. En particulier, nous faisons le lien entre la
th\'eorie de Kimura \cite{ki} et la question de l'alg\'ebricit\'e des 
projecteurs de K\"unneth pairs (th\'eor\`eme \ref{tkim.1}). Notre r\'esultat le 
plus abouti est le th\'eor\`eme suivant: si on d\'efinit une \emph{ca\-t\'e\-go\-rie 
de Kimura} comme \'etant une ca\-t\'e\-go\-rie $K$-lin\'eaire pseudo-ab\'elienne
mo\-no\-\"{\i}\-dale sy\-m\'e\-tri\-que rigide, o\`u $K$ est un corps de
caract\'eristique z\'ero, dont tout objet est de dimension finie au sens
de Kimura (voir section \ref{kimura}), on a (th\'eor\`eme \ref{tctki}):

\begin{Thm}\label{Tkim} Toute ca\-t\'e\-go\-rie de Kimura $\sA$ est de
Wedderburn. Son radical $\sR$ est mo\-no\-\"{\i}\-dal, et la ca\-t\'e\-go\-rie
quotient $\sA/\sR$ est semi-simple tannakienne, apr\`es changement
convenable de la contrainte de commutativit\'e.
\end{Thm}

\medskip La troisi\`eme partie d\'ebute sur le rappel d'une version
cat\'egorique, due \`a Mitchell, de la cohomologie de Hochschild; elle
nous sert \`a prouver le point a) de l'analogue cat\'egorique suivant du
th\'eor\`eme de Wedderburn-Malcev (\cf th\'eor\`emes \ref{T1}, \ref{T2}
et \ref{t4}).

\begin{Thm}\label{th1} a) Soit $\sA$ une $K$-ca\-t\'e\-go\-rie de Wedderburn, et
soit $\pi: \sA\to \sA/\rad(\sA)$ le foncteur de projection. Alors $\pi$ admet 
une section fonctorielle. Deux telles sections sont conjugu\'ees.\\
b) Si $\sA$ est mo\-no\-\"{\i}\-dale avec $End({\un})=K$, et si le
radical est compatible \`a la structure mo\-no\-\"{\i}\-dale (de sorte que $\pi$
est un foncteur mo\-no\-\"{\i}\-dal), $\pi$ admet une section mo\-no\-\"{\i}\-dale.
Deux telles sections sont conjugu\'ees par un isomor\-phis\-me mo\-no\-\"{\i}\-dal.\\
c) Si enfin $\sA$ est sy\-m\'e\-tri\-que, toute section mo\-no\-\"{\i}\-dale est 
sy\-m\'e\-tri\-que (c'est-\`a-dire respecte les tressages) \`a condition que $\car 
K\ne 2$.
\end{Thm}

La preuve de b), tr\`es technique, se trouve au \S
\ref{s6}. Elle utilise aussi la cohomologie de Hochschild pour
l'alg\`ebre libre sur les objets de la ca\-t\'e\-go\-rie. Ainsi, a) repose
essentiellement sur le fait qu'une $K$-alg\`ebre s\'eparable est de
dimension $0$ (au sens de \cite[ch. IX, \S 7]{ce}), tandis que b) repose
sur le fait qu'une $K$-alg\`ebre libre est de dimension $1$\dots

La compatibilit\'e des sections mo\-no\-\"{\i}\-dales
\`a des tressages, dont il s'agit au point c) du th\'eor\`eme
\ref{th1}, est
\'etudi\'ee en d\'etail au
\S \ref{tre}. Elle n'est pas automatique; toutefois, en caract\'eristique
dif\-f\'e\-ren\-te de $2$, si le tressage r\'esiduel est sy\-m\'e\-tri\-que, son
image par toute section mo\-no\-\"{\i}\-dale $s$ donne lieu \`a un tressage
sy\-m\'e\-tri\-que canonique sur $\sA$, ind\'ependant du choix de $s$ et qui ne
co\"{\i}ncide avec le tressage originel que si celui-ci \'etait
sy\-m\'e\-tri\-que. Ceci
s'applique notamment, en caract\'eristique z\'ero,
\`a la quantification de Drinfeld-Cartier d'une ca\-t\'e\-go\-rie mo\-no\-\"{\i}\-dale
sy\-m\'e\-tri\-que munie d'un tressage infinit\'esimal (exemple \ref{ex4}).

Le th\'eor\`eme \ref{th1} s'applique en particulier aux ca\-t\'e\-go\-ries de
Kimura, en vertu du th\'eor\`eme \ref{Tkim} (\cf th\'eor\`eme
\ref{kimsplit}).

Au \S \ref{repr}, on conf\`ere une structure g\'eom\'etrique aux
constructions  pr\'e\-c\'e\-den\-tes: groupo\"{\i}des pro-unipotents des
sections (\resp des sections mo\-no\-\"{\i}\-da\-les).

\medskip La derni\`ere partie explore les cons\'equences des r\'esultats
pr\'ec\'edents en th\'eorie des re\-pr\'e\-sen\-ta\-tions.

Le \S \ref{radrep} est consacr\'e \`a l'\'etude de
$\sA/\rad(\sA)$, lorsque $\sA$ est la ca\-t\'e\-go\-rie des re\-pr\'e\-sen\-ta\-tions
d'un sch\'ema en groupes affine $G$ (par exemple un groupe alg\'ebrique
lin\'eaire) sur un corps
$K$ de ca\-rac\-t\'e\-ris\-ti\-que nulle. On peut conclure des
r\'esultats pr\'ec\'edents que
$\sA/\rad(\sA)$ est \'equi\-va\-len\-te \`a ca\-t\'e\-go\-rie des
re\-pr\'e\-sen\-ta\-tions
d'un sch\'e\-ma en groupes pro-r\'eductif $\Pred(G)$ contenant $G$.

Il est plus technique de rendre cette construction canonique. La
question se simplifie si l'on travaille, comme au \S \ref{jm}, dans la
ca\-t\'e\-go\-rie
$\bGaff_K$ dont les objets sont les $K$-groupes affines et les
mor\-phis\-mes $G\to H$ sont donn\'es par les ensembles quotients de
$Hom_K(G,H)$ par la relation d'\'equivalence $\sim$ telle que $f\sim g$
s'il existe $h\in H_K$ tel que $g=hfh^{-1}$. Soit $\bGred_K$ la
sous-ca\-t\'e\-go\-rie pleine form\'ee des
groupes pro-r\'eductifs. Alors (\cf th\'eor\`eme \ref{t1}):

\begin{Thm} Le foncteur d'inclusion $ \bGred_K\to
\bGaff_K$ admet un adjoint (et quasi-inverse) \`a gauche:
$G\mapsto \Pred(G)$.
\end{Thm}

On en d\'eduit une s\'erie de r\'esultats concrets concernant d'une part
la structure des groupes alg\'ebriques, et d'autre part les
re\-pr\'e\-sen\-ta\-tions in\-d\'e\-com\-po\-sa\-bles des groupes alg\'ebriques.

L'une de ces applications concerne la notion d'\emph{enveloppe
r\'eductive} d'un sous-groupe ferm\'e d'un groupe r\'eductif, \ie de
sous-groupe r\'eductif interm\'ediaire minimal. On a (\cf th\'eor\`eme
\ref{T3})

\begin{Thm} Deux enveloppes r\'eductives de $G$ dans $H$
sont toujours conjugu\'ees par un \'el\'ement de $h \in H(K)$ commutant
\`a $G$.
 \end{Thm}

Dans un paragraphe ant\'erieur, on examine l'avatar non mo\-no\-\"{\i}\-dal,
plus simple, de ces constructions, ce qui m\`ene \`a la notion
d'enveloppe semi-simple d'une alg\`ebre profinie (sur un corps parfait).
Ces enveloppes sont intimement li\'ees aux alg\`ebres d'Auslander. Nous
les calculons dans le cas des alg\`ebres h\'er\'editaires de dimension
finie.

\medskip
Enfin, dans un appendice, nous donnons pour m\'emoire diverses
ca\-rac\-t\'e\-ri\-sa\-tions des ca\-t\'e\-go\-ries semi-simples, dont la
plupart sont sans doute bien connues des sp\'ecia\-lis\-tes.

\medskip
Les auteurs se sont beaucoup amus\'es \`a \'ecrire cet article, dont
l'\'e\-la\-bo\-ra\-tion s'est r\'ev\'el\'ee jusqu'au bout pleine de
rebondissements. Ils craignent que ce c\^ot\'e ludique n'\'echappe aux
lecteurs. Ils esp\`erent n\'eanmoins que ceux-ci prendront plaisir \`a
lire l'\'enonc\'e de certains th\'eor\`emes.

\subsection*{Post-scriptum} Apr\`es la soumission de cet article,
nous avons appris que Peter J. O'Sullivan avait introduit et \'etudi\'e
les ca\-t\'e\-go\-ries que nous avons baptis\'ees ca\-t\'e\-go\-ries de
Kimura au \S
\ref{kimura} de mani\`ere ind\'ependante, sous le nom de
\emph{ca\-t\'e\-go\-ries semi-positives}. Il a obtenu la plupart des
r\'esultats du \S \ref{kimura}, et bien plus, \`a l'exception
toutefois du th\'eor\`eme de nilpotence \ref{nouvellekim}.
N\'e\-an\-moins, il a pu utiliser ses r\'esultats pour obtenir une
d\'emonstration du th\'eor\`eme \ref{th1} ci-dessus dans le cas
particulier des ca\-t\'e\-go\-ries de Kimura (qu'il faudrait maintenant
rebaptiser ca\-t\'e\-go\-ries de Kimura-O'Sullivan), totalement
diff\'erente de la n\^otre et beaucoup plus g\'eom\'etrique. Il a ainsi
d\'emontr\'e l'existence des enveloppes pro-r\'eductives de mani\`ere
ind\'ependante de notre travail \cite{os}.

Par ailleurs, O'Sullivan a apport\'e une r\'eponse compl\`ete \`a notre
question \ref{q19.1}: voir son appendice dans cet article.

\newpage

\addtocontents{toc}{{\bf I. Radicaux et nilpotence}\hfill\thepage}
\
\bigskip
\begin{center}
\large\bf I. Radicaux et nilpotence
\end{center}
\bigskip
L'objectif de cette partie est l'\'etude des ca\-t\'e\-go\-ries qui sont
extensions d'une ca\-t\'e\-go\-rie additive semi-simple par un ``id\'eal"
(localement) nilpotent. Com\-me en alg\`ebre non commutative, la notion
de radical joue ici un r\^ole de premier plan. On \'etudie aussi ce qui
lui arrive par extension des scalaires.

\section{Id\'eaux et radicaux}\label{s1}

Soit $K$ un anneau commutatif unitaire. Pour pr\'evenir les paradoxes
ensemblistes, il est utile de fixer d\`es
\`a pr\'esent un univers $\sU$ (contenant $K$). Les ensembles
d'objets et de fl\`eches d'une \emph{petite} ca\-t\'e\-go\-rie se trouvent donc
dans $\sU$.

\subsection{$K$-ca\-t\'e\-go\-ries, ca\-t\'e\-go\-ries $K$-li\-n\'e\-aires et
$K$-ca\-t\'e\-go\-ries \allowbreak pseu\-do-a\-b\'e\-lien\-nes}\label{biprod}
Par
\emph{$K$-ca\-t\'e\-go\-rie}, nous entendons une ca\-t\'e\-go\-rie
telle que pour tout couple d'objets
$(A,B)$, les mor\-phis\-mes de $A$ vers $B$ forment un $K$-module $\sA(A,B)$,
et que la composition des mor\-phis\-mes soit $K$-lin\'eaire (pour $K=\Z$, on
dit aussi ca\-t\'e\-go\-rie pr\'e-additive).

Un \emph{$K$-foncteur} entre deux $K$-ca\-t\'e\-go\-ries est un foncteur
$K$-lin\'eaire (sur les $K$-modules de mor\-phis\-mes).

\'Etant donn\'e deux objets $A,A'\in \sA$, un \emph{biproduit} de
$(A,A')$ est un syst\`eme $(C,i,i',p,p')$, avec $C\in \sA$, $i\in
\sA(A,C)$, $i'\in \sA(A',C)$, $p\in \sA(C,A)$, $p'\in \sA(C,A')$, le tout
v\'erifiant les identit\'es
\begin{equation}\label{eqbiprod}
pi=1_A,\quad p'i'=1_B,\quad ip+i'p'=1_C.
\end{equation}

Un tel biproduit, s'il existe, est \`a la fois un produit et un
coproduit et est d\'etermin\'e \`a isomor\-phis\-me unique pr\`es. On dit que
$\sA$ est \emph{$K$-lin\'eaire} si tout couple $(A,A')$ d'objets de $\sA$
admet un biproduit.  Le choix d'un tel biproduit pour chaque couple
d'objets $(A,A')$ d'une petite ca\-t\'e\-go\-rie
$K$-lin\'eaire $\sA$ d\'efinit alors un $K$-foncteur
$\oplus:\sA\times\sA\to \sA$, d\'etermin\'e \`a isomor\-phis\-me naturel
\emph{unique} pr\`es, qui fait de $\sA$ une ca\-t\'e\-go\-rie mo\-no\-\"{\i}\-dale
sy\-m\'e\-tri\-que (l'unit\'e \'etant l'objet nul).

D'apr\`es une variante de \cite[ch. VIII, prop. 4]{maclane},  un
foncteur $T:\sA\to \sB$ entre deux ca\-t\'e\-go\-ries $K$-lin\'eaires est un
$K$-foncteur si et seulement s'il transforme un biproduit en un
biproduit. Si on s'est donn\'e des biproduits $\oplus$ sur
$\sA$ et $\sB$, il existe alors un isomor\-phis\-me naturel canonique
$\oplus \circ (T,T)\cong T\circ \oplus$.

Enfin, une $K$-ca\-t\'e\-go\-rie $\sA$ est
\emph{pseudo-ab\'elienne} si tout projecteur a un
noyau (et une image): pour tout objet
$A\in
\sA$ et tout
\'el\'ement
$e\in
\sA(A,A)$ v\'erifiant $e^2=e$, il existe un objet $B\in \sA$ et un
mor\-phis\-me
$f:B\to A$ tel que $f$ soit un noyau de $(0,e)$. Si $(C,g)$ est un noyau
de $1-e$, le mor\-phis\-me $f\coprod g:B\oplus C\to A$ est alors un
isomor\-phis\-me.

Toute $K$-ca\-t\'e\-go\-rie ab\'elienne est $K$-lin\'eaire et
pseudo-ab\'elienne.

\begin{rem} La terminologie varie suivant les auteurs. Ainsi, Saavedra
\cite[I.0.1.2]{saavedra} appelle ca\-t\'e\-go\-rie $K$-lin\'eaire ce que nous
appelons
$K$-ca\-t\'e\-go\-rie. De m\^eme, une ca\-t\'e\-go\-rie pseudo-ab\'elienne est
souvent appel\'ee \emph{karoubienne}\footnote{Cette notion n'est pas
sp\'ecifiquement additive, \cf par exemple
\protect\cite[VI.4.1.2.1]{saavedra}.}.
\end{rem}

On a les d\'efinitions suivantes, \'etroitement li\'ees aux
pr\'ec\'edentes:

\begin{defn}\label{ep-loc} Une sous-ca\-t\'e\-go\-rie pleine $\sB$ d'une
$K$-ca\-t\'e\-go\-rie $\sA$ est
\begin{thlist}
\item \emph{\'epaisse} si elle est stable par facteur directs
(re\-pr\'e\-sen\-ta\-bles dans $\sA$);
\item \emph{localisante} si elle est \'epaisse et stable par sommes
directes quelconques (re\-pr\'e\-sen\-ta\-bles dans $\sA$).
\end{thlist}
\end{defn}

\begin{sorite}\label{so0} Soient $\sA,\sB$ deux $K$-ca\-t\'e\-go\-ries et $T,T'$
deux
$K$-foncteurs de $\sA$ vers $\sB$. Soit $u:T\Rightarrow T'$ une
transformation naturelle. Alors, pour tout biproduit $A\oplus B$ dans
$\sA$, on a
\[u_{A\oplus B}=u_A\oplus u_B\]
modulo les isomor\-phis\-mes naturels $T(A\oplus B)\simeq T(A)\oplus
T(B)$.\qed
\end{sorite}

\subsection{Compl\'etions}\label{rap} \'Etant
donn\'e une $K$-ca\-t\'e\-go\-rie $\sA$, on
peut former son enveloppe
$K$-lin\'eaire
$\sA^\oplus$ (adjonction de sommes directes finies,
cf.
\cite[\S 2]{kelly}), son enveloppe pseudo-ab\'elienne
$\sA^\natural$
(scindage
d'idempotents, cf.
\cite[\S 1,\S
11]{mitchell}), et son enveloppe $K$-lin\'eaire
pseudo-ab\'elienne
$\sA^{\oplus\natural}$ (parfois
appel\'ee compl\'etion projective, ou
de Cauchy, ou de Karoubi\dots),
formant un carr\'e de
sous-ca\-t\'e\-go\-ries
pleines
\begin{equation}\label{eq1}
\begin{matrix}
\;&\;& \sA^\oplus &\;&\;\\
\;&\nearrow\;&\;&\searrow &\;\\
      \sA &\;&\;&\;& \sA^{\oplus\natural}. \\
      \;&\searrow\;&\;&\nearrow &\;\\
                        \;&\;&
\sA^\natural &\;&\; \end{matrix}\end{equation}

Pour r\'ef\'erence, rappelons les
constructions de $\sA^\oplus$ et de $\sA^\natural$:

\subsubsection{$\sA^\oplus$} Les objets sont des
suites finies
$(A_1,A_2,\dots, A_n)$ d'objets de $\sA$, not\'ees $A_1\oplus A_2
\oplus \dots \oplus A_n $ \cite[\S 2]{kelly}. Les mor\-phis\-mes entre
$A_1\oplus A_2
\oplus \dots \oplus A_n $ et $A'_1\oplus A'_2
\oplus \dots \oplus A'_m $ sont donn\'es par les matrices dont le
coefficient d'ordre $(i,j)$ est dans
$\sA(A_i,A'_j)$ (se composant selon la r\`egle usuelle).

Cette construction est $2$-covariante en $\sA$: un foncteur $T: \sA \to
\sB$ \'etant donn\'e, on lui associe
foncteur $T^\oplus: \sA^\oplus \to \sB^\oplus$ construit sur $T$
composante par composante. Il est donc
``strictement additif" par d\'efinition: l'isomor\-phis\-me canonique
$T^\oplus(X)\oplus T^\oplus(Y)=
T^\oplus(X\oplus Y)$ est l'identit\'e pour tout couple d'objets 
$(X,Y)$ de $\sA^\oplus$.
De m\^eme, \`a une transformation naturelle $u:T\Rightarrow T'$
correspond une transformation naturelle $u^\oplus:T^\oplus\Rightarrow
{T'}^\oplus$ sur le mod\`ele du sorite \ref{so0}.

\subsubsection{$\sA^\natural$} Les objets sont les couples $(A,e)$ avec
$A\in\sA$ et $e\in \sA(A,A)$ v\'erifiant $e^2=e$. Pour un autre tel
couple $(B,f)$, on pose
\[\sA^\natural((A,e),(B,f))=f\sA(A,B)e\subset \sA(A,B).\]

\begin{sloppypar}
Les mor\-phis\-mes se composent de mani\`ere \'evidente. Cette construction
est \'egalement $2$-covariante en $\sA$ ($T^\natural(A,e)=(T(A),T(e))$,
$u^\natural_{(A,e)}=T'(e)u_A T(e)$).
\end{sloppypar}

\subsubsection{Propri\'et\'es universelles} Les
deux constructions ci-dessus ont des propri\'et\'es universelles. Quitte
\`a introduire des $2$-ca\-t\'e\-go\-ries, on peut les interpr\'eter  comme des
adjoints \`a gauche. Soit $\{K\}$ la $2$-ca\-t\'e\-go\-rie dont les objets sont
les petites $K$-ca\-t\'e\-go\-ries, les $1$-mor\-phis\-mes les
$K$-foncteurs et les $2$-mor\-phis\-mes les transformations naturelles. Soient
$\{K\}^\oplus,\{K\}^\natural,\{K\}^{\oplus\natural}$ les
sous-$2$-ca\-t\'e\-go\-ries pleines de $\{K\}$ (m\^emes $1$-mor\-phis\-mes et
m\^emes $2$-mor\-phis\-mes) form\'ees des ca\-t\'e\-go\-ries
$K$-li\-n\'e\-ai\-res, pseudo-ab\'eliennes et $K$-li\-n\'e\-ai\-res
pseudo-ab\'eliennes. Alors:

\begin{sorite}\label{so01} a) Pour tout couple de petites $K$-ca\-t\'e\-go\-ries
$(\sA,\sB)$, avec $\sB$ $K$-lin\'eaire, les foncteurs de restriction
\begin{align*}
\{K\}^\oplus(\sA^\oplus,\sB)&\to \{K\}(\sA,\sB)\\
\{K\}^\natural(\sA^\natural,\sB)&\to \{K\}(\sA,\sB)
\end{align*}
sont pleinement fid\`eles et surjectifs. De plus, pour tout $T\in
\{K\}(\sA,\sB)$, la fibre de $T$ est un groupo\"{\i}de \`a groupes
d'automor\-phis\-mes triviaux.\\
b) Pour tout couple de petites $K$-ca\-t\'e\-go\-ries
$(\sA,\sB)$, les foncteurs
\begin{align*}
\{K\}(\sA,\sB)&\to \{K\}^\oplus(\sA^\oplus,\sB^\oplus)\\
\{K\}(\sA,\sB)&\to \{K\}^\natural(\sA^\natural,\sB^\natural)
\end{align*}
sont pleinement fid\`eles (mais pas essentiellement surjectifs en
g\'en\'eral).\qed
\end{sorite}

(Le fait que les foncteurs consid\'er\'es en a) ne soient pas
bijectifs sur les objets provient des choix possibles de biproduits et
de noyaux.)

\subsection{Alg\`ebres \`a plusieurs objets}\label{1.3}
Dans cet article, nous adoptons le point de vue de B. Mitchell
\cite{mitchell} et
consid\'erons les $K$-ca\-t\'e\-go\-ries comme
des ``$K$-alg\`ebres (associatives) \`a plusieurs objets".

Suivant ce point de vue, un id\'eal (bilat\`ere) $\sI$
d'une
$K$-ca\-t\'e\-go\-rie $\sA$ est
la donn\'ee, pour tout
couple
d'objets $(A,B)$, d'un sous-$K$-module de $\sI(A,B)\allowbreak \subset
\sA(A,B)$
tel que
pour tout couple de mor\-phis\-mes $(f \in \sA(A,A'), g \in
\sA(B,B'))$, on
ait
$g\sI(A',B)f \subset \sI(A,B')$.

On peut alors
former la $K$-ca\-t\'e\-go\-rie quotient $\sA /\sI$ (qui a les
m\^emes
objets que $\sA$). Si $A=\oplus A_i, B= \oplus B_j$, alors
$\sI(A,B)$
s'identifie
ca\-no\-ni\-que\-ment \`a $\oplus_{i,j} \sI (A_i,B_j)$.

\begin{ex}\label{ex1.3.1} Si $X$ est un ensemble d'objets de $\sA$, stable 
par
biproduit, la famille d'ensembles
\[\sI_X(A,B)=\{f\in \sA(A,B)\mid f \text{ se factorise \`a travers un objet
de } X\}\]
est un id\'eal de $\sA$. Le quotient $\sA/\sI$ est souvent not\'e $\sA/X$.
\end{ex}

On a une
notion \'evidente d'id\'eal produit (\resp somme, \resp intersection) de
deux id\'eaux
$\sI\cdot\sJ$ (\resp $\sI+\sJ$, \resp $\sI\cap \sJ$). On prendra toutefois
garde de ne pas confondre les
id\'eaux $(\sI\cdot\sJ)(A,A)=\{\sum f\circ
g,B \in Ob\sA, f\in \sA(B,A),
g\in \sA(A,B)\}$ et
$\sI(A,A)\cdot\sJ(A,A)$ de la $K$-alg\`ebre
$\sA(A,A)$.

Le noyau d'un $K$-foncteur est l'id\'eal $\Ker T$ de $\sA$ form\'e des
mor\-phis\-mes que $T$ annule; $T$ induit une $K$-\'equivalence entre le
quotient $\sA/\Ker T$ et une sous-ca\-t\'e\-go\-rie non pleine de $\sB$.

\begin{lemme}\label{l1} Supposons $\sA$ $K$-lin\'eaire.
Soient $\sI,\sJ$
deux id\'eaux de $\sA$. Supposons que
$\sI(A,A)\subset \sJ(A,A)$ pour tout
objet $A$ de $\sA$. Alors
$\sI\subset \sJ$.
\end{lemme}

\prf Soient $A,B\in \sA$. On a
$\sA(A\oplus B,A\oplus B)=\sA(A,A)\oplus
\sA(A,B)\oplus
\sA(B,A)\oplus \sA(B,B)$. En utilisant les injections de
$A$ et $B$
dans $A\oplus B$ et les projections de $A\oplus B$ sur $A$ et
$B$, on
voit que $\sA(A,B)\cap \sI(A\oplus B,A\oplus B)=\sI(A,B)$, et
de
m\^eme pour $\sJ$. Le lemme en r\'esulte.\qed

Notons $K$-$Mod$ la
ca\-t\'e\-go\-rie $K$-lin\'eaire des $K$-modules.

\begin{defn}\label{D2 1/2} Un {\it $\sA$-module} (\`a gauche) est un
$K$-foncteur $M:\sA\to
K$-$Mod$.\\
Il est dit \emph{fini} si tous les $M(A)$ sont de type fini sur
$K$.
\end{defn}

Les $\sA$-modules (\resp les $\sA$-modules finis) 
forment une ca\-t\'e\-go\-rie
a\-b\'e\-lien\-ne $K$-lin\'eaire,
not\'ee 
$\sA\hbox{--}Mod$ \footnote{Notation
compatible \`a la pr\'ec\'edente 
si
    l'on consid\`ere $K$ comme
ca\-t\'e\-go\-rie \`a un seul objet.} 
(\resp $\sA\hbox{--}Modf$). Tout
$K$-foncteur $F:\sA\to \sB$
induit 
par
composition des $K$-foncteurs en sens 
inverse
$F^\ast:\sB\hbox{--}Mod\to \sA\hbox{--}Mod$ et 
$F^\ast:\sB\hbox{--}Modf\to
\sA\hbox{--}Modf$.

Nous aurons aussi besoin de la notion suivante:

\begin{defn}\label{d1} Soit $A\in \sA$. 
Un \emph{$A$-id\'eal \`a gauche}
de $\sA$ est la 
donn\'ee,
pour tout
$B\in \sA$, d'un sous-$K$-module $\sI(B)$ de 
$\sA(A,B)$,
telle que la
famille des
$\sI(B)$ soit stable par 
composition \`a
gauche.
\end{defn} 

On a une notion \'evidente 
d'id\'eal \`a gauche somme (\resp
intersection) de deux id\'eaux \`a 
gauche. 
 
\medskip Pour tout $A\in \sA$, notons ${}_A \sA$ le 
foncteur $B\mapsto
\sA(A,B)$. Ceci d\'efinit un $K$-foncteur 
contravariant
\begin{align*}
Y_\sA:\sA^{\rm o}&\to 
\sA\hbox{--}Mod\\
A&\mapsto  {}_A \sA.
\end{align*}

Pour que l'image 
de $Y_\sA$
soit contenue dans $\sA\hbox{--}Modf$, il faut et il 
suffit que
$\sA(A,B)$ soit un $K$-module de type fini pour tous 
$A,B\in \sA$. De ce
point de vue, un $A$-id\'eal \`a gauche n'est 
autre qu'un sous-objet de
${}_A \sA$.

On a dualement les 
ca\-t\'e\-go\-ries $Mod\hbox{--}\sA$ et $Modf\hbox{--}\sA$
des 
$\sA$-modules \`a droite\footnote{C'est celle consid\'er\'ee 
par
Street
\cite{street}.} et des $\sA$-modules \`a droite finis 
(contravariants),
ainsi que les
$A$-id\'eaux \`a droite. On a le
$K$-foncteur covariant
\begin{align*}
{}_\sA Y:\sA&\to Mod\hbox{--}\sA\\
A&\mapsto   \sA{}_A:B\mapsto \sA(B,A).
\end{align*}

\begin{defn}\label{d1.1} Un objet $X$ d'une ca\-t\'e\-go\-rie $\sC$ est
\emph{compact} si le foncteur $Y\mapsto \sC(X,Y)$ commute aux limites
inductives quelconques (re\-pr\'e\-sen\-ta\-bles dans $\sC$).
\end{defn}

Le r\'esultat suivant est bien connu (par exemple \cite[III.2 et
III.7]{maclane}, \cite{street}):

\begin{prop}\label{p1.1}
a) Pour tout $A\in \sA$ et tout $M\in
Mod\hbox{--}\sA$, l'homomor\-phis\-me
\begin{align*}
Mod\hbox{--}\sA(\sA{}_A,M)&\to M(A)\\
f&\mapsto f(1_A)
\end{align*}
est bijectif.\\
b) Le foncteur ${}_\sA Y$ est pleinement fid\`ele. Un
module appartenant \`a son i\-ma\-ge essentielle est dit
\emph{re\-pr\'e\-sen\-ta\-ble}.\\
c) La ca\-t\'e\-go\-rie
$Mod\hbox{--}\sA$ est stable par limites inductives et projectives
quelconques.\\
d) Tout objet $M$ de $Mod\hbox{--}\sA$ est limite inductive d'objets
re\-pr\'e\-sen\-ta\-bles.\\
e) Si $\sA$ est stable par limites inductives quelconques, le foncteur
${}_\sA Y$ admet un adjoint \`a  gauche
$\colim:Mod\hbox{--}\sA\to \sA$, qui en est aussi un inverse \`a gauche.\\
f) Pour tout objet $M$
de \allowbreak $Mod\hbox{--}\sA$, les conditions suivantes sont
\'equi\-va\-len\-tes:
\begin{thlist}
\item $M$ est compact.
\item Pour tout syst\`eme inductif $(N_i)_{i\in I}$ d'objets de
$Mod\hbox{--}\sA$ index\'e par une petite
ca\-t\'e\-go\-rie $I$, l'homomor\-phis\-me
\[\colim Mod\hbox{--}\sA(M,N_i)\to Mod\hbox{--}\sA(M,\colim N_i)\]
est surjectif.
\item $M$ est facteur direct d'une somme directe d'objets
re\-pr\'e\-sen\-ta\-bles\footnote{Si $\sA$ est $K$-lin\'eaire
pseudo-ab\'elienne, cette condition \'equivaut au fait que $M$ soit
re\-pr\'e\-sen\-ta\-ble.}.
\end{thlist}
De plus, $M$ est projectif.
\end{prop} 

Rappelons rapidement la 
d\'emonstration, pour la commodit\'e du lecteur.
Dans a)
(lemme de 
Yoneda enrichi), l'application inverse est donn\'ee par
$m\mapsto 
(f\mapsto M(f)m)$. b) en r\'esulte imm\'ediatement. c) est
clair, puisque $K\hbox{--}Mod$ a la m\^eme propri\'et\'e: les limites se
calculent ``argument par argument".

Pour voir d), on consid\`ere la
ca\-t\'e\-go\-rie $\sC$ dont les objets sont les couples $(A,x)$ avec $A\in \sA$
et $x\in M(A)$, un mor\-phis\-me $(A,x)\to (B,y)$ \'etant un mor\-phis\-me $f\in
\sA(A,B)$ tel que $M(f)(x)=y$: par a), le syst\`eme inductif
\begin{align*}
T:\sC&\to Mod\hbox{--}\sA\\
(A,x)&\mapsto \sA_A
\end{align*}
s'envoie vers $M$ et on v\'erifie que $\colim T\to M$ est un isomor\-phis\-me.

On d\'eduit e) de d) comme suit: soit $M\in Mod\hbox{--}\sA$. \'Ecrivons
$M=\colim \sA_{A_i}$: on d\'efinit alors $\colim M=\colim A_i$. On
utilise b) pour montrer que ce foncteur est bien d\'efini et pour
v\'erifier sa propri\'et\'e d'adjonction.

On d\'eduit \'egalement f) de a) et d) comme suit:
(iii) $\If$ (i) r\'esulte imm\'ediatement de a). (i) $\If$ (ii) est
\'evident. Pour voir que (ii) $\If$ (iii), soit
$M$ un objet v\'erifiant (ii): \'ecrivons $M=\colim_{i\in I} \sA_{A_i}$,
o\`u
$I$ est une petite ca\-t\'e\-go\-rie. Il existe alors un ensemble fini $I_0$
d'objets de $I$ tel que l'identit\'e: $M\to M$ se factorise \`a travers
$\displaystyle\bigoplus_{i\in I_0}\sA_{A_i}$.

Enfin, la derni\`ere assertion r\'esulte du lemme \ref{lA.2}.\qed

\begin{rems}\ \label{r1.1}
\begin{itemize}
\item[a)] On laisse au lecteur le soin d'obtenir les \'enonc\'es
correspondants pour $\sA\hbox{--}Mod$, en rempla\c cant $\sA$ par
$\sA^{\rm o}$.
\item[b)] On dit qu'un objet $M\in Mod\hbox{--}\sA$ est \emph{de
pr\'esentation finie} si le foncteur ${}_M Mod\hbox{--}\sA$ commute aux
limites inductives \emph{filtrantes}. (Comme il est additif, il commute
alors aussi aux sommes directes quelconques.) On se gardera de confondre
cette notion avec celle d'objet compact, qui est plus restrictive.
On peut montrer que $M$ est de pr\'esentation finie si et seulement s'il
s'ins\`ere dans une suite exacte
\[P_1\to P_0\to M\to 0\]
o\`u $P_1$ et $P_0$ sont compacts.
\end{itemize}
\end{rems}

Voici une r\'eciproque de la proposition \ref{p1.1}:

\begin{prop} Soit $\sA$ une ca\-t\'e\-go\-rie ab\'elienne stable par limites
inductives quelconques. Notons $\sA_\comp$ la sous-ca\-t\'e\-go\-rie pleine de
ses objets compacts. Supposons $\sA_\comp$ dense dans $\sA$,
c'est-\`a-dire que tout objet de $\sA$ est limite inductive d'objets de
$\sA_\comp$. Alors le foncteur ``de Yoneda"
\begin{align*}
\sA&\to Mod\hbox{--}\sA_\comp\\
A&\mapsto (\sA_A:C\mapsto \sA(C,A))
\end{align*}
est une \'equivalence de ca\-t\'e\-go\-ries.
\end{prop}

\prf La pleine fid\'elit\'e se d\'emontre comme dans le lemme de Yoneda
enrichi: Si $A\in \sA$ est limite inductive d'objets compacts $C_i$ et
$M\in Mod\hbox{--}\sA_\comp$, l'homomor\-phis\-me
\begin{align*}
Mod\hbox{--}\sA_\comp(\sA_A,M)&\to \lim M(C_i)\\
f&\mapsto (f(\iota_i))
\end{align*}
o\`u $\iota_i$ est le mor\-phis\-me $C_i\to A$, est un isomor\-phis\-me d'inverse
$(m_i)\mapsto (g\mapsto M(g_i)m_i)$ pour $C\in \sA_\comp$ et $g\in
\sA(C,A)=\colim \sA(C,C_i)$ \'ecrit sous la forme $\iota_i\circ g_i$
pour $i$ assez grand. On en d\'eduit que, pour $A,B\in \sA$,
\[Mod\hbox{--}\sA_\comp(\sA_A,\sA_B)\stackrel{\sim}{\to} \lim
\sA_B(C_i)=\lim \sA(C_i,B)=\sA(A,B).\]

\noindent Pour l'essentielle surjectivit\'e, on utilise le fait que $\sA_\comp$ est
dense dans $Mod\hbox{--}\sA_\comp$ d'apr\`es la proposition 1.3.6 d).\qed

La construction $\sA\mapsto Y_\sA$ \emph{n'est pas} $2$-fonctorielle en
$\sA$. Par contre:

\begin{cor}[\cf \protect{\cite[p. 140]{street}}]\label{morita} 
Soit
$F:\sA\to
\sB$ un
$K$-foncteur. Alors $F^*:\sB\hbox{--}Mod\to 
\sA\hbox{--}Mod$ est une
\'equivalence de ca\-t\'e\-go\-ries si et 
seulement si $F^{\oplus\natural}$ est
une \'equivalence de 
ca\-t\'e\-go\-ries. On dit alors que $F$ est 
une
\emph{\'e\-qui\-va\-len\-ce de Morita}.\qed
\end{cor}

\prf Cela 
r\'esulte imm\'ediatement de la proposition \ref{p1.1} 
f).\qed

\begin{lemme}\label{corresp} Soit $A\in \sA$; notons encore 
$A$ son image
dans
$\sA^\oplus$, $\sA^\natural$ et 
$\sA^{\oplus\natural}$. Alors 
\\a) Il
y a bijection entre les 
$A$-id\'eaux \`a gauche dans $\sA$,
$\sA^\oplus$,
$\sA^\natural$ 
et
$\sA^{\oplus\natural}$.
\\ b) Il
y a bijection entre les id\'eaux 
(bilat\`eres) dans $\sA$,
$\sA^\oplus$,
$\sA^\natural$ 
et
$\sA^{\oplus\natural}$.
\\c) Pour tout id\'eal $\sI$ de $\sA$, 
notons $\sI^\oplus$ l'id\'eal
correspondant de $\sA^\oplus$. Alors on 
a un isomor\-phis\-me canonique
$\sA^\oplus/\sI^\oplus= (\sA/\sI)^\oplus$.
\\d) Pour tout id\'eal $\sI$ de $\sA$, notons $\sI^\natural$ l'id\'eal
correspondant de $\sA^\natural$. Alors on a un foncteur 
canonique
pleinement fid\`ele
 $\sA^\natural/\sI^\natural\to 
(\sA/\sI)^\natural$. C'est un
isomor\-phis\-me si et seulement si les 
idempotents des alg\`ebres
d'endomor\-phis\-mes $(\sA/\sI)(A,A)$ se 
rel\`event en des idempotents
des al\-g\`e\-bres d'endomor\-phis\-mes 
$\sA(A,A)$.
 \end{lemme}

\prf a): Cela r\'esulte imm\'ediatement de 
l'interpr\'etation des
$A$-i\-d\'e\-aux \`a gauche
comme \'etant les 
sous-modules du $\sA$-module \`a
gauche
${}_A\sA:B\mapsto \sA(A,B)$ 
et du fait que le carr\'e
\eqref{eq1}
induit un carr\'e 
d'\'equivalences de
ca\-t\'e\-go\-ries (corollaire 
\ref{morita}):
\[
\begin{matrix}   \;&\;& 
\sA^\oplus\hbox{--}Mod
&\;&\;\\
\;&\nearrow\;&\;&\searrow 
&\;\\

\sA\hbox{--}Mod &\;&\;&\;& \sA^{\oplus\natural}\hbox{--}Mod. 
\\

\;&\searrow\;&\;&\nearrow &\;\\
 
\;&\;&
\sA^\natural\hbox{--}Mod &\;&\;
\end{matrix}\]
b): voir 
\cite[prop.
2]{street}

c) est imm\'ediat.

d): pour tout couple 
d'objets $(A,e), (A,e')$ de $\sA^\natural$, on 
a
\[\sA^\natural((A,e),(A',e'))=e'\sA(A,A')e\]
et
\[\sI^\natural((A,e) 
,(A',e'))=e'\sI(A,A')e.\]
 
On a  donc
\[\sA^\natural/\sI^\natural((A,e),(A',e'))= 
\bar
e'((\sA/\sI)^\natural(A,A'))\bar
e=(\sA/\sI)^\natural((A,\bar 
e),(A',\bar e'))
\]
o\`u $\bar e,\bar e'$ d\'esignent 
les images de $e$ et
$e'$ modulo $\sI$. Le foncteur $(A,e)\mapsto 
(A,\bar e)$ est donc
pleinement fid\`ele. L'hypoth\`ese de 
rel\`evement des idempotents
\'equivaut \`a dire qu'il est bijectif sur les objets.
\qed

\subsection{Radical} On d\'emontre \cite{kelly}
que
\begin{multline*}
\rad(\sA)(A,B)=\\
\{f\in \sA(A,B)\mid
\text{pour tout}\ g\in
\sA(B,A), 1_A-gf \text{  est
inversible}\}
\end{multline*}
d\'efinit un id\'eal $\rad(\sA)$ de
$\sA$.

\begin{defn}\label{D1} Cet id\'eal s'appelle le
\emph{radical} de $\sA$.
\end{defn}
Si $\sA$ est une $K$-alg\`ebre,
vue comme $K$-ca\-t\'e\-go\-rie \`a un seul
objet, on retrouve une
d\'efinition du
radical de Jacobson.

\begin{rem} \label{gabriel} Un autre radical a \'et\'e 
introduit dans la
th\`ese de Gabriel \cite{gabriel} (c'est celui qui est \'evoqu\'e dans
\cite{bru2}):
\[\rad'(\sA)(A,B)=
\{f\in \sA(A,B)\mid
\text{pour tout}\ g\in
\sA(B,A), gf \text{  est
nilpotent}\}.\]

On pourrait donc parler de radical de Gabriel et de radical de Kelly (le
dernier \'etant une version cat\'egorique du radical de Jacobson\dots) On
a \'evidemment $\rad'(\sA)\subset \rad(\sA)$, avec \'egalit\'e si et
seulement si, pour tout $A\in \sA$, $\rad(\sA)(A,A)$ est un nilid\'eal.
\end{rem}

La caract\'erisation suivante de $\rad'$ est parfois utile:

\begin{lemme}\label{rad'} Soit $H:\sA\to Projf_L$ un foncteur 
$K$-lin\'eaire de $\sA$ vers la ca\-t\'e\-go\-rie des $L$-modules 
projectifs de type fini, o\`u $L$ est une $K$-alg\`ebre. Alors, pour 
$(A,B)\in \sA\times\sA$, on a
\[\rad'\sA(A,B)\subset \{f\in \sA(A,B)\mid \forall g\in \sA(B,A), tr 
H(gf)=0\},\]
avec \'egalit\'e si $H$ est fid\`ele et $L$ est de caract\'eristique $0$.
\end{lemme}

\prf L'inclusion $\subset$ est claire, puisque la trace d'un 
endomor\-phis\-me nilpotent est nulle. Pour l'inclusion oppos\'ee, soit $f$ 
un \'el\'ement du second membre. Alors, pour tout $g\in \sA(B,A)$ et tout 
$n\ge 1$, on a
\[tr H(gf)^n= tr H((gf)^n) = tr((gf)^{n-1}gf) =0.\]

Comme $L$ est de caract\'eristique $0$, cela implique que $H(gf)$ est 
nilpotent, donc \'egalement $gf$ par fid\'elit\'e.\qed

\begin{sloppypar}
\begin{prop}\label{P1} $a)$ Le
radical est le plus grand id\'eal
$\rad(\sA)$ tel que
$\rad(\sA)(A,A)$
soit contenu (\resp co\"{\i}ncide) avec le radical
de Jacobson de
$\sA(A,A)$ pour tout objet $A$.\\
$b)$ C'est aussi le
plus grand id\'eal de $\sA$
tel que le foncteur quotient $\sA\to \sA 
/ \rad(\sA)$ soit
conservatif, i.e. refl\`ete les isomor\-phis\-mes.
Ce 
foncteur refl\`ete aussi r\'etractions 
et
cor\'etractions\footnote{Rappelons qu'une r\'etraction 
(\resp
cor\'etraction) est un mor\-phis\-me inversible \`a droite (\resp 
\`a
gauche); on dit aussi \'epimor\-phis\-me scind\'e (\resp monomor\-phis\-me
scind\'e).}.
\end{prop}
\end{sloppypar}

Voir \cite[th. 1]{kelly}, \cite[\S 4]{mitchell}, \cite[prop.
6]{street}. On voit tout de suite que
la notion de radical est auto-duale.

\begin{cor}\label{corresp1} Les radicaux des ca\-t\'e\-go\-ries
$\sA,\sA^\oplus,\sA^\natural$ et $\sA^{\oplus\natural}$ se
cor\-res\-pon\-dent par les bijections du lemme
\ref{corresp} b).
\end{cor}

\prf Cela d\'ecoule du point a) de la proposition \ref{P1}.\qed

\begin{defn}\label{D2} Soit $T:\sA\to \sB$ un $K$-foncteur
entre deux $K$-ca\-t\'e\-go\-ries. Alors $T$ est dit \emph{radiciel} si
$T(\rad(\sA))\subset
\rad(\sB)$.
\end{defn}

\begin{sloppypar}
\begin{lemme} \label{L1} Si $T$ est plein, il est radiciel. D'autre
part, si $T$ est conservatif, on a
$T^{-1}(\rad(\sB))\subset \rad(\sA)$.
\\
Si $T$ est radiciel et conservatif, on a $T^{-1}(\rad(\sB))=\rad(\sA)$.
Si $T$ est m\^eme plein,
$T(\rad(\sA))$ est \'egal \`a la restriction de $\rad(\sB)$ \`a $T(\sA)$.
\end{lemme}

\prf Le premier point d\'ecoule imm\'ediatement du point $a)$ de la 
proposition
\ref{P1}.
Prouvons le second: soit
$f\in
\sA(A,B)$ tel que
$T(f)\in \rad(\sB)(T(A),T(B))$. Alors pour tout
$g\in
\sA(B,A)$, $1_{T(A)}-T(g)\circ T(f)$ est inversible. Donc $1_A-g\circ
f$ est inversible puisque $T$ est
conservatif.
Prouvons la derni\`ere assertion: $T$
radiciel $\If$
$\sA\subset T^{-1}T(\rad(\sA))\subset T^{-1}(\rad(\sB))$, et $T$
conservatif
$\If$
$T^{-1}(\rad(\sB))\subset \rad(\sA)$, d'o\`u $T^{-1}(\rad(\sB))=
\rad(\sA)$. Si $T$ est plein, la derni\`ere assertion en d\'ecoule.\qed

\begin{meg}\label{meg1} a) Si $T$ est conservatif, 
l'inclusion
$T^{-1}(\rad(\sB))\subset \rad(\sA)$  devient fausse en 
g\'en\'eral si
l'on remplace les radicaux par leurs puissances 
$n$-i\`emes, $n>1$
(m\^eme si $T$ est pleinement fid\`ele). 
\\
b) Il 
d\'ecoule du premier point du lemme que pour tout id\'eal $\sI$ de 
$\sA$, $\rad(\sA/\sI)\,$
contient l'image de $\,\rad(\sA)$ dans 
$\sA/\sI$ (la projection $T:\sA\to \sA/\sI$ est un foncteur
plein); 
mais il en est en g\'en\'eral distinct si $\sI \not\subset \rad(\sA)$ 
(exemple: $\Z\to
\Z/4\Z$).  
 
\end{meg}
\end{sloppypar}

\begin{lemme} \label{L2} Soit $\sA$ une 
$K$-ca\-t\'e\-go\-rie, et
soit $S$ un objet de $\sA$ tel que $\sA(S,S)$ 
soit un corps (par exemple
un objet simple dans un ca\-t\'e\-go\-rie 
ab\'elienne). Soit
$A$ un objet quelconque de $\sA$. Alors tout 
mor\-phis\-me
$f \in \sA(S,A)$ qui n'est pas dans le radical est une 
cor\'etraction, et
tout mor\-phis\-me $g\in \sA(A,S)$ qui n'est pas dans 
le ra\-di\-cal est une
r\'etraction.
         Si en outre 
$\sA(S,S)\allowbreak=K$, alors
\[\rad({\sA})(S,A)= \{f \in \sA(S,A 
),
\forall g\in
\sA(A ,S),
         g\circ f =0\} . 
\]
\end{lemme}

\prf Par d\'efinition du radical, il existe $h \in 
\sA(A,S)$ tel que $1_S
- hf$ ne soit pas inversible dans $\sA(S,S)$. 
Comme $\sA(S,S)$ est un
corps, ceci entra\^{\i}ne que $1_S = hf $. 
Pour $g$ on raisonne
dualement.

Prouvons la derni\`ere assertion. 
Tout  $f \in \sA(S,A )$ tel que
$g\circ f =0$ pour tout $g\in
\sA(A 
,S)$ n'est pas une cor\'etraction, donc appartient au
radical. 
R\'eciproquement, si $f\in
\rad({\sA})(S,A)$, alors, par 
d\'efinition, l'\'el\'ement $1-g\circ f$
de $K$ est inversible, 
c'est-\`a-dire
non nul, pour tout $g \in \sA(A ,S)$ par d\'efinition. 
Comme $g$
peut
\^etre multipli\'e par
toute constante de 
$K$,
ceci
entra\^{\i}ne $g\circ f=0$. \qed

\section{ca\-t\'e\-go\-ries de Wedderburn}\label{s2}

Dans ce paragraphe, on suppose que \emph{$K$ est un corps.}

\subsection{Ca\-t\'e\-go\-ries semi-simples}

\begin{defn}\label{D3} Une $K$-ca\-t\'e\-go\-rie est dite
\emph{semi-simple} si tout $\sA$-module \`a gauche est 
semi-simple. 
\end{defn}

Cette notion est manifestement Morita-invariante.  En
particulier, elle est stable par passage \`a
$\sA^{\oplus},\sA^{\natural},\sA^{\oplus\natural}$. Elle 
a
\'et\'e
\'etudi\'ee dans \cite{leduc}, \cite[\S
4]{mitchell} 
et
\cite{street}\footnote{Prendre garde \`a la terminologie: les 
notions de
semi-simplicit\'e utilis\'ees dans
\cite{leduc} et dans 
\cite{street}
sont plus faibles.}.
Pour la commodit\'e
du lecteur, 
nous avons donn\'e en
appendice un expos\'e des 
nombreuses
caract\'erisations
des ca\-t\'e\-go\-ries semi-simples, et 
quelques
contre-exemples. Nous en extrayons la proposition suivante:

\begin{prop} \label{P4/3} Soit
$\sA$ une (petite) ca\-t\'e\-go\-rie $K$-lin\'eaire.
Les conditions sui\-van\-tes sont
\'equi\-va\-len\-tes: \\
a) $\;\sA$ est semi-simple,\\
b) $\;\rad(\sA)=0$
et pour tout
$A\in \sA$, $\sA(A,A)$ est
un anneau artinien,\\
c) pour tout objet
$A$,
$\sA(A,A)$ est une
$K$-alg\`ebre semi-simple\footnote{Selon la terminologie de Bourbaki
\cite{alg}; selon la terminologie
anglo-saxonne, on dirait plut\^ot ``semisimple artinian", \cf
\cite{rowen}.}.\\
Sous ces conditions, $\sA$ est ab\'elienne si et seulement si elle est
pseudo-ab\'elienne.
\end{prop}

En particulier, la notion de semi-simplicit\'e est
\emph{auto-duale}, et {\it stable par quotient par tout id\'eal}. Voir
\cite[\S 4]{mitchell} ou l'appendice de ce travail.

En particulier, si $\sA$ est
semi-simple (non n\'e\-ces\-sai\-re\-ment $K$-lin\'eaire ni pseudo-ab\'elienne),
on a $\Ext^i_\sA(M,N)=0$ pour tout couple
$(M,N)$ de
$\sA$-modules et tout $i>0$.

\begin{lemme} \label{l2.1} Toute sous-ca\-t\'e\-go\-rie pleine $\sB$ d'une
$K$-ca\-t\'e\-go\-rie semi-sim\-ple $\sA$ est semi-simple. De plus, si $\sA$ et
$\sB$ sont
$K$-lin\'eaires pseudo-a\-b\'e\-lien\-nes, il existe une unique
sous-ca\-t\'e\-go\-rie pleine $\sC$ de $\sA$ telle que $\sA=\sB\coprod \sC$
(\cf d\'efinition \ref{A.D7} c)).
\end{lemme}

\prf Par invariance de Morita, on se ram\`ene imm\'ediatement au cas o\`u
$\sA$ et $\sB$ sont $K$-lin\'eaires pseudo-ab\'eliennes. La premi\`ere
assertion r\'esulte alors de la caract\'erisation c) de la proposition
\ref{P4/3}. Pour voir la deuxi\`eme, notons $S(\sA)$ (\resp $S(\sB)$)
l'ensemble des classes d'isomor\-phis\-mes (types) d'objets simples de $\sA$
(\resp de $\sB$). Il est clair que les objets de $\sB$ sont les sommes
directes d'objets simples de type appartenant \`a $\sB$. La ca\-t\'e\-go\-rie
$\sC$ est alors la ca\-t\'e\-go\-rie des sommes directes d'objets simples de
$\sA$ de type appartenant \`a $S(\sA)-S(\sB)$.\qed

\begin{meg} Dans le lemme \ref{l2.1}, si $\sA$ est
mo\-no\-\"{\i}\-dale et 
si $\sB$ est une sous-ca\-t\'e\-go\-rie mo\-no\-\"{\i}\-dale de
$\sA$, il n'est 
pas vrai en g\'en\'eral que $\sC$ soit mo\-no\-\"{\i}\-dale.
\end{meg}

Le lemme suivant est \`a mettre en regard de \ref{meg1}b.

\begin{lemme}\label{NewL1} Soit $\sA$ une 
$K$-ca\-t\'e\-go\-rie telle que
$\sA/\rad(\sA)$ soit semi-sim\-ple.  Pour 
tout id\'eal $\sI$, le radical
de
$\sA/\sI$ est l'image de 
$\rad(\sA)$ dans $\sA/\sI$. \end{lemme}

\prf Le point est 
l'inclusion $\rad(\sA/\sI)\subset \sA/\rad(\sA)+\sI$. Il d\'ecoule de 
ce que le foncteur
de projection $ \sA/\sI\to \sA/\rad(\sA)+\sI$ est 
plein, donc radiciel, et de ce que le radical de $
\sA/\rad(\sA)+\sI$ 
est nul puisque c'est une ca\-t\'e\-go\-rie semi-simple (comme quotient de 
$\sA/\rad(\sA)$).  \qed

Voici enfin deux compl\'ements au lemme \ref{l2.1}, qui montrent que la 
structure des ca\-t\'e\-go\-ries semi-simples est ``triviale":

\begin{prop} Soit $\sA$ une petite ca\-t\'e\-go\-rie ab\'elienne semi-simple, o\`u 
$K$ est un corps commutatif. Notons $S(\sA)$ (resp $\bI(\sA)$, $C(\sA)$) 
l'ensemble des classes d'isomor\-phis\-mes d'objets simples de $\sA$ (resp. 
l'ensemble des id\'eaux de $\sA$, l'ensemble des sous-ca\-t\'e\-go\-ries 
strictement \'epaisses de $\sA$, c'est-\`a-dire strictement pleines et 
stables par somme directe et facteur direct). Consid\'erons les applications 
suivantes:
\begin{itemize}
\item $SC:C(\sA)\to S(\sA)$: $\sB\mapsto S(\sB)$.
\item $CS:S(\sA)\to C(\sA)$: $X\mapsto\{\bigoplus S_i\mid S_i\in X\}$.
\item $IC:C(\sA)\to \bI(\sA)$: $\sB\mapsto \sI_\sB$ (\cf exemple 
\ref{ex1.3.1}).
\item $SI: \bI(\sA)\to S(\sA)$: $\sI\mapsto \{S\in S(\sA)\mid 1_S\in \sI\}$.
\end{itemize}
Alors $SC,CS,IC$ et $SI$ sont bijectives, $SC$ et $CS$ sont inverses l'une 
de l'autre, $SI\circ IC = SC$ et $CS\circ SI = IC^{-1}$.
\end{prop}

\begin{sloppypar}
\prf Le fait que $SC\circ CS =Id$ est \'evident. Si $\sB\in C(\sA)$, on a 
$CS\circ SC(\sB)\subset \sB$ parce que $\sB$ est $K$-lin\'eaire, avec 
\'egalit\'e parce que $\sB$ est pseudo-ab\'elienne. De plus, $SI\circ 
IC(\sB)=\{S\mid 1_S \text{ se factorise \`a travers un objet de } \sB\}$ 
contient clairement $S(\sB)$. R\'eciproquement, pour un $S\in SI\circ 
IC(\sB)$, il existe $a:S\to X$ et $b:X\to S$ tels que $ba=1_S$. Soit 
$S'=a(S)$: restreignant $b$ \`a $S'$, on se ram\`ene au cas o\`u $X$ est 
simple. Ceci montre que $S\in S(\sB)$ puisque $\sB$ est strictement pleine. 
\end{sloppypar}

Il reste \`a montrer que $IC\circ CS\circ SI = Id$. Pour $\sI\in\bI$, 
$IC\circ CS\circ SI(\sI)$ est l'id\'eal des mor\-phis\-mes qui se factorisent 
\`a travers un objet de la forme $\bigoplus S_i$, o\`u les $S_i$ sont 
simples et tels que $1_{S_i}\in\sI$ pour tout $i$. Cet id\'eal est 
\'evidemment contenu dans $\sI$. R\'eciproquement, soit $f\in\sI$: 
\'ecrivons $f$ comme somme directe d'un mor\-phis\-me nul et d'un isomor\-phis\-me 
$u$ (lemme \ref{lA.3}). Alors $u\in\sI$. Si $S$ est un facteur simple de la 
source de $u$, $1_S$ se factorise \`a travers $u$, donc $1_S\in\sI$. \qed

\begin{prop} Soit $\sI$ un id\'eal de $\sA$. Alors:
\begin{enumerate}
\item Le foncteur de projection $\pi:\sA\to \sA/\sI$ est essentiellement 
surjectif (i.e. $\sA/\sI$ est pseudo-ab\'elienne).
\item $\pi$ a une section $i$, qui en est un adjoint \`a gauche et \`a 
droite.
\item Posons $\sB=CS\circ SI(\sI)$. Alors on a $\sA=\sB\coprod i(\sA/\sI)$.
\item La suite
\[0\to K_0(\sB)\to K_0(\sA)\to K_0(\sA/\sI)\to 0\]
est exacte scind\'ee.
\end{enumerate}
\end{prop}

\begin{sloppypar}
\prf Pour tout $A\in \sA$, $\sI(A,A)$ est un id\'eal bilat\`ere de la 
$K$-alg\`ebre semi-simple $\sA(A,A)$; l'homomor\-phis\-me d'anneaux $\sA(A,A)\to 
\sA(A,A)/\sI(A,A)=\sA/\sI(\pi(A),\pi(A))$ a donc une unique section $s_A$. 
En particulier, tout idempotent de $\sA/\sI(\pi(A),\pi(A))$ se rel\`eve dans 
$\sA(A,A)$, ce qui prouve 1). De m\^eme, on construit 2) \`a l'aide des 
$s_A$; ses propri\'et\'es sont imm\'ediates. 3) en r\'esulte gr\^ace \`a la 
proposition pr\'ec\'edente, et 4) r\'esulte de 3).\qed
\end{sloppypar}

\subsection{Ca\-t\'e\-go\-ries s\'eparables}

\begin{prop}\label{Psep} Soit $A$ une $K$-alg\`ebre 
(associative
unitaire). Les conditions suivantes sont 
\'equi\-va\-len\-tes:
\begin{enumerate}
\item \label{Psep1} $A$ est un 
$A$-bimodule projectif.
\item \label{Psep2} Toute d\'erivation de $A$ 
vers un bimodule est
int\'erieure. 
\item \label{Psep3} $A$ est de 
dimension finie sur $K$, est semi-simple,
et le reste apr\`es toute 
extension des scalaires.
\item\label{Psep4} $A$ est semi-simple, et 
le reste apr\`es toute
extension des scalaires.
\item\label{Psep5} 
$A$ est de dimension finie sur $K$, sans radical, et
le reste apr\`es 
toute extension des scalaires.  
\item\label{Psep6} L'alg\`ebre 
$A^{\rm o}\otimes_K A$ est semi-simple.
\end{enumerate}
 
\end{prop}

Pour \ref{Psep1} $\iff$ \ref{Psep2}, voir \cite[4]{cq}. 
Pour  \ref{Psep1}
$\iff$ \ref{Psep3} $\iff$ \ref{Psep4}, on renvoie 
\`a \cite[ch. IX, th.
7.10]{ce}. Pour \ref{Psep4} $\iff$ \ref{Psep5}, 
voir
\cite[7.5, cor.]{alg}. Pour \ref{Psep1} $\iff$ \ref{Psep6}, 
voir
\cite[ch. IX, th. 7.9]{ce}.

En particulier, l'hypoth\`ese de 
finitude dans \cite[ch. IX, th.
7.10]{ce} n'est pas n\'ecessaire. 
Pour la commodit\'e du lecteur, nous
reproduisons l'argument 
de
\cite[4]{cq} qui montre que \ref{Psep1} $\If$
$\dim_K A<\infty$. 
Soit $s$ une section de l'application de
multiplication $A\otimes_K 
A\to A$, comme homomor\-phis\-me de $A$-bimodules.
\'Ecrivons 
$s(1)=\sum_1^n x_i\otimes y_i$, avec $n$ minimal; en
particulier, les 
$y_i$ sont lin\'eairement ind\'ependants sur $K$, donc
s\'epar\'es 
par les formes $K$-lin\'eaires sur $A$. Notons
$I$ le $K$-espace de 
dimension finie $\sum_1^n Kx_i \subset A$. Alors $I$
est un id\'eal 
\`a gauche de $A$: en effet, pour tout $a\in A$, on a
$as(1)=s(1)a$, 
d'o\`u $\sum ax_ia^\ast(y_i)=\sum x_ia^\ast(y_ia)$ pour
tout $a^\ast 
\in Hom_K(A,K)$, et le r\'esultat. En outre, l'action de $A$
par 
multiplication \`a gauche sur
$I$ est fid\`ele du fait que $\sum x_i 
y_i= 1$. Par cons\'equent, $A$
s'injecte dans 
$End_K(I)$.
\qed

\begin{defn} \label{D4sep} Une $K$-alg\`ebre $A$ 
est dite
\emph{s\'eparable} (ou \emph{absolument semi-simple}) si 
elle v\'erifie les conditions de la
proposition 
\ref{Psep}.
\end{defn}

\begin{rem} Nous nous
\'ecartons ici 
l\'eg\`erement de la terminologie de Bourbaki
\cite[7.5]{alg}, pour 
qui une alg\`ebre est s\'eparable si elle est sans
radical et le 
reste apr\`es toute extension des scalaires. Dans le cas
commutatif, 
on dit d'ailleurs plut\^ot
\emph{\'etale} que 
s\'eparable.
\end{rem}

\begin{lemme}\label{l2.2} Toute alg\`ebre de 
dimension finie sans radical
sur un corps parfait $K$ est 
s\'eparable.
\end{lemme}

\prf Voir \cite[7.5,7.6]{alg}. 
\qed

\begin{lemme}\label{l2.3} Le produit tensoriel de deux 
$K$-alg\`ebres
s\'eparables est s\'e\-pa\-rab\-le.
\end{lemme}

\prf 
Voir \cite[ch. IX, prop. 7.4]{ce}.\qed

\begin{defn}\label{D4sep1} 
Une ca\-t\'e\-go\-rie $K$-lin\'eaire $\sA$ est dite
\emph{s\'eparable} (ou 
\emph{ab\-so\-lu\-ment semi-simple}) si la
ca\-t\'e\-go\-rie $\sA_L$ 
qui s'en d\'eduit\footnote{Voir la
d\'efinition \ref{d3} ci-dessous 
pour plus de d\'etails.} en tensorisant
les mor\-phis\-mes par $L$ est 
semi-simple pour toute extension
$L/K$.
\end{defn} 

La notion de 
ca\-t\'e\-go\-rie s\'eparable est Morita-invariante (\cf
proposition 
\ref{P4/3}), auto-duale, et stable par passage au quotient
par un 
id\'eal. 

D'apr\`es le th\'eor\`eme
\ref{tsep}, on en a les autres 
caract\'erisations sui\-van\-tes:

\begin{itemize}
\item (si $\sA$ 
est $K$-lin\'eaire) Pour tout objet
$A\in
\sA$,
$\sA(A,A)$ est une 
$K$-alg\`ebre s\'eparable.
\item La ca\-t\'e\-go\-rie 
enveloppante\footnote{Voir le \S
\ref{prodtens} ci-dessous pour plus 
de d\'etails.}
$\sA^e=\sA^{\rm o}\boxtimes_K\sA$ est 
semi-simple.
\end{itemize}

\begin{prop}\label{p2.1} Si $K$ est 
parfait, une $K$-ca\-t\'e\-go\-rie $\sA$ est
s\'eparable si et seulement si 
elle est semi-simple et $\dim_K
\sA(A,B)<\infty$ pour tout couple 
d'objets $(A,B)$.
\end{prop}

\prf Cela r\'esulte du lemme 
\ref{l2.2}.\qed

\subsection{Ca\-t\'e\-go\-ries 
semi-primaires}

\begin{defn}\label{D4} Une $K$-ca\-t\'e\-go\-rie $\sA$ 
est
\emph{semi-primaire}
si
\begin{thlist}
\item pour tout objet 
$A\in \sA$, le
radical
$\rad(\sA(A,A))$ est nilpotent;
\item 
$\sA/\rad(\sA)$ 
est
semi-simple.
\end{thlist}
\end{defn}

\begin{rems}\label{R1}\
\begin{itemize}
\item[$a)$]
Dans le cas d'une ca\-t\'e\-go\-rie \`a un 
seul
objet, on retrou\-ve une
notion connue en alg\`ebre non 
commutative sous
le nom d'``an\-neau semi-pri\-mai\-re"
\cite{rowen} 
\footnote{Ces anneaux
sont caract\'eris\'es par l'e\-xis\-tence d'une 
borne uniforme pour
la longueur de toute 
cha\^{\i}ne
d\'ecrois\-san\-te de mo\-du\-les
cycliques, \cf 
\cite[2.7.7]{rowen}.}.
Ainsi une ca\-t\'e\-go\-rie $K$-lin\'eaire 
est
semi-pri\-mai\-re si et seulement si tous ses
anneaux 
d'endomor\-phis\-mes sont semi-pri\-mai\-res.
\item[$b)$] Cette notion 
est auto-duale.

 \item[$c)$] Une ca\-t\'e\-go\-rie
$K$-lin\'eaire 
pseudo-ab\'elienne
semi-primaire n'est 
pas
n\'e\-ces\-sai\-re\-ment
ab\'elienne (exemple: la ca\-t\'e\-go\-rie des 
modules
projectifs de type
fini sur une $K$-alg\`ebre artinienne 
non
semi-simple). 

\item[$d)$] La condition (i) entra\^{\i}ne que 
les radicaux de Kelly
et de Gabriel co\-\"{\i}n\-ci\-dent.

\end{itemize}
\end{rems}  

L'un des aspects classiques des anneaux 
semi-primaires est la propri\'et\'e de rel\`eve\-ment des 
idempotents
modulo le radical, qui suit du lemme classique 
suivant

\begin{lemme}[\cf 
\protect{\cite[1.1.28]{rowen}}]\label{idemp} Soit $N$ un
nil-id\'eal 
d'un anneau $A$. Alors tout idempotent de $A/N$ se rel\`eve
en un 
idempotent de $A$.   \qed  
\end{lemme}

\begin{prop}\label{P3/2} 
a) Si $\sA$ est semi-primaire, il en est de
m\^eme de $\sA/\sI$ pour 
tout id\'eal
$\sI$, et le radical de $\sA/\sI$ est l'image du radical 
de $\sA$.\\ 
b) Supposons $\sA$ semi-primaire. Alors
$\sA/\rad(\sA)$ 
est pseudo-ab\'elienne si et seulement si $\sA$ l'est
(et m\^eme 
ab\'elienne si $\sA$ est
$K$-lin\'eaire). \\ 
Plus g\'en\'eralement, 
le foncteur canonique (\cf 
lemme \ref{corresp} d)) 
$\sA^\natural/\rad(\sA^\natural)\allowbreak\to
(\sA/\rad(\sA))^\natural$ est un isomor\-phis\-me.\\ 
c) La notion de $K$-ca\-t\'e\-go\-rie 
semi-primaire est Morita-invariante.
En par\-ti\-cu\-lier, si
$\sA$ 
est semi-primaire, 
alors
$\sA^{\oplus},\sA^{\natural},\sA^{\oplus\natural}$ le sont 
(et
r\'e\-ci\-pro\-que\-ment).\\ 
d) Si $\sA$ est $K$-lin\'eaire et 
si $\sA(A,A)$ est un anneau
artinien pour tout
$A\in \sA$, alors 
$\sA$ est semi-primaire. \\ 
e) Si $\sA$ est
semi-primaire,
$\sB$ est 
semi-simple et
$T:\sA\to\sB$ est un $K$-foncteur plein, alors il 
existe une unique
factorisation de $T$ en un $K$-fonc\-teur 
plein
$\bar T:\sA/\rad(\sA)\to \sB.$
\\ f) Si $\sA$ est
semi-primaire 
et pseudo-ab\'elienne,
et $T:\sA\to\sB$ est un $K$-foncteur radiciel 
qui n'envoie aucun objet
non nul de $\sA$ sur l'objet nul de $\sB$, 
alors  $T$  est conservatif.
\end{prop}   

\prf a) Cela r\'esulte 
ais\'ement du lemme \ref{NewL1}.

b) suit des lemmes \ref{corresp} d) 
et \ref{idemp} (et du fait qu'une
ca\-t\'e\-go\-rie $K$-lin\'eaire 
pseudo-ab\'elienne semi-simple est
ab\'elienne, \cf \ref{A.P4}). 

c) 
Pour l'invariance par \'equivalence de Morita, il s'agit d'apr\`es 
le
corollaire \ref{morita} de voir que $\sA$ est semi-primaire si 
et
seulement si
$\sA^{\oplus\natural}$ l'est. D'apr\`es le corollaire 
\ref{corresp1},
la trace sur $\sA$ du radical de 
$\sA^{\oplus\natural}$ est le radical de
$\sA$. Il suit de ceci, du 
lemme \ref{corresp} d) et du lemme \ref{l2.1}
que
$\sA$ est 
semi-primaire si
$\sA^{\oplus\natural}$ l'est. 

Supposons 
r\'eciproquement $\sA$ semi-primaire. Le point b)
entra\^{\i}ne que 
$\sA^{\natural}$ l'est. On peut donc supposer 
$\sA$
pseudo-ab\'elienne. En vertu du lemme \ref{corresp} c), il 
suffit de 
voir que si les radicaux de
$\sA^\oplus(A,A)$ 
et
$\sA^\oplus(B,B)$ sont nilpotents, il en est de m\^eme du radical 
de
$\sA^\oplus(A\oplus B,A\oplus B)$. 

Supposons donc que 
$\sR(A,A)^m=\sR(B,B)^n=0$, et montrons qu'alors
$\sR(A\oplus 
B,A\oplus B)^{m+n+2}=0$. Par additivit\'e, on se ram\`ene 
\`a
supposer que le compos\'e de toute cha\^{\i}ne de
$N\ge 2n+2$ 
mor\-phis\-mes $f_i\; (i=1,\dots, N)$ composables appartenant 
\`a
$\sR(A,A),\sR(A,B),\sR(B,A)$ ou $\sR(B,B)$ est nul.

Consid\'erons par exemple une cha\^{\i}ne partant de $A$ et 
aboutissant en
$A$. La composition correspondante appartient \`a un 
produit de la 
forme
\begin{multline*}
\sR(A,A)^{m_{k+1}}\sR(B,A)\sR(B,B)^{n_k}\sR(A, 
B)\sR(A,A)^{m_{k}}\cdots\\
\cdots\sR(A,A)^{m_2}\sR(B,A)\sR(B,B)^{n_1}\ 
sR(A,B)\sR(A,A)^{m_1}
\end{multline*}
avec 
$m_1+\dots+m_{k+1}+n_1+\dots+n_k+2k+1=N$. Ce produit est contenu 
dans
$\sR(A,A)^a$, avec $a=m_1+\dots+m_{k+1}+k$. D'autre part, le 
produit
interm\'ediaire
\[ 
\sR(B,B)^{n_k}\sR(A,B)\sR(A,A)^{m_{k}}
\cdots\sR(A,A)^{m_2}\sR(B,A)\sR 
(B,B)^{n_1}
\]
est contenu dans $\sR(B,B)^b$, avec 
$b=n_1+\dots+n_k+k-1$. On a donc
$a\ge m$ ou $b\ge n$ et le mor\-phis\-me 
compos\'e
est nul. Les autres cas se traitent de la m\^eme mani\`ere.

d) Un raisonnement analogue \`a celui
fait dans le cas des 
anneaux
montre que, pour toute $\sA$, le radical de $\sA/\rad(\sA)$ 
est r\'eduit
\`a $0$. En particulier,  pour tout objet $A$ de $\sA$, 
la $K$-alg\`ebre
$(\sA/\rad(\sA))(A,A)$ est semi-simple. Comme
$\sA$ 
est suppos\'ee $K$-lin\'eaire, on conclut par la
proposition 
\ref{P4/3} que $\sA$ est semi-simple.   

e) suit de la premi\`ere 
assertion du lemme \ref{L1}.

f) Le $K$-foncteur radiciel $T$ induit 
un $K$-foncteur
\[\bar T: (\sA/\rad(\sA))^\oplus\to 
(\sB/\rad(\sB))^\oplus .\] 

Tout comme $T$, $\bar T$ n'envoie aucun 
objet non nul de
$(\sA/\rad(\sA))^\oplus$ sur l'objet nul de 
$(\sB/\rad(\sB))^\oplus$.
D'autre part, $(\sA/\rad(\sA))^\oplus$ est 
ab\'elienne semi-simple
puisque $\sA$ est semi-primaire 
pseudo-ab\'elienne (utiliser le point b)
ci-dessus). Il en d\'ecoule 
que
$\bar T$ est conservatif (lemme \ref{utile}), et il en est de 
m\^eme de
$T$ par
\ref{P1} b.
\qed 

\begin{prop}\label{ablfspr} Soit 
$\sA$ une ca\-t\'e\-go\-rie ab\'elienne dont
tout objet est de lon\-gueur 
finie. Alors
$\sA$ est semi-primaire. 
\end{prop} 

\prf Le plus 
court est d'envoyer $\sA$ par un foncteur pleinement
fid\`ele dans 
une ca\-t\'e\-go\-rie de modules (th\'eor\`eme de plongement 
de
Freyd-Mitchell). Dans le cas d'une ca\-t\'e\-go\-rie de modules, 
l'assertion
revient \`a ceci: l'anneau d'endomor\-phis\-mes de tout 
module $M$ de
longueur finie est semi-primaire. Ce r\'esultat, qui 
repose sur le lemme
de Fitting, est bien 
connu,
\cf
\cite[2.9.10]{rowen} (la longueur du module est une borne 
pour
l'\'echelon de nilpotence du radical).\qed

\subsection{Ca\-t\'e\-go\-ries de Wedderburn}  

\begin{defn} 
\label{D4wed} Une $K$-ca\-t\'e\-go\-rie $\sA$ est dite de
\emph{Wedderburn} 
si 
\begin{thlist}
\item[(i)] pour tout objet $A\in \sA$, le radical 
$\rad(\sA(A,A))$ est
nilpotent;
\item[(ii')] $\sA/\rad(\sA)$ est une 
$K$-ca\-t\'e\-go\-rie s\'eparable.  
\end{thlist}
\end{defn}

\begin{meg} 
Il est clair qu'une $K$-ca\-t\'e\-go\-rie de
Wedderburn est semi-primaire. 
Nous verrons en \ref{C3.1.} qu'elle le
reste apr\`es toute extension 
des scalaires. La r\'eciproque est
\'evidemment \emph{fausse} au 
moins quand $K$ est imparfait, comme le
montre l'exemple de la 
$K$-ca\-t\'e\-go\-rie \`a un seul objet donn\'e par
une extension finie 
ins\'eparable de $K$.
\end{meg}

\begin{rem}\label{R1wed}
Dans le cas 
d'un seul objet, les alg\`ebres correspondantes
sont celles 
aux\-quelles s'appliquent le th\'eor\`eme de scindage de
Wedderburn 
dont  il sera question au \S \ref{s5}, d'o\`u la terminologie.
Nous 
les appellerons \emph{alg\`ebres de 
Wedderburn}.
\end{rem}

\begin{prop}\label{Pwed} a) La notion de 
$K$-ca\-t\'e\-go\-rie de Wedderburn
est stable par passage au quotient par 
un id\'eal.\\
b) La notion de $K$-ca\-t\'e\-go\-rie de Wedderburn
est 
Morita-invariante. En particulier, si
$\sA$ est de Wedderburn, 
alors
$\sA^{\oplus},\sA^{\natural},\sA^{\oplus\natural}$ le sont 
(et
r\'e\-ci\-pro\-que\-ment).\\ 
c) Si
$K$ est parfait et 
si
$\sA(A,B)$ est de dimension finie sur $K$ pour tout 
couple
$(A,B)\in Ob(\sA)$, alors
$\sA$ est de Wedderburn. En 
particulier, toute ca\-t\'e\-go\-rie tannakienne sur $K$ est de Wedderburn. 

        \end{prop}

\prf a) d\'ecoule de \ref{P3/2} a): si $\sA$ est 
de Wedderburn, $\sA/\sI$ est semi-primaire, et le quotient
par son 
radical (qui n'est autre que l'image du radical de
$\sA\,$) est 
s\'eparable, en tant que quotient de $\sA/\rad(\sA)$.
 
b) On sait 
d\'ej\`a que
$\sA^{\oplus},\sA^{\natural},\sA^{\oplus\natural}$ 
sont
semi-primaires. Il est alors
facile de conclure, compte-tenu de 
ce que, par
le lemme \ref{corresp} et la proposition \ref{P3/2} b), 
le quotient de
$\sA^{\oplus\natural}$ par son radical est
isomorphe 
\`a $(\sA/\rad(\sA))^{\oplus\natural} $.  

c) Il revient au m\^eme 
de dire que toute alg\`ebre d'endomor\-phis\-mes de
$\sA^\oplus$ est de 
dimension finie. D'apr\`es le point d) de la
proposition
\ref{P3/2}, 
on voit que $\sA^\oplus$, donc $\sA$, est semi-primaire. On
conclut 
par le fait que toute alg\`ebre de dimension finie sur un
corps 
parfait est s\'eparable.
\qed

\section{Radical infini et nilpotence 
renforc\'ee}\label{radrenf} 

\subsection{Le radical infini}

Soit 
$\sA$ une $K$-ca\-t\'e\-go\-rie.

\begin{defn}\label{d3.1} Le radical 
infini de $\sA$ est l'id\'eal
$\displaystyle
\rad^\omega(\sA)= 
\bigcap_{n\geq 1}
\rad^n(\sA)$.
\end{defn}

Cet id\'eal intervient dans diverses questions, \eg dans la th\'eorie
du type de re\-pr\'e\-sen\-ta\-tion des
$K$-alg\`ebres de dimension finie (pour $\sA = A$-$Modf$), et en liaison
avec les conjectures de Bloch-Beilinson-Murre sur les motifs (pour $\sA=$
la ca\-t\'e\-go\-rie des motifs de Gro\-then\-dieck pour l'\'equivalence
rationnelle - motifs `de Chow'- \`a coefficients dans
$K$, \cf \ref{BBM})

  \begin{lemme}\label{radinf} Si $T: \sA \to \sB$ est un $K$-foncteur
radiciel entre $K$-ca\-t\'e\-go\-ries, alors $T(\rad^\omega(\sA))\subset
\rad^\omega(\sB)$.
\end{lemme}

C'est clair. \qed

\begin{sloppypar}
\begin{defn} \label{D4swed} Une
$K$-ca\-t\'e\-go\-rie $\sA$ est dite
\emph{strictement semi-primaire (\resp strictement de Wedderburn)} si
\begin{thlist}
\item[(i')] pour
tout couple d'objets $(A,B)$, il existe un
entier $n>0$ tel que $(\rad(\sA))^n(A,B)=0$;
\item[(ii')] $\sA/\rad(\sA)$ est une $K$-ca\-t\'e\-go\-rie semi-simple
(\resp s\'eparable).
         \end{thlist}
\end{defn}
\end{sloppypar}

Il est clair qu'une ca\-t\'e\-go\-rie strictement semi-primaire (\resp
strictement de Wedderburn) est semi-primaire (\resp de
Wedderburn) et que son radical infini est nul. La r\'eciproque est vraie
  si tous les anneaux d'endomor\-phis\-mes sont artiniens.

Par ailleurs, il est clair qu'une ca\-t\'e\-go\-rie semi-primaire n'ayant 
qu'un nombre fini d'objets est strictement
semi-primaire.

\begin{contrex}\label{Cssp} Consid\'erons la ca\-t\'e\-go\-rie $\sA$
suivante:
\begin{itemize}
\item $Ob (\sA)=E$, o\`u $E$ est un ensemble
ordonn\'e dense quelconque
(dense signifie qu'entre deux \'el\'ements
distincts on peut toujours en trouver
un troisi\`eme distinct);
\item $\sA(x,y)=K$ si $x\leq y$, $0$ sinon;
\end{itemize}
la composition \'etant induite par la multiplication dans
$K$. On a alors
\[(\rad\,\sA)(x,y) =
\begin{cases} K &\text{si 
$x<y$}\\
0&\text{sinon.} \end{cases}\]
Ainsi $(\rad\,\sA)(x,x)=0$ 
pour tout $x$, donc $\sA$ est semi-primaire.
Mais si $x<y$, 
$(\rad\,\sA)(x,y)=(\rad\, \sA )^n(x,y)$ pour tout
$n>0$ (gr\^ace \`a 
la densit\'e de $E$). On a donc $\rad^\omega(\sA)=\rad(\sA)\neq 
0$.

Si $E$ est un groupe commutatif ordonn\'e, il est facile de 
munir $\sA$
d'une  structure mo\-no\-\"{\i}\-dale sy\-m\'e\-tri\-que rigide. 
Quitte \`a passer
\`a $\sA^{\oplus\natural}$, on peut m\^eme 
construire un exemple
$K$-lin\'eaire pseudo-ab\'elien (mais non 
ab\'elien).  
 \end{contrex} 

\subsection{Cas des ca\-t\'e\-go\-ries de 
modules}\label{ccm} Prenons \`a
pr\'esent pour
$\sA$ la ca\-t\'e\-go\-rie 
des {\it modules de type fini sur une alg\`ebre
$A$ de dimension 
finie sur un corps $K$}. Dans ce cas, le th\'eor\`eme
d'existence des 
suites d'Auslander-Reiten montre que le radical $\sR$ 
de
$A\hbox{--}Modf$ est engendr\'e (\`a droite ou \`a gauche) par 
les
mor\-phis\-mes dits irr\'eductibles (\ie dans $\sR\setminus \sR^2$), 
\cf
\cite{sisk} ou \cite{ks} pour plus de pr\'ecisions. Ce point, 
et
l'analyse d\'etaill\'ee des carquois d'Auslander-Reiten, sont \`a 
la base
du r\'esultat suivant, qui lie le type de re\-pr\'e\-sen\-ta\-tion 
de
$A$ au radical infini 
$\sR^\omega$.

\begin{thm}[\protect{\cite{sisk}, 
\cite{ks}}]\label{t3.1}  a) $A$ est de
type de re\-pr\'e\-sen\-ta\-tion fini 
(\ie il n'y a qu'un nombre fini de classes
d'isomor\-phis\-mes 
d'ind\'e\-com\-po\-sa\-bles dans $A\hbox{--}Modf$) 
$ \iff \sR^\omega= 0 \iff 
\exists n\;$ tel que
$\;\sR^n=0$. \\ 
b) (Supposons $K$ 
alg\'ebriquement clos\footnote{Il est probable que
cette hypoth\`ese, 
qui figure dans
\cite{ks}, peut \^etre affaiblie en: $K$ 
parfait.})
$A$ est de type de re\-pr\'e\-sen\-ta\-tion infini sauvage
$ \iff 
\sR^\omega$ n'est pas nilpotent.   
 \end{thm}

Un cas 
interm\'ediaire est celui de l'alg\`ebre $\bar\F_2[(\Z/2\Z)^2]$
(de 
type de repr\'e\-sen\-ta\-tion infini mod\'er\'e); $\sR^\omega$ 
est
non nul mais nilpotent.    
 
Un exemple standard d'alg\`ebre de 
type de re\-pr\'e\-sen\-ta\-tion infini
sauvage est donn\'e par le quotient 
de
$K[[T_1,T_2]]$ par le cube de l'id\'eal maximal, ou m\^eme par 
l'id\'eal
$(T_1^2, T_1T_2^2, T_2^3)$ (Drozd, \cf
\cite{gabrielrev}, 
\cite{nathanson}, \cite[1]{ringel}); il s'ensuit que
$A\hbox{--}Modf$ 
n'est pas strictement semi-primaire. 

Nous verrons plus loin qu'une 
ca\-t\'e\-go\-rie tannakienne alg\'ebrique sur
un corps de 
caract\'eristique nulle n'est pas n\'e\-ces\-sai\-re\-ment
strictement 
semi-primaire. 

Tout ceci indique que la notion de ca\-t\'e\-go\-rie 
strictement semi-primaire
est beaucoup trop restrictive pour nos 
besoins, et motive notre choix de
traiter syst\'ematique\-ment des 
ca\-t\'e\-go\-ries semi-primaires (ou de
Wedderburn).

\begin{sloppypar}
\begin{rem} Un exemple remarquable de 
ca\-t\'e\-go\-rie strictement de
Wedderburn devrait \^etre fourni par la 
ca\-t\'e\-go\-rie des motifs purs pour
une
\'e\-qui\-va\-le\-nce ad\'equate 
quelconque; cet exemple conjectural est
bri\`evement discut\'e en 
\ref{ex2} ci-dessous.  \end{rem}
\end{sloppypar}

\section{Radical et 
extension des scalaires}\label{s3}

\subsection{\ } Il est connu que 
le radical d'une alg\`ebre se comporte
``mal" par extension des 
scalaires en g\'en\'eral, \cf
\cite{alg}, \cite{rowen}. La situation 
s'am\'eliore toutefois dans le
cas ``de 
Wedderburn".

\begin{prop}\label{L5} 1) Soit $L/K$ une extension de 
corps. Soit $A$
une $K$-al\-g\`e\-bre de Wedderburn, \ie
$\rad(A)$ 
est nilpotent et $A/\rad(A)$ est s\'e\-pa\-ra\-ble. Alors \\
a) 
$\rad(A\otimes_K L)\cap A= \rad(A)$,\\
b) $\rad(A)\otimes_K L= 
\rad(A\otimes_K L)$,\\
c) $A\otimes_K L$ est de Wedderburn.

2) Si 
$A\otimes_K L$ est semi-primaire et si $\rad(A \otimes_K L)\subset 
\rad(A)\otimes_K L$ pour toute
extension $L/K$, alors
$A$ est de 
Wedderburn. 
               \end{prop}

\prf 1) Il est clair que 
$\rad(A)\otimes_K L$ est un id\'eal nilpotent
de $A\otimes_K L$, donc 
contenu dans
$\rad(A\otimes_K L)$. D'autre part $A/\rad(A)$ est 
s\'eparable, donc
aussi $(A/\rad(A))\otimes_K L= (A\otimes_K 
L)/(\rad(A)\otimes_K L)$. En particulier, cette derni\`ere
est 
semi-simple, donc
$\rad(A)\otimes_K L$ contient $\rad(A\otimes_K L)$. 
D'o\`u $b)$ et
$c)$. Le point $a)$ en d\'ecoule
imm\'ediatement.

2) 
Supposons $A\otimes_K L$ semi-primaire pour toute extension $L/K$. En 
particulier $\rad(A)$ est nilpotent,
donc $\rad(A \otimes_K L)\supset 
\rad(A)\otimes_K L$, et on a finalement \'egalit\'e compte tenu 
de
l'hypoth\`ese. Il suit que $(A/\rad(A))\otimes_K L\cong 
(A\otimes_K L)/\rad(A \otimes_K L)$ est semi-simple
pour toute 
extension $L/K$, donc que $A/\rad(A)$ est s\'eparable. Ainsi $A$ est 
de Wedderburn. 
\qed

\begin{rem}\label{imparf} La condition que 
$\rad(A \otimes_K L)\subset \rad(A)\otimes_K L$ est
n\'ecessaire: 
si
$K$ est imparfait, toute extension finie ins\'eparable $A/K$ non 
triviale fournit un exemple de $K$-alg\`ebre
``absolument 
semi-primaire" mais pas de Wedderburn.
\end{rem}

\medskip
Si $\sA$ 
est une $K$-ca\-t\'e\-go\-rie, notons $\sA_L$ la $L$-ca\-t\'e\-go\-rie qui
s'en 
d\'eduit en tensorisant les mor\-phis\-mes par $L$.

\begin{meg} Si $\sA$ 
est pseudo-ab\'elienne, il n'en est pas de m\^eme de
$\sA_L$ en 
g\'en\'eral \footnote{Voir \S\ref{ext} ci-dessous pour plus
de 
d\'etails.}. Si $\sA$ est ab\'elienne, il n'en est pas de m\^eme 
de
$(\sA_L)^\natural$ en g\'en\'eral. Un contre-exemple est fourni 
par la
ca\-t\'e\-go\-rie des re\-pr\'e\-sen\-ta\-tions de dimension finie 
de
$\bG_a\times \bG_a$ (\cf fin du \S 
\ref{jm}.)
\end{meg}

\begin{sloppypar}
\begin{cor}\label{C3.1.} Soit 
$\sA$ une ca\-t\'e\-go\-rie de Wedderburn. Alors,
 pour toute extension 
$L/K$, $\sA_L$ est de Wedderburn. De plus
\begin{thlist}
 \item $ 
\rad(\sA_L)= \rad(\sA)\otimes_K L$.
\item Le foncteur d'extension des 
scalaires $\sA\to \sA_L$ est radiciel.
\item Le 
carr\'e
\[\begin{CD}
\sA@>>> \sA_L\\
@VVV @VVV\\
\sA/\rad(\sA)@>>> 
\sA_L/\rad(\sA_L)
\end{CD}\]
est naturellement cocart\'esien: le 
foncteur naturel
$\sA_L/\rad(\sA_L)\to (\sA/\rad(\sA))_L$ est une 
\'equivalence 
de
ca\-t\'e\-go\-ries.
\qed
\end{thlist}
\end{cor}
\end{sloppypar}

\begin{thm} \label{P2} Soient $K$ un corps, $L$ une extension 
de
$K$,
$V$ un $L$-espace vectoriel de dimension finie, et
$A$ 
une
sous-$K$-alg\`ebre de $End_L(V)$. On suppose que, pour tout 
$a\in
A$, le poly\-n\^o\-me
caract\'eristique de $a$ dans $V$ est \`a 
coefficients dans $K$.
Alors \\
a) Le radical $R$ de $A$ est 
nilpotent d'\'echelon
$\le n= \dim V$.\\ 
b) Si $K$ est infini, $A/R$ 
est semi-simple, produit d'au plus $n$
composants simples. Si de plus 
$K$ est parfait, $A/R$ est s\'eparable.
\\
c) Supposons $K$ infini et 
parfait. Alors le radical de la
sous-$L$-alg\`ebre $AL$ de
$End_L(V)$ 
engendr\'ee par $A$ est $RL$, on a
$RL\cap A=R$, et l'application 
canonique $(A/R)\otimes_K L\to AL/RL$ est
bijective.
\end{thm}

\prf 
a) D'apr\`es \cite[\S 11, ex. 1 a)]{alg}, $R$ est \emph{a priori} 
un
nil-id\'eal. Montrons m\^eme que tout $u\in R$ v\'erifie
$u^n=0$. 
Pour cela, il suffit de voir que les valeurs propres de
$u$ dans $V$ 
sont toutes nulles. Supposons le contraire et soit
$\lambda\ne 0$ 
une
telle valeur propre. Comme $\lambda$ est alg\'ebrique sur $K$, il 
existe
un poly\-n\^o\-me $Q\in K[T]$ tel que $\lambda Q(\lambda)=1$. Mais 
$1-uQ(u)$
est inversible, donc $1-\lambda Q(\lambda)$ est inversible, 
contradiction.

Comme l'image de $R$ dans $End_L(V)$ est un 
semi-groupe
multiplicatif form\'e d'\'el\'ements nilpotents, un 
th\'eor\`eme de N.
Jacobson permet de conclure que $R$
est nilpotent 
d'\'echelon $\le \dim V$, \cf \cite[2.6.30]{rowen}.

b) Cela 
r\'esulte de \cite[\S 11, ex. 1 d)]{alg}.

c) Notons provisoirement 
$R_L$ le radical de $AL$. D'apr\`es la
proposition
\ref{L5} b), on 
a
$\rad(A\otimes_K L)=R\otimes_K L$, donc l'homomor\-phis\-me 
surjectif
$A\otimes_K L\to AL$ induit un
homomor\-phis\-me surjectif 
$(A/R)\otimes_K L\to AL/R_L$. Comme $A/R$
est s\'eparable et $R$ 
nilpotent, le th\'eor\`eme de Wedderburn
(\cf th\'eor\`eme \ref{T1} 
ci-dessous, dans le cas d'un seul objet)
permet de relever $A/R$ dans 
$A$. Notons
ce relev\'e $B$, et soit $K'$ le centre d'un facteur 
simple de $B$: $K'$ est
une extension finie de $K$. Soit $x$ un 
\'el\'ement primitif de $K'$: comme le
poly\-n\^o\-me caract\'eristique 
de $x$ op\'erant sur $V$ est \`a
coefficients dans $K$, c'est 
n\'e\-ces\-sai\-re\-ment une puissance du poly\-n\^o\-me
caract\'eristique de $x$ 
op\'erant sur $K'$. Ceci implique que la
surjection $K'\otimes_K 
L\to
K'L$ est bijective, et donc que $B\otimes_K L\to BL$ est 
bijective. Il en
r\'esulte que $A/R\otimes_K L\to AL/R_L$ est 
bijective, d'o\`u l'on
d\'eduit $R_L=RL$. Enfin, $(R_L\cap A=) RL\cap 
A = R$ d\'ecoule
imm\'ediatement de la nilpotence de $R$. \qed

\begin{rems} \label{R2}\ 
\begin{itemize}
\item[a)]  L'hypoth\`ese 
que $K$ est parfait n'est pas superflue dans b)
et c). On obtient un 
contre-exemple en prenant pour
 $A=L$ une extension finie 
ins\'eparable de $K$, $V= L\otimes_K L$ avec
la structure gauche de 
$L$-espace, et
$A\inj End_L(V)$ donn\'e par $\ell\mapsto 
(\ell_1\otimes \ell_2 \mapsto
\ell_1\otimes \ell\ell_2)$.
\item[b)] 
Si l'on suppose que $V=\bigoplus V_i$ est un $L$-espace
vectoriel 
gradu\'e et que $A$ est une sous-alg\`ebre de $\prod End(V_i)$,
la 
borne sur l'\'echelon de nilpotence $n$ de $R$ dans le point a) 
du
th\'eor\`eme \ref{P2} se raffine en:
$n\le max\dim V_i$ (m\^eme 
d\'emonstration). 
\item[c)] Le point a) du th\'eor\`eme \ref{P2} 
reste valable si on
remplace l'hypoth\`ese que pour tout $a\in
A$, le 
poly\-n\^o\-me caract\'eristique de $a$ est \`a coefficients dans 
$K$
par: $K$ est de
cardinal strictement sup\'erieur \`a $\dim_K A$, 
\cf
\cite[no. 12, ex. 
16]{alg}.
\end{itemize}
\end{rems}

\begin{sloppypar}
\begin{cor}\label{c6}  Soient $K$ un corps parfait infini, $L$ une
extension de
$K$, 
$\sA$ une
ca\-t\'e\-go\-rie $K$-lin\'eaire et $H:\sA\to Vec_L$ un 
$K$-foncteur
fid\`ele de $\sA$ vers la ca\-t\'e\-go\-rie des $L$-espaces 
vectoriels
de dimension finie. On suppose que, pour tout 
endomor\-phis\-me $a\in
\sA(A,A)$, le 
poly\-n\^o\-me
ca\-rac\-t\'e\-ris\-ti\-que de $H(a)$ est
\`a 
coefficients dans $K$. 
Alors\\
a) $\sA$ est une $K$-ca\-t\'e\-go\-rie de 
Wedderburn.\\
b) Soit $H(\sA) L$ la sous-ca\-t\'e\-go\-rie
(non pleine) 
$L$-lin\'eaire de $Vec_L$ engendr\'ee par $H(\sA)$: les
objets de 
$H(\sA) L$ sont les objets de $H(\sA)$ et, pour deux tels
objets 
$X,Y$, 
$\big(H(\sA)L\big)(X,Y)=\allowbreak\big(H(\sA)(X,Y)\big)L$.
Alors:\\ 
i) $\rad(H(\sA)L)=\rad(H(\sA))L$.\\
ii) $\rad(H(\sA)L)\cap 
H(\sA)=\rad(H(\sA))$.\\
iii) Le foncteur $\sA/\rad(\sA)_L\to 
H(\sA)L/\rad(H(\sA)L)$ induit par $H$
est une \'equivalence de 
ca\-t\'e\-go\-ries. 
\\
En particulier, $H(\sA) L$ est une $L$-ca\-t\'e\-go\-rie 
de Wedderburn. \qed 
\end{cor}
\end{sloppypar}

\begin{rem}\label{rem 
en fam} La conclusion $a)$ vaut encore si, au lieu
d'un $L$ et d'un 
$H$, on a une famille finie d'extensions de corps
$L_i/L$ et une 
famille fid\`ele de foncteurs $H_i:
\sA\to Vec_{L_i}$ et si l'on 
suppose que pour tout $a\in
\sA(A,A)$, le 
poly\-n\^o\-me
ca\-rac\-t\'e\-ris\-ti\-que de $H_i(a)$ est
\`a 
coefficients dans $K$. En effet, plongeant les $L_i$ dans un 
corps
commun $L$, et posant
$H(A)= \bigoplus_i H_i(A)\otimes_{L_i} 
L$, on se ram\`ene au cas du
corollaire.   
\end{rem}

 \begin{ex} 
Comme nous le verrons en \ref{wedtr}, le corollaire pr\'ec\'edent 
s'applique notamment au cas des
ca\-t\'e\-go\-ries tannakiennes sur un 
corps de caract\'eristique nulle. 
\end{ex}

\section{Extensions des 
scalaires na\"{\i}ve et non na\"{\i}ve}\label{ext}

Dans ce 
paragraphe, $K$ est un corps commutatif. Nous allons 
comparer
l'extension des scalaires ``na\"{\i}ve" \'etudi\'ee dans le 
paragraphe
pr\'ec\'edent \`a une extension des scalaires plus 
sophistiqu\'ee, due
\`a Saavedra, dont nous aurons besoin dans les 
paragraphes \ref{compl} et 
\ref{jm}.

\subsection{Deux types 
d'extensions des scalaires}\label{extsc}

\begin{defn}\label{d3} 
Soient $\sA$ une $K$-ca\-t\'e\-go\-rie et $L$ une
$K$-alg\`ebre 
(commutative unitaire).\\ 
a) La ca\-t\'e\-go\-rie $\sA_L$ a pour objets 
les objets de $\sA$; si $A,B\in
\sA_L$, $\sA_L(A,B)\allowbreak 
=L\otimes_K \sA(A,B)$. On appelle $\sA_L$
l'\emph{extension des 
scalaires na\"{\i}ve} de $\sA$ de $K$ \`a $L$.\\
b) On note 
$\sA_{(L)}=Rep_K(L,\sA)$ la ca\-t\'e\-go\-rie des
$K$-re\-pr\'e\-sen\-ta\-tions de 
$L$ dans $\sA$: un objet de $\sA_{(L)}$ est un
couple $(A,\rho)$, 
o\`u $A\in \sA$ et $\rho$ est un homomor\-phis\-me
de $K$-alg\`ebres 
unitaires $L\to
\sA(A,A)$; un mor\-phis\-me $f:(A,\rho)\to (A',\rho')$ 
est un \'el\'ement de
$\sA(A,A')$ commutant \`a $\rho$ et $\rho'$. On 
appelle $\sA_{(L)}$
l'\emph{extension des scalaires \`a la Saavedra 
de $\sA$ de $K$ \`a $L$}
(\cf 
\cite[II.1.5]{saavedra}).
\end{defn}

\begin{lemme} \label{l10} a) La 
$K$-structure de $\sA_{(L)}$ s'enrichit
en une $L$-structure par la 
loi
\[\lambda f=\rho'(\lambda) f=f\rho(\lambda)\]
pour $\lambda\in L$ 
et $f\in \sA_{(L)}((A,\rho),(A',\rho'))$.\\
b) Si $\sA$ est additive (\resp pseudo-ab\'elienne, ab\'elienne, stable
par li\-mi\-tes inductives quelconques), il en est de m\^eme de
$\sA_{(L)}$.\qed
\end{lemme}

\begin{rems}\label{r4}\
\begin{itemize}
\item[a)] Pour que cette construction soit raisonnable, il faut que $\sA$
soit ``assez grosse" ou que $L$ soit ``assez petit". Par exemple, si
$\sA(A,A)$ est de
$K$-dimension finie pour tout objet $A$, alors $\sA_{(L)}=0$ d\`es que
$L$ est un corps infini contenant $K$.
\item[b)] Si $\sA$ est semi-simple, il ne suit pas que
$\sA_{(L)}$ le soit: prendre pour $L/K$ une extension finie ins\'eparable
de corps, et pour
$\sA$ la $K$-ca\-t\'e\-go\-rie \`a un objet $L$ (ou, si on pr\'ef\`ere $Vec_L$).
C'est toutefois le cas (trivialement) si $L$ est une extension infinie de
$K$, car dans ce cas $\sA_{(L)}=0$ (voir a))
\item[c)]
Soit
$\sA= Rep^\infty_K G$ la Ind-ca\-t\'e\-go\-rie attach\'ee \`a la ca\-t\'e\-go\-rie
$Rep_K G$ des
$K$-re\-pr\'e\-sen\-ta\-tions de dimension finie d'un $K$-groupe affine.
Alors
$\sA_{(L)}$ s'identifie canoniquement \`a $Rep^\infty_L(G)$, \cf
\cite[III.1]{saavedra}. M\^eme remarque avec $Rep_K(G)$ quand $L/K$ est
finie.
\item[d)] Soit $\sA$ une $K$-ca\-t\'e\-go\-rie quelconque. Alors
$(Mod\hbox{--}\sA)_{(L)}$ s'identifie \`a la ca\-t\'e\-go\-rie des
$K$-foncteurs de $\sA$ vers $Mod\hbox{--}L$. Cette derni\`ere s'identifie
canoniquement \`a $Mod\hbox{--}\sA_{L}$ via l'extension des scalaires
na\"{\i}ve. D'apr\`es la proposition \ref{p1.1} f), les objets compacts
de cette ca\-t\'e\-go\-rie s'identifient \`a $(\sA_L)^{\oplus\natural}$.
\end{itemize}
\end{rems}

\subsection{Sur l'extension des scalaires \`a la Saavedra}\label{s.saa}
Notre but est de comparer $\sA_L$ \`a $\sA_{(L)}$. On a un ``foncteur
fibre"
\'evident
\begin{equation}\label{eqB.2}
\omega_L:\sA_{(L)}\to \sA
\end{equation}
qui est fid\`ele et conservatif, commutant aux limites projectives et
inductives quelconques, exact si
$\sA$ est ab\'elienne.

La premi\`ere chose
\`a faire est de d\'efinir un foncteur en sens inverse, de
$\sA$ vers $\sA_{(L)}$. Ce n'est possible que sous des hypoth\`eses
restrictives.

\begin{lemme}[\cf \protect{\cite[prop. II.1.5.1.1]{saavedra}}]\label{l7}
a) Soient
$A\in
\sA$ et
$V\in Vec_K$. Le foncteur
\[B\mapsto Hom_K(V,\sA(A,B))\]
est re\-pr\'e\-sen\-ta\-ble dans les cas suivants:
\begin{thlist}
\item $\sA$ est $K$-lin\'eaire et $\dim_K V<\infty$.
\item $\sA$ est stable par limites inductives quelconques.
\end{thlist}
On note un repr\'esentant $V\otimes_K A$. \\
b) La construction $(V,A)\mapsto
V\otimes_K A$ est bifonctorielle et bilin\'eaire. Si
$W$ est un autre
$K$-espace vectoriel, on a un isomor\-phis\-me canonique
\[W\otimes_K(V\otimes_K A)\simeq (W\otimes_K V)\otimes_K A\]
pourvu que les deux membres aient un sens. (Nous nous autoriserons de la
canonicit\'e de cet isomor\-phis\-me pour passer sous silence tout tel
pa\-ren\-th\'e\-sa\-ge.)\\
c) Si $\sA$ est $K$-lin\'eaire
et que $\dim_K V<\infty$, alors pour tout
$B\in \sA$, le foncteur
\[A\mapsto \sA(V\otimes_K A,B)\]
est re\-pr\'e\-sen\-ta\-ble par $V^*\otimes_K B$, o\`u $V^*$ est le dual de $V$.
\qed
\end{lemme}

Soient $A,B\in \sA$ et $V\in Vec_K$. Supposons que les objets $V\otimes_K
A$ et $V\otimes_K B$ soient d\'efinis. L'homomor\-phis\-me de fonctorialit\'e
\[\sA(A,B)\to \sA(V\otimes_K
A,V\otimes_K B)\simeq Hom_K(V,\sA(A,V\otimes_K B))\]
fournit un autre homomor\-phis\-me
\begin{equation}\label{eqB.1}
V\otimes_K \sA(A,B)\to \sA(A,V\otimes_K B).
\end{equation}

\begin{lemme}\label{l11} L'homomor\-phis\-me \eqref{eqB.1} est un
isomor\-phis\-me dans les cas suivants:
\begin{thlist}
\item $\dim_K V<\infty$;
\item $A$ est compact (\ie le foncteur $B\mapsto \sA(A,B)$ commute aux
limites inductives quelconques)
\end{thlist}
\end{lemme}

\prf (i) \'Ecrire $V\simeq K^n$ et se ramener \`a $n=1$ par additivit\'e.
(ii) passer \`a la limite sur (i).\qed

\begin{hyps}\label{h1}
\`A partir de maintenant, on suppose que $\sA$ est additive
(c'est-\`a-dire $K$-lin\'eaire). On se donne une $K$-alg\`ebre $L$
(commutative unitaire). On suppose:
\begin{itemize}
\item soit que $\dim_K L<\infty$;
\item soit que $\sA$ est stable par limites inductives quelconques.
\end{itemize}
\end{hyps}

\begin{lemme}\label{l8} Sous les
hypoth\`eses \ref{h1},\\
a) La donn\'ee
d'une re\-pr\'e\-sen\-ta\-tion $\rho:L\to \sA(A,A)$ comme dans la d\'efinition
\ref{d3} b) \'equivaut \`a la donn\'ee d'un mor\-phis\-me
$\rho^\vee:L\otimes_K A\to A$ tel que le diagramme
\[\begin{CD}
L\otimes_K L\otimes_K A@>1\otimes \rho^\vee>> L\otimes_K A\\
@V{\mu_L\otimes 1}VV @V{\rho^\vee}VV\\
L\otimes_K A@>\rho^\vee>> A,
\end{CD}\]
o\`u $\mu_L$ est la multiplication de $L$, soit commutatif.\\
b) Le mor\-phis\-me
\[\begin{CD} \rho_A^\vee:L\otimes_K L\otimes_K A@>\mu_L\otimes 1>>
L\otimes_K A
\end{CD}\]
v\'erifie la condition de a), donc d\'efinit une re\-pr\'e\-sen\-ta\-tion
canonique $\rho_A:L\to \sA(L\otimes_K A,L\otimes_K A)$.\qed
\end{lemme}

Sous les hypoth\`eses \ref{h1}, le lemme \ref{l8} d\'efinit un
$K$-foncteur
\begin{align}
\sA&\to \sA_{(L)}\label{eqB.3}\\
A&\mapsto (L\otimes_K A,\rho_A).\notag
\end{align}

\begin{lemme}[\protect{\cite[II.1.5.2]{saavedra}}]\label{l12} Le foncteur
\eqref{eqB.3} est adjoint \`a gauche au foncteur
$\omega_L$ de \eqref{eqB.2}. Il envoie objet compact
de $\sA$ sur objet compact de $\sA_{(L)}$.
\end{lemme}

(La seconde assertion provient du fait que $\omega_L$ commute aux limites
inductives quelconques.)\qed

Le lemme \ref{l10} b) montre que le foncteur \eqref{eqB.3} s'\'etend alors
en un $L$-foncteur
\begin{equation}\label{eqB.4}
\Xi_L:\sA_L\to \sA_{(L)}.
\end{equation}

\subsection{Comparaison de $\sA_L$ et
$\sA_{(L)}$}

\begin{prop}\label{p6} Soient $A,B$ deux objets de $\sA$. Si
$\dim_K L=\infty$, supposons que $A$ soit compact. Alors l'homomor\-phis\-me
\[L\otimes_K\sA(A,B)\to \sA_{(L)}(L\otimes_K A,L\otimes_K B)\]
est bijectif. En particulier, $\Xi_L$ est pleinement fid\`ele si
$\dim_K L<\infty$, et sa restriction \`a la sous-ca\-t\'e\-go\-rie pleine des
objets compacts de $\sA_L$ est pleinement fid\`ele en g\'en\'eral (et
prend ses valeurs dans les objets compacts de $\sA_{(L)}$ d'apr\`es le
lemme \ref{l12}).
\end{prop}

\prf Cela r\'esulte des lemmes \ref{l11} et \ref{l12}.\qed

\begin{thm}\label{t2} Le
foncteur
\[\Xi_L^\natural:(\sA_L)^\natural\to \sA_{(L)}\]
induit par $\Xi_L$ est une \'equivalence de ca\-t\'e\-go\-ries pour toute 
$K$-ca\-t\'e\-go\-rie $\sA$
pseudo-ab\'elienne si et seulement si $L$ est une $K$-alg\`ebre
\emph{\'etale}. Dans ce cas, $\Xi_L^\natural$ induit une
\'e\-qui\-va\-le\-nce de ca\-t\'e\-go\-ries sur les sous-ca\-t\'e\-go\-ries pleines
form\'ees des objets compacts.
\end{thm}

\prf  Vu la proposition \ref{p6}, le foncteur
$\Xi_L^\natural$ est une \'equivalence de ca\-t\'e\-go\-ries si et
seulement s'il est essentiellement surjectif. D'autre part, d'apr\`es
la proposition \ref{Psep}, $L$ est \'etale sur $K$ si et seulement si le
$L$-bimodule $L$ est projectif.

1) Supposons d'abord $L$ projectif sur $L\otimes_K L$. La suite exacte de
$L$-bimodules $0\to
\Ker \mu \to L\otimes_K L \stackrel{\mu}{\longrightarrow} L\to 0$ est
scind\'ee; on note
$s$ un scindage.

Dire que $\Xi_L^\natural$ est essentiellement
surjectif est dire que tout objet $B$ dans
$\sA_{(L)}$ est facteur direct de
$L\otimes_K A$ pour un objet $A$ de $\sA$ convenable (comme
$\Xi_L^\natural$ est pleinement fid\`ele, l'expression ``facteur
direct" n'est pas ambig\"ue).  Nous allons montrer
que $A =\omega_L(B)$ convient, ce qui montrera du m\^eme coup la seconde
assertion.

Notons $\iota_2$ l'homomor\-phis\-me de $L$-alg\`ebres $L\to L\otimes_K L$
donn\'e par $\ell\mapsto 1\otimes \ell$. Dans
$\sA_{(L)}$, on a un isomor\-phis\-me
$L\otimes_K
\omega_L(B)\cong (L\otimes_K L)\otimes_{L,\iota_2} B$, et un mor\-phis\-me
$\mu\otimes_L 1_B:(L\otimes_K L)\otimes_{L,\iota_2} B \to
L\otimes_{L,\iota_2} B = B$ dont
$s\otimes_L 1_B$ est une section. Ceci montre bien que $B$ est facteur
direct de $L\otimes \omega_L(B)$.

2) Supposons maintenant que $L$ ne soit pas $L\otimes_K L$-projectif. Pour
obtenir un contre-exemple, prenons
$\sA= Vec_L$ (vue comme $K$-ca\-t\'e\-go\-rie a\-b\'e\-lien\-ne): $\sA_{(L)}$
est la  ca\-t\'e\-go\-rie des $(L\otimes_K L)$-modules
finis. Alors $L$, vu comme objet de
$\sA_{(L)}$, n'est pas facteur direct d'un objet de la forme $L\otimes_K
V$, $V \in Ob(\sA)$, puisqu'un tel objet est un
$(L\otimes_K L)$-module libre. \qed

\begin{sloppypar}
\begin{ex} Prenons pour $\sA$ la ca\-t\'e\-go\-rie $Rep_K G$ des
$K$-re\-pr\'e\-sen\-ta\-tions de dimension finie d'un $K$-groupe affine.
On note comme ci-dessus
$ Rep^\infty_K G$ la Ind-ca\-t\'e\-go\-rie attach\'ee \`a $\sA$, de sorte
qu'on a une \'equivalence naturelle (mo\-no\-\"{\i}\-da\-le)
$(Rep^\infty_K G)_{(L)}\cong Rep^\infty_L G_L$. En outre, si $L/K$ est
(finie) s\'eparable, elle induit une \'equivalence naturelle
(mo\-no\-\"{\i}\-dale) entre objets compacts
$(Rep_K G)_{(L)}\cong Rep_L G_L$, et aussi (en vertu du th\'eor\`eme
ci-dessus) $\cong \allowbreak(Rep_K G)_L^\natural$.
\end{ex}
\end{sloppypar}

\newpage

  \addtocontents{toc}{{\bf II. Radical et rigidit\'e}\hfill\thepage}
\
\bigskip
\begin{center}
\large\bf II. Radical et rigidit\'e
\end{center}
\bigskip

L'objectif de cette partie est l'\'etude du radical
$\sR$ en pr\'esence d'une structure mo\-no\-\"{\i}\-dale sur $\sA$. Il 
s'av\`ere que la
question de la compatibilit\'e du radical \`a la structure mo\-no\-\"{\i}\-dale
est fort d\'elicate, et n'a pas toujours une r\'eponse positive. Nous
comparons chemin faisant $\sR$ au plus grand id\'eal propre $\sN$
compatible \`a la structure mo\-no\-\"{\i}\-dale, et analysons le quotient 
$\sA/\sR$.

\medskip
Pour simplifier, on suppose $\sA$ \emph{stricte} (la contrainte
d'associativit\'e et les contraintes d'unit\'e sont l'identit\'e).

\section{G\'en\'eralit\'es}\label{s3.5}

\subsection{Id\'eaux mono\"{\i}daux, duaux.}\label{6.2} Soit $K$ un anneau
commutatif unitaire, et soit $(\sA,\bullet)$ une
$K$-ca\-t\'e\-go\-rie mo\-no\-\"{\i}\-dale.

\begin{defn}\label{d6.1} Un id\'eal $\sJ$ de $\sA$ est dit
\emph{mo\-no\-\"{\i}\-dal}\footnote{Jannsen \cite{jannsen2} dit plut\^ot
\emph{tensoriel}.} s'il est stable par les transformations
$1_C\bullet -$ et $- \bullet 1_C$ pour tout objet $C$ de $\sA$.
\end{defn}

\begin{lemme}\label{l6.1} a) Si $\sJ$ est mo\-no\-\"{\i}\-dal, il est stable par
produit mo\-no\-\"{\i}\-dal \`a gauche et \`a
droite par un mor\-phis\-me arbitraire. En particulier, on a une
structure mo\-no\-\"{\i}\-dale
induite sur le quotient $\sA/\sJ$.\\
b) La notion d'id\'eal mo\-no\-\"{\i}\-dal est stable par somme,
intersection, produit.
\end{lemme}

\prf a) Soient $f\in \sJ(A,B)$ et $g\in \sA(C,D)$. On a alors $g\bullet
f=(g\bullet 1_B)\circ
(1_C\bullet f)\in\sJ(C\bullet A,D\bullet B)$. On raisonne de m\^eme
pour les produits \`a
droite.

b) est imm\'ediat.
\qed

   Rappelons qu'un
objet $A$ d'une ca\-t\'e\-go\-rie mo\-no\-\"{\i}\-dale $\sA$ \emph{admet
un dual} (\`a droite) s'il existe un objet $A^\vee\in \sA$ et des mor\-phis\-mes
d'\'evaluation
\[\epsilon_A:A\bullet A^\vee\to \un\]
et de co\'evaluation
\[\eta_A:\un\to A^\vee\bullet A\]
tels que les compos\'es
\[\begin{CD}
A@> 1_A\bullet\eta_A>> A\bullet A^\vee\bullet A@> \epsilon_A\bullet 1_A>>
A\\
A^\vee@>\eta_A\bullet 1_{A^\vee}>> A^\vee\bullet A\bullet
A^\vee@>1_{A^\vee}\bullet\epsilon_A >> A^\vee
\end{CD}\]
soient \'egaux \`a l'identit\'e. Le triplet
$(A^\vee,\epsilon_A,\eta_A)$ est d\'e\-ter\-mi\-n\'e \`a isomor\-phis\-me
unique pr\`es.

Il en r\'esulte que le foncteur $- \bullet A^\vee$ est adjoint \`a 
gauche au foncteur
$-\bullet A$ et que le foncteur $ A^\vee \bullet -$ est adjoint \`a 
droite au foncteur
$A\bullet -$. Plus pr\'ecis\'ement, pour tous objets
$B,C$ de $\sA$, l'homomor\-phis\-me compos\'e
\[\begin{CD}
\sA(A\bullet B,C)@>1_{A^\vee}\bullet ->> \sA(A^\vee\bullet A\bullet B,
A^\vee\bullet C) @>(\eta_A\bullet 1_B)^*>>\sA(B,A^\vee\bullet C)
\end{CD}\]
est un isomor\-phis\-me, d'inverse le compos\'e
\[\begin{CD}
\sA(B,A^\vee\bullet C)@>1_{A}\bullet ->>
\sA(A\bullet B,A\bullet A^\vee\bullet C) @>(\epsilon_A\bullet
1_C)_*>>\sA(A\bullet B,C).
\end{CD}\]
et de m\^eme pour l'autre adjonction (\cf \eg \cite[1]{bru1}).

Si $A$ et $B$ ont pour duaux \`a droite $A^\vee$ et $B^\vee$, alors
$A\bullet B$ a pour dual \`a droite
\[(A\bullet B)^\vee =B^\vee\bullet A^\vee\]
avec
\begin{equation}\label{eq6.2.0} \varepsilon_{A\bullet B}=
\varepsilon_{A}\circ (1_A\bullet \varepsilon_B \bullet
1_{A^\vee}).\end{equation}

\medskip
On d\'efinit de mani\`ere duale la notion de dual \`a gauche ${}^\vee
A$. Si $A$ admet un dual
\`a droite et un dual \`a gauche, on a $({}^\vee A)^\vee={}^\vee(
A^\vee)=A$. En g\'en\'eral, ${}^\vee A$ et $A^\vee$ ne co\"{\i}ncident
pas.

\begin{defn}
On dit que
$(\sA,\bullet)$ est
\emph{rigide} si tout objet a un dual \`a droite et est un dual \`a
droite (ou de mani\`ere \'equi\-va\-len\-te, si tout objet a un dual \`a droite
et un dual \`a gauche).
\end{defn}

\begin{sorite}\label{sorrig} Si $\sA$ est rigide, il en est de m\^eme de
$\sA^\oplus,\sA^\natural$ et de $\sA^{\oplus\natural}$ (notations du \S
\ref{s1}). Il en est aussi de m\^eme du quotient de $\sA$ par tout 
id\'eal mo\-no\-\"{\i}\-dal. \qed
\end{sorite}

Soient $A$ et $B$ deux objets de ${\sA}$, $A$ ayant un dual \`a droite.
Par adjonction, on a un isomor\-phis\-me canonique de $K$-modules:
\begin{align}\label{eq6.2.1}\iota_{AB}:&{\sA}(\un , A^\vee\bullet B)\iso
{\sA}(A,B)\notag\\
\iota_{AB}(\phi)&=(\epsilon_A\bullet 1_B)\circ (1_A\bullet \phi),
\end{align}
d'inverse
  \begin{equation}\label{eq6.2.2}
\iota_{AB}^{-1}(f)= (1_{A^\vee} \bullet f)\circ  \eta_A.
\end{equation}

Si $B$ a aussi un dual \`a droite, la composition de deux mor\-phis\-mes $f\in
\sA(A,B)$ et $g\in \sA(B,C)$ peut se calculer comme suit (``composition
des correspondances"):
\begin{equation}\label{eq6.2.3}g\circ f=
\iota_{AC}\left((1_{A^\vee}\bullet
\epsilon_B \bullet 1_C)(\iota_{AB}^{-1} f\bullet \iota_{BC}^{-1}
g)\right). \end{equation}

Appliquant $\iota_{AC}^{-1}$, cette formule fondamentale exprime la
commutativit\'e du diagramme
\[\begin{CD}
\un @>{\eta_A}>>A^\vee A @>1_{A^\vee}\bullet f>> A^\vee B@= A^\vee B \\
@V{\eta_B}VV@V{1_{A^\vee A}\bullet \eta_B}VV@V{1_{A^\vee
B}\bullet\eta_B}VV\left|\right|\\
B^\vee B@>\eta_A\bullet 1_{B^\vee B}>> A^\vee AB^\vee
B@>{1_{A^\vee}\bullet\epsilon_B\bullet 1_B}>>A^\vee B B^\vee B
@>{1_{A^\vee}\bullet\epsilon_B\bullet 1_B}>>  A^\vee B\\
@V{1_{B^\vee}\bullet g}VV@V{1_{A^\vee AB^\vee}\bullet g}VV@V{1_{A^\vee
BB^\vee}\bullet g}VV @V{1_{A^\vee}\bullet g}VV\\
B^\vee C@>\eta_A\bullet 1_{B^\vee C}>> A^\vee AB^\vee C@>1_{A^\vee}\bullet
f\bullet 1_{B^\vee C}>>A^\vee B B^\vee C
@>{1_{A^\vee}\bullet\epsilon_B\bullet 1_C}>> A^\vee C.
\end{CD}\]

On a aussi les variantes
\begin{equation}\label{eq6.2.4} g\circ \iota_{AB}(\phi)=
\iota_{AC}((1_{A^\vee}\bullet g)\circ \phi)
\end{equation}
\begin{equation}\label{eq6.2.5}  \iota_{BC}(\psi) \circ f=
\iota_{AC}((\varepsilon_B\bullet 1_C)\circ (f\bullet 1_{B^\vee\bullet
C})\circ (1_A\bullet \psi ) )
\end{equation}
qui se lisent sur le m\^eme diagramme.

\medskip
Si $A, B$ ont des duaux \`a droite et que $f\in \sA(A,B)$, on d\'efinit
(\loccit) le transpos\'e
$f^t$ (ou dual \`a droite) de $f$ comme \'etant le compos\'e
\begin{equation}\label{eq6.2.6}
B^\vee =\un \bullet B^\vee \stackrel{\eta_A\bullet 1}{\To} A^\vee\bullet A
\bullet B^\vee\stackrel{1\bullet f\bullet 1}{\To} A^\vee \bullet B
\bullet B^\vee \stackrel{1\bullet \varepsilon_B }{\To} A^\vee \bullet \un
=A^\vee, \end{equation}
  ce qui s'\'ecrit aussi \begin{equation}\label{eq6.2.6,5} 
f^t=(1_{A^\vee}\bullet
\varepsilon_B )\circ (\iota_{AB}^{-1}(f) \bullet 1_{B^\vee}) , 
\end{equation}
  et on a
alors la formule
\begin{equation}\label{eq6.2.7}  (g\circ f)^t= f^t \circ
g^t.\end{equation}

On a de m\^eme une notion de transpos\'e ${}^tf$ \`a gauche, et
$({}^tf)^t={}^t(f^t)=f$. En g\'en\'eral, ${}^tf$ et $f^t$ ne
co\"{\i}ncident pas.

\medskip {\it On suppose d\'esormais que $\sA$ est rigide.}

\begin{lemme}\label{L4.1} Tout id\'eal mo\-no\-\"{\i}\-dal $\sJ$ v\'erifie
\[{\sJ} (A,B)=\iota_{AB}\big({\sJ}(\un , A^\vee\bullet B)\big) \]
et
\[\big({\sJ}(A,B)\big)^t = {\sJ}(B^\vee ,A^\vee).  \]
          \end{lemme}

\prf La premi\`ere formule d\'ecoule des formules \eqref{eq6.2.1},
\eqref{eq6.2.2}. Pour la seconde,
\begin{multline*}
\big({\sJ}(A,B)\big)^t = (1_{A^\vee}\bullet \varepsilon_B )\circ ( 
{\sJ}(\un , A^\vee\bullet
B)  \bullet 1_{B^\vee}) \\
\subset(1_{A^\vee}\bullet \varepsilon_B )\circ   {\sJ}( {B^\vee} , 
A^\vee\bullet
B\bullet  {B^\vee})  \subset {\sJ}(B^\vee ,A^\vee) .
\end{multline*}

\begin{sloppypar}
L'analogue pour les duaux \`a gauche de ces deux formules conduit \`a
${}^t\big({\sJ}(C,D)\big) = {\sJ}({}^\vee D ,{}^\vee C)$. En posant
$B={}^\vee C, \,A= {}^\vee D$ (ce qui est loisible puisque $\sA$ est
rigide), et en appliquant ${}^t$, on obtient finalement l'inclusion
oppos\'ee $ {\sJ}(B^\vee ,A^\vee)\subset \big({\sJ}(A,B)\big)^t$.
\qed
\end{sloppypar}

On a une r\'eciproque partielle du lemme \ref{L4.1} (mais voir aussi le
lemme \ref{L4.2} dans le cas sy\-m\'e\-tri\-que):

\begin{lemme}\label{L4.2bemol} Soit $\sI$  un  $\un$-id\'eal \`a gauche de
$\sA$ (d\'efinition \ref{d1}). La formule
\[{\sI}^{(\bullet)}(A,B)=\iota_{AB}\big({\sI}(A^\vee\bullet B)\big)\]
d\'efinit un id\'eal de $\sA$ contenant  $\sI$. Cet id\'eal est
mo\-no\-\"{\i}\-dal \`a gauche (stable par produit $\bullet$ \`a gauche).
\end{lemme}

\prf Soit $g\in
{\sA}(B,C)$ un mor\-phis\-me. En utilisant la formule \eqref{eq6.2.4}, on
calcule
\begin{multline*}
g\circ {\sI}^{(\bullet)}(A,B) = g\circ  \iota_{AB}\big({\sI}(A^\vee\bullet
B)\big) = \iota_{AC}\big((1_{A^\vee}\bullet g)\circ
{\sI}(A^\vee\bullet B)\big) \\
\subset \iota_{AC}\big(  {\sI}(A^\vee\bullet
C)\big) = {\sI}^{(\bullet)}(A,C).
\end{multline*}

La stabilit\'e par composition \`a droite se d\'emontre de
m\^eme en utilisant \eqref{eq6.2.5}.

Enfin, soit $C$ un objet de ${\sA}$. Observons que
l'homomor\-phis\-me
\[1_C\bullet - :  {\sA}(A,B)\to {\sA}(C\bullet A,C\bullet B)\]
n'est autre que
\[f\mapsto \iota_{C\bullet A, C\bullet B}[ (1_{A^\vee}\bullet \eta_C
\bullet 1_B)\circ \iota_{AB}^{-1}(f)]\]
o\`u $\eta_C: \un \to C^\vee\bullet C$ est la co\'evaluation. On a donc
\begin{multline*}
1_C\bullet {\sI}^{(\bullet)}(A,B) = \iota_{C\bullet A, C\bullet B}[
(1_{A^\vee}\bullet \eta_C \bullet 1_B)\circ
          \big({\sI}(A^\vee\bullet  B)\big)\\
\subset \iota_{C\bullet A, C\bullet B}  \big({\sI}(A^\vee\bullet
C^\vee\bullet C\bullet B)\big) = {\sI}^{(\bullet)}(C\bullet A, C\bullet B).
\end{multline*}

On conclut en appliquant (la moiti\'e du) lemme \ref{l6.1}.\qed

\subsection{Le cas sy\-m\'e\-tri\-que} On suppose d\'esormais que $\sA$ (qui
rappelons-le, est suppos\'ee rigide) est munie d'un tressage sy\-m\'e\-tri\-que
$R$ (une contrainte de commutativit\'e, dans le langage de Saavedra). On
renvoie \`a \ref{Tre} pour plus de d\'etails sur les tressages.

Sous cette hypoth\`ese, ${}^\vee A= A^\vee$ (duaux \`a droite et \`a
gauche co\"{\i}ncident), et on a
\begin{equation}\label{eq6.2.7,3} \epsilon_{A^\vee}=\epsilon_A\circ
R_{A^\vee, A},\eta_{A^\vee}=R_{A^\vee, A}\circ
\eta_A , A^{\vee\vee}=A ,\end{equation}
\begin{equation}\label{eq6.2.7,6}  {}^tf= f^t, \;\;\;\;(f^t)^t=f .
\end{equation}

On a aussi l'identit\'e
\begin{equation}\label{eq6.2.11}
\iota_{C^\vee,C^\vee}(R_{C^\vee,C}\circ \phi)=\iota_{C,C}(\phi)^t
\end{equation}
(en effet, par \eqref{eq6.2.6,5} et \eqref{eq6.2.2}, l'image par
$\iota_{C^\vee,C^\vee}^{-1}$ du membre de droite s'\'ecrit
$(1_{C^\vee}\bullet \varepsilon_C)\circ (\phi\bullet 1_{C^\vee})$, tandis 
que
celle du membre de gauche s'\'ecrit $(\varepsilon_{C^\vee} \bullet 
1_{C^\vee})\circ (1_{C^\vee}\bullet (R_{C^\vee,C}\circ
\phi))= ((\varepsilon_{C } \circ R_{C^\vee,C})\bullet 
1_{C^\vee})\circ (1_{C^\vee}\bullet  \phi)$, et l'\'egalit\'e des deux
membres s'ensuit par permutation des facteurs $C^\vee$.)

\medskip Par ailleurs, si $A,B,C,D$ sont quatre objets de $\sA$ et $\phi\in
\sA(\un,A^\vee\bullet B)$, $\psi\in \sA(\un, C^\vee\bullet D)$, on a la
formule
\begin{equation}\label{eq6.2.8}\iota_{A,B}(\phi)\bullet\iota_{C,D}(\psi)=
\iota_{A\bullet C,B\bullet D}((R_{A^\vee\bullet B,C^\vee} \bullet
1_D)\circ(\phi\bullet\psi)).\end{equation}

\medskip Le lemme suivant compl\`ete \ref{L4.2bemol}.

\begin{lemme}\label{L4.2} Soit $\sI$  un  $\un$-id\'eal \`a gauche de
$\sA$ (d\'efinition \ref{d1}). La formule
\[{\sI}^{(\bullet)}(A,B)=\iota_{AB}\big({\sI}(A^\vee\bullet B)\big)\]
d\'efinit un id\'eal mo\-no\-\"{\i}\-dal ${\sI}^{(\bullet)}$; c'est le plus
petit id\'eal mo\-no\-\"{\i}\-dal
de $\sA$ contenant 
$\sI$.
\end{lemme}

\prf On a vu que ${\sI}^{(\bullet)}$ est un 
id\'eal contenant
$\sI$, stable par produit \`a gauche.  Par 
tressage, on a aussi
${\sI}^{(\bullet)}(A,B) \bullet 1_C \subset 
{\sI}^{(\bullet)}(A\bullet C,
B\bullet C)$ pour trois objets $A,B,C$. 
La deuxi\`eme assertion est
imm\'ediate  compte tenu du lemme 
\ref{L4.1}.   \qed

\subsection{Interpr\'etation} Pour expliciter le 
sens des lemmes
\ref{L4.1} et
\ref{L4.2}, 
notons:

\begin{itemize}
\item $\bG$ l'ensemble des $\un$-id\'eaux 
\`a gauche de $\sA$;
\item $\bI$ l'ensemble des id\'eaux 
[bilat\`eres] de $\sA$;
\item $\bT$ l'ensemble des id\'eaux 
[bilat\`eres] mono\"{\i}daux de
$\sA$.
\end{itemize}

Ces trois 
ensembles sont ordonn\'es par inclusion. On a le diagramme
sui\-vant 
d'applications
croissantes:
\[\begin{CD}
\bT&& 
\begin{smallmatrix}i\\\longrightarrow\\\longleftarrow\\j
\end{smallmatrix}&&\bI\\
&\scriptstyle t\displaystyle\nwarrow&& 
\scriptstyle
r\displaystyle\swarrow\\ &&\bG
\end{CD}\]
o\`u

\begin{sloppypar}
\begin{itemize}
\item $i$ associe \`a un id\' 
eal mo\-no\-\"{\i}\-dal l'id\'eal sous-jacent;
\item $j$ associe \`a un 
id\'eal le plus petit id\'eal mo\-no\-\"{\i}\-dal le
contenant;
\item $r$ 
associe \`a un id\'eal le $\un$-id\'eal \`a gauche 
sous-jacent
($r(\sI)(A)=\sI(\un,A)$);
\item $t$ associe \`a un 
$\un$-id\'eal \`a gauche $\sI$ l'id\'eal
mo\-no\-\"{\i}\-dal 

$\sI^{(\bullet)}$ construit dans le lemme 
\ref{L4.2}.
\end{itemize}
\end{sloppypar}

On a \'evidemment $ji=Id$. 
Le lemme \ref{L4.1} dit que $tri=Id$,
tandis que le lemme
\ref{L4.2} 
dit que $rit=Id$. En particulier, \emph{les applications
$t$ et $ri$ 
sont des
bijections inverses l'une de l'autre}. D'autre part, soit 
$\sI\in
\bI$. Comme $j(\sI)$
contient $r(\sI)$, on a $j(\sI)\supset 
tr(\sI)$, et $j(\sI)=tr(\sI)$
si et seulement si
$tr(\sI)\supset 
\sI$. 

Par abus de notation, on \'ecrira encore $\sI^{(\bullet)}$ 
pour $tr(\sI)$
si $\sI\in \bI$.

\begin{sloppypar}
\medskip Par 
ailleurs, les ensembles $\,\bG,\,\bI,\,\bT$ sont munis
d'op\'erations 
internes ``som\-me", ``intersection". Il est clair que 
les
applications $i,t,r$ respectent ces 
op\'erations.
\end{sloppypar}

Le cas du produit est plus 
int\'eressant. {\it A priori}, seuls $\bI$ et
$\bT$ sont munis d'une 
op\'eration ``produit" \'evidente. On peut alors
munir $ \bG $ d'un 
tel produit en posant: $\sI\cdot\sJ :=
ri(t(\sI)\cdot t(\sJ))$ (en 
formule: 
\begin{equation}\label{eq15}
(\sI\cdot\sJ)(B)   = \sum_A 
\iota_{AB}(\sI(A^\vee \bullet
B))\circ \sJ(A) \;\; 
).
\end{equation}      

Les applications $i,t,r$ respectent alors aussi le 
produit. On a aussi la formule 
\begin{equation}\label{eq15sym}
(\sI\cdot\sJ)(B)   = \sum_{u: 
A\bullet A'\to B} 
 Im(\sJ(A)  \bullet \sI(A') \to \sA(\un, B)) 
\end{equation}
qui montre que le produit des id\'eaux mono\"{\i}daux (dans le cas rigide)
est {\it sy\-m\'e\-tri\-que}. 
 
Pour $A'=A^\vee\bullet B$ et 
$u=
\varepsilon_A
\bullet 1_B$, le terme du second membre 
de
(\ref{eq15sym}) s'\'ecrit en effet $  (\varepsilon_A \bullet 
1_B)\circ (\sJ(A)\bullet \sI(A^\vee \bullet
B))$, qui co\"{\i}ncide 
d'apr\`es (\ref{eq6.2.5}) avec le terme courant du second membre de 
(\ref{eq15}) (donc
avec
$(\sI\cdot\sJ)(B) $). 
R\'eciproquement,
$\sJ(A) 
\bullet \sI(A')\subset  (t(\sJ) \cdot 
t(\sI))(A\bullet A') $ puisque $t(\sJ)$ et $t(\sI)$ sont des 
id\'eaux
mono\"{\i}daux, d'o\`u il d\'ecoule que le second membre de 
(\ref{eq15sym}) est contenu dans $(\sI\cdot\sJ)(B)
$.

\begin{ex} 
\label{ex2} Si $\sA$ est la ca\-t\'e\-go\-rie mo\-no\-\"{\i}\-dale des
motifs de 
Gro\-then\-dieck pour l'\'equivalence rationnelle (motifs `de
Chow') 
\`a coefficients dans
$K$ et si $K$ est un corps, il y a bijection 
entre les id\'eaux
mono\"{\i}daux de $\sA$ et les relations 
d'\'equivalence ad\'equates sur
les cycles alg\'ebriques \`a 
coefficients dans $K$.

Il est plus facile de d\'ecrire la bijection 
entre $\un$-id\'eaux \`a
gauche de $\sA$ et relations d'\'equivalence 
ad\'equates (\cf
\cite[4]{jannsen2}). Soit $\sI$ un tel 
$\un$-id\'eal, c'est-\`a-dire la
donn\'ee d'une famille de 
sous-groupes $\sI(M)$ de $\sA(\un,M)$ stable
par l'action \`a gauche 
des correspondances, o\`u $M$ d\'ecrit les objets
de $\sA$. Pour le 
motif $M$ d'une vari\'et\'e $X$ tordu $m$ fois
\`a la Tate, on a 
$\sA(\un,M)=CH^m(X)_K$, donc $\sI$ d\'efinit une
relation 
d'\'equivalence ad\'equate. Inversement, toute 
relation
d'\'equivalence ad\'equate d\'efinit un $\un$-id\'eal \`a 
gauche de la
ca\-t\'e\-go\-rie des correspondances de Chow (non 
effectives),
 et on v\'erifie que les $\un$-id\'eaux \`a gauche de 
cette ca\-t\'e\-go\-rie
sont en correspondance bijective avec ceux de son 
enveloppe
pseudo-ab\'elienne $\sA$. 

On d\'eduit de (\ref{eq15}) la 
formule pour le produit de relations
d'\'equivalence a\-d\'e\-qua\-tes 
(correspondant au produit d'id\'eaux, et
not\'e $\ast$ dans \loccit). 

\end{ex}

\section{Traces}

\medskip
{\it Dans tout ce paragraphe, 
on suppose que $(\sA,\bullet)$ est $K$-lin\'eaire mo\-no\-\"{\i}\-dale 
rigide,
sy\-m\'e\-tri\-que, et que
$End(\un) =K$} (en revanche, on ne 
suppose pas $\sA$ ab\'elienne).

\subsection{Traces et id\'eal 
$\sN$}\label{traces}  
Rappelons que la {\it
trace} d'un 
endomor\-phis\-me
$h\in
\sA(C,C)$ est
l'\'e\-l\'e\-ment $tr(h) \in K= End(\un) $ d\'efini par
\[\epsilon_C \circ R_{C^\vee, C} \circ({\iota_{CC}}^{-1}(h) ): {\bf
1 }\to \un ,\]
o\`u $\epsilon_C:   C\bullet C^\vee \to  \un $ est l'\'evaluation. Dans
la situation sy\-m\'e\-tri\-que dans laquelle nous nous pla\c cons, la formule
s'\'ecrit aussi
\begin{equation}\label{eq6.2.8,5} tr(h)= \epsilon_{C^\vee}
\circ({\iota_{CC}}^{-1}(h) ) \end{equation}
  La {\it dimension} (ou rang) de l'objet
$A$ est la trace de $1_A$.

   Soient $f\in \sA(A,B)$ et $g\in \sA(B,A)$. Des formules \eqref{eq6.2.0}
et \eqref{eq6.2.3}, on d\'eduit que
\[tr(g\circ f)=\epsilon_{A^\vee}(1_{A^\vee}\bullet \epsilon_B\bullet
1_A)(\iota_{AB}^{-1}f\bullet \iota_{BA}^{-1}g)=\epsilon_{A^\vee B}\circ
(\iota_{AB}^{-1}f\bullet \iota_{BA}^{-1}g)\]
d'o\`u
\begin{equation}\label{eq6.2.9}tr(g\circ f)= tr(f\circ g)\end{equation}
et, en appliquant \eqref{eq6.2.6}
\begin{equation}\label{eq6.4.1}
tr(g\circ f)=(\iota_{AB}^{-1}f)^t\circ
\iota_{BA}^{-1}g \;\; (=tr((\iota_{AB}^{-1}f)^t\circ \iota_{BA}^{-1}g)).
\end{equation}

En appliquant \eqref{eq6.2.11}, on obtient aussi l'identit\'e
\[tr(h^t)=tr(h).\]

Voici une autre formule utile (\cf
\cite[7.2]{de}): si $f$  et $g$ sont des endomor\-phis\-mes de $A$ et $B$, on
a
\begin{equation}\label{eq6.2.10}tr(f\bullet g)= tr(f) tr(g).
\end{equation}

Introduisons \`a pr\'esent l'id\'eal $\sN$, protagoniste principal de ce
paragraphe.

\begin{lemme} La formule
\[{\sN}(A,B)=\{f \in \sA(A,B), \forall g\in \sA(B,A), tr(g\circ f)=0\} \]
d\'efinit un id\'eal mo\-no\-\"{\i}\-dal de $\sA$ \footnote{Dans \cite{bru2}
cet id\'eal appara\^{\i}t, dans le cadre plus g\'en\'eral (tressage non
n\'e\-ces\-sai\-re\-ment sy\-m\'e\-tri\-que)
des \emph{tortils} ou \emph{ca\-t\'e\-go\-ries enrubann\'ees}, sous le nom
d'id\'eal des mor\-phis\-mes n\'egligeables.}. On a
\[{\sN}(A,B)=\iota_{AB}\{f \in \sA(\un ,A^\vee\bullet B),
\forall g\in\sA(A^\vee\bullet B,\un ),  g\circ f =0\} . \]
\end{lemme}

\prf La stabilit\'e par composition \`a gauche est claire et la
stabilit\'e par composition \`a
droite s'en d\'eduit en vertu de la formule \eqref{eq6.2.9}. Il 
suffit
alors de prouver la seconde assertion, qui dit que 
$\sN=\sN^{(\bullet)}$.
C'est imm\'ediat \`a partir de la formule 
\eqref{eq6.4.1}.  \qed

\begin{ex}\label{ex1} Dans le cas de 
l'exemple \ref{ex2} (motifs),
l'id\'eal mo\-no\-\"{\i}\-dal
$\sN$ 
correspond \`a l'\'equivalence num\'erique des cycles 
alg\'ebriques.
\end{ex}

\begin{sloppypar}
\begin{lemme}\label{Ldim} On suppose que 
$K$ est un corps. Soit
$f_1,\dots,f_n
\in
\sA(A,A)$ une famille de 
$n$ \'el\'ements li\-n\'e\-ai\-re\-ment
ind\'ependants 
modulo
$\sN(A,A)$. Alors $(f_1,\dots,f_n)$ correspond \`a une 
cor\'etraction
${\un}^n\inj A^\vee\bullet A$.   
\end{lemme} 
\end{sloppypar}

\prf 
La donn\'ee des $f_i$ fournit un mor\-phis\-me
${\un}^n\rightarrow 
A^\vee\bullet A$ (de $i$-\`eme composante
$\iota_{AA}^{-1}(f_i)$). 
Remarquons
que, par d\'efinition de $\sN$, la trace induit une forme 
bilin\'eaire
(sy\-m\'e\-tri\-que) non 
d\'eg\'en\'er\'ee
\[\sA(A,A)/\sN(A,A)\times \sA(A,A)/\sN(A,A)\to 
K.\]

Il exi\-ste donc des \'el\'ements $g_1,\dots, g_n$ 
de
$\sA(A,A)$ tels que $tr(g_i\circ f_j) = \delta_{ij}$ (symbole 
de
Kronecker). Le mor\-phis\-me compos\'e
\[\begin{CD}
A^\vee\bullet A 
@>\circ g_1,\dots ,\circ
g_n>>  (A^\vee\bullet A)^n
@>R_{A^\vee 
A},\dots, R_{A^\vee A}>>
(A\bullet A^\vee )^n  @>{ev,\dots 
,ev}>>{\un}^n
\end{CD}\]
est un inverse \`a gauche  de 
${\un}^n\inj
A^\vee\bullet A$. R\'eciproquement, toute cor\'etraction 
comme dans
l'\'enonc\'e s'obtient clairement 
ainsi.
\qed

\begin{prop}\label{new p1} On suppose que $K$ est un 
corps. Soit $\sA$
un ca\-t\'e\-go\-rie $K$-lin\'eaire mo\-no\-\"{\i}\-dale 
sy\-m\'e\-tri\-que rigide, avec
$K=End(\un )$. On note $\sR$ son radical. 
Alors\\
a) $\sN=\sR^{(\bullet)}$: $\sN(A,B)= \iota_{AB}(\sR(\un, 
A^\vee \bullet
B))$ pour tout couple $(A,B)$. \\
b) $\sN$ est le plus 
grand id\'eal mo\-no\-\"{\i}\-dal
distinct de $\sA$.\\
c) Soit $\sJ$ un 
id\'eal mo\-no\-\"{\i}\-dal distinct de $\sA$. Si
$\sA/\sJ$ est 
semi-simple, alors $\sJ=\sN$. 
\end{prop}

\prf $a)$ Que $\sN$ soit 
distinct de $\sA$ se voit sur le couple
$(\un ,\un)$. Pour 
l'\'egalit\'e $\sN=\sR^{(\bullet)}$, il s'agit de
voir que pour tout 
objet
$A$,
$\sN(\un ,A)=\sR(\un ,A)$. C'est bien le cas, car 
d'apr\`es le lemme \ref{L2},
\[{\sR}(\un ,A)= \{f \in \sA(\un ,A 
),
\forall g\in
\sA(A ,\un ),
          g\circ f =0\} . \]
$b)$ Prouvons que tout id\'eal mo\-no\-\"{\i}\-dal $\sJ$ distinct de $\sA$ est
contenu dans $\sN$: il suffit de faire voir
que pour tout objet $A$, $\sJ(\un ,A)\subset \sN(\un ,A)$,
c'est-\`a-dire que pour tout $f\in \sJ(\un ,A)$ et tout $g\in
\sA(A,\un )$, $g\circ f=0 \in K$. Si tel
n'\'etait pas le cas, on pourrait supposer
$g\circ f=1$, donc $1\in \sJ(\un ,\un )$, d'o\`u $\sJ(\un ,A)=\sA(\un
,A)$. Par ailleurs, $\sN(\un,\un)=0$, donc $\sN\neq \sA$.

$c)$ Si $\sA/\sJ$ est semi-simple, alors $\sR\subset \sJ$. Donc
$\sR^{(\bullet)}=\sN$ est contenu dans $\sJ^{(\bullet)}=\sJ$. On a
l'inclusion oppos\'ee
d'apr\`es $i)$.\qed

\begin{cor}  L'id\'eal $\sN$ ne d\'epend pas du choix du
tressage (sy\-m\'e\-tri\-que) $R$.\qed
\end{cor}

\begin{prop}\label{new p3/2} Sous les m\^emes hypoth\`eses, les 
conditions suivantes sont
\'equi\-va\-len\-tes:
\begin{thlist}
  \item $\sR \supset \sN $,
\item  le foncteur de projection $\sA\to \sA/\sN$ est conservatif,
\item  pour tous objets $A,B$ de $\sA$, l'identification
\[\iota_{AB}:\;\sA({\bf 1}, A^{\vee}\bullet B)\iso \sA(A,B)\]
induit une inclusion
\[\sR(A,B)\supset \sR({\bf 1},A^{\vee}\bullet B).\]
  \end{thlist}
\noindent En outre,  $\sR$ est un id\'eal
mo\-no\-\"{\i}\-dal si et seulement si $\;\sR =\sN $.
\end{prop}

\prf (i) $\iff$ (ii)  r\'esulte de \ref{P1} b); (i) $\iff$ (iii) n'est
autre qu'une reformulation de la proposition \ref{new p1} a). La
derni\`ere assertion suit de ce que $\sN=\sR^{(\bullet)}$ (ou de (i)
$\iff$ (iii)).\qed

\begin{cor}\label{c2} Si $\sA$ est semi-simple, $\sN=0$.
     Si $\;\sR\ne \sN$, $\sR\ne 0$.
\end{cor}

En effet, la premi\`ere (\resp deuxi\`eme) affirmation
r\'esulte 
de la proposition \ref{new p1}  (\resp  \ref{new p1} et
\ref{new p3/2}) 
compte tenu de ce que le radical d'une ca\-t\'e\-go\-rie
semi-simple est nul 
(\resp de ce que
$0$ est un id\'eal mo\-no\-\"{\i}\-dal).\qed

\begin{ex}\label{BBM} Revenons \`a l'exemple de la ca\-t\'e\-go\-rie
mo\-no\-\"{\i}\-dale sym\'etri\-que rigide $\sA$ des motifs de Chow
\`a 
coefficients rationnels. Des conjectures tr\`es fortes dues \`a A.
Beilinson 
\cite{beilinson} et de J. Murre \cite{murre} (qui sont en fait
\'equi\-va\-len\-tes, \cf \cite{jannsen2}) pr\'edisent l'existence d'un
filtration 
finie fonctorielle sur les objets de $\sA$. Chaque cran de la
filtration est 
suppos\'e provenir d'une \'equivalence ad\'equate, donc
d\'efinit un id\'eal 
mo\-no\-\"{\i}\-dal $F^p$. D'apr\`es  \cite{jannsen2}, il
r\'esulterait de ces 
conjectures  et des conjectures standard que $F^p$
est la puissance
$p$-i\`eme de l'id\'eal $F^1$ (correspondant \`a l'\'equivalence
homologique) et que
$F^1=\sN$. La finitude conjecturale de la filtration a 
pour
cons\'equence que $F^1 \subset \sR$, donc $F^1= \sR = \sN$ ($\sR$ est 
mo\-no\-\"{\i}\-dal),
et implique alors:
$\sA$ est strictement de Wedderburn (noter que
$\sR^\omega=
\cap F^p= 0$). Ces propri\'et\'es se transf\'erent aux ca\-t\'e\-go\-ries de
motifs  purs pour une \'equivalence
ad\'equate quelconque.

Cet exemple conjectural est l'une des principales sources d'inspiration
de ce travail.
  \end{ex}

\begin{rem}\label{r6.1} La compatibilit\'e de $\sN$  vis-\`a-vis de
     l'extension des sca\-lai\-res (na\"{\i}ve) ne pose pas de probl\`eme
(contrairement
\`a celle du radical $\sR$, \cf \S 3.) C'est prouv\'e
dans \cite[1.4.1]{bru2}, dans un contexte beaucoup plus large.
\end{rem}

\subsection{Traces et puissances ext\'erieures}\label{tr et pext}

Cette partie est inspir\'ee par \cite[7.2]{de} et par le r\'ecent
travail de Kimura (\cite{ki}, voir aussi \cite{gu-pe}). Ici la sym\'etrie
du tressage est essentielle.

Soit $n$ un entier $>0$. On fait op\'erer le groupe sy\-m\'e\-tri\-que
$\mathfrak{S}_n$ sur
$A^{\bullet n}$ au moyen du tressage [sy\-m\'e\-tri\-que!] $R$, en \'ecrivant
tout \'el\'ement $\sigma$ de $\mathfrak{S}_n$ comme produit de
transpositions de la forme $(i\;i+1)$ (le r\'esultat ne d\'epend pas de
cette \'ecriture,
\cf
\cite[7.2]{de}). Plus g\'en\'eralement, si $A_1,\dots,A_n\in \sA$ et
$\beta\in \mathfrak{S}_n$, on a un isomor\-phis\-me bien d\'etermin\'e
\[\beta:A_1\bullet\dots\bullet A_n\to A_{\beta(1)}\bullet \dots\bullet
A_{\beta(n)}\]
covariant en $\beta$. Si $B_1,\dots,B_n$ sont d'autres objets, pour toute
famille $f_i\in
\sA(A_i,B_i)$, on a
\begin{equation}\label{eq6.5.1} \beta \circ
(f_1\bullet\dots\bullet f_n)=(f_{\beta(1)}\bullet\dots \bullet
f_{\beta(n)})\circ \beta.\end{equation}

On a aussi les formules suivantes qui nous seront utiles
\begin{equation}\label{eq6.5.1,4/3} \iota_{A_{\beta(1)}\bullet \ldots
\bullet A_{\beta(n)},A_{\beta(1)}\bullet
\ldots \bullet A_{\beta(n)}}((\beta^{-1},\beta)\circ\phi)=\beta\circ
   \iota_{A_1\bullet \ldots \bullet A_n,A_1\bullet \ldots \bullet
A_n}(\phi) \circ \beta^{-1} 
\end{equation}
\begin{equation}\label{eq6.5.1,5/3} 
\varepsilon_{A_{\beta(n)}^\vee\bullet \ldots \bullet 
A_{\beta(1)}^\vee}\circ(\beta^{-1},\beta)
=\varepsilon_{A_n^\vee\bullet
\ldots
\bullet A_1^\vee },
\end{equation}
o\`u 
$(\beta^{-1},\beta)$ d\'esigne la ``permutation" 
\[A_n^\vee\bullet 
\ldots \bullet A_1^\vee \bullet A_1\bullet \ldots
\bullet A_n  \to 
A_{\beta(n)}^\vee\bullet \ldots \bullet 
A_{\beta(1)}^\vee \bullet 
A_{\beta(1)}\bullet \ldots \bullet 
A_{\beta(n)}.\]

{\it Cas 
particulier de \eqref{eq6.5.1}}: soit $\sigma\in 
\mathfrak{S}_n$, et 
supposons que
$B_i=A_{\sigma(i)}$ pour tout $i$. On peut alors 
consid\'erer
l'endomor\-phis\-me suivant de $A_1\bullet\dots\bullet A_n$ 
(\cf
\cite[7.2]{de}):
\[F=\sigma^{-1}\circ (f_1\bullet\dots\bullet 
f_n).\]

Nous nous proposons de calculer la \emph{trace de 
$F$}.

\'Ecrivons $\sigma$ comme produit de cycles 
disjoints
$\sigma_1\circ\dots\circ \sigma_p$, et soit $I_k\subset 
\{1,\dots,n\}$ le
support de $\sigma_k$. Pour tout $i\in I_k$, on 
peut consid\'erer
l'endomor\-phis\-me de 
$A_i$
\[F_{k,i}=f_{\sigma^{\ell_k-1}(i)}\circ\dots\circ 
f_{\sigma(i)}\circ f_i\]
o\`u $\ell_k=|I_k|$ est la longueur du cycle 
$\sigma_k$. D'apr\`es
\eqref{eq6.2.9}, la trace $tr(F_{k,i})$ ne 
d\'epend pas du choix de $i$:
notons-l\`a 
$t_k$.

\begin{prop}\label{p6.4.1} On a
\[tr(F)=\prod_{k=1}^p 
t_k.\]
\end{prop}

\prf Soit $\beta\in \mathfrak{S}_n$. D'apr\`es 
\eqref{eq6.5.1} et
\eqref{eq6.2.9}, on 
a
\begin{multline*}
tr(\beta\sigma^{-1}
\beta^{-1}\circ(f_{\beta(1)}\bullet\dots\bullet 
f_{\beta(n)}))\\
=tr(\beta\circ\sigma^{-1}\circ(f_1\bullet\dots\bullet 
f_n)\circ
\beta^{-1})\\
=tr(\sigma^{-1}\circ(f_1\bullet\dots\bullet 
f_n)).
\end{multline*}

Cette formule permet de se ramener au cas 
o\`u les $I_k$ forment une
partition croissante de 
$\{1,\dots,n\}$:
\[I_k=\{\ell_1+\dots+\ell_{k-1}+1,\dots,\ell_1+\dots+ 
\ell_k\}\]
et o\`u, de plus, pour tout $k$, $\sigma_k$ est le cycle 
``croissant"
\[\sigma_k=(\ell_1+\dots+\ell_{k-1}+1,\cdots 
,\ell_1+\dots+\ell_k).\]

Posant 
alors
\[F_k=\sigma_k^{-1}\circ
(f_{\ell_1+\dots+\ell_{k-1}+1}\bullet\d 
ots\bullet f_{\ell_1+\dots+\ell_k})\]
une it\'eration de la formule 
\eqref{eq6.2.10} donne
\[tr(F)=\prod_{k=1}^p tr(F_k).\]

Ceci 
ram\`ene la d\'emonstration de la proposition \ref{p6.4.1} au
cas 
o\`u $\sigma$ est la permutation cyclique de $\{1,\dots,n\}$. 
Par
it\'eration de la formule \eqref{eq6.2.3}, et compte tenu 
de
\eqref{eq6.2.0} et \eqref{eq6.2.8}, on 
trouve
\begin{multline}\label{eq6.5.3}
\iota^{-1}_{A_1A_1}(f_n\circ 
\dots \circ
f_1)\\
= (1_{A_1^\vee}\bullet\epsilon_{A_2}\dots 
\bullet\epsilon_{A_{n}}
\bullet 1_{A_1})\circ (\iota_{A_1A_2}^{-1} 
f_1\bullet \dots \bullet
\iota_{A_nA_1}^{-1} f_n)  \\
= 
(1_{A_1^\vee}\bullet\epsilon_{A^\vee_n\bullet\dots \bullet 
A^\vee_2}
\bullet 1_{A_1})\circ (\iota_{A_2\bullet\dots \bullet A_n 
\bullet
A_1,A_2\bullet\dots \bullet A_n \bullet A_1}^{-1} F 
).
\end{multline}

En appliquant
$\epsilon_{A^\vee}$ ``aux deux 
facteurs extr\^emes $A^\vee$ et $A$" (ce
qui fait sens gr\^ace au 
tressage), on obtient
\begin{multline}\label{eq6.5.4}
tr(f_n\circ \dots \circ f_1)=
\epsilon_{A^\vee_1\bullet A^\vee_n\bullet\dots \bullet A^\vee_2}
(\iota_{A_2\bullet\dots \bullet A_n \bullet A_1,A_2\bullet\dots \bullet
A_n \bullet A_1}^{-1} F )\\
=\epsilon_{A^\vee_1\bullet
A^\vee_n\bullet\dots\bullet A^\vee_2}  (\iota_{A_2\bullet\dots \bullet
A_n \bullet A_1,A_2\bullet\dots \bullet A_n \bullet A_1}^{-1} F )\\
=\epsilon_{A^\vee_n\bullet\dots \bullet A^\vee_1}
(\iota_{A_1\bullet\dots \bullet A_n  ,A_1\bullet\dots \bullet A_n}^{-1}
(\sigma \circ  (f_1\bullet\dots\bullet f_n)))
= tr(\sigma\circ(f_1\bullet\dots\bullet f_n))
\end{multline}
(\cf aussi \cite[7.2]{de}). \qed

\begin{cor}\label{c6.5.1} Supposons $A_1=\dots=A_n=A$ et
$f_1=\dots=f_n=f$. Alors
\[tr(F)=\prod_{k=1}^p tr(f^{\circ \ell_k}).\]
En particulier, pour $f=1$, on a
$tr(\sigma)=\dim(A)^p$.\qed
\end{cor}

{\it Dans la suite de ce paragraphe (jusqu'\`a \ref{tracesnil}), on
suppose que $\sA$ est pseudo-ab\'elienne.}

\begin{defn}\label{alt, sym} Soit $n$ un entier naturel tel que $n!$
soit inversible dans
$K$, et soit $A\in \sA$. La \emph{puissance ext\'erieure $n$-i\`eme}
$\Lambda^n A$ de
$A$ est l'image du projecteur d'anti\-sym\'e\-tri\-sa\-tion sur
$A^{\bullet n}$
\begin{equation}\label{eq6.5.2} a_n=\frac{1}{n!}\sum_{\sigma \in \mathfrak
S_n}\,sgn(\sigma)
\sigma.\end{equation}
  De m\^eme, la puissance sy\-m\'e\-tri\-que  $n$-i\`eme  ${\bf S}^n
A$ d'un objet
$A$ de $\sA$ est l'image du projecteur de sym\'e\-trisation
    $s_n=\frac{1}{n!}\sum_{\sigma \in \mathfrak S_n}\, \sigma.$
\end{defn}

   Ces constructions sont {\it fonctorielles} en $A$.

\begin{prop}[\cf \protect{\cite[7.2]{de}}]\label{a la Deligne} Soit $A$ un
objet de $\sA$, et soit $f$ un endomor\-phis\-me de $A$. Alors on a
\[tr(\Lambda^n(f))=\frac{1}{n!}
\sum_{\sigma\in \mathfrak{S}_n} sgn(\sigma) \prod_{k=1}^p tr(f^{\circ
\ell_k})\]
o\`u, pour tout $\sigma\in \mathfrak{S}_n$, $\ell_1,\dots,\ell_p$ sont les
longueurs des cycles disjoints constitutant $\sigma$. En particulier, si
$d=\dim A$, la dimension de
$\Lambda^n A$ est
\[\binom{d}{n} = \frac{d(d-1)\dots (d-n+1)}{n!}.\]
De m\^eme, on a
\[tr(\mathbf{S}^n(f))=\frac{1}{n!}
\sum_{\sigma\in \mathfrak{S}_n} \prod_{k=1}^p tr(f^{\circ
\ell_k})\]
\[\dim \mathbf{S}^n(A)=\binom{d+n-1}{n} = \frac{d(d+1)\dots
(d+n-1)}{n!}.\]
\end{prop}

\prf  La premi\`ere assertion r\'esulte directement du corollaire
\ref{c6.5.1}. La seconde r\'esulte du cas particulier $f=1$ et de
l'identit\'e bien connue
\[\binom{d}{n}=\frac{1}{n!}\sum_{\sigma\in \mathfrak{S}_n} sgn(\sigma)
d^{\alpha(\sigma)}\]
o\`u $\alpha(\sigma)$ est le nombre de cycles intervenant dans $\sigma$. Les
assertions pour les puissances sy\-m\'e\-tri\-ques se d\'emontrent de m\^eme.
\qed

\begin{cor} Sous les hypoth\`eses de la proposition \ref{a la Deligne},\\
a) $tr(f^{\circ n})$ s'exprime comme un poly\-n\^o\-me universel en les
$tr(\Lambda^i(f))$ pour $i=1,\dots,n$ \`a coefficients dans $\Z[1/n!]$.\\
b) On a
\begin{equation}\label{eq6.6.4}tr(\Lambda^n (1_A+f))= \sum_{i=0}^n
tr(\Lambda^i f ).
\end{equation}
\end{cor}

\prf a) Posons $t_i = tr(f^{\circ i})$: le second
membre de la proposition \ref{a la Deligne} est donc un poly\-n\^o\-me
universel en les $t_i$. On remarque que $t_n$ intervient dans le terme
$\prod_{k=1}^p tr(f^{\circ\ell_k})$ si et seulement si $\sigma$ est un
cycle de longueur $n$, et qu'alors ce terme est
\'egal \`a $t_n$. Tous ces cycles ont la m\^eme signature $(-1)^{n-1}$,
et ils sont au nombre de $(n-1)!$. Ainsi, le second membre de la
proposition \ref{a la Deligne} est un poly\-n\^o\-me du premier degr\'e en
$t_n$, et le coefficient de $t_n$ est \'egal \`a $(-1)^{n-1}/n$.
L'assertion en r\'esulte, par r\'ecurrence sur $n$.

b) Gr\^ace \`a a), on voit que   $tr(\Lambda^n(1+f))$   s'\'ecrit
comme un poly\-n\^o\-me universel \`a coefficients dans $\Z[\frac{1}{n!}]$ en
les $tr(\Lambda^i f )$ pour $i\leq n$. Prenant pour $\sA$ la ca\-t\'e\-go\-rie
mo\-no\-\"{\i}\-dale des modules libres de type fini sur $K$, on obtient la
formule voulue.\qed

   En nous inspirant de Kimura \cite{ki} (qui se pla\c cait
dans le cadre des motifs de Chow, voir aussi
\cite[3.15]{gu-pe}), nous allons d\'emontrer
un r\'esultat plus pr\'ecis que \ref{p6.4.1} et \ref{a la Deligne}, qu'on
retrouve en appliquant $\varepsilon_{A_1^\vee}$ aux facteurs extr\^emes.
On reprend les notations de \ref{p6.4.1}.

\begin{prop}\label{a la Kimura} Pour tout $\sigma\in  \mathfrak S_n$,
notons
$\ell_1(\sigma)$ la longueur du cycle $\sigma_1$ de $\sigma$ ayant $1$
dans son support. On a
\begin{multline}\label{eq6.6.3}\iota_{A_1
A_1}\big((1_{A_1^\vee}\bullet\epsilon_{A_n^\vee
\bullet\dots
\bullet A_2^\vee} \bullet 1_{A_1  })( \iota_{A_1\bullet \dots \bullet
A_n,A_1\bullet \dots \bullet
A_n}^{-1}( \sigma^{-1}\circ  f_1\bullet\dots\bullet f_n))\big)\\
=t_\sigma.f_{\sigma_1^{\ell_1-1}(1)}\circ\dots\circ f_{\sigma_1(1)}\circ
f_1,
\end{multline}
o\`u $t_\sigma= 1$ si $\sigma$ est cyclique
($\ell_1=n$), et $t_\sigma  \in K $ est un produit de traces de
mon\^omes non triviaux en les $f_i$ si $\sigma$ n'est pas cyclique.

A fortiori, si $A_1=\ldots = A_n =A$ et $f_1=\ldots = f_n =f \in
\sA(A,A)$, et si $n!$ est inversible dans $K$, on obtient
\begin{equation}\label{eq6.6.4a}
\iota_{AA}\big((1_{A^\vee}\bullet\epsilon_{A^{\bullet (n-1)}}
\bullet 1_{A })\circ (\iota_{A^{\bullet n},A^{\bullet n}}^{-1}(\Lambda^n
f ))\big)= \sum_{\sigma\in \mathfrak
S_n}  t'_\sigma f^{\ell_1(\sigma)}\end{equation}
\begin{equation}\label{eq6.6.4b}
\iota_{AA}\big((1_{A^\vee}\bullet\epsilon_{A^{\bullet (n-1)}}
\bullet 1_{A })\circ (\iota_{A^{\bullet n},A^{\bullet
n}}^{-1}(\mathbf{S}^n f ))\big)=
\sum_{\sigma\in \mathfrak
S_n}  t"_\sigma f^{\ell_1(\sigma)}\end{equation}
o\`u $t'_\sigma= \frac{(-1)^{n-1}}{n!}$
si $\sigma$ est cyclique, et
$t'_\sigma=\pm\frac{1}{n!}\times$ un produit de traces de puissances
non nulles de $f$ sinon (\resp m\^emes formules pour $t"_\sigma$, sans
les signes).
\end{prop}

\prf Abr\'egeons provisoirement $\iota_{AA}$ en $\iota_A$, et
\[(1_{A_1^\vee}\bullet\epsilon_{A_n^\vee
\bullet\dots
\bullet A_2^\vee} \bullet 1_{A_1  })\circ   \iota_{A_1\bullet \dots
\bullet A_n,A_1\bullet \dots \bullet
A_n}^{-1}(-) \]
en $\Theta_{A_1\ldots A_n}(-)$.

   Soit $\beta\in \mathfrak{S}_n$ une permutation telle que $\beta(1)=1$.
Compte tenu de \eqref{eq6.5.1}, \eqref{eq6.5.1,4/3} et
\eqref{eq6.5.1,5/3}, on a
   \begin{multline*}
\iota_{A_1 A_1} \big(   \Theta_{ A_{1}   A_{\beta (1)} \ldots
A_{\beta (n)}  }(\beta \sigma^{-1}
\beta^{-1}\circ (f_1\bullet f_{\beta (2)} \bullet \dots \bullet 
f_{\beta (n)} ))\big) \\=
\iota_{A_1 A_1} \big(\Theta_{A_1A_{\beta (2)}\ldots A_{\beta (n)} }
(\beta\circ \sigma^{-1}
\circ (f_1\bullet\dots\bullet f_n)\circ \beta^{-1})\big)\\
=\iota_{A_1 A_1}\big(\Theta_{A_1\ldots A_n}(\sigma^{-1}\circ
(f_1\bullet\dots\bullet f_n))\big).
   \end{multline*}
   Comme dans la preuve de \ref{p6.4.1}, et parce que $\sigma_1$ est
choisi de telle sorte que $1$ soit dans son support, cette formule
permet de se ramener au
cas o\`u les
$I_k$ forment une partition croissante de $\{1,\dots,n\}$ et o\`u, de
plus, pour tout $k$, $\sigma_k$ est le cycle ``croissant"
\[\sigma_k=(\ell_1+\dots+\ell_{k-1}+1,\cdots ,\ell_1+\dots+\ell_k),\]
    ce qui permet d'\'ecrire
   \[   \sigma^{-1}\circ f_1\bullet\dots\bullet f_n   =
(\sigma_1^{-1}\circ f_1\bullet\dots\bullet f_{\ell_1})\bullet \dots
\bullet  (\sigma_p^{-1}\circ f_{n-{\ell_p}}\bullet\dots\bullet f_n ). \]
      Compte tenu de \eqref{eq6.2.0} et
\eqref{eq6.2.8}, on en tire
   \begin{multline*} \iota_{A_1} \big(\Theta_{A_1\ldots A_n}(
f_1\bullet\dots\bullet f_n)\big)
     \\ =\iota_{A_1  }\big( (1_{A_1^\vee} \bullet
\epsilon_{A^\vee_{\ell_1} \bullet \dots \bullet A^\vee_2 }\bullet
1_{A_1  })
   (\iota_{A_1 \bullet \ldots \bullet A_{\ell_1} }^{-1}(\sigma_1^{-1}\circ
f_1\bullet\dots\bullet f_{\ell_1}) )\bullet \\
   \epsilon_{A^\vee_{\ell_1+\ell_2} \bullet \dots \bullet
A^\vee_{\ell_1+1} }
   (\iota_{A_{\ell_1 +1} \bullet \ldots \bullet A_{\ell_1+\ell_2}
}^{-1}(\sigma_2^{-1}\circ
f_{\ell_1+1}\bullet\dots\bullet f_{\ell_1+\ell_2}) )\bullet \\ \ldots
\bullet
   \epsilon_{A^\vee_{n } \bullet \dots \bullet A^\vee_{n-\ell_p} }
   (\iota_{A_{n-\ell_p} \bullet \ldots \bullet
A_{n}}^{-1}(\sigma_p^{-1}\circ f_{n-\ell_p}\bullet\dots\bullet f_{n})
)\big)
\end{multline*}
ce qui n'est autre que \[(f_{\ell_1}\circ \cdots \circ
f_1).tr(f_{\ell_1+\ell_2}\circ
\dots \circ f_{\ell_1 +1}).\cdots .tr(f_n \circ \cdots \circ
f_{n-\ell_p})\] par  \eqref{eq6.5.3}  et  \eqref{eq6.5.4}. Ceci
d\'e\-mon\-tre la premi\`ere assertion, et la seconde s'ensuit
imm\'edia\-tement.
   \qed

\begin{prop}\label{nouvelle}
Soit $A\ne 0$ un objet de $\sA$. On suppose qu'il existe un entier $n>0$
tel que $n!$ soit inversible dans $K$ et tel que
tel que $\Lambda^n A= 0$ (\resp que $\mathbf{S}^nA=0$). Alors\\
i) $\sN(A,A)$ est un nil-id\'eal de $\sA(A,A)$ d'\'echelon de nilpotence
$\leq n$. C'est m\^eme un id\'eal nilpotent de $\sA(A,A)$, d'\'echelon de
nilpotence born\'e en fonction de $n$.  En particulier, l'image de $A$
n'est pas nulle dans $\sA/\sN$.\\  ii) Si $K$ est int\`egre de
caract\'eristique nulle, la dimension de $A$ (\resp l'oppos\'e de la
dimension de
$A$) est un entier naturel $\le n$.
\end{prop}

\prf i) la formule \eqref{eq6.6.4a}
nous dit que  \[\iota_{A  }(1_{A ^\vee}\bullet\epsilon_{A^{\bullet
(n-1)}}\bullet 1_{A })(
\iota_{A^{\bullet n}}^{-1}(\Lambda^n f ))= \sum_{\sigma\in \mathfrak
S_n}  t'_\sigma f^{\ell_1(\sigma)},\]
o\`u $t'_\sigma= \frac{({-1)}^{n-1}}{n!}$ si $\sigma$ est cyclique, et
$t'_\sigma=\pm\frac{1}{n!}\times$  un produit de traces de puissances non
nulles de $f$ sinon. Pour
$f\in \sN$, ces derniers termes s'annulent, et on a donc
\[\iota_{A}(1_{A ^\vee}\bullet\epsilon_{A^{\bullet (n-1)}}
\bullet 1_{A })(\iota_{A^{\bullet n}}^{-1}(\Lambda^n f ))=
\frac{({-1)}^{n-1}}{n}f^n.\]

Mais $\Lambda^n f = 0$ si $\Lambda^n A= 0$, et on
trouve finalement que $f^n=0$. De m\^eme avec $S^n$ au lieu de
$\Lambda^n$. La nilpotence de $\sN(A,A)$ r\'esulte alors du th\'eor\`eme
de Nagata-Higman \cite{nagata,higman}, rappel\'e ci-dessous pour la 
commodit\'e du lecteur.

La troisi\`eme assertion de i) est imm\'ediate.

ii) D'apr\`es la proposition \ref{a la Deligne}, la dimension de
$\Lambda^n A$ est donn\'ee par le coefficient bin\^omial $\binom{\dim
A}{n}$, qui ne s'annule que si $\dim A$ est un entier naturel $\le n$.

Le cas d'une puissance sy\-m\'e\-tri\-que se traite de m\^eme.\qed

  \begin{lemme}[Nagata-Higman]\label{naghig} Soient $n$ un entier $>0$ 
et $K$ un anneau
commutatif unitaire dans lequel $n!$ est inversible. Soit
$R$ une $K$-alg\`ebre associative non unitaire telle que tout 
\'el\'ement $x\in R$ v\'erifie $x^n=0$. Alors
$R^{2^n-1}=0$.
\end{lemme}

\prf (d'apr\`es P. Higgins) On raisonne par r\'ecurrence sur $n$, en
posant $R= R_n$. Posons
\[\xi (x,y) = \sum_{i=0}^{n-1} x^i y x^{n-1-i}
\text{ et }\eta (x,y,z)=\sum_{i,j= 0}^{n-1} x^i z y^j x^{n-i-1}
y^{n-j-1}.\]

En \'ecrivant $(x+iy)^n= 0$ pour $i=0,\ldots n-1$, on
obtient par Van der Monde (compte tenu de ce que $(n-1)!$ est inversible
dans $K$) que $\xi(x,y)$ est identiquement nul sur $R_n\times R_n$. Donc
$\eta(x,y,z)=\sum_{j=0}^{n-1} \xi(x, zy^j)y^{n-1-j}=0$, mais
$\eta(x,y,z)$ se r\'e\'ecrit $nx^{n-1}zy^{n-1} +
\sum_{j=0}^{n-2}\, x^j z \xi (y, x^{n-1-j}) $, d'o\`u $x^{n-1} zy^{n-1}
= 0$ (compte tenu de ce que $n$ est inversible dans $K$). Soit $I_{n-1}$
l'id\'eal de $R_n$ engendr\'e par les puissances $(n-1)$-i\`emes.  La
$K$-alg\`ebre $R_{n-1}= R_n/I_{n-1}$ est nil d'\'echelon $n-1$. Par
r\'ecurrence, $R_{n-1}^{2^{n-1}-1}=0$,
\ie $R_n^{2^{n-1}-1}\subset I_{n-1}$, et on conclut de ce qui
pr\'ec\`ede que
\[R_n^{2^n-1}= R_n^{2^{n-1}-1}.R_n.R_n^{2^{n-1}-1}\subset I_{n-1}\cdot
R_n\cdot I_{n-1}=0.\]
  \qed

\subsection{Trace des endomor\-phis\-mes nilpotents}\label{tracesnil}

Nous ignorons si un endomor\-phis\-me nilpotent est toujours de trace
nilpotente, sans hypoth\`ese suppl\'e\-men\-taire. Nous allons n\'eanmoins 
le
d\'e\-mon\-trer dans deux cas particuliers.

\begin{prop}\label{nil, tr} Supposons que $\sA$ soit pseudo-ab\'elienne
et que la dimension de tout objet de
$\sA$ soit un entier naturel. Soit $A$ un objet de dimension $n$ telle
que $n!$ soit inversible dans $K$. Alors tout endomor\-phis\-me nilpotent de
$A$ est de trace nilpotente dans $K$.
\end{prop}

\prf Par le th\'eor\`eme de Krull et par extension des scalaires
(na\"ive), on se ram\`ene au cas o\`u $K$ est un corps. Par extension des
scalaires de $K$
\`a $K[t]$, la formule \eqref{eq6.6.4} donne
$tr(\Lambda^n (1+t f ))= 
\sum_{i=0}^n t^i tr(\Lambda^i f)$. Pour
montrer que $tr(f)=0$ pour 
tout $f$ nilpotent, il suffit donc d'\'etablir
que
$tr(\Lambda^n (1+f 
))=1$ pour tout tel $f$. Or par la proposition
\ref{a la Deligne}, le 
rang de $\Lambda^n A$ est $1$. Il suffit donc de
d\'emontrer que pour 
tout objet $B$ de rang $1$, $(\sA/\sN)(B,B)=K.1_B
\stackrel{tr 
\,\cong}{\To} K$, compte tenu de ce que $K$ est un corps.
Cela 
d\'ecoule du lemme suivant, qui est une cons\'equence imm\'ediate 
de
\ref{Ldim}:

\begin{lemme}\label{dim 0 ou 1} Supposons que $K$ 
soit un corps, que
$\sA$ soit pseudo-a\-b\'e\-lien\-ne et que la 
dimension de tout objet
soit un entier naturel. Alors tout objet de dimension $0$ de
$\sA/\sN$ est nul, et pour tout objet $B$ de $\sA/\sN$ de dimension
$1$, le mor\-phis\-me canonique
$\eta_B:\un\to B^\vee\bullet B$ est un isomor\-phis\-me. \qed
\end{lemme}

\begin{prop}[\cf \protect{\cite[1.4.3]{bru2}}]\label{a la bru} Supposons
qu'il existe un foncteur mo\-no\-\"{\i}\-dal sy\-m\'e\-tri\-que
$H$ de $\sA$ vers une ca\-t\'e\-go\-rie ab\'elienne mo\-no\-\"{\i}\-dale
sy\-m\'e\-tri\-que rigide
$\sV$; supposons en outre que l'application $K\to  \sV(\un, \un)$ soit
injective. Alors tout endomor\-phis\-me nilpotent est de trace nulle.
\end{prop}

\prf On se ram\`ene au cas $\sA=\sV$, donc \`a supposer $\sA$
ab\'elienne. Cela permet de factoriser
$f=m\circ e$ en
\'epi suivi de mono. Si $f^n=0$, on a donc
\[0=(me)^n=m(em)^{n-1}e\]
d'o\`u $(em)^{n-1}=0$ puisque $m$ (\resp $e$) est simplifiable \`a gauche
(\resp \`a droite). On \'ecrit alors $tr(me)=tr(em)$, et
on conclut par r\'ecurrence sur l'\'echelon de nilpotence.
\qed

\begin{rem} Dans \ref{nil, tr}, contrairement \`a \ref{a la bru},
on ne peut remplacer ``trace nilpotente" par ``trace nulle" si $K$ n'est
pas r\'eduit: la ca\-t\'e\-go\-rie $\sA$ des modules libres de type fini
sur l'alg\`ebre des nombres duaux $D=K[\epsilon]/(\epsilon^2)$ est
mo\-no\-\"{\i}\-dale,
$D$-lin\'eaire, sy\-m\'e\-tri\-que, rigide, et v\'erifie $End\,\un= D$ (mais
n'est pas ab\'elienne). Or
$\epsilon.1_\un$ est un endomor\-phis\-me nilpotent de trace non nulle.
\end{rem}

\subsection{$\otimes$-Radicaux}\label{tens-rad} On suppose toujours
$\sA$ mo\-no\-\"{\i}\-dale sy\-m\'e\-trique rigide.

\begin{defn}\label{d6.2} Soit $\sJ$ un id\'eal mo\-no\-\"{\i}\-dal de $\sA$. Le
{\rm
$\otimes$-radical} de $\sJ$ est d\'efini par:
\[\sqrt[ \otimes]{ \sJ}(A,B)=\{f\in \sA(A,B), \exists n \in {\bf N},
f^{\bullet n} \in \sJ(A^{\bullet
n},B^{\bullet n})\}.\]
\end{defn}

\begin{lemme}\label{l8.1}
i) $\sqrt[ \otimes]{ \sJ}$ est un id\'eal mo\-no\-\"{\i}\-dal.\\
ii) Pour tout objet $A$, l'id\'eal $\sqrt[\otimes]{0}(A,A)$ est
nil. Qui plus est, l'id\'eal bilat\`ere engendr\'e par un \'el\'ement
   quelconque de $\sqrt[\otimes]{0}(A,A)$ est nilpotent.\\
iii) $\sqrt[ \otimes]{0}$ est contenu dans
$\sR\cap \sN$ (rappelons que $\sR$ d\'esigne le radical de $\sA$.)
\end{lemme}

\prf $i)$ Pour montrer que $\sqrt[ \otimes]{ \sJ}$ est un id\'eal
mo\-no\-\"{\i}\-dal, le seul point non trivial est de
v\'erifier que $ \sqrt[ \otimes]{ \sJ}(A,B)$ est stable par addition.
Soient $f$ et $g$ deux \'el\'ements de
cet ensemble. Il existe donc $n \in {\bf N}$ tel que $f^{\bullet n}
=g^{\bullet n}\in \sJ ,$ et l'on voit alors
que $(f+g)^{\bullet 2n}\in \sJ $ puisque $\sJ $ est un id\'eal
mo\-no\-\"{\i}\-dal.

$ii)$ (\cf \cite[d\'em. de la prop. 2.16]{ki}). Soit $f\in \sqrt[
\otimes]{0}(A,A)$. Pour tous
$g_1,
\ldots ,  g_{n+1}\in \sA(A,A)$, on a $g_{n+1}\bullet f
\bullet g_n \bullet \ldots \bullet f \bullet g_{1}=0$ (appliquer le
tressage), d'o\`u, par
\eqref{eq6.5.3},
\begin{multline*}
\iota_{AA}^{-1}(g_1\circ f \circ g_2 \circ \ldots \circ f
\circ g_{n+1})=\\
(1_{A^\vee}\bullet\epsilon_{A^\vee \bullet\dots\bullet A^\vee }
\bullet 1_A)\circ  \iota_{A^{\bullet 2n+1},A^{\bullet 2n+1}}^{-1}
(\sigma^{-1} \circ(g_{n+1}\bullet f
\bullet g_n \bullet \ldots \bullet f \bullet g_{1}))
\end{multline*}
o\`u $\sigma$ est la permutation cyclique de
$\{1,\dots,2n+1\}$.  On obtient donc que  $g_1\circ f \circ g_2\circ
\ldots\circ f \circ g_{n+1}=0$.

$iii)$ $\sqrt[ \otimes]{0}\subset \sR$ suit imm\'ediatement de $ii)$.
Pour $\sqrt[ \otimes]{0}\subset \sN$, on utilise la
formule \eqref{eq6.2.10}: $tr(f^{\bullet n})= tr(f)^n$.
\qed

\begin{ex}\label{ex1'} Dans le cas de l'exemple \ref{ex2} (motifs),
l'id\'eal mo\-no\-\"{\i}\-dal $\sqrt[ \otimes]{0}$ correspond \`a
l'\'equivalence de smash-nilpotence sur les cycles al\-g\'e\-bri\-ques,
introduite par V. Voevodsky \cite{voevodsky}. Sa conjecture principale
(loc. cit.) revient \`a dire que $\sqrt[ \otimes]{0}=\sN$ (dans le cas
o\`u le corps de coefficients est de caract\'eristique nulle).

Par ailleurs, dans \cite{thomason}, Thomason calcule $\sqrt[ 
\otimes]{0}$ dans le cas
des ca\-t\'e\-go\-ries d\'eriv\'ees de ca\-t\'e\-go\-ries de faisceaux coh\'erents.
\end{ex}

\begin{rem}\label{r6.1a} La compatibilit\'e de 
$\sqrt[\otimes]{0}$
vis-\`a-vis de l'extension des sca\-lai\-res 
(na\"{\i}ve) est
imm\'ediate. 
\end{rem}

\section{Th\'eor\`emes de 
structure}\label{nouveau par}

Dans ce paragraphe, nous tirons les 
fruits des calculs pr\'ec\'edents,
pour obtenir des th\'eor\`emes de 
structure: ceux-ci sont rassembl\'es
dans \ref{rtss}. 
Dans tout le 
paragraphe, {\it $K$ est un corps et $\sA$ est une 
ca\-t\'e\-go\-rie 
$K$-lin\'eaire mo\-no\-\"{\i}\-da\-le sy\-m\'e\-tri\-que 
rigide,
avec
$End(\un )=K$}. 

On note $\sR$ le radical, et $\sN$ le 
plus grand id\'eal mo\-no\-\"{\i}\-dal
distinct de $\sA$. (\cf \ref{new 
p1}).    

\subsection{Structure de $\sA/\sN$}\label{a sur n}\

\begin{prop}\label{new p2} On suppose qu'il existe un foncteur 
$F:
\sA\to Vec_L$ vers la ca\-t\'e\-go\-rie des espaces vectoriels de 
dimension
finie sur une extension $L$ de $K$, v\'erifiant 
$\;F(\un)\neq 0$.\\
a) Alors les mor\-phis\-mes de $\sA/\sN$ forment des 
$K$-espaces de
dimension finie. Plus pr\'ecis\'ement,
\[\dim_K 
(\sA/\sN)(A,B)\leq \dim_LF(A^\vee\otimes B).\dim_LF(\un).\]
En 
particulier, $\sA/\sN$ est semi-primaire.\\
b) Si $F$ est la 
compos\'ee $G\circ H$ d'un
foncteur mo\-no\-\"{\i}\-dal sy\-m\'e\-tri\-que 
$H$ de $\sA$ vers une
ca\-t\'e\-go\-rie ab\'elienne mo\-no\-\"{\i}\-dale 
sy\-m\'e\-tri\-que rigide
$\sV$, et d'un foncteur $G:
\sV\to Vec_L$, 
alors
$\sR \subset\sN$ et $\sA/\sN$ est semi-simple.
 
\end{prop}

\prf $a)$ On va donner deux preuves, la premi\`ere 
valable seulement si
$F$ est $K$-lin\'eaire.

1) Comme 
$(\sA/\sN)(A,B)\otimes_K L\simeq (\sA_L/\sN_L)(A,B)$ 
(\cf
remarque
\ref{r6.1}), on peut supposer $L=K$. D'autre part, par 
rigidit\'e, on peut
supposer que $A=\un$. Consid\'erons alors 
l'application 
$K$-lin\'eaire d\'efinie par $F$:
\[F:\sA(\un,B)\to 
Hom_K(F(\un),F(B)).\]

Soit $f\in \sA(\un,B)$ tel que $F(f)=0$. Pour 
tout $g\in \sA(B,\un)$, on
a $F(gf)=0$, donc $gf=0$: en effet, 
$F:K=End(\un)\to End_K(F(\un))$ est
injectif puisque $F(\un)\ne 0$. 
Ainsi $\Ker F\subset \sN(\un,B)$. On a
donc
\[\dim \sA/\sN(\un,B)\le 
\dim \sA(\un,B)/\Ker F\le \dim
Hom(F(\un),F(B)).\]

2)  Soit $V$ le 
sous-$L$-espace de $  Hom_L(F(A^\vee\bullet  B),F(\un)) $
engendr\'e 
par $F(\sA(A^\vee\bullet  B,\un))$. On peut donc trouver 
des
\'el\'ements $b_1, \dots, b_n$ de
$\sA(A^\vee\bullet  B,\un)$, 
avec $n\leq \dim_L F(A^\vee\bullet
B).\dim_LF(\un)$, qui engendrent 
$V$.  Soit $a\in
\sA(A,B)\cong
\sA(\un,A^\vee\bullet  B)$. Dire que 
$a\in \sN(A,B)\cong
\sN(\un,A^\vee\bullet  B)$, c'est dire que 
$b\circ a=0$ dans $K=End(\un)$
pour tout
$b\in \sA(A^\vee\bullet 
B,\un)$. Comme $F$ est un foncteur, cela
revient
\`a dire que 
$V\circ F(a)=0\in L$. Comme les $F(b_i)$ engendrent $V$
sur $L$, cela 
revient finalement \`a dire que $b_i\circ a=0$ pour
$i=1,\dots n$. 
Ainsi $a\mapsto (b_i\circ a)_i$ d\'efinit
un plongement 
$K$-lin\'eaire de $(\sA/\sN)(A,B) \cong (\sA/\sN)(\un,
A^\vee\bullet 
B)$ dans $K^n$.

$b)$ Comme $\sA/\sN $ est semi-primaire, toutes 
les
$K$-alg\`ebres d'endo\-mor\-phis\-mes $\sA/\sN(C,C)$
sont 
semi-primaires, donc leurs radicaux sont nil\-po\-tents.
D'apr\`es 
\ref{a la bru}, tout \'el\'ement de ce radical est donc de
trace 
nulle.
  En appliquant ceci
\`a $C=A\oplus B$, on voit 
que
$\sR(A,B)\subset \sN(A,B)$. Donc $\sR\subset
\sN$, et $\sA/\sN$ 
est semi-simple puisque semi-primaire et
de radical 
nul.\qed

\begin{rems} \label{rmot} 1) La preuve de $a)$ est une 
variante
abstraite de la preuve classique de la finitude des cycles 
alg\'ebriques
modulo l'\'e\-qui\-va\-len\-ce nu\-m\'e\-ri\-que. 
L'hypoth\`ese est tr\`es
faible, mais non superflue: on peut montrer 
que l'exemple de
\cite[2.19]{de} (il s'agit d'une ind-ca\-t\'e\-go\-rie 
obtenue \`a partir de la
ca\-t\'e\-go\-rie mo\-no\-\"{\i}\-dale sy\-m\'e\-tri\-que 
rigide
$\sA$ librement engendr\'ee par un objet $X$ de dimension 
transcendante
sur
$\Q$) v\'erifie nos conditions et est m\^eme 
ab\'elienne, mais que
$\sA/\sN$ n'est pas semi-primaire (en fait, il 
d\'ecoule du calcul de
$\sA(X,X)$ fait dans
\loccit que 
$\sN(X,X)=0$).

2) Le m\^eme \'enonc\'e (sans les in\'egalit\'es de 
dimensions)
vaut si l'on suppose seulement que $L$ est une 
$K$-alg\`ebre semi-simple
commutative plut\^ot qu'un corps. En effet, 
on se ram\`ene \`a la
proposition \ref{new p2} en projetant sur la 
ca\-t\'e\-go\-rie $Vec_{L_i}$, o\`u
$L_i$ est un facteur simple de $L$ tel 
que $F(\un)\otimes_L L_i\ne 0$.
\end{rems}

  Voici un autre 
\'enonc\'e qui va dans le m\^eme sens:

\begin{prop}\label{new p5/2} 
On suppose que $\sA$ est 
pseudo-ab\'elienne et que les
dimensions 
des objets sont des entiers naturels.\\
$a)$ Alors les mor\-phis\-mes de 
$\sA/\sN$ forment des
$K$-espaces de dimension finie. On a plus 
pr\'ecis\'ement, pour tout
objet $A$:
\[\dim_K (\sA/\sN)(A,A)\leq 
(\dim A )^2.\]
$b)$ En fait  $\sR\subset
\sN$, et $\sA/\sN$ est 
semi-simple.
\end{prop}

\prf $a)$ Comme $\sA$ est $K$-lin\'eaire, on 
se ram\`ene au cas des 
endomor\-phis\-mes. Choisissons 
$n$
endomor\-phis\-mes $f_1,\dots, f_n$ dont les i\-ma\-ges dans 

$\sA(A,A)/\sN(A,A)$ sont
lin\'eairement in\-d\'e\-pen\-dan\-tes. Il 
s'agit de voir que $n\leq
(\dim\,A)^2=
\dim\,A^\vee\bullet A$. Au 
lemme \ref{Ldim}, on a vu que les 
$f_1,\dots, f_n$ donnent lieu 
\`a
un facteur direct ${\un}^n$ dans $A^\vee\bullet
A$ dans la 
ca\-t\'e\-go\-rie mo\-no\-\"{\i}\-dale pseudo-ab\'elienne $\sA/\sN$. 
Tout
suppl\'ementaire sera de dimension un entier
naturel par 
hypoth\`ese, d'o\`u l'in\'egalit\'e voulue.

$b)$: comme dans la 
proposition pr\'ec\'edente, mais en invoquant 
\ref{nil, tr} au lieu 
de \ref{a la bru}.\qed

\subsection{Radical, traces, 
semi-simplicit\'e. Trois th\'eor\`emes 
de
structure}\label{rtss}\

Comme dans le sous-paragraphe 
pr\'ec\'edent, on suppose que {\it $K$ est
un corps, et que $\sA$ est 
une  ca\-t\'e\-go\-rie $K$-lin\'eaire mo\-no\-\"{\i}\-dale
sy\-m\'e\-tri\-que rigide, 
avec
$End(\un )=K$}.

   \begin{thm}\label{wedtr} On suppose $K$ de caract\'eristique nulle et
$\sA$ munie d'un foncteur mo\-no\-\"{\i}\-dal sy\-m\'e\-tri\-que fid\`ele $H:\sA\to
L$-$Modf$, o\`u $L$ est une $K$-alg\`ebre commutative semi-simple. Alors
$\sA$ est de Wedderburn et
$H$ a les propri\'et\'es
\'enum\'er\'ees dans le corollaire \ref{c6} b).
\end{thm}

\prf Comme dans la remarque \ref{rem en fam}, on se ram\`ene au cas o\`u
$L$ est un corps. D'apr\`es le corollaire
\ref{c6}, il suffit alors de  v\'erifier que, pour tout $a\in
\sA(A,A)$, le poly\-n\^o\-me
caract\'eristique de $H(a)$ est \`a coefficients dans $K$. Or, comme
$\car K=0$, ces coefficients s'expriment comme po\-ly\-n\^o\-mes
\`a coefficients rationnels en les $tr(H(a^i))$. Puisque $H$ est
mo\-no\-\"{\i}\-dal sy\-m\'e\-tri\-que, on a $tr(H(a^i))=H(tr(a^i))\in
H(End(\un))=K$.
\qed

Au vu de l'exemple \ref{ex1}, le point $a)$ de l'\'enonc\'e
suivant peut
\^etre con\-si\-d\'e\-r\'e comme une version abstraite des
r\'esultats de Jann\-sen \cite{jannsen} (\cf \cite{ak(note)}).

\begin{thm}\label{absJannsen} a) On suppose qu'il existe un
foncteur mo\-no\-\"{\i}\-dal sy\-m\'e\-tri\-que $H$ de $\sA$ vers une
ca\-t\'e\-go\-rie ab\'elienne mo\-no\-\"{\i}\-dale sy\-m\'e\-tri\-que $\sV$, et
un foncteur $\,G:\sV\to L$-$Modf$ v\'erifiant
$\,G(\un)\neq 0$, o\`u $L$ est une $K$-alg\`ebre commutative
semi-simple.  Alors
\begin{thlist}
\item $\sR\subset \sN$,
\item $\sA/\sN$ est semi-simple,
\item le seul id\'eal mo\-no\-\"{\i}\-dal $\sI$ de $\sA$ tel que $\sA/\sI$ soit
semi-simple est $\sI=\sN$.
\end{thlist}
b) Supposons de plus que
$K$ soit de caract\'eristique nulle, que
$G=Id$ et que $H$ soit fid\`ele. Alors
\begin{thlist}
\item $\sR=\sN$,
\item $\sA$ est de
Wedderburn,
\item $\sA$ est pseudo-ab\'elienne si et seulement si
$\sA/\sN$ est ab\'elienne.
\end{thlist}
\end{thm}

\prf a) (i) et (ii) suivent de la proposition \ref{new p2} b) et de
la remarque \ref{rmot} 2); (iii), indiqu\'e pour m\'emoire, a d\'ej\`a
\'et\'e vu dans la proposition \ref{new p1} c).

b) Pour (i), il faut
voir que
$\sN \subset \sR$. On remarque que, via $H$, la trace
mo\-no\-\"{\i}\-dale $tr$ s'interpr\`ete comme une  ``vraie" trace au sens de
l'alg\`ebre lin\'eaire. Comme $K$ est suppos\'e de caract\'eristique
nulle, un endomor\-phis\-me d'un $L$-module de type fini dont toutes les
puissances strictement positives sont de trace nulle est nilpotent, et on
d\'eduit de l\`a que $\sN(A,A)$ est un nil-id\'eal, donc contenu dans
$\sR(A,A)$, pour tout objet $A$ de $\sA$. Le point (ii)
d\'ecoule du th\'eor\`eme pr\'ec\'edent. Enfin, (iii) r\'esulte de ceci et
de la proposition \ref{P3/2} b).
\qed

\begin{rems}\label{nouv rem} 1) Si $\sA$ est tannakienne et que $K$
est de caract\'eristique
nulle, on verra dans la proposition \ref{c4} que $\sA/\sR$ est encore
tannakienne (neutre si $\sA$ l'est). Par contre, on verra aussi que
l'hypoth\`ese de ca\-rac\-t\'e\-ris\-ti\-que nulle est n\'ecessaire
(contre-exemple \ref{Carp}, ou remarque \ref{r1} 1).

2) Par ailleurs, il ne suffirait pas dans le th\'eor\`eme
\ref{absJannsen} de supposer que
$H$ est un foncteur fibre $\Z/2$-gradu\'e, \cf contre-exemple \ref{super}
ci-dessous. 
   \end{rems}

   \begin{thm}\label{nouveau} Supposons $K$ de
caract\'eristique 
nulle et $\sA$ pseudo-a\-b\'e\-lien\-ne. Consid\'erons 
les 
propri\'et\'es suivantes: \\
a) pour tout objet $A$ de $\sA$, il 
existe un entier $\,n>0\,$ tel que
$\Lambda^n A= 0$,\\
b) $\sA/\sN$ 
est semi-simple, et pour tout objet $A$ de $\sA$, $\sN(A,A)$ est un 
id\'eal nilpotent,\\
c) $\sR=\sN$,\\
d) $\sA$ ne contient pas d'objet 
fant\^ome, \ie
d'objet non nul dont l'image dans $\sA/\sN$ est 
nulle,\\
e) $\sA$ est de Wedderburn,\\
f) $\sA/\sN$ est tannakienne 
semi-simple. \\
g) la dimension de tout objet de $\sA$ est un entier 
naturel,\\
Alors
\begin{thlist}
\item $a)\If b) \If c) \If d) + e)$ 
et $a)\If f) + g)$,
\item sous $g)$, on a $a)\iff b)  \iff c)\iff 
d)$. 
\end{thlist}
\end{thm}

\begin{rems} 1) La propri\'et\'e $a)$ 
pr\'esente l'avantage sur les
suivantes d'\^etre testable sur des 
``g\'en\'erateurs tensoriels" de
$\sA$: soit
$(A_\alpha)$ une famille 
d'objets $A_\alpha$ de $\sA$ tels que tout
objet de $\sA$ soit 
isomorphe \`a un facteur direct d'une somme directe
de 
$\bullet$-mon\^omes en les $A_\alpha$; on peut montrer que $a)$ 
est
v\'erifi\'e d\`es que pour tout $\alpha$, il existe $n_\alpha$ 
tel que
$\Lambda^{n_\alpha} A_\alpha= 0$.\\
2) Les propri\'et\'es 
$a), f)$ et $g)$ d\'ependent du tressage, alors
que $b), c), d)$ et 
$e)$ n'en d\'ependent pas.\\ 
3) Nous verrons en \ref{super} que  $g) 
+ e)$ ne suffit pas \`a impliquer $
\sR=\sN$).
\\ 4) Les 
propri\'et\'es du th\'eor\`eme sont par exemple v\'erifi\'ees dans le 
cas de la ca\-t\'e\-go\-rie
mo\-no\-\"{\i}\-dale des fibr\'es vectoriels sur une 
vari\'et\'e projective g\'eo\-m\'e\-tri\-quement connexe sur
$K$. 

\end{rems}

\prf S\'eparons $b)$ en $b_1)$: $\sA/\sN$ est 
semi-simple, et $b_2)$: $\forall \sA$, $\sN(A,A)$ est un 
id\'eal
nilpotent. Nous allons d\'emontrer que 
$a)\If b_2)$, 
$b_1) 
+b_2)\If c)\If d)$,
$a)\If g)\If b_1)$,   $ b_1)+b_2)+c)\If e)$, $a)+ 
b_1)+b_2)\If f)$ et $d)+ g) \If a)$.

$a)\If b_1)$ suit de la 
proposition \ref{nouvelle} i).

$b_1)+b_2)\If c)$: l'inclusion $\sR 
\supset \sN$ vient gr\^ace au fait que le
radical d'une $K$-alg\`ebre 
contient tout nilid\'eal. L'inclusion
oppos\'ee de la 
semi-simplicit\'e de $\sA/\sN$ (\ref{L1}).

$c)\If d)$: en effet, la 
projection $\sA \to \sA/\sR$ \'etant
conservative, elle ne peut 
envoyer un objet non nul sur un objet nul.

$a)\If g)$: cela 
r\'esulte de la proposition \ref{nouvelle} ii).

  $ b_1)+b_2)+c)\If 
e)$: imm\'ediat.

  Il reste \`a montrer que $a)+ b_1)+b_2)\If f)$. 
Par 
$b_2)$, $\sA/\sN$ est
pseudo-ab\'elienne tout comme $\sA$ (lemme 
\ref{idemp}). Par $b_1)$, 
elle est donc ab\'elienne (\cf 
\ref{A.P4}).
Comme elle v\'erifie aussi $a)$, le th\'eor\`eme de 
Deligne 
\cite[7.1]{de} montre qu'elle est tannakienne.

$g)\If b)$: 
cela r\'esulte de la proposition \ref{new
p5/2}.

$d)+ g) \If a)$: 
Soit $A$ un objet de dimension $d$. Notons
$\bar A$ son image dans 
$\sA/\sN$ et posons $B:=\Lambda^{d+1}A$.
Montrons que $B=0$. Par $d)$, il suffit de montrer que son image
$\bar B$ dans $\sA/\sN$, qui n'est autre que $\Lambda^{d+1}\bar A$, est
nulle. Notons que $\bar B  $ est de dimension
$\binom{d}{d+1}=0$. Par $g)$, et d'apr\`es \ref{dim 0 ou 1},
on a bien
$\bar B=0$ puisque $\sA$ est pseudo-ab\'elienne.
  \qed

\subsection{Variante $\Z/2$-gradu\'ee}\label{Z/2}

\begin{sloppypar}
En th\'eorie des motifs, la ca\-t\'e\-go\-rie mo\-no\-\"{\i}\-dale des motifs purs
modulo l'\'equivalence num\'erique est parfois appel\'ee la ``fausse"
ca\-t\'e\-go\-rie des motifs. La raison en est que cette ca\-t\'e\-go\-rie ne peut
pas \^etre tannakienne car la dimension d'un objet $M$ est un entier
\'eventuellement n\'egatif: si $H$ est une cohomologie de Weil et si
$\tilde M$ est un rel\`evement de $M$ \`a la ca\-t\'e\-go\-rie des motifs
modulo l'\'equivalence $H$-homologique, alors $\dim M = \dim \tilde M
=\sum_{i\in \Z} (-1)^i H^i(M)$. En particulier, la condition g) du
th\'eor\`eme \ref{nouveau} ne s'applique pas dans ce contexte.

Pour corriger ce d\'efaut, on fait la conjecture standard que les
projecteurs de K\"unneth de tout motif $M$ (relatifs \`a $H$) sont
alg\'ebriques, ce qui permet de modifier la contrainte de
commutativit\'e de fa\c con \`a obtenir la formule $\dim M =\sum_{i\in
\Z} H^i(M)\ge 0$. Pour y parvenir, il suffit en fait que la somme des
projecteurs de K\"unneth pairs (ou impairs, de mani\`ere
\'equi\-va\-len\-te), soit alg\'ebrique.
\end{sloppypar}

Faute de savoir d\'emontrer cette conjecture standard en g\'en\'eral, on
peut se limiter
\`a la sous-ca\-t\'e\-go\-rie pleine form\'ee des motifs $M$ pour lesquelles
elle est vraie: c'est le point de vue adopt\'e dans \cite{ak(note)}.
(C'est une sous-ca\-t\'e\-go\-rie mo\-no\-\"{\i}\-dale qui contient au moins les
motifs attach\'es aux courbes, aux surfaces et aux vari\'et\'es
ab\'eliennes, et qui contient tous les motifs si le corps de base est un
corps fini.) Le formalisme de la modification de la contrainte de
commutativit\'e est abstrait dans la proposition suivante:

\begin{prop}\label{new p4} Supposons $K$ de caract\'eristique nulle, et 
soit $L$ une $K$-alg\`ebre semi-simple commutative. Supposons $\sA$
munie  d'un  foncteur mo\-no\-\"{\i}\-dal sy\-m\'e\-tri\-que fid\`ele $H$ vers la   
ca\-t\'e\-go\-rie mo\-no\-\"{\i}\-dale sy\-m\'e\-tri\-que $Modf_L^\pm$ des $L$-modules de
type fini $\Z/2$-gradu\'es munie de la r\`egle de Koszul.
Consid\'erons, pour tout $A\in \sA$, la projection $p^+_A\in End(H(A))$
sur le facteur
$H^+(A)$. Soit
$\sA^\pm$ la sous-ca\-t\'e\-go\-rie pleine de $\sA$ form\'ee des objets $A$ tels
que $p^+_A= H(\pi^+_A)$ pour un $\pi^+_A\in \sA(A,A)$. On note $\sR^\pm,
\sN^\pm$ la restriction de $\sR, \sN$
\`a $\sA^\pm$. Alors
\begin{thlist}
\item  $\sA^\pm$ est une sous-ca\-t\'e\-go\-rie $K$-lin\'eaire mo\-no\-\"{\i}\-dale
de $\sA$, stable par facteurs directs et rigide.
\item Il existe sur $\sA^\pm$ une contrainte de commutativit\'e
telle 
que le foncteur mo\-no\-\"{\i}\-dal compos\'e
\[\sA^\pm\to 
\sA\stackrel{H}{\to} Vec_L^\pm\to Vec_L,\]
o\`u le dernier foncteur 
est le foncteur d'oubli, soit
sy\-m\'e\-tri\-que.
\end{thlist}
En 
particulier, $\sA^\pm$ est de Wedderburn, $\sR^\pm=\sN^\pm$, 
et
$\sA^\pm/\sN^\pm$ est semi-simple. 
 \end{prop} 

\prf 
V\'erification imm\'ediate pour (i): si $A,B\in 
\sA^\pm$,
$\pi^+_{A\bullet B}$ est donn\'e par 
$\pi^+_A\bullet\pi^+_B
+(1-\pi^+_A)\bullet(1-\pi^+_B)$. De m\^eme, si 
$A^\vee$ est le dual de
$A\in
\sA^\pm$, $\pi^+_{A^\vee}$ est donn\'e 
par le transpos\'e $(\pi^+_A)^t$.
Pour
$A\in \sA^\pm$, 
posons
$\epsilon_A=2\pi^+_A-1$: on a $\epsilon_A^2=1$. Pour $A,B\in 
\sA$,
soit $R_{A,B}:A\bullet B\to B\bullet A$ la contrainte 
de
commutativit\'e. Si $A,B\in \sA^\pm$, 
posons
\begin{equation}\label{kunneth}
R'_{A,B}=R_{A,B}\epsilon_A\bullet\epsilon_B.
\end{equation}

On voit tout de suite que $R'$ est une 
nouvelle contrainte de
commutativit\'e, v\'erifiant la condition 
(ii): pour la naturalit\'e,
cela r\'esulte du fait que $\pi_A^+$, 
donc $\epsilon_A$, est naturel en
$A$, et la v\'erification des 
identit\'es de tressage sy\-m\'e\-tri\-que est
imm\'ediate. La derni\`ere 
affirmation r\'esulte du th\'eor\`eme
\ref{absJannsen}.\qed

Si l'on 
est dans la situation de la proposition \ref{new p4}, on 
peut
appliquer le th\'eor\`eme \ref{nouveau} \`a $\sA^\pm$ munie de 
la
nouvelle contrainte de commutativit\'e. En exprimant tout en 
termes de
l'ancienne contrainte, cela revient \`a utiliser au lieu 
des puissances
ext\'erieures les ``super-puissances ext\'erieures" 
d'un objet
$A=A_+\oplus A_-$ de
$\sA^\pm$, d\'efinies 
par
\[s\Lambda^n A= \bigoplus_{i+j=n}\Lambda^i A_+\bullet 
\mathbf{S}^j A_-.\]

Nous reviendrons sur ces super-puissances 
ext\'erieures dans
la section suivante \ref{kimura}, o\`u nous 
d\'emontrerons que la
ca\-t\'e\-go\-rie $\sA^\pm$ de la proposition 
\ref{new p4} est
\emph{intrins\`eque} (\ie ind\'ependante de 
$H$).

\section{Le Pair et 
l'Impair}\label{kimura}

\begin{sloppypar}
La ca\-t\'e\-go\-rie $\sA$ est 
toujours $K$-lin\'eaire mo\-no\-\"{\i}\-dale
sy\-m\'e\-tri\-que rigide, qu'on suppose en outre {\it pseudo-ab\'elienne};
$K$ est un corps de caract\'eristique nulle.
\end{sloppypar}

\subsection{Fondements}

\begin{defn}\label{dkim.1} a) Un objet $A\in \sA$ est \emph{pairement de
dimension finie au sens de Kimura} s'il existe un entier $n>0$ tel que
$\Lambda^n(A)=0$. \\
b) Un objet $A\in \sA$ est \emph{impairement de
dimension finie au sens de Kimura} s'il existe un entier $n>0$ tel que
$\mathbf{S}^n(A)=0$. \\
c) Un objet $A\in \sA$ est \emph{de
dimension finie au sens de Kimura} s'il admet une d\'ecomposition en somme
directe $A=A_+\oplus A_-$ telle que $A_+$ soit pairement de dimension
finie et $A_-$ soit impairement de dimension finie.
\end{defn}

Cette d\'efinition a \'et\'e introduite par Kimura dans \cite{ki} pour les
motifs de Chow. Un grand nombre de ses r\'esultats sont valables, avec
leur d\'e\-mons\-tra\-tion, dans notre cadre g\'en\'eral. Pour les
d\'emonstrations de ce qui suit, nous renvoyons donc \`a Kimura, sauf pour
les \'enonc\'es qui ne figurent pas chez lui.

\begin{nota} Si $A$ est pairement (\resp impairement) de dimension finie
au sens de Kimura,  on note $\kim(A)$ le plus petit entier $n>0$ tel que
$\Lambda^{n+1}(A)=0$ (\resp tel que
$\mathbf{S}^{n+1}(A)=0$).
\end{nota}

Cette notation abusive ne pr\^etera pas \`a confusion car nous
pr\'eciserons \`a chaque fois si l'on est dans le cas pair ou
impair \cite{platon}\footnote{Du reste, on verra ult\'erieurement qu'un
objet
\`a la fois pairement et impairement de dimension finie est nul.}.

\begin{ex} L'objet unit\'e $\un$ est pairement de dimension finie; on a
$\kim(\un)=1$.
\end{ex}

\begin{prop}\label{pkim.1} a) Toute somme directe, tout facteur direct
d'un objet pairement (\resp impairement) de dimension finie est pairement
(\resp impairement) de dimension finie.\\
b) Soient $A,B\in \sA$, pairement ou impairement de
dimension finie. Alors $A\bullet B$ est pairement de dimension
finie si $A$ et $B$ sont de m\^eme parit\'e, impairement de dimension
finie sinon. De plus, $\kim(A\bullet B)\le\kim(A)\kim(B)$.\\
c) Le dual d'un objet pairement (\resp impairement) de dimension finie
est pairement (\resp impairement) de dimension finie.\\
d) Si $A$ est
pairement de dimension finie, il en est de m\^eme de 
toutes ses puissances
ext\'erieures.\\
e) Si $A$ est impairement de 
dimension finie, il en est de m\^eme de ses
puissances sy\-m\'e\-tri\-ques 
d'ordre impair, tandis que ses puissances
sy\-m\'e\-tri\-ques d'ordre pair 
sont pairement de dimension finie.
\end{prop}

\prf a) r\'esulte des 
isomor\-phis\-mes
\[\Lambda^n(A\oplus 
B)\simeq\bigoplus_{p+q=n}
\Lambda^p(A)\bullet\Lambda^q(B)\]
\[\mathbf{ 
S}^n(A\oplus 
B)\simeq\bigoplus_{p+q=n}
\mathbf{S}^p(A)\bullet\mathbf{S}^q(B).\]

b): \cf \cite[prop. 5.10]{ki}.

c) r\'esulte des isomor\-phis\-mes 
$\mathfrak{S}_n$-(anti-)\'equivariants
\[(A^\vee)^{\bullet n}\simeq 
(A^{\bullet n})^\vee.\]

d) et e) r\'esultent imm\'ediatement de a) 
et
b).\qed

\begin{lemme}\label{lkim.1} a) Soit $A\in \sA$, pairement 
de dimension
finie. Alors $\dim A$ est un entiernaturel $\le \kim 
A$.\\
b) Soit $A\in \sA$, impairement de dimension
finie. Alors 
$-\dim A$ est un entier naturel $\le \kim A$.
\end{lemme}

\prf C'est 
une reformulation de la proposition \ref{nouvelle} 
ii).\qed

\begin{prop} \label{lkim.2} Soit $A\in\sA$, pairement ou 
impairement de
dimension finie. Supposons que $\dim A=0$. Alors 
$A=0$.
\end{prop}

\prf Soit $\sB$ la plus petite
sous-ca\-t\'e\-go\-rie 
pleine rigide \'epaisse (\ie stable par facteurs directs)
de $\sA$ 
contenant $A\bullet A^\vee$. D'apr\`es la proposition
\ref{pkim.1}, 
tout objet de $\sB$ est pairement de dimension finie.
D'apr\`es le 
lemme \ref{lkim.1}, cela implique que la dimension de tout
objet de 
$\sB$ est un entier naturel. De plus, 
$\dim(A\bullet
A^\vee)=\dim(A)^2=0$. D'apr\`es la 
proposition
\ref{nouvelle} i), $\sN(A,A)$ est donc un id\'eal 
nilpotent de
$\sA(A,A)=\sB(A,A)$, et d'apr\`es la proposition 
\ref{new p5/2} a),
$\sB/\sN(A,A)=0$. Il en r\'esulte bien que 
$A=0$.\qed

\begin{thm}\label{pkim.2} a) Soit $A\in \sA$ pairement de 
dimension
finie. Alors $\kim A\allowbreak=\dim A$.\\
b) Soit 
$A\in\sA$ impairement de dimension finie. Alors $\kim A = 
-\dim
A$.
\end{thm} 

\begin{sloppypar}
\prf a) Soit $m'=\dim A$. Vu 
le lemme \ref{lkim.1} a), il suffit de
montrer que 
$\Lambda^{m'+1}(A)=0$. Mais d'apr\`es la proposition \ref{a 
la
Deligne}, on a $\dim\Lambda^{m'+1}(A)=\binom{m'}{m'+1}=0$. 
L'assertion
r\'esulte donc de la proposition pr\'ec\'edente. 

\end{sloppypar}

b) M\^eme d\'emonstration, en utilisant une 
puissance sy\-m\'e\-tri\-que.\qed

Le corollaire suivant n'est pas 
d\'emontr\'e chez Kimura.

\begin{cor}\label{ckim.1}
a) Si $A$ et $B$ 
sont pairement de dimension finie, alors $\kim(A\oplus
B)=\kim A+\kim 
B$.\\
b) Si $A$ et $B$ sont impairement de dimension finie, alors 
$\kim(A\oplus
B)=\kim A+\kim 
B$.\qed
\end{cor}

\begin{prop}\label{kimnil} Soient $A_+$ (\resp 
$A_{-}$) un objet de
$\sA$ pairement (\resp impairement) de dimension 
finie. Soient $m=\kim
A_+$ et $n=\kim A_{-}$. Alors pour tout
$f\in 
\sA(A_+,A_{-})$ et tout $g\in \sA(A_{-}, A_+)$, on a
$f^{\bullet 
(mn+1)}=0$,  $g^{\bullet (mn+1)}=0$.
\end{prop}

\prf Voici une 
variante de \cite[d\'em de la prop. 6.1]{ki}, un peu plus
simple: 
gr\^ace \`a la proposition \ref{pkim.1} b) et par dualit\'e, on 
se
ram\`ene \`a traiter le cas de $f$ avec $A_+=\un$. Le 
mor\-phis\-me
$f^{\bullet (n+1)}:\un=\un^{\bullet (n+1)}\to A_-^{\bullet 
(n+1)}$ est
$\mathfrak{S}_n$-\'equivariant, ce qui implique qu'il se 
factorise \`a
travers $\mathbf{S}^{n+1}(A_-)=0$.\qed

\begin{prop} 
\label{pkim.3} Soit $A$ un objet de dimension finie au sens
de 
Kimura, et supposons donn\'ees deux d\'ecompositions 
$A\simeq
A_+\oplus A_-\simeq B_+\oplus B_-$, avec $A_+,B_+$ pairement 
de dimension
finie et $A_-,B_-$ impairement de dimension finie. Alors 
$A_+\simeq B_+$
et $A_-\simeq B_-$.
\end{prop}

\prf \cf \cite[d\'em. 
de la prop. 6.3]{ki}. On remarquera que cette
d\'emonstration utilise 
le corollaire \ref{ckim.1} a).\qed

\begin{lemme}\label{lkim.3} Soit 
$A\in\sA$, et soit $A=A_+\oplus A_-$ une
d\'ecomposition en somme 
directe. Posons
\[s\Lambda^n(A_+,A_-)= \bigoplus_{i+j=n}\Lambda^i 
A_+\bullet \mathbf{S}^j
A_-.\]
Alors les conditions suivantes sont 
\'equi\-va\-len\-tes:
\begin{thlist}
\item Il existe un entier $n>0$ tel 
que $s\Lambda^n(A_+,A_-)=0$.
\item $A_+$ est pairement de dimension 
finie et $A_-$ est impairement de
dimension finie.
\end{thlist}
De plus, dans (i), on peut prendre $n=\kim A_+ +\kim A_-+1$, et ceci est
le plus petit choix possible.
\end{lemme}

\prf Laiss\'ee au lecteur (pour voir que $\kim A_+ +\kim A_-+1$ est
le plus petit choix possible, observer que \[\dim s\Lambda^{\kim
A_++\kim A_-}(A_+,A_-)=\dim s\Lambda^{\dim
A_+-\dim A_-}(A_+,A_-)=1,\] en utilisant le corollaire \ref{ckim.1}).\qed

\begin{prop}[\protect{\cite[prop. 6.9]{ki}}]\label{pkim.4} Tout facteur
direct d'un objet de dimension finie est de dimension finie.\qed
\end{prop}

\begin{defn}\label{dkim.2} Soit $A$ un objet de dimension finie au sens
de Kimura. On d\'efinit
\[\kim A =\kim A_+ +\kim A_-=\dim A_+ -\dim A_-\]
\[s\Lambda^n A= \bigoplus_{i+j=n}\Lambda^i A_+\bullet \mathbf{S}^j A_-\]
o\`u $A=A_+\oplus A_-$ est une d\'ecomposition de $A$ en somme d'un objet
pairement de dimension finie et d'un objet impairement de dimension
finie.
\end{defn}

La proposition \ref{pkim.3} montre que $\kim A$ {\it ne d\'epend que de
$A$ et que $s\Lambda^n A$ ne d\'epend que de $A$ \`a isomor\-phis\-me pr\`es}.

On a maintenant
la g\'en\'eralisation suivante de la proposition
\ref{nouvelle}:

\begin{prop}\label{nouvellekim} Soit $A\in \sA$ un objet de dimension
finie au sens de Kimura. Alors $\sN(A,A)$ est un id\'eal nilpotent de
$\sA(A,A)$, d'\'echelon born\'e en termes de $\kim(A)$.  En particulier,
pour $A\ne 0$, l'image de $A$ n'est pas nulle dans $\sA/\sN$.
\end{prop}

\prf D'apr\`es \ref{naghig}, il suffit de faire voir que $\sN(A,A)$ 
est un nilid\'eal d'\'echelon born\'e en
termes de $\kim(A)$. L'argument est analogue \`a celui de \cite[prop. 
7.5]{ki}: tout
\'el\'ement $f$ de $\sN(A,A)$ peut s'\'ecrire comme une matrice $2\times
2$ sur la d\'ecomposition $A=A_+\oplus A_-$. On \'ecrit $f=f_+ + f_- +
f_\pm$, o\`u $f_+$ est le terme pr\'eservant $A_+$, $f_-$ le terme
pr\'eservant $A_-$ et $f_\pm$ est la somme des deux termes antidiagonaux.
Comme $f_+f_-=f_-f_+=0$, un mon\^ome typique intervenant dans le
d\'eveloppement de $f^n$ est de la forme
\[m=f_{\epsilon_1}^{k_1}\circ f_\pm\circ f_{\epsilon_2}^{k_2}\circ
f_\pm\circ\dots\circ f_\pm\circ f_{\epsilon_r}^{k_r}\]
avec $\epsilon_i\in \{+,-\}$.

D'apr\`es la proposition \ref{nouvelle}, on a $m=0$ d\`es que l'un des
$k_i$ est assez grand. D'autre part, d'apr\`es la proposition
\ref{kimnil} et le lemme \ref{l8.1} i), on a $f_\pm^{\bullet N}=0$ pour
$N$ assez grand (ind\'ependamment de $f$),  ce qui implique (lemme 
\ref{l8.1} ii)) que $m=0$ d\`es
que le nombre de fois o\`u
$f_\pm$ appara\^{\i}t est $\ge N$. Ceci conclut la d\'emonstration.\qed

\begin{rem} On pourrait aussi utiliser le raisonnement de la
d\'e\-mons\-tra\-tion de la proposition \ref{P3/2} c) \`a partir de la
proposition \ref{nouvelle}. Ceci donne l'estimation suivante pour
l'\'echelon de nilpotence $N$ de $\sN(A,A)$:
\[N\le 2^{\kim A_++1}+2^{\kim A_-+1}\le 2^{\kim A+1}+1\]
bien meilleure que celle fournie par la d\'emonstration ci-dessus.
\end{rem}

\begin{defn}\label{defki} Une {\it ca\-t\'e\-go\-rie de Kimura}\footnote{Cette
notion a \'et\'e introduite ind\'ependamment par P. O'Sullivan sous le
nom de ``semi-positive category", comme
nous l'avons appris apr\`es la soumission de ce texte \`a publication 
(\cf postscriptum
\`a l'introduction g\'en\'erale).} sur un corps de
ca\-rac\-t\'e\-ris\-ti\-que nulle
$K$ est une 
ca\-t\'e\-go\-rie
$K$-lin\'eaire mo\-no\-\"{\i}\-dale sy\-m\'e\-tri\-que rigide, 
v\'erifiant $End(\un)=K$,
pseudo-ab\'elienne, et dont tout objet est 
de dimension finie au sens de
Kimura. 

\end{defn}

\subsection{Ca\-t\'e\-go\-ries de Kimura et ``projecteurs de 
K\"unneth"}\

Soient $K$ un corps de caract\'eristique nulle et $\sA$ 
une ca\-t\'e\-go\-rie
$K$-lin\'eaire mo\-no\-\"{\i}\-dale sy\-m\'e\-tri\-que rigide, 
v\'erifiant
$End(\un)=K$, pseudo-ab\'elienne. Soit $\sA_\kim$ la plus 
grande
sous-ca\-t\'e\-go\-rie de Kimura de $\sA$, \ie la sous-cat\'ego\-rie 
pleine de
$\sA$ form\'ee des objets de dimension finie au sens de 
Kimura: c'est
une sous-cat\'ego\-rie \'epaisse et rigide de $\sA$, 
d'apr\`es les
propositions \ref{pkim.1} et \ref{pkim.4}.

\begin{thm}\label{tkim.1} Supposons que $\sqrt[\otimes]{0}=0$. 
Alors
\\ a) La ca\-t\'e\-go\-rie $\sA_{\kim}$ est munie 
d'une
$\Z/2$-graduation compatible \`a la structure 
mo\-no\-\"{\i}\-dale.\\
b) Pour tout $A\in \sA_{\kim}$, notons $\pi^+_A$ le 
projecteur de $A$
d\'efinissant $A_+$ et $\epsilon_A=2\pi_A^+-1$: on 
a $\epsilon_A^2=1$.
D\'efinissons une nouvelle contrainte de 
commutativit\'e $R'$ sur
$\sA_{\kim}$ par la formule \eqref{kunneth}. 
Alors
\begin{thlist}
\item Pour tout $n\ge 0$, l'objet 
$s\Lambda^n(A)$ est la $n$-i\`eme
puissance sy\-m\'e\-tri\-que de $A$ 
relative \`a $R'$.
\item La dimension de $A$ relative \`a $R'$ est 
\'egale \`a $\kim A$; en
particulier, c'est un entier $\ge 
0$.

\noindent Cette propri\'et\'e caract\'erise d'ailleurs la 
graduation dont il est question en $a)$.
\end{thlist}
c) Soient $L$ 
une $K$-alg\`ebre semi-simple commutative et $H:\sA\to
Modf_L^\pm$ un 
foncteur v\'erifiant les hypoth\`eses de la proposition
\ref{new p4}. 
Alors avec les notations de \loccit, on a
$\sA_{\kim}=\sA^\pm$. Munie 
de la contrainte induite par $R'$, le
quotient de $\sA_{\kim}$ par 
son radical est un ca\-t\'e\-go\-rie tannakienne 
semi-simple.
\end{thm}

\prf a) Puisque $\sqrt[\otimes]{0}=0$, la 
proposition \ref{kimnil}
montre que si $A$ et $B$ 
sont
$\epsilon$-ment de dimension finie \`a la Kimura avec des 
parit\'es
diff\'erentes, alors $\sA(A,B)=0$. Ceci permet de renforcer 
la
proposition \ref{pkim.3} en l'\'enonc\'e suivant: si 
$A\in\sA_{\kim}$, la
d\'ecomposition de $A$ en somme directe d'un 
objet pair et d'un objet
impair est \emph{unique}. Il en r\'esulte 
imm\'ediatement que cette
d\'ecomposition est fonctorielle en $A$; 
les autres propri\'et\'es
r\'esultent de la proposition 
\ref{pkim.1}.

b) ne pr\'esente pas de difficult\'e.

c) Pour voir 
l'inclusion $\sA_{\kim}\subset\sA^\pm$, il suffit de voir
que 
$H(\pi_A^+)=p_A^+$ (avec les notations de la proposition 
\ref{new
p4}). En effet, on voit imm\'ediatement que $H(A_+)$ est 
pair et que
$H(A_-)$ est impair (\cf \cite[d\'em. de la prop. 
3.9]{ki}). Pour voir
l'inclusion oppos\'ee, il suffit de voir que si 
$A=A_+\oplus A_-\in
\sA^\pm$, alors $A_+$ est pairement de dimension 
finie et $A_-$ est
impairement de dimension finie. Choisissons $n$ 
tel que  $H(\Lambda^n(A_+))=\Lambda^n(H(A_+))=0$. Comme $H$ 
est
fid\`ele, cela entra\^{\i}ne que le mor\-phis\-me identit\'e de 

$\Lambda^n(A_+)$ est nul, donc
 $\Lambda^n(A_+)=0$; de m\^eme pour 
$A_-$.
  Le changement de contrainte de commutativit\'e comme 
en
\ref{new p4} ram\`ene \`a la situation du th\'eor\`eme 
\ref{nouveau},
qui permet de conclure que la ca\-t\'e\-go\-rie quotient par $\sR = \sN$ munie
de la contrainte induite par $R'$ est tannakienne semi-simple.\qed

\begin{thm}\label{tctki} Toute ca\-t\'e\-go\-rie de Kimura $\sA$ est
de Wedderburn. Son radical $\sR$ est le plus grand id\'eal mo\-no\-\"{\i}\-dal
distinct de $\sA$. Le quotient de $\sA$ par $\sR$ est une
ca\-t\'e\-go\-rie tannakienne semi-simple (apr\`es
changement appropri\'e de la contrainte de commutativit\'e).
\end{thm}

\prf On a $\sA=\sA_\kim$. La premi\`ere affirmation r\'esulte de la
proposition \ref{nouvellekim}. Pour le reste, on se ram\`ene au
th\'eor\`eme
\ref{tkim.1} en quotientant par $\sqrt[\otimes]{0}$ (\cf le lemme
\ref{l8.1} iii))\qed

\begin{sloppypar}
\begin{rem}\label{nilplusplus}
\'Etant donn\'e une ca\-t\'e\-go\-rie de Kimura $\sA$, le th\'eor\`eme
\ref{tctki} nous permettra d'appliquer les th\'eor\`emes de scindage
mo\-no\-\"{\i}\-dal \ref{T2} et \ref{t4} \`a la projection
$\sA\to \sA/\sR$. En particulier, \emph{on pourra munir $\sA$
d'une $\Z/2$-graduation compatible \`a celle de
$\sA/\sqrt[\otimes]{0}$, unique \`a \'equivalence mo\-no\-\"{\i}\-dale
pr\`es}. Voir le th\'eor\`eme \ref{kimsplit}.
  \end{rem}
\end{sloppypar}

\begin{ex}\label{cj kimura} Kimura \cite[7.1]{ki} conjecture que
tout motif de Chow est de dimension finie en
son sens (c'est-\`a-dire, avec la terminologie ci-dessus, que la
ca\-t\'e\-go\-rie 
des motifs de Chow est une
$\Q$-ca\-t\'e\-go\-rie de Kimura). Cette conjecture 
implique donc que tout 
motif ``fant\^ome",
\ie  num\'eriquement nul, est 
nul.
R\'eciproquement, modulo l'existence pour chaque motif de Chow $A$
d'un 
projecteur $\pi^+_A$ qui s'envoie sur
le projecteur de K\"unneth pair 
$p^+_A$, la non-existence de motifs
fant\^omes non nuls implique la 
conjecture
de Kimura (consid\'erer le motif $s\Lambda^n A$ pour $n$ assez 
grand).

Par ailleurs, la conjecture de Murre
mentionn\'ee en \ref{BBM} 
implique l'e\-xis\-ten\-ce d'un tel $p^+_A$,
et,  jointe aux conjectures 
standard de
Grothendieck (\'equivalence homologique = \'equivalence 
num\'erique),
implique la non-exis\-ten\-ce de motifs
fant\^omes non nuls. 
On voit ainsi que la conjecture de Murre (ou, ce qui
est
\'equivalent, la 
conjecture de Bloch-Beilinson) jointe aux
conjectures standard implique la 
conjecture de
Kimura.
\end{ex}

\section{Exemples et compl\'ements}

\subsection{Trois contre-exemples}\label{3contrex} Voici trois
e\-xem\-ples de ca\-t\'e\-go\-ries mo\-no\-\"{\i}da\-les
$K$-lin\'eaires sy\-m\'e\-tri\-ques et rigides, v\'erifiant $End(\un )=K$,
ab\'e\-lien\-nes, de Wedderburn, mais
\emph{dont le  radical n'est pas un
id\'eal mo\-no\-\"{\i}\-dal}.

\begin{contrex}\label{super}\begin{sloppypar}
Partons de $D=K[\epsilon],\epsilon^2=0$, vue comme
super-alg\`ebre ($\epsilon$ plac\'e
en degr\'e impair). On consid\`ere le super-groupe alg\'ebrique
$GL(1\vert 1, D)$. Son
groupe de points est donn\'e par les matrices
$2\times 2$ inversibles du type
\[\begin{pmatrix} a &\epsilon b \\
\epsilon c & d
\end{pmatrix}, a,b,c,d \in K, a\neq 0, d\neq 0 \]
\cf \cite[ch. 1]{qfs}. On prend pour $\sA$ la ca\-t\'e\-go\-rie des
re\-pr\'e\-sen\-ta\-tions de
$GL(1\vert 1, D)$. C'est une sous-ca\-t\'e\-go\-rie $K$-lin\'eaire 
mo\-no\-\"{\i}\-dale rigide non
pleine de la ca\-t\'e\-go\-rie
mo\-no\-\"{\i}\-dale sy\-m\'e\-tri\-que rigide des super-$D$-modules.
\end{sloppypar}

La re\-pr\'e\-sen\-ta\-tion standard de $GL(1\vert 1, D)$ est le
super-$D$-module $D^{1\vert 1}=D\times \Pi D$ (o\`u
$\Pi$ est le foncteur changement de parit\'e). Un endomor\-phis\-me de
$D^{1\vert 1}$ est
donn\'e par une matrice
         $\begin{pmatrix} a &\epsilon b \\
\epsilon c & d
\end{pmatrix}$ comme ci-dessus; s'il commute \`a tout \'el\'ement
$\begin{pmatrix} a'
&\epsilon b' \\
\epsilon c' & d'
\end{pmatrix}\in GL(1\vert 1, D)$, on trouve que $ba'+db'=ab'+bd',
ca'+dc'=ac'+cd'$ pour tous $ a',b',c',d' \in K, a'\neq 0, d'\neq 0$,
d'o\`u, $b=c=0,a=d$. Donc $\sA(D^{1\vert
1},D^{1\vert 1})\cong K$ est r\'eduit aux homoth\'eties, et son
radical est nul.

\begin{sloppypar}
La trace mo\-no\-\"{\i}\-dale de $\begin{pmatrix} a &\epsilon b \\
\epsilon c & d
\end{pmatrix}$ est la supertrace $a-d$. Elle s'annule donc pour toute
homoth\'etie, ce qui
montre que
$\sN(D^{1\vert 1},D^{1\vert 1})=\sA(D^{1\vert 1},D^{1\vert 1})$. En
d'autres termes, $D^{1\vert 1}$ est un objet
``fant\^ome": il s'envoie sur l'objet nul dans $\sA/\sN$.
\end{sloppypar}

On a donc $0\subset \sR \subset \sN \subset \sA$ et toutes les
inclusions sont strictes  (\cf th\'eor\`eme \ref{absJannsen} et corollaire
\ref{c2}); en particulier
$\sA$ n'est ni
semi-simple, ni de Kimura. Elle est en revanche de Wedderburn 
(finitude des espaces d'homomor\-phis\-mes), et la
dimension de tout objet est un entier naturel ($\sA/\sN\sim Vec_K$).
\end{contrex}

\begin{contrex}\label{Carp} $K$ est maintenant suppos\'e de
caract\'eristique
$p>0$, et on prend pour $\sA = Rep_K(GL_p)$
la ca\-t\'e\-go\-rie des re\-pr\'e\-sen\-ta\-tions de dimension finie de $GL_p$.
Soit $V$ la re\-pr\'e\-sen\-ta\-tion
standard. On a $\sA(V,V)\allowbreak=K$, $\sR(V,V)=0$ et, comme toute
homoth\'etie de $V$ est de trace nulle, $\sN(V,V)=K$.
Ici encore, $V$ est un objet ``fant\^ome", \ie nul modulo $\sN$.
\end{contrex}

\begin{contrex}\label{Car3} $K$ est maintenant suppos\'e de
caract\'eristique
$3$, et on prend $\sA = Rep_K(SL_2)$.
  Soit $V$ la re\-pr\'e\-sen\-ta\-tion
standard. On a $V^{\otimes 2}\cong \un \oplus S^2V$, et on v\'erifie
comme ci-dessus que $S^2V$ est un objet ``fant\^ome". On en d\'eduit que
la ca\-t\'e\-go\-rie mo\-no\-\"{\i}\-dale sy\-m\'e\-tri\-que rigide $\sA/\sN$ est
\'e\-qui\-va\-len\-te
\`a celle des $K$-super-espaces de dimension finie.
\end{contrex}

\subsection{Compl\'ement: l'alg\`ebre commutative
$HH_0(\sA)$}\label{h0} Soient $\sA$ une petite $K$-ca\-t\'e\-go\-rie et $\sI$
un id\'eal de $\sA$.

On note $H_0(\sA,\sI)$ le $K$-module\footnote{L'explication de
cette notation appara\^itra au paragraphe suivant.} engendr\'e par les
classes d'isomor\-phis\-mes de couples
$(A,h)$ o\`u
$A$ est un objet de $\sA$ et $h\in \sI(A,A)$, quotient\'e par les
relations suivantes:
\begin{itemize}
\item l'application canonique $\sI(A,A)\to H_0(\sA,\sI)$ est
$K$-lin\'eaire,
\item pour tout $f\in \sI(A,B)$ et tout $g\in \sA(B,A)$, les classes de
$(A,gf)$ et $(B,fg)$ co\"{\i}ncident dans
$H_0(\sA,\sI)$.
\end{itemize}

Cette construction est fonctorielle en $(\sA,\sI)$.

On note $[A,h]$ la classe de $(A,h)$.  Notons que si $\sA$ est
$K$-lin\'eaire, et si $A,A'$ sont deux objets, on a $[A,1_A= pi] =
[A\oplus A', ip]$ avec les notations $i,p$ de
\eqref{eqbiprod}.

On note encore  $H_0(\sA,\sA)=HH_0(\sA)$. Ce $K$-module est le
r\'eceptacle universel des traces abstraites, \ie des applications
$t: A\to M$ v\'erifiant $t(fg)=t(gf)$ pour tout couple de mor\-phis\-mes
allant en directions oppos\'ees. Il appara\^{\i}t en th\'eorie des
n\oe uds, \cf \cite[2]{vogel}, dans le cadre mo\-no\-\"{\i}\-dal comme
ci-dessous.

\begin{sloppypar} \medskip
Supposons maintenant $\sA$ mo\-no\-\"{\i}\-dale. Par
fonctorialit\'e, la loi
$\bullet$ induit sur
$HH_0(\sA)$ une structure de
$K$-alg\`ebre, commutative si $\sA$ est munie d'un tressage
sy\-m\'e\-tri\-que\footnote{L\`a encore, il suffirait que
$\sA$ soit munie d'une structure de tortil, la sym\'etrie du tressage
n'\'etant pas vraiment n\'ecessaire.}. Si $\sI$ est mo\-no\-\"{\i}\-dal,
$H_0(\sA,\sI)$ est un id\'eal de $HH_0(\sA)$, et le quotient
$HH_0(\sA)/H_0(\sA,\sI)$ s'identifie \`a $HH_0(\sA/\sI)$. Par exemple,
$H_0(\sA,\sqrt[\otimes]{0})$ est un nil-id\'eal de $HH_0(\sA)$ (par
d\'efinition, tout \'el\'ement de $\sqrt[\otimes]{0}(\sA,\sA)$ est
$\bullet$-nilpotent).
\end{sloppypar}

Les
traces abstraites \`a valeurs dans une $K$-alg\`ebre commutative $R$
v\'e\-ri\-fiant de plus
$t(f\bullet g)=t(f)t(g)$ pour tout couple d'endomor\-phis\-mes $(f,g)$ sont
en bijection avec les
$R$-points de
$\Spec\, HH_0(\sA)$. Si $\sA$ est rigide, la transposition d\'efinit une
involution de $HH_0(\sA)$; la trace mo\-no\-\"{\i}\-dale
$tr$ d\'efinit un \emph{$K$-point canonique} de $\Spec HH_0(\sA)$ ou
encore une \emph{augmentation} de la $K$-alg\`ebre $HH_0(\sA)$,
compatible avec l'involution.

  \medskip Supposons de plus que $\sA$ v\'erifie les hypoth\`eses du 
th\'eor\`eme \ref{absJannsen}. On a donc
$\sR=\sN$.  Posons
$\bar\sA= \sA/\sR = \sA/\sN$. Le raisonnement de la d\'emonstration de la
proposition
\ref{a la bru} montre que toute trace abstraite s'annule sur tout
endomor\-phis\-me nilpotent, donc sur les $\sR(A,A)$. Donc le $K$-dual de
$HH_0(\sA)$ co\"{\i}ncide avec le
$K$-dual de $HH_0(\bar\sA)$, et finalement $HH_0(\sA)\cong
HH_0(\bar\sA)$.

\begin{ex}\label{hh0} Supposons $\sA$ tannakienne semi-simple neutre 
sur un corps $K$
alg\'e\-bri\-que\-ment clos de caract\'eristique nulle, de sorte que $\sA$
est \'e\-qui\-va\-len\-te \`a la ca\-t\'e\-go\-rie des $K$-re\-pr\'e\-sen\-ta\-tions de
dimension finie d'un $K$-groupe pro-r\'eductif $G$.

Par lin\'earit\'e, les \'el\'ements
$HH_0(\sA)$ sont des $K$-combinaisons d'objets
$[A,h]$ avec $A$ simple; par le lemme de Schur, on peut m\^eme
supposer $h=1_A$. On en d\'eduit que
\[HH_0(\sA)\cong R_K(G)\otimes_\Z K,\]
o\`u $R_K(G)$ d\'esigne l'anneau des re\-pr\'e\-sen\-ta\-tions de $G$.
\end{ex}

\newpage
\addtocontents{toc}{{\bf III. Sections}\hfill\thepage}
\
\bigskip
\begin{center}\large\bf III. Sections
\end{center}
\bigskip

L'objectif de cette partie est d'\'etablir des analogues cat\'egoriques
du th\'e\-o\-r\`e\-me de Wedderburn, qui affirme l'existence de sections
pour la projection d'une alg\`ebre de dimension finie sur un corps parfait
vers son quotient par le radical. On \'etablit de tels analogues tant
dans le cadre lin\'eaire g\'en\'eral que dans le cadre mo\-no\-\"{\i}\-dal ou
mo\-no\-\"{\i}\-dal tress\'e.

\section{Cohomologie de Hochschild-Mitchell}\label{hochschild}\label{s4}

\subsection{Produit tensoriel de $K$-ca\-t\'e\-go\-ries}\label{prodtens} Ici,
$K$ d\'esigne de nouveau un anneau commutatif unitaire quelconque.

\begin{defn}\label{D5} (\cf \cite[\S 11]{mitchell}). Soient $\sA,\sB$
deux $K$-ca\-t\'e\-go\-ries. Le
\emph{produit tensoriel} de $\sA$ par $\sB$ est la $K$-ca\-t\'e\-go\-rie
$\sA\boxtimes_K\sB$ dont les objets sont les couples $(A,B)$ ($A\in \sA,
B\in
\sB$), avec \[(\sA\boxtimes_K \sB)((A, B),(A', B'))=\sA(A,A')\otimes_K
\sB(B,B').\]
\end{defn}
On a des isomor\-phis\-mes canoniques
\[\sA\boxtimes_K\sB\cong
\sB\boxtimes_K\sA,\quad(\sA\hbox{-}Mod)^\sB\cong (\sB\hbox{-}Mod)^\sA
\cong
(\sA\boxtimes_K \sB)\hbox{-}Mod.\]

Par ailleurs, on note $ \sA^o$ la ca\-t\'e\-go\-rie oppos\'ee de $\sA$ et 
on pose ${\sA}^e= \sA^o  \boxtimes_K\sA$. Les
${\sA}^e$-modules (\`a gauche) sont parfois
appel\'es
$\sA$-bimodules. Alors $\sA$ d\'efinit de mani\`ere \'evidente un
$\sA$-bimodule (par la loi $(X,Y)\mapsto \sA(X,Y)$) qu'on identifie
\`a
$\sA$; un id\'eal
(bilat\`ere) de $\sA$ (\cf \S 1) n'est autre qu'un sous-bimodule de $\sA$.

\begin{lemme}\label{L6} Supposons que $K$ soit un corps. Si $\sA$ et $\sB$ 
sont
s\'eparables (\cf d\'efinition \ref{D4sep1}), alors il en est de m\^eme
de $\sA\boxtimes_K \sB$. A fortiori $\sA\boxtimes_K \sB$ (et en
particulier ${\sA}^e$) est semi-simple.
\end{lemme}

Cela r\'esulte du lemme \ref{l2.3}.\qed

\begin{lemme}\label{L6'} Supposons seulement que les $K$-ca\-t\'e\-go\-ries
$\sA$ et $\sB$ soient de Wedderburn (\cf d\'efinition
\ref{D4wed}). Alors il en est de m\^eme de $\sA\boxtimes_K \sB$.
\end{lemme}

  \prf En effet, il d\'ecoule du lemme pr\'ec\'edent que pour tout
objet $A$ de $\sA$ et tout objet $B$ de
$\sB$, le quotient de la $K$-alg\`ebre $(\sA\boxtimes_K
\sB)((A,B),(A,B))$ par l'id\'eal bilat\`ere
\[\rad(\sA)(A,A)\otimes_K
\sB(B,B)+\sA(A,A)\otimes_K \rad(\sB)(B,B)\]
est s\'eparable, donc
semi-simple. Il en d\'ecoule que cet id\'eal
contient le ra\-di\-cal. Par ailleurs, il est clair que c'est un
nil-id\'eal, donc contenu dans le radical, donc
finalement \'egal \`a $\rad(\sA\boxtimes_K \sB)$. D'o\`u le
r\'esultat.\qed

La caract\'erisation suivante des ca\-t\'e\-go\-ries de Wedderburn est celle
qui nous servira effectivement dans cette partie.

\begin{prop}\label{Wedcaract} On suppose que $K$
  est un corps. Une $K$-ca\-t\'e\-go\-rie $\sA$ est de Wedderburn si et 
seulement si $\sA^e$ est
semi-primaire et $\sA^e/\rad(\sA^e)=(\sA/\rad(\sA))^e $.
\end{prop}

\prf Le lemme pr\'ec\'edent montre que $\sA^e$ est semi-primaire et
v\'erifie
$\sA^e/\rad(\sA^e)=(\sA/\rad(\sA))^e $ d\`es lors que $\sA$ est de 
Wedderburn. R\'eciproquement, si $\sA^e$
est semi-primaire et $\sA^e/\rad(\sA^e)=(\sA/\rad(\sA))^e $, on a
$\rad(\sA^e)\supset id\otimes_K
\rad(\sA)$, donc la condition de nilpotence (dans la d\'efinition des 
cat\'ego\-ries semi-primaires) passe
de $\sA^e$ \`a $\sA$; d'autre part, $\sA^e/\rad(\sA^e)$ est 
semi-simple par hypoth\`ese, et on conclut que
$(\sA/\rad(\sA))^e $ est semi-simple, donc $\sA/\rad(\sA)$ est s\'epa\-ra\-ble, 
\cf
\ref{D4sep1} .\qed

\begin{sloppypar}
Le m\^eme exemple que dans la remarque \ref{imparf} montre que la
condition
$\sA^e/\rad(\sA^e)=(\sA/\rad(\sA))^e $ est n\'ecessaire.
\end{sloppypar}

\subsection{Complexe de Hochschild-Mitchell}\label{7.2} Soit $\sA$ une
petite
$K$-ca\-t\'e\-go\-rie, telle que pour tout
$(A,B)\in
\sA\times \sA$, $\sA(A,B)$
soit un $K$-module projectif.

Notons $C_*(\sA)$ le complexe de Hochschild de $\sA$ \cite[\S
17]{mitchell}: c'est un complexe de
${\sA}^e$-modules tel que
\[\scriptstyle 
C_m(\sA)(A,B)=\displaystyle\bigoplus_{A_0,\dots,A_{m}}
\scriptstyle
\s 
A(A,A_0)\otimes_K\sA(A_0,A_1)\otimes_K\dots\otimes_K\sA(A_{m-1}
,A_m) 
\otimes_K \sA(A_m,B),\quad m\geq 0\]
o\`u $(A_0,\dots,A_m)$ d\'ecrit 
les $n$-uplets d'objets de $\sA$. La
diff\'erentielle 
est
\[
d_m(f_0\otimes
f_1\otimes\dots\otimes 
f_{m+1})=
\sum_{i=0}^m(-1)^if_0\otimes
f_1\otimes\dots \otimes 
f_if_{i+1}\otimes\dots \otimes f_{m+1}.
\]

Il est connu 
que $C_*(\sA)$ est une r\'esolution
projective du ${\sA}^e$-module 
${\sA}$.

\medskip
\begin{sloppypar}
Pour tout ${\sA}^e$-module $M$, 
on note
$C^*(\sA,M)$ le complexe de 
$K$-modules
$Hom_{\sA^e}(C_*(\sA), M)$.  Plus concr\`etement, pour 
$m>0$, un
\'el\'ement de $C^m(\sA,M)(A,B)$ {\it
\'equivaut}
\`a la 
donn\'ee d'un \'el\'ement de
$M(A,B)$ pour chaque $m$-uplet de 
fl\`eches composables $(f_1,\dots,
f_{m})$, $f_1$ \'etant de source 
$A$ et
$f_m$ de but $B$, cette donn\'ee \'etant $K$-multilin\'eaire 
en
$(f_1,\dots, f_{m})$; pour $m=0$,
c'est la donn\'ee d'un 
\'el\'ement de $M(A,A)$ pour tout $A$.

Voici comment: on obtient une 
telle donn\'ee en
consid\'erant, dans
\allowbreak 
$Hom_{\sA^e}(C_*(\sA), M)(A_0,A_m)$ le
facteur correspondant \`a 
$A=A_0, f_0=1_{A_0}, B=A_{m},
f_{m}=1_{A_{m}}$. R\'eciproquement, une 
telle
donn\'ee fournit une cocha\^{\i}ne de Hochschild gr\^ace \`a 
la
commutativit\'e du carr\'e (bivariance 
en
$(A,B)$):
\[\begin{CD}
Hom_{\sA^e}(C_*(\sA), M)(A,B) 
@>>>M(A,B)\\
@V{(f_0,f_{m+1})}VV 
@V{(f_0,f_{m+1})}VV\\
Hom_{\sA^e}(C_*(\sA), M)(A_0,A_m) @>>> 
M(A_0,A_m).
\end{CD}\]
\end{sloppypar}

Avec cette interpr\'etation, 
la diff\'erentielle de Hochschild prend
la forme 
habituelle
\begin{multline*}(d^m\phi)(f_1,\dots,
f_{m+1})=f_1\phi(f_2, 
\dots, 
f_{m+1})\\
+\sum_1^m
(-1)^i\phi(f_1,\dots,f_if_{i+1},\dots,
f_{m+1})\\ 
+(-1)^{m+1}\phi(f_1,\dots, f_m)f_{m+1}.
\end{multline*}

Puisque 
$C_*(\sA)$ est une r\'esolution projective du
${\sA}^e$-module 
${\sA}$, on a
\[H^i(\sA, M):= 
H^i(C^*(\sA,M))=\Ext^i_{\sA^e}({\sA},M).\]

Si $\sA$ est s\'eparable, 
$\sA^e$ est semi-simple, donc ces groupes de cohomologie sont donc
 
nuls pour tout $M$ et tout
$i>0$ (\cf lemme \ref{L6} et proposition 
\ref{P4/3}). 
 
\begin{sloppypar}
\medskip
La formation de $C_*(\sA)$ 
est naturelle en $\sA$. Plus
pr\'ecis\'ement, tout
$K$-foncteur $F: 
\sA\to \sB$ donne lieu \`a des applications
$K$-lin\'eaires 
$C_n(\sA)(A,B)\allowbreak\to C_n(\sB)(FA,FB)$,
et, pour tout 
$\sB^e$-module $N$, \`a un homomor\-phis\-me de complexes
de $K$-modules 
dans l'autre sens
\[ C^*(\sB ,N)\to  C^*(\sA,F^\ast 
N).\]
\end{sloppypar}
  
\begin{rem} \label{rhoch} (Homologie de 
Hochschild) Si $M$ et $N$ sont
deux
$\sA$-bimodules, on peut voir $M$ 
comme
${\sA^e}$-module \`a droite, et $N$ comme ${\sA^e}$-module \`a 
gauche,
ce qui permet de former le $K$-module 
$M\otimes_{\sA^e} N$ 
(\cf \cite[pp. 29, 71]{mitchell}). Cette
construction \'etant 
fonctorielle en $N$, on peut b\^atir un complexe
$M\otimes_{\sA^e} 
C_*(\sA)$, et on a 
\[H_i(\sA, M):= H_i(M\otimes_{\sA^e} 
C^*(\sA))=\Tor_i^{\sA^e}(M,{\sA}).\]
En particulier, pour $i=0$, on 
retrouve la notion du \S  \ref{h0}.
\end{rem}

\section{Un 
``th\'eor\`eme de Wedderburn \`a plusieurs objets"}\label{s5}

Dans 
ce paragraphe, on se donne un corps $K$. On renvoie \`a 
la
d\'efinition \ref{D4wed} pour la notion de $K$-ca\-t\'e\-go\-rie
de 
Wedderburn.

\subsection{Existence et ``unicit\'e" de sections} Le 
th\'eor\`eme de
structure fondamental pour les
$K$-alg\`ebres de 
Wedderburn (th\'eor\`eme de Wedderburn) dit
que toute extension d'une 
$K$-alg\`ebre s\'eparable de dimension finie par un id\'eal
nilpotent 
se scinde. En suivant la preuve
co\-ho\-mo\-lo\-gi\-que classique de 
Hochschild et Whitehead, nous
g\'en\'eralisons ce r\'esultat aux 
alg\`ebres ``\`a
plusieurs objets".

\begin{thm}\label{T1}
Soit $\sA$ 
une petite $K$-ca\-t\'e\-go\-rie de Wedderburn, de
radi\-cal $\sR$. Posons 
$\bar\sA = \sA/\sR$. Alors:\\
a) Le $K$-foncteur de projection 
$\pi_{\sA}:\sA\to \bar\sA$
admet une section $s$.\\
b) Si $s'$ est 
une autre section, il existe une famille $(u_A)_{A\in
\sA}, u_A \in 
1_A+\sR(A,A)$, telle que pour tous $A,B\in
\sA$ et tout $f\in 
\bar\sA(A,B)$, on ait $s'(f)=
u_Bs(f)(u_A)^{-1}$.
\end{thm}

\prf 
$1)$ Supposons d'abord $\sR^2=0$. Alors $\sR$
h\'erite d'une 
structure de $\bar\sA^e$-module par passage au quotient de
l'action 
de $\sA$.

Pour tout $A,B\in \sA$, choisissons une 
section
$K$-lin\'eaire
\[\sigma_{A,B}: \bar\sA(A,B)\to \sA(A,B)\]
 
de la projection naturelle. Si $C$
est un troisi\`eme objet de $\sA$ 
et si $(f,g)\in \bar\sA (A,B)\times
               \bar\sA (B,C)$, 
notons
\[\gamma(g,f)=\sigma_{A,C}(gf)-\sigma_{B,C}(g)\sigma_{A,B}(f).\]

C'est un \'el\'ement de $\sR(A,C)$, et, si $D$ est un quatri\`eme 
objet
et
$h\in \bar\sA (C,D)$, on a la 
relation
\[\gamma(h,gf)+h\gamma(g,f)=\gamma(hg,f)+\gamma(h,g)f.\]

Ain 
si $\gamma$ d\'efinit un \emph{$2$-cocycle de Hochschild} 
de
$\bar\sA$ \`a valeurs dans $\sR$. D'a\-pr\`es le \S 
\ref{hochschild}, ce
$2$-cocycle est un $2$-cobord; autrement dit, il 
existe une fonction
$\rho$ telle que, pour tous $f,g$ composables, on 
ait
\[\gamma(g,f)=\rho(gf)-g\rho(f)-\rho(g)f.\]

La nouvelle section 
$s(f)=\sigma(f)-\rho(f)$ v\'erifie alors
$s(gf)=s(g)s(f)$ pour tout 
couple $(f,g)$ composable.

$2)$ Posons $\delta(f)=s'(f)-s(f)\in 
\sR(A,B)$. On a, pour tout couple
$(f,g)$ 
composable,
\[\delta(gf)=\delta(g)f+\delta(f)g\]
donc $\delta$ est un 
$1$-cocycle de Hochschild. C'est donc un cobord
$f\mapsto n_B f-fn_A$ 
pour un choix convenable d'\'el\'ements $n_A\in
\sR(A,A)$. En posant 
$u_A= 1_A+n_A$, on a bien $s'(f)=
u_Bs(f)(u_A)^{-1}$.

$3)$ Supposons 
maintenant $\sR$ nilpotent d'\'echelon $r+1$. On 
raisonne
par
r\'ecurrence sur $r$, le cas $r=1$ \'etant trait\'e 
ci-dessus. Soit
$\sA^{(r)}=\sA/\sR^{r}$: le radical $ \sR^{(r)} 
=
\sR/\sR^r$ de $\sA^{(r)}$ est nilpotent d'\'echelon $r$, donc il 
existe
une
section $\bar s:\bar\sA\to \sA^{(r)}$ comme
dans 
l'\'enonc\'e, et toute autre section s'en d\'eduit par
``conjugaison" 
par une famille de
$\bar u_A\in 1_A+\sR^{(r)} (A,A)$.

Soit $\bar 
s(\bar\sA)$ l'image d'une telle section $\bar
s$ dans $\sA^{(r)}$. 
Pour tout
mor\-phis\-me $f$ dans $\sA$, notons $ f^{(r)}$ son image dans 
$\sA^{(r)}$.
Posons
\begin{align*}
{\tilde\sA}(A,B)&=\{f\in 
\sA(A,B)\mid  f^{(r)}\in 
\bar
s(\bar\sA)(A,B)\}\\
{\tilde\sR}(A,B)&=\sR(A,B)\cap 
{\tilde\sA}(A,B).
\end{align*}

Alors ${\tilde\sA}$ est une 
sous-ca\-t\'e\-go\-rie (non pleine) de $\sA$,
${\tilde\sR}$ est 
son
radical, ${\tilde\sR}^2=0$, et ${\tilde\sA}/{\tilde\sR}=\bar\sA$. 
Il
existe donc une section $\tilde s:\bar
s(\bar\sA)\to {\tilde\sA}$; 
la compos\'ee $s=\iota\circ\tilde s\circ
\bar s:\bar\sA\to 
\bar
s(\bar\sA)\to{\tilde\sA}\to
\sA$, o\`u $\iota$ est l'inclusion 
canonique, est la section cherch\'ee.

Si $s'$ est une autre section, 
notons $\bar s'$ la composition
$\bar\sA\overset{s'}{\to}\sA\to 
\sA^{(r)}$, et soit $(\bar u_A)$ une
famille
d'\'el\'ements de 
$1_A+\sA^{(r)}(A,A)$ tels que $\bar s'(f)=\bar u_B
\bar s(f) (\bar 
u_A)^{-1}$ pour tout $f\in
\sA^{(r)}(A,B)$. Choisissons des relev\'es 
$u_A$ des $\bar u_A$ dans
$\sA(A,A)$. Quitte \`a changer $s'$ en 
la
section conjugu\'ee par les $(u_A)^{-1}$, on se ram\`ene au cas 
o\`u
$\bar s'=\bar s$. On conclut par l'\'etape
$2)$.

  $4)$ 
\'Etablissons maintenant l'assertion $a)$ du
th\'eor\`eme dans le cas 
g\'en\'eral. Si l'on supposait $\sA$ strictement
de Wedderburn, il 
suffirait de la voir comme limite projective des
$\sA^{(r)}$ du 
num\'ero pr\'ec\'edent, cette limite \'etant
alors localement 
stationnaire sur les $Hom$s, et d'appliquer directement
ce num\'ero. 
Mais nous supposons $\sA$ seulement de Wedderburn, ce qui
complique 
la d\'emonstration.

Consid\'erons l'ensemble des couples
$(\sB, s)$ 
form\'es d'une sous-ca\-t\'e\-go\-rie pleine $\sB$ de $\sA$ et 
d'une
section fonctorielle $s$ de la projection $\sB\to \bar \sB$ 
($\bar \sB$
peut \^etre vue comme sous-ca\-t\'e\-go\-rie
pleine de 
$\bar\sA$). Cet ensemble est non vide: d'apr\`es les 
pas
pr\'ec\'edents, il contient par exemple des
couples $(\sB, s)$ 
d\`es que $Ob(\sB)$ est fini.

On ordonne cet ensemble en 
d\'ecr\'etant que $(\sB', s')\prec (\sB,
s)$ si $\sB'\subset \sB$ et 
si
$s'$ est la restriction de $s$ \`a $\bar\sB'$. Il est clair que 
cet
ensemble ordonn\'e est inductif,
donc admet un \'el\'ement 
maximal $(\sB, s)$ d'apr\`es Zorn.
D\'emontrons par l'absurde que 
$\sB=\sA$.
Sinon, soit $X$ un objet de $\sA$ qui n'est pas dans 
$\sB$. On peut
remplacer $\sA$ par sa sous-ca\-t\'e\-go\-rie
pleine dont 
les objets sont $X$ et ceux de $\sB$, et il s'agit de
construire une 
extension de
$s$ \`a $\bar \sA$. Pour cela, on peut derechef 
remplacer $\sA$ par
toute sous-ca\-t\'e\-go\-rie non pleine $\sA'$
ayant 
m\^emes objets, et telle que $\sA'/(\sA'\cap \sR)=\bar\sA$.

Il y a 
un choix minimal d'une telle $\sA'$: celle dont
les mor\-phis\-mes sont 
des sommes finies de
compositions de fl\`eches de $s(\sB)$ et de 
fl\`eches de $\sA$ de
source ou de but $X$. Le radical de 
cette
$K$-ca\-t\'e\-go\-rie $\sA'$ n'est autre que $\sA'\cap \sR$: en 
effet, tous
les $(\sA'\cap \sR)(A,A)$ sont des
nil-id\'eaux, ce qui 
entra\^{\i}ne que $(\sA'\cap \sR)\subset
\rad(\sA')$; 
r\'eciproquement, $\sA'/(\sA'\cap
\sR)=\bar\sA$ est semi-simple, donc 
$ \rad(\sA')\subset \sA'\cap \sR$.

\begin{sloppypar}
Le point est 
que $\sA'\cap \sR$ est (globalement) nilpotent, ce qui
permet, par le 
pas 3), d'\'etendre $s$
comme souhait\'e. En fait, nous allons voir 
que si $n$ est
l'\'echelon de nilpotence de l'id\'eal
$\sR(X,X)$ de 
$\sA(X,X)$, alors $ \rad(\sA')^{(2n+1)}=0$. En effet,
comme le 
radical de $s(\sB)$ est nul, on
observe que pour tout couple d'objets 
$(A,B)$ de $\sA$, tout
\'el\'ement de $\rad(\sA')(A,B)$ s'\'ecrit 
comme combinaison lin\'eaire
finie $\sum f^i_{XB}\circ g^i_{AX}$, 
avec
$f^i_{XB}\in \sA(X,B), g^i_{AX}\in \sA(A,X)$. Tout \'el\'ement 
de
$\rad(\sA')^{(2n+1)}(A,B)$ s'\'ecrit comme
somme finie de termes 
de la 
forme
\end{sloppypar}
\begin{multline*}\sum_{i_1,i_3,\dots,i_{2n+1}} 
f^{i_{2n+1}}_{XB}
g^{i_{2n+1}}_{A_{2n}X}\left( \sum_{i_{2n}} f^{i_{2n 
}}_{XA_{2n}}
g^{i_{2n }}_{A_{2n-1}X}\right)f^{i_{2n-1}}_{XA_{2n-1}} 
\dots\\  \left(
\sum_{i_{2}} f^{i_{2 
}}_{XA_{2}}
g^{i_2}_{A_1X}\right)f^{i_1}_{XA_1} 
g^{i_1}_{AX}.
\end{multline*}

Or chaque fragment 
$g^{i_{2k+1}}_{A_{2k}X}( \sum_{i_{2k}} f^{i_{2k
}}_{XA_{2k}}
g^{i_{2k 
}}_{A_{2k-1}X})f^{i_{2k-1}}_{XA_{2k-1}}$ est 
dans
\[\sA(A_{2k},X)\circ 
\rad(\sA')(A_{2k-1},A_{2k})\circ
\sA(X,A_{2k-1})\subset 
\sR(X,X).\]

Comme ces fragments apparaissent $\;n\;$ 
fois cons\'ecutives, on trouve bien que
$\;\rad(\sA')^{(2n+1)}(A,B)=0\;$. 
Ceci ach\`eve la preuve de l'assertion $a)$.

\medskip $5)$ Un 
argument analogue s'applique pour l'assertion $b)$.
Soient
$s$ et 
$s'$ deux sections fix\'ees.
Par Zorn, il existe une sous-ca\-t\'e\-go\-rie 
pleine maxi\-male $\sB$ de
$\sA$ telle que les restrictions de $s$ et 
$s'$ \`a $\bar\sB$ soient
conjugu\'ees au sens de l'assertion 
$b)$.
Raisonnant par l'absurde comme ci-dessus, on peut encore 
supposer que
$Ob(\sA)=Ob(\sB)\cup\{X\}$. Quitte \`a
conjuguer $s$ 
au-dessus de $\bar\sB$, on peut supposer que les
restrictions de $s$ 
et $s'$ \`a $\bar\sB$
co\"{\i}ncident. On observe alors que $s$ et 
$s'$ prennent alors
leurs valeurs dans la sous-ca\-t\'e\-go\-rie
$\sA'$ non 
pleine de $\sA$ d\'efinie comme ci-dessus. Comme on l'a
vu, son 
radical est globalement nilpotent,
ce qui permet d'appliquer le pas 
3) et de conclure que $s$ et $s'$
sont conjugu\'ees via une 
famille
d'\'el\'ements $u_A\in 1_A+\sR(A,A)$.\qed

Le lemme 
sui\-vant, souvent utile, r\'esulte
imm\'ediatement de la proposition 
\ref{P3/2} b) et du fait que $\bar \sA$
est 
semi-simple:

\begin{lemme}\label{l2} Supposons $\sA$ 
pseudo-ab\'elienne. Soit $s$ une
section de
$\pi_\sA$. Les foncteurs 
$\pi_\sA$ et $s$ induisent des bijections
inverses l'une de l'autre 
entre les objets in\-d\'e\-com\-po\-sa\-bles
de $\sA$ et les objets 
irr\'eductibles de $\bar\sA$.\qed
\end{lemme}

\begin{rem}\label{r3} 
 Supposons $K$ al\-g\'e\-bri\-que\-ment clos. 
Soit $A$ un objet
in\-d\'e\-com\-po\-sa\-ble. Alors par Schur, $\bar\sA(A,A)=K$, donc il y a
un unique relev\'e dans la
$K$-alg\`ebre locale $\sA(A,A)$. Si $B$ est un ind\'e\-com\-po\-sa\-ble non
isomorphe \`a $A$, alors $\pi_\sA(A)$ et
$\pi_\sA(B)$ sont des irr\'eductibles non isomorphes dans $\bar\sA$,
donc $\bar\sA(A, B)=0$. Il en d\'ecoule que $s$ existe et est unique si
les objets de $\sA$ sont tous ind\'e\-com\-po\-sa\-bles et deux \`a deux non
isomorphes.

Voici comment tirer de l\`a une preuve \'el\'ementaire (non
co\-ho\-mo\-lo\-gi\-que) de l'existence d'une section fonctorielle de
$\pi_\sB$ pour une petite ca\-t\'e\-go\-rie semi-primaire $\sB$
sur un corps alg\'ebriquement clos $K$, dont les
$K$-alg\`ebres d'endomor\-phis\-mes sont de dimension finie. Supposons
d'abord $\sA$ $K$-li\-n\'e\-ai\-re et pseudo-ab\'elienne. Il suit de cette
hypoth\`ese que tout objet de
$\sB$ est somme directe finie d'ind\'e\-com\-po\-sa\-bles. Il r\'esulte du \S
\ref{i1} ci-dessous que l'\'enonc\'e du th\'eor\`eme \ref{T1} est
invariant par
\'equivalence $K$-lin\'eaire de ca\-t\'e\-go\-ries et par passage aux
compl\'etions $K$-lin\'eaires et pseudo-ab\'eliennes. Il suffit alors
d'appliquer le raisonnement pr\'ec\'edent \`a
une sous-ca\-t\'e\-go\-rie pleine $\sA$ de
$\sB$ dont les objets forment un syst\`eme de repr\'esentants des classes
d'isomorphie d'objets ind\'e\-com\-po\-sa\-bles de $\sA$.
\end{rem}

\subsection{Variantes relatives}

\begin{sloppypar}
\begin{prop} \label{P3} Soit $T:\sB\to \sA$ un $K$-foncteur radiciel entre
deux
petites $K$-ca\-t\'e\-go\-ries de Wedderburn. Notons
$\bar \sA=\sA/\rad(\sA)$, $\bar \sB=\sB/\rad(\sB)$ et $\bar T:\bar
\sB\to \bar \sA$ le foncteur
induit. Soit $s_\sA:\bar\sA\to \sA$ (\resp $s_\sB:\bar\sB\to \sB$ )
une section de la
projection canonique. Alors il existe un
syst\`eme $(u_X)$ d'\'el\'ements de $1_{T(X)}+ \rad(\sA)(T(X),T(X))$
tel que $s_\sA\bar
T(h)=u_{X'}Ts_\sB(h)(u_X)^{-1}$ pour
tout $h\in \bar \sB(X,X')$.
\end{prop}
\end{sloppypar}

\prf 1) Posons $\sR=\rad(\sA)$. Commen\c cons par le cas o\`u
$\sR^2=0$. On a $Ts_\sB(h)-s_\sA\bar T(h)\in \rad(\sA)$
pour tout $h\in \bar \sB$. Ceci d\'efinit une d\'erivation de $\bar\sB$
\`a valeurs dans le ${\bar\sB}^e$-module $\rad(\sA)$ ($1$-cocycle de
Hoch\-schild). Cette
d\'erivation est int\'erieure ($1$-cobord de 
Hochschild): il existe une famille
$n_{X}=u_{X}-1_{T(X)}\in\sR(T(X), T(X))$ v\'erifiant,
pour tout $h\in \bar 
\sB(X,X')$,
\[Ts_\sB(h)-s_\sA\bar T(h) =n_{X'} Ts_\sB(h) - 
Ts_\sB(h)n_{X},\]
               d'o\`u l'assertion.

\begin{sloppypar}
2) Supposons 
maintenant seulement que $\sR=\rad(\sA)$ soit nilpotent
d'\'echelon 
$r+1$, et raisonnons par
r\'ecurrence sur $r$ (le cas $r=1$ \'etant 
acquis). Consid\'erons la
composition
$\sB\overset{T}{\to}\sA\to 
\sA^{(r)}=\sA/\sR$. En appliquant
l'hypoth\`ese de r\'ecurrence au 
foncteur
compos\'e, on trouve un syst\`eme $(u_X)$ d'\'el\'ements de 
$1_{T(X)}+\sR(T(X),T(X))$ tel que  
\[s_\sA\bar 
T(h)\equiv u_{X'}Ts_\sB(h)(u_X)^{-1}\pmod{\sR^{r}}.\]
\end{sloppypar}

Quitte
\`a 
remplacer $s_\sA$ par la section conjugu\'ee
par les $(u_A)^{-1}$, on 
se ram\`ene au cas o\`u
$s_\sA\bar
T(h)\equiv Ts_\sB(h)$ modulo 
$\sR^{r}$. On peut alors remplacer $\sB$
par $\bar\sB$ et $\sA$ par 
la
ca\-t\'e\-go\-rie $\tilde \sA$ introduite au pas $3)$ de la preuve 
du
th\'eor\`eme \ref{T1}. Comme le radical
de cette derni\`ere est de 
carr\'e nul, on conclut par le pas pr\'ec\'edent.

3) Dans le cas 
g\'en\'eral, on applique un raisonnement \`a la Zorn
comme au dernier 
pas de la preuve
du th\'eor\`eme \ref{T1}: bri\`evement, il existe 
une
sous-ca\-t\'e\-go\-rie pleine maximale $\sC$ de $\sB$ telle
que la 
proposition vaille pour $\sC$. On peut alors supposer que
$\sB$ 
contient un seul objet hors de
$\sC$, et de m\^eme que $A$ contient 
un seul objet hors de
$T(Ob(\sC))$. On constate comme
ci-dessus que 
le radical de cette nouvelle ca\-t\'e\-go\-rie $\sA$ est
alors 
(globalement) nilpotent.
\qed

Dans le paragraphe suivant, nous 
appliquerons la proposition \ref{P3} dans
le cas o\`u $\sA$ est munie 
d'une structure mo\-no\-\"{\i}\-dale $\bullet$ et
o\`u $T$ est le foncteur 
canonique
$\sA\boxtimes_K \sA\to \sA: X=(A, B)\mapsto A\bullet B$, 
suppos\'e
radiciel (ce qui revient \`a supposer que le radical est 
un
id\'eal mo\-no\-\"{\i}\-dal). Pour
$h=f\otimes g,s=s_\sA$, on a alors
 
\[s_\sA\bar
T(h) =s( f\bullet  g),Ts_\sB(h)=s(f)\bullet s( g) .\]

Le 
compl\'ement de Malcev au th\'eor\`eme de Wedderburn dit qu'on
peut 
choisir un
scindage de Wedderburn dont l'image contient une 
sous-alg\`ebre
semi-simple donn\'ee d'une $K$-alg\`ebre
de 
Wedderburn. En voici une variante \`a plusieurs 
objets.

\begin{cor}\label{cmalcev} Sous les
hypoth\`eses du 
th\'eor\`eme \ref{T1},
supposons don\-n\'ee en outre une 
sous-ca\-t\'e\-go\-rie
   $\sB$ {\rm s\'eparable} de $\sA$ 
(non
n\'ecessai\-re\-ment pleine), avec $Ob(\sB)=Ob(\sA)$. Alors il 
existe
une section $s$ de
$\pi_\sA$ v\'erifiant $(s\circ \pi_\sA 
)\vert \sB = id_{\sB}$.
\end{cor}

\prf Partant d'une
section de 
$\pi_\sA$ (dont l'existence est
assur\'ee par le th\'eor\`eme 
\ref{T1}), on la modifie en appliquant
la proposition 
pr\'e\-c\'e\-den\-te pour
obtenir une section fixant 
$\sB$.
\qed

\section{Sections 
mo\-no\-\"{\i}\-dales}\label{s6}

\subsection{Foncteurs 
mono\"{\i}daux}\label{9.1} Soient $(\sB,\top)$ et
$(\sA,\bullet)$ 
deux ca\-t\'e\-go\-ries
mo\-no\-\"{\i}\-da\-les. Rappelons qu'un 
\emph{foncteur mo\-no\-\"{\i}\-dal} est
la donn\'ee $(s,(\tilde s_{AB}))$ 
d'un foncteur
$s: \sB\to \sA$ et
d'isomor\-phis\-mes fonctoriels $\tilde 
s_{AB}: s(A)\bullet s(B)\allowbreak\to
s(A\top B)$ v\'erifiant 
une
compatibilit\'e aux contraintes d'associati\-vit\'e et d'unit\'e 
(on
\'ecrit souvent par abus $s$ au lieu de
$(s,(\tilde s_{AB}))$ 
pour abr\'eger). Les foncteurs mono\"{\i}daux
se composent, avec la 
r\`egle
$\widetilde{s's}_{AB}=s'(\tilde s_{AB})\circ \tilde 
s'_{s(A)s(B)}$.

Un \emph{mor\-phis\-me de foncteurs mono\"{\i}daux} est 
une
transformation naturelle des foncteurs sous-jacents,
compatible 
aux donn\'ees $\tilde s_{AB}, \tilde s'_{AB}$, et aux
unit\'es, \cf 
\cite[I.4.1.1]{saavedra}.
En d'autres termes, un tel mor\-phis\-me est la 
donn\'ee pour tout $A\in \sB$
d'un mor\-phis\-me $u_A:s(A)\to s'(A)$, 
naturel en $A$, tel que 
les
diagrammes
\begin{equation}\begin{CD}\label{eq4}
s(A)\bullet 
s(B)@>\tilde s_{A,B}>> s(A\top B)\\
@Vu_A\bullet u_B VV @Vu_{A\top 
B}VV\\
s'(A)\bullet s'(B)@>\tilde s'_{A,B}>> s'(A\top 
B)
\end{CD}\end{equation}
soient tous commutatifs, et tel que le 
mor\-phis\-me
\begin{equation}\begin{CD}\label{eq4'}
s({\un})={\un} 
@>u_{\un}>> s'({\un})={\un} \\\end{CD}\end{equation}
   soit 
l'identit\'e dans l'anneau commutatif $End({\un})$.

Un foncteur 
mo\-no\-\"{\i}\-dal $\pi: \sA\to \sB$ \'etant donn\'e, une {\it
section} de 
$\pi$ est un foncteur
mo\-no\-\"{\i}\-dal $s: \sB\to \sA$ tel que $\pi\circ 
s=id_\sB$
(\'egalit\'e de foncteurs mono\"{\i}daux). Deux
sections 
$s,s'$ sont dites {\it isomorphes} s'il existe un
isomor\-phis\-me de 
foncteurs mono\"{\i}daux
$s\Rightarrow s'$ tel que $\pi\circ 
s\Rightarrow \pi \circ s'$ soit
le foncteur identique.

Dans une 
$K$-ca\-t\'e\-go\-rie mo\-no\-\"{\i}\-dale, le bifoncteur $\bullet$ est 
 $K$-bilin\'eaire. Nous supposerons
toujours nos ca\-t\'e\-go\-ries 
mo\-no\-\"{\i}\-dales non nulles (\ie
$End(\un) \neq 
0$).

\subsection{Existence et ``unicit\'e" de sections 
mo\-no\-\"{\i}\-dales} Apr\`es
ces pr\'elimi\-nai\-res, nous 
pouvons
\'enoncer la version mo\-no\-\"{\i}\-dale du th\'eor\`eme 
de
Wedderburn \`a plusieurs objets. On suppose de nouveau que
$K$ est 
un corps.

\begin{sloppypar}
\begin{thm}\label{T2} Soit $(\sA, \bullet)$ une 
petite
$K$-ca\-t\'e\-go\-rie de Wedderburn mo\-no\-\"{\i}\-da\-le. 
Supposons que
les $End({\un})$-bimodules $\sA(A,B)$ soient 
commu\-tatifs,
\ie l'action \`a gauche co\"{\i}ncide avec 
l'action
\`a droite\footnote{Une ca\-t\'e\-go\-rie mo\-no\-\"{\i}\-dale 
v\'erifiant cette
hypoth\`ese est appel\'ee \emph{ca\-t\'e\-go\-rie
de 
Penrose} dans \cite{bru1}. C'est par exemple le cas 
si
$End({\un})=K$, ou en pr\'esence d'un tressage.}.
\\
Supposons 
que
$\rad(\sA)$ soit mo\-no\-\"{\i}\-dal, d'o\`u une structure 
mo\-no\-\"{\i}\-dale sur
$\bar
\sA=\sA/\rad(\sA)$. Notons $\pi_\sA:
\sA\to 
\bar \sA$ la projection canonique, vue comme foncteur
mo\-no\-\"{\i}\-dal. 
Alors:\\
a) Il existe une section mo\-no\-\"{\i}\-dale $s$ de $\pi_A$. 
\\
b) Toute autre section mo\-no\-\"{\i}\-dale de $\pi_\sA$ est
isomorphe 
\`a $s$.
\end{thm}
\end{sloppypar}

\subsection{Sorites}\label{sorital} Nous allons d\'emontrer le 
th\'eor\`eme \ref{T2} par
r\'eduction au cas strict. Pour cela, 
quelques sorites sont n\'ecessaires.

\begin{sorite}\label{so1} Soit 
$(s,\tilde s):(\sB,\top)\to (\sA,\bullet)$
un foncteur mo\-no\-\"{\i}\-dal, 
et soit $u:s\Rightarrow t$ un isomor\-phis\-me
naturel de $s$ sur un 
autre foncteur $t$. Alors il existe une unique
structure 
mo\-no\-\"{\i}\-dale $\tilde t$ sur $t$ faisant de $u$ un mor\-phis\-me
de 
foncteurs mono\"{\i}daux.
\end{sorite}

\prf D\'efinissons
\[\tilde 
t_{A,B}= u_{A\top B}\tilde s_{A,B} (u_A\bullet u_B)^{-1}\]
(\cf 
\eqref{eq4}). Il faut voir que $\tilde t$ est compatible avec 
les
contraintes d'associativit\'e et d'unit\'e. Cela
r\'esulte de 
calculs triviaux (d'o\`u le terme ``sorite"), ce que nous
laissons au 
lecteur sceptique le soin de v\'erifier.\qed

\begin{sorite}\label{so2} Soient $(\sA,\bullet),\bar\sA,\pi_\sA$ 
comme
dans le th\'eor\`eme \ref{T2}. Alors le th\'e\-o\-r\`e\-me 
\ref{T2}
r\'esulte de l'\'enonc\'e plus faible obtenu en rempla\c 
cant dans a) et
b) le mot ``section" par ``quasi-section", o\`u une 
quasi-section
mo\-no\-\"{\i}\-dale de
$\pi_\sA$ est par d\'efinition un 
foncteur mo\-no\-\"{\i}\-dal $s_0:\bar
\sA\to \sA$ tel que $\pi_\sA\circ 
s_0$ soit \emph{naturellement
isomorphe} \`a 
$Id_{\bar\sA}$.
\end{sorite}

\begin{sloppypar}
\prf Supposons 
donn\'ee une telle quasi-section mo\-no\-\"{\i}\-da\-le
$(s_0,\tilde 
s_0)$. Soit $\bar u:Id_{\bar \sA}\iso \pi_\sA s_0$ un
isomor\-phis\-me 
naturel. Pour tout $A\in \sA$, choisissons un \'el\'ement
$u_A\in 
\sA(A,s_0(A))$ se projetant sur $\bar u_A$, avec
$u_{\un}=1_{\un}$. 
D\'efinissons un
foncteur $s:\bar\sA\to \sA$ par les 
formules:
\begin{align*}
s(A)&=A\\
s(f)&=u_B^{-1} s_0(f) 
u_A
\end{align*}
pour $f\in \bar\sA(A,B)$. Alors $s$ est une section 
de $\pi_\sA$
naturellement isomorphe \`a $s_0$ via $u=(u_A)$, et le 
sorite \ref{so1}
montre qu'il lui est associ\'e une unique structure 
mo\-no\-\"{\i}\-dale faisant
de $u$ un mor\-phis\-me de foncteurs 
mono\"{\i}daux. Ceci fournit la partie a) du th\'eor\`eme
\ref{T2}, 
et le point concernant la partie b) est trivial.\qed 

\end{sloppypar}

\begin{sorite}\label{so3} Soient 
$(\sA,\bullet),\bar\sA,\pi_\sA$ et
$(\sB,\top),\bar\sB,\pi_\sB$ comme 
dans le th\'e\-o\-r\`e\-me \ref{T2}.
Soit
$(f,\tilde f):(\sB,\top)\to 
(\sA,\bullet)$ un $K$-foncteur mo\-no\-\"{\i}\-dal
tel que $f$ soit une 
\'equivalence de ca\-t\'e\-go\-ries. Alors $f$ est
radiciel, et les parties 
a) et b) du th\'eor\`eme \ref{T2} sont vraies
pour $\sA$ si et 
seulement si elles sont vraies pour $\sB$.
\end{sorite}

\prf Le fait 
que $f$ est radiciel est clair (et ne d\'epend pas des
structures 
mo\-no\-\"{\i}\-dales). Le reste de l'\'enonc\'e r\'esulte du 
sorite
\ref{so2}, en remarquant que les conditions dudit sorite 
sont
manifestement invariantes par \'equivalence mo\-no\-\"{\i}\-dale de 
ca\-t\'e\-go\-ries
mo\-no\-\"{\i}\-dales.\qed

\begin{sloppypar}
\subsection{Cat\' 
egories mo\-no\-\"{\i}\-dales strictes}\label{ml} Une 
petite 
ca\-t\'e\-go\-rie
mono\"{\i}\-dale
$(\sB,\top)$ est dite {\it stricte} si 
l'ensemble des objets muni de
la loi $\top$ est un
mono\"{\i}de (\ie 
si $\top$ est strictement associative).
\end{sloppypar}

Un foncteur 
mo\-no\-\"{\i}\-dal $(s,(\tilde s_{AB})): (\sA,\bullet)\to 
(\sB, \top)$ 
est dit \emph{strict} si
$s(A)\bullet s(B)= s(A\top
B)$ et si
$\tilde 
s_{AB}$ est l'identit\'e pour tout couple $(A,B)$. C'est par
exemple 
le cas du foncteur identique.

Dans le cas strict, la situation se 
simplifie pour les
mor\-phis\-mes de foncteurs mono\"{\i}daux: si 
$u:s\Rightarrow
s'$ est un tel mor\-phis\-me, alors il v\'erifie 
l'identit\'e
\[u_{A\bullet B}=u_A\top u_B\]
pour tout couple d'objets 
$(A,B)$.

\begin{rem} Soit $\sA^\sA$ la ca\-t\'e\-go\-rie des endofoncteurs 
additifs de $\sA$. Munie de la
composition, c'est une ca\-t\'e\-go\-rie 
mo\-no\-\"{\i}\-dale stricte. Pour toute ca\-t\'e\-go\-rie mo\-no\-\"{\i}\-dale, le 
foncteur
$A\mapsto A\bullet ?$ de $\sA $ vers $\sA^\sA$ admet une 
structure mo\-no\-\"{\i}\-dale d\'efinie par les 
contraintes
d'associativit\'e et d'unit\'e. {\it Ce foncteur est 
strict si et seulement si $(\sA, \bullet)$ est
mo\-no\-\"{\i}\-dale 
stricte}.
\end{rem} 

\medskip Si $\sB$ est mo\-no\-\"{\i}\-dale stricte, 
il en est
alors de m\^eme de $\sB^\oplus$: en effet,
rappelons que 
les objets de
$\sB^\oplus$ sont des suites finies $\underline A= 
(A_1,A_2,\dots,
A_r)$ d'objets de $\sB$, not\'ees $A_1\oplus 
A_2
\oplus \dots \oplus A_r $ (\cf \S\ref{rap}); les composantes 
de
$\underline A \top \underline A'$ sont
par d\'efinition les 
$A_i\top A'_j$, num\'erot\'ees selon l'ordre
lexi\-cographique. Il 
est clair qu'on
obtient ainsi une ca\-t\'e\-go\-rie mo\-no\-\"{\i}\-dale stricte. 
Si $\sB$ est $K$-lin\'eaire, les
foncteurs \'evidents $\displaystyle 
\sB^\oplus
\begin{smallmatrix}u^\oplus\\\rightarrow\\\leftarrow\\v^\oplus
\end{smallmatrix}
\sB$ (formation du mot \`a une lettre, \resp 
parenth\'esage
canonique) sont sous-jacents \`a des 
foncteurs
mono\"{\i}daux stricts (vis-\`a-vis de $\bullet$).

Tout 
foncteur mo\-no\-\"{\i}\-dal $s:\,\sB\to \sA$, il induit, composante par 
 composante, un foncteur mo\-no\-\"{\i}\-dal
$s^\oplus: \sB^\oplus\to 
\sA^\oplus$, strict si $s$ l'est.

\medskip Une petite ca\-t\'e\-go\-rie 
$K$-lin\'eaire $\sC$ sera dite {\it 
strictement $K$-lin\'eaire} si 
son biproduit
$\oplus$ (et l'objet nul) fait de $\sC$ une ca\-t\'e\-go\-rie 
 mono\"{\i}\-dale stricte. C'est le cas de
$\sB^\oplus$, pour toute 
petite $K$-ca\-t\'e\-go\-rie $\sB$. Dans une 
ca\-t\'e\-go\-rie strictement 
$K$-lin\'eaire, les
sommes directes finies $A_1\oplus \dots \oplus 
A_n$ sont d\'efinies 
sans 
ambig\"uit\'e.

\begin{sloppypar}
\subsubsection{La construction de 
Mac Lane \protect\cite[XI 3]{maclane}}\label{constrM} 
\`A toute 
petite ca\-t\'e\-go\-rie mono\"{\i}\-dale $(\sA,\bullet)$, on
associe 
une petite ca\-t\'e\-go\-rie mono\"{\i}\-dale stricte
$(\sA^{str},\top)$, 
et deux foncteurs mono\"{\i}daux $\displaystyle 
\sB
\begin{smallmatrix}u\\\rightarrow\\\leftarrow\\v\end{smallmatrix}
\sA$, avec $u$ strict, tels que le compos\'e $uv$ soit 
le
foncteur mo\-no\-\"{\i}\-dal identique de
$\sA$ (les foncteurs 
sous-jacents \`a $u$ et $v$ sont 
donc
quasi-inverses).
\end{sloppypar}

Rappelons cette construction: 
on prend pour objets de $\sA^{str}$
les mots (finis) dont les lettres 
sont les objets de $\sA$. Les mor\-phis\-mes
entre mots sont les 
mor\-phis\-mes entre les objets de $\sA$
correspondants obtenus par 
pa\-ren\-th\'e\-sa\-ge
canonique (toutes les pa\-ren\-th\`e\-ses 
commencent au d\'ebut du mot). Le
produit $\top$, souvent omis de la 
notation, est donn\'e au niveau des
objets par la concat\'enation des 
mots; $v$ est donn\'e par la
``formation de mots d'une lettre", 
tandis
que $u$ est donn\'e par le ``pa\-ren\-th\'e\-sa\-ge canonique" 
({\it
e.g.} $u(ABC)= (A\bullet B)\bullet
C)$). L'unit\'e est le mot 
vide: $u(\emptyset)=\un$.

C'est une $K$-ca\-t\'e\-go\-rie si $\sA$ 
l'\'etait, et $u,v$ sont
alors des $K$-foncteurs.

Noter que lorsque 
$A$ est $K$-lin\'eaire, la construction 
$\sA^\oplus$ est un cas 
parti\-cu\-lier de cette
construction, en prenant 
$\bullet=\oplus$.

\medskip Si $s: \sB\to \sA$ est un foncteur, on 
lui
associe canoni\-quement comme suit un foncteur mo\-no\-\"{\i}\-dal 
strict
$s^{str}: \sB^{str}\to \sA^{str}$: on associe au 
mot
$w=B_1\dots B_n$ de $\sB^{str}$ le mot
$s^{str}(w)=s(B_1)\dots 
s(B_n)$ de
$\sA^{str}$ ($s^{str}(\emptyset)=\emptyset$).

On associe 
\`a $f\in
\sB^{str}(w,w')= \sB(u(w) ,u(w'))$ le 
mor\-phis\-me
$s^{str}(f)=  s(f) 
\in
\sA(u(s^{str}(w)),u(s^{str}(w')))=\sA^{str}(s^{str}(w),s^{str}(w') 
)$.

Si $s$ est sous-jacent \`a un foncteur mo\-no\-\"{\i}\-dal strict, on 
a
$us^{str}=su$ (donc $s= u_\sA \circ s^{str}\circ 
v_\sB$).

\medskip
Vu la construction de Mac Lane et le sorite 
\ref{so3}, pour d\'emontrer
le th\'e\-o\-r\`e\-me \ref{T2}, on peut 
supposer $\sA$ stricte. Le
th\'eor\`eme
\ref{T2} r\'esulte alors de 
la proposition plus pr\'ecise 
suivante:

\begin{prop}\label{p5}
Soient 
$(\sA,\bullet),\bar\sA,\pi_\sA$ comme
dans le th\'eor\`eme \ref{T2}. 
Si en outre $\sA$ est mo\-no\-\"{\i}\-dale stricte,\\
a) Il existe une 
section mo\-no\-\"{\i}\-dale stricte $s$ de $\pi_\sA$.\\
b) Toute autre 
section mo\-no\-\"{\i}\-dale stricte de $\pi_\sA$ est
isomorphe \`a 
$s$.\\
c) Toute section mo\-no\-\"{\i}\-dale de $\pi_\sA$ est naturellement 
isomorphe
\`a une section mo\-no\-\"{\i}\-dale stricte de 
$\pi_\sA$.
\end{prop}

\subsection{R\'eductions} Commen\c cons par 
une r\'eduction du probl\`eme
en deux 
lemmes.

\begin{lemme}\label{L8'} La proposition \ref{p5} d\'ecoule 
de l'\'enonc\'e
correspondant pour la ca\-t\'e\-go\-rie mo\-no\-\"{\i}\-dale 
stricte $K$-lin\'eaire
$\sA^\oplus$. \end{lemme}

\prf Soit $S$ une 
section
mo\-no\-\"{\i}\-dale stricte de $\pi_{\sA^\oplus}$. Le foncteur 
mo\-no\-\"{\i}\-dal
strict (pleinement) fid\`ele 
$\sA\overset{v^\oplus}{\to}\sA^\oplus$
induit un diagramme 
commutatif
\[\begin{CD}
\bar\sA@>\bar v>>\bar\sA^\oplus \\
@VsVV 
@VSVV\\
\sA @>v>> \sA^\oplus .
\end{CD} \]
R\'eciproquement, si $s$ 
est donn\'ee, on peut prendre $S= s^\oplus$. \qed

En identifiant 
$(\bar\sA^{str})^\oplus$ et 
$\overline{(\sA^{str})^\oplus}$, ce 
lemme nous ram\`ene \`a prouver la
proposition \ref{p5} dans le 
cas
strictement $K$-lin\'eaire, strictement 
mo\-no\-\"{\i}\-dal.

\begin{lemme}\label{L8''} Supposons que $\sA$ soit 
strictement 
$K$-lin\'eaire, et mo\-no\-\"{\i}\-da\-le stricte.
 
\noindent Soit $s$ une section fonctorielle
$K$-lin\'eaire de $\pi_\sA$ (vue 
comme simple $K$-foncteur).

Supposons que $s( 
f\bullet
g)=s(f)\bullet s( g)$ pour tout couple d'endomor\-phis\-mes 
$(f\in
\bar\sA^\oplus(A,A),g\in \bar\sA^\oplus(B,B))$, et que 
 $s(1_{\un})=1_{\un}$. Alors $s^\oplus$
d\'efinit une unique section 
mo\-no\-\"{\i}\-dale stricte
    de $\pi_{\sA^\oplus}$.
\end{lemme}

\prf 
Il s'agit de v\'erifier l'identit\'e $s^\oplus( 
f\bullet
g)=s^\oplus(f)\bullet s^\oplus( g)$ pour tout couple 
$(f\in
\bar\sA^\oplus(w_1,w'_1),g\in \bar\sA^\oplus(w_2,w'_2))$, et 

$s^\oplus(\emptyset)=1_{\emptyset}$ (ce second point
est 
imm\'ediat).

Nous utilisons transitoirement la notation $w$ des mots 
de \ref{ml}, 
dans le cas o\`u la structure mo\-no\-\"{\i}\-dale
est 
donn\'ee par $\oplus$; au lieu de $u(w)$, la somme directe (dans 
la 
ca\-t\'e\-go\-rie strictement $K$-lin\'eaire
$\sA$) des lettres de $w$ sera 
not\'ee $\pointplus{w}$ 
($\pointplus{\emptyset}=\un$). Via les 
identifications
canoniques
\[\sA^\oplus(w,w')= 
\sA(\pointplus{w},
\pointplus{w}'),\]
on peut \'ecrire 
$s^\oplus(f)=s(f), s^\oplus(g)=s(g), 
s^\oplus(f\bullet g)=s(f\bullet 
g)$.

Les identifications ci-dessus sont compatibles \`a $\bullet$ 
(
$\pointplus{\widehat{{w_1}\bullet{w_2}}}=\pointplus{w}\bullet
\point 
plus{w}'$).

On peut donc \'ecrire $s^\oplus(f)\bullet s^\oplus( 
g)=s(f)\bullet s( g)$
(en revanche, on prendra garde que $s$ n'est 
pas suppos\'ee stricte 
vis-\`a-vis de $\oplus$,\; \ie on ne 
suppose
pas que le mor\-phis\-me canonique 
$\pointplus{\widehat{s^\oplus(w)}}\to 
s(\pointplus{w})\;$  soit 
l'\'egalit\'e).

L'hypoth\`ese du lemme se traduit donc par 
l'identit\'e $s^\oplus( f\bullet
g)=s^\oplus(f)\bullet s^\oplus( g)$ 
dans \[
\sA^\oplus(w_1\bullet w_2,w_1\bullet w_2)= 
\sA(\pointplus{w_1}\bullet 
\pointplus{w_2},\pointplus{w_1}\bullet
\pointplus{w_2})\] pour tout 
couple d'endomor\-phis\-mes $(f\in
\bar\sA^\oplus(w_1,w_1),g\in 
\bar\sA^\oplus(w_2,w_2))$.

\begin{sloppypar}
Or $s^\oplus$ est strict vis-\`a-vis de 
$\oplus$ par 
construction (par contraste avec $s$).
Cela permet de 
ramener l'identit\'e ci-dessus pour tout couple de 
mor\-phis\-mes 
$(f\in
\bar\sA^\oplus(w_1,w'_1),g\in \bar\sA^\oplus(w_2,w'_2))$ au 
cas d'un couple
d'endomor\-phis\-mes en rempla\c cant $w_1$ et $w'_1$ par 
$w_1\oplus 
w'_1$, et $w_2$ et $w'_2$ par $w_2\oplus w'_2$.
\qed
\end{sloppypar}

Ce 
lemme nous ram\`ene \`a traiter d'endomor\-phis\-mes, plut\^ot 
que
d'homomor\-phis\-mes.

\begin{sloppypar}
\medskip\noindent
{\bf 
Notation abr\'eg\'ee.} Nous \'ecrirons $\sA(A)$ (\resp $\sR(A)$)
au 
lieu de
$\sA(A,A)$ (\resp $\sR(A,A)$) dans ce
paragraphe. Pour un 
objet $A\in\sA$ et $a,b\in \sA(A)$, nous noterons
comme d'habitude 
$[a,b]=ab-ba$ le commutateur de $a$ et $b$, et $ad(a)$
l'application 
$b\mapsto [a,b]$.

\subsection{Au c\oe ur du probl\`eme: 
d\'emonstration de la proposition
\protect\ref{p5}}  \allowbreak 
Comp\-te tenu des lem\-mes
pr\'e\-c\'e\-dents, on suppose d\'esormais 
que
$\sA$ est une petite ca\-t\'e\-go\-rie, strictement 
$K$-lin\'eaire,
mo\-no\-\"{\i}\-dale stricte (notons que $\bar\sA$ a la 
m\^eme
propri\'et\'e). On pro\-c\`e\-de par \'etapes. Nous 
commencerons par
d\'e\-mon\-trer c), cette d\'e\-mons\-tra\-tion 
\'etant la plus simple des
trois et servant de prototype \`a celles 
de a) et b).
\end{sloppypar} 

\subsubsection{D\'ebut de la preuve de 
c)} Soit $s$ une section
$K$-fonctorielle mo\-no\-\"{\i}\-da\-le 
de
$\pi_A$. On a donc des isomor\-phis\-mes
\[\begin{CD}
A\bullet 
B=s(A)\bullet s(B)@>\tilde s_{A,B}>> s(A\bullet B)=A\bullet 
B.
\end{CD}\]

Les $\tilde s_{A,B}$ appartiennent \`a $1+\sR(A\bullet 
B)$ et
v\'erifient la relation de coh\'erence (dans $1+\sR(A\bullet 
B\bullet C)$)
\begin{equation}\label{coh1}
\tilde s_{A\bullet 
B,C}\circ (\tilde s_{A,B}\bullet
1_C)=\tilde s_{A,B\bullet C}\circ 
(1_A\bullet \tilde s_{B,C}).
\end{equation}

Nous allons 
montrer:

\begin{lemme}\label{l4} Il existe une famille 
double
$(u_A^{(r)})_{A\in\sA,r\ge 1}$, avec $u_A^{(r)}\in1+\sR(A)$, 
ayant
les propri\'et\'es suivantes:
\begin{thlist}
\item On a 
$u_A^{(r+1)}\equiv u_A^{(r)}\pmod{\sR(A)^r}$ pour tout
$A\in\sA$ et 
tout $r\ge 1$.
\item En posant
\[s_r(f)=u_B^{(r)}\circ s(f)\circ 
(u_A^{(r)})^{-1}\]
pour $A,B\in \sA$ et $f\in
\bar\sA(A,B)$, on 
a
\[(\tilde s_r)_{A,B}=u_{A\bullet B}^{(r)}\circ\tilde 
s_{A,B}\circ
(u_A^{(r)}\bullet u_B^{(r)})^{-1}\in 1+\sR(A\bullet 
B)^r\]
pour tout couple d'objets 
$(A,B)$.
\end{thlist}
\end{lemme}

Supposons le lemme \ref{l4} 
d\'emontr\'e. Pour tout $A$, choisissons un
entier $r(A)$ tel que 
$\sR(A)^{r(A)}=0$. Posons
\[u_A=u_A^{(r(A))}=u_A^{(r(A)+1)}=\dots\in 
1+\sR(A)\]
et
\[s'(f)=u_B\circ s(f)\circ (u_A)^{-1}\]
pour $A\in \sA$ 
et $f\in \bar\sA(A)$. Pour $A,B\in \sA$, on a donc
\[\tilde 
s'_{A,B}=(\tilde s_r)_{A,B}\]
d\`es que 
$r\ge
\sup(r(A),r(B),r(A\bullet B))$. De plus, le membre de droite de 
cette
\'egalit\'e est \'egal \`a $1_{A\bullet B}$. Ainsi, la
section 
$s'$ v\'erifie les conclusions de c).

Nous allons d\'emontrer le 
lemme \ref{l4} par r\'ecurrence sur $r$,
partant du cas trivial $r=1$ 
o\`u l'on prend $u_A^{(1)}=1_A$ pour tout
$A$.

\subsubsection{Une 
relation de cocycle, I} Supposons $r\ge
1$ et trouv\'ee une famille 
$u_A^{(r)}$ v\'erifiant la condition (ii) du
lemme \ref{l4}. Quitte 
\`a remplacer $s$ par $s_r$, on peut supposer
$s_r=s$, $u^{(r)}=1$. 
On a donc $\tilde s_{A,B}\in 1+\sR(A\bullet B)^r$
pour tout couple 
d'objets $(A,B)$. Posons $n_{A,B}=1_{A\bullet 
B}-\tilde
s_{A,B}\in\sR(A\bullet B)^r$. La relation \eqref{coh1} se 
r\'e\'ecrit
alors
\begin{multline*}
(1_{A\bullet B\bullet 
C}-n_{A\bullet B,C})\circ ((1_{A\bullet
B}-n_{A,B})\bullet 
1_C)\\
=(1_{A\bullet B\bullet C}-n_{A,B\bullet C})\circ 
((1_A\bullet
(1_{B\bullet 
C}-n_{B,C}))
\end{multline*}
d'o\`u
\begin{equation}\label{coc1}
n_{A\bullet B,C}+n_{A,B}\bullet 1_C\equiv n_{A,B\bullet C} 
+1_A\bullet
n_{B,C}\pmod{\sR(A\bullet B\bullet 
C)^{r+1}}.
\end{equation}

Nous cherchons une 
famille
$(u_A^{(r+1)})_{A\in \sA}$ v\'erifiant:
\begin{thlist}
\item 
$u_A^{(r+1)}\equiv 1\pmod{\sR(A)^r}$ pour tout
$A\in\sA$.
\item 
$u_{A\bullet B}^{(r+1)}\circ\tilde s_{A,B}\circ
(u_A^{(r+1)}\bullet 
u_B^{(r+1)})^{-1}\in 1+\sR(A\bullet B)^{r+1}$
pour tout couple 
d'objets $(A,B)$.
\end{thlist}

Supposons le probl\`eme r\'esolu et 
posons $m_A=1_A-u_A^{(r+1)}\in
\sR(A)^r$. On a alors 
l'identit\'e
\begin{multline*}
(1_{A\bullet B}-m_{A\bullet B})\circ 
(1_{A\bullet
B}-n_{A,B})\circ ((1_A-m_A)\bullet 
(1_B-m_B))^{-1}\\
\equiv 1\pmod{\sR(A\bullet 
B)^{r+1}}
\end{multline*}
soit 
encore
\begin{equation}\label{cob1}
n_{A,B}\equiv m_A\bullet 1_B 
+1_A\bullet m_B -m_{A\bullet
B}\pmod{\sR(A\bullet 
B)^{r+1}}.
\end{equation}

R\'eciproquement, l'existence d'une 
famille $(m_A)$ v\'erifiant la
congruence \eqref{cob1} est clairement 
\'equi\-va\-len\-te\`a l'existence
d'une famille $(u_A^{(r+1)})$ comme 
ci-dessus.

\subsubsection{Une relation de cocycle, II} La congruence 
\eqref{coc1}
fait penser \`a une relation de $2$-cocycle, et la 
congruence \eqref{cob1}
\`a une relation de $2$-cobord. Pour donner 
corps \`a ces impressions,
introduisons encore quelques 
notations.

Si $w=A_1A_2\dots A_m$ est un mot form\'e d'objets de 
$\sA$, notons
$\point{w}$ l'objet $A_1\bullet A_2\bullet\dots \bullet 
A_m$ de $\sA$ ou
$\bar\sA$ (notation plus commode que celle $u(w)$ 
utilis\'ee plus haut);
notons les 
formules
\begin{align*}
\point{(\widehat{\point{w_1}\point{w_2}})}&
=\point{\widehat{w_1w_2}}\\
\point{\emptyset}&=\un.
\end{align*}

Pour 
tout $A\in\sA$, posons $\sR(A)^{[r]}=\sR(A)^r/\sR(A)^{r+1}$. On 
peut
faire op\'erer les mots \`a gauche et \`a droite sur le 
$K$-espace
vectoriel
\[\sR^{[r]} = \prod_w \sR(\point{w})^{[r]}\]
 
de la mani\`ere suivante:
\[(w_0\cdot 
m)(w)=\begin{cases}1_{\point{w_0}}\bullet
m(w')&\text{si $w$ est de 
la forme $w_0w'$}\\
0& \text{sinon}
\end{cases}\]
\[(m\cdot 
w_1)(w)=
\begin{cases}
m(w')\bullet 1_{\point{w_1}} & \text{si $w$ 
est de la forme 
$w'w_1$}\\
0&\text{sinon.}
\end{cases}\]

(L'expression $m(w)$ 
d\'esigne la $w$-composante de $m$).

\begin{sloppypar}
Il n'est pas 
difficile de voir que cette
op\'eration est bien d\'efinie, et 
qu'\'etendue par
$K$-lin\'earit\'e, elle munit $\sR^{[r]}$ d'une 
structure de
$K\langle Ob(\sA)\rangle$-bimodule (o\`u
$K\langle 
Ob(\sA)\rangle$ d\'esigne la $K$-alg\`ebre libre de base
les objets 
de
$\sA$).
\end{sloppypar}

\subsubsection{Une relation de cocycle, 
III -- fin de la preuve de c)}
Nous som\-mes maintenant en mesure 
d'interpr\'eter \eqref{coc1} comme
relation de $2$-cocycle de 
Hochschild pour l'al\-g\`e\-bre $K\langle
Ob(\sA)\rangle$.

En effet, 
reprenons la famille $(n_{A,B})$ ci-dessus. Avec
les notations 
pr\'e\-c\'e\-den\-tes, on a une famille induite
\[\bar n_{AB}\in 
\sR(AB)^{[r]}.\]

On d\'efinit alors une $2$-cocha\^{\i}ne 
de
Hochschild $\bar n$ pour
$K\langle Ob(\sA)\rangle$ \`a va\-leurs 
dans le $K\langle
Ob(\sA)\rangle$-bimodule $\sR^{[r]}$,
        en 
associant \`a tout couple
de mots non vides $(w_0,w_1)$ 
l'\'el\'ement
$\bar n_{w_0,w_1}\in \sR^{[r]}$ d\'efini par
\[\bar 
n_{w_0,w_1}(w) =\begin{cases}
\bar n_{\point{w_0},\point{w_1}}& 
\text{si $w_0
w_1=w$,}\\
0&\text{sinon.}
\end{cases}\]

Calculons le 
bord de Hochschild de $\bar n$:
\[(d^2(\bar n))_{w_0,w_1,w_2}= 
w_0\cdot \bar n_{w_1,w_2} -
\bar n_{w_0w_1,w_2} + \bar n_{w_0,w_1w_2} 
-\bar n_{w_0,w_1}\cdot w_2.\]

Sa $w$-composante est nulle, par 
d\'efinition, si $w_0w_1w_2\neq w$.
Si $ w_0w_1w_2= w$, elle est 
\'egale \`a
\[
1_{\point{w_0}}\bullet
{\bar 
n}_{\point{w_1},
\point{w_2}}-{\bar 
n}_{\point{w_0}\bullet
\point{w_1},
\point{w_2}}
+{\bar 
n}_{\point{w_0},\point{w_1}\bullet
\point{w_2}}
-   {\bar 
n}_{\point{w_0},\point{w_1}}\bullet
1_{\point{w_2}}
=0\text{ 
d'apr\`es \eqref{coc1}.}
\]

Donc $\bar n$ est un 
$2$-cocycle. Comme $K\langle Ob(\sA)\rangle$ est
une $K$-alg\`ebre 
\emph{libre}, c'est un \emph{$2$-cobord}, \cf
\cite[ch. IX, ex. 2 ou 
ch. X.5]{ce}. En d'autres termes, il existe une
famille 
d'\'el\'ements $\bar m_w\in \sR^{[r]}$ telle que l'on 
ait
l'identit\'e
\[\bar n_{w_0,w_1} =w_0\cdot \bar m_{w_1}+\bar 
m_{w_0}\cdot w_1 -\bar
m_{w_0w_1}.\]

Soient $A,B\in\sA$. Pour 
$w_0=A$, $w_1=B$, la composante
de cette identit\'e dans 
$\sR(w_0w_1)^{[r]}$ n'est autre que
\[{\bar n}_{AB}=1_A\bullet{\bar 
m}_B  -{\bar m}_{A\bullet B}
+{\bar m}_A\bullet 1_B.\]

En relevant 
les $\bar m_A$ en des $m_A\in \sR(A)^r$, on
obtient la famille 
$(u_A^{(r+1)}=1_A-m_A)$ cherch\'ee.

\subsubsection{D\'ebut de la 
preuve de a)} Soit $s$ une section
$K$-fonctorielle de
$\pi_A$. 
L'existence en est garantie par le th\'eor\`eme \ref{T1}.

   Comme 
les $End({\un})$-bimodules $\sA(A,B)$ sont
commutatifs par 
hypoth\`ese, il est loisible de remplacer
chaque $s(f)$ par 
$s(1_{\un})^{-1}\circ s(f)$, et donc de
supposer 
$s(1_{\un})=1_{\un}$.

En outre, d'apr\`es la proposition \ref{P3} 
appliqu\'ee \`a
$\sB=\sA\boxtimes_K
\sA$, il existe une famille $(u_{A,B}\in
               1_{A\bullet B}+\sR(A\bullet B))$ v\'erifiant 
,
pour tout $(f,g)\in  \bar\sA(A)\times \bar\sA(B)$:
\[s( f\bullet 
g)= u_{A,B}(s(f)\bullet s(g))u_{A,B}^{-1}.\]
    Nous allons 
montrer:

\begin{lemme}\label{l5} Il existe une famille 
double
$(u_A^{(r)})_{A\in\sA,r\ge 1}$, avec $u_A^{(r)}\in1+\sR(A)$, 
ayant
les propri\'et\'es suivantes:
\begin{thlist}
\item On a 
$u_A^{(r+1)}\equiv u_A^{(r)}\pmod{\sR(A)^r}$ pour tout
$A\in\sA$ et 
tout $r\ge 1$.
\item En posant
\[s_r(f)=u_B^{(r)}\circ s(f)\circ 
(u_A^{(r)})^{-1}\]
pour $A,B\in \sA$ et $f\in
\bar\sA(A,B)$, on 
a
\[s_r(f\bullet g)\equiv s_r(f)\bullet 
s_r(g)\pmod{\sR(A\bullet
B)^{r}}\]
pour tout couple d'objets $(A,B)$ 
et 
tout
$(f,g)\in\bar\sA(A)\times\bar\sA(B)$.
\end{thlist}
\end{lemme}

L 
e lemme \ref{l5} implique a) de la m\^eme mani\`ere que le 
lemme
\ref{l4} impliquait c). On proc\`ede de nouveau par 
r\'ecurrence sur $r$.

\subsubsection{Une relation de cocycle, I} 
Supposons $r\ge
1$ et trouv\'ee une famille $u_A^{(r)}$ v\'erifiant 
la condition (ii) du
lemme \ref{l5}. Quitte \`a remplacer $s$ par 
$s_r$, on peut supposer
$s_r=s$, $u^{(r)}=1$. 
Posons
\[n_{A,B}=u_{A,B}-1_{A\bullet B}.\]

On a 
alors
\[ad(n_{A,B})(s(f)\bullet s(g))\in \sR(A\bullet B)^r\]
pour 
tout $(f,g)\in \bar\sA(A)\times\bar\sA(B)$.

Soit $h \in \bar\sA(C)$. 
Nous
allons calculer de deux fa\c cons
$s( f\bullet  g\bullet 
h)$
mo\-du\-lo $\sR(A\bullet B\bullet C)^{r+1}$, compte-tenu 
de
l'associativit\'e stricte de $\bullet$.

On a, modulo 
$\sR(A\bullet B\bullet C)^{r+1}$,
\begin{multline*}
s(( f\bullet 
g)\bullet  h)
= s( f\bullet  g)\bullet s( h) +[n_{A\bullet B ,C 
},(s(f)\bullet
s(g))\bullet s(h)] \\
= (s( f)\bullet s( g))\bullet s( 
h) +[n_{A ,B },s(f)\bullet
s(g))]\bullet
s(h) \\+ [n_{A \bullet B ,C 
},(s(f)\bullet s(g))\bullet s(h)]\\
              = (s( f)\bullet s( 
g))\bullet s( h)+[n_{A , B }\bullet
1_{C}+n_{A\bullet B ,C 
},(s(f)\bullet
s(g))\bullet s(h)]
\end{multline*}
et de 
m\^eme
\begin{multline*} s_r( f\bullet ( g\bullet  h))= \\
s( 
f)\bullet (s( g)\bullet
s( h))+[1_A \bullet n_{B , C }+n_{A , B 
\bullet C },s(f)\bullet
(s(g)\bullet s(h))].
\end{multline*}

Puisque 
$\sA$ est mo\-no\-\"{\i}\-dale stricte, on en
d\'eduit 
que
\begin{multline}\label{coc2}
ad(1_A\bullet
{n}_{B,C}-{n}_{A\bullet 
B, C}  +{n}_{A,B\bullet C} - {n}_{A,B}\bullet 1_C)
(s( f\bullet 
g\bullet  h))\\
\in\sR(A\bullet B\bullet C)^{r+1}
\end{multline}
pour 
tous endomor\-phis\-mes $f,g,h$ comme ci-dessus.

D'autre part, si 
$(u^{(r+1)}_A=1+m_A)_{A\in\sA}$ est une famille
d'\'el\'ements comme 
dans la conclusion du lemme \ref{l5} (donc $m_A\in
\sR(A)^r$ pour 
tout $A$), un calcul analogue donne 
l'identit\'e
\begin{equation}\label{cob2}
ad(n_{A,B}-m_A\bullet 
1_B-1_A\bullet
m_B+m_{A\bullet B})(s(f)\bullet s(g))\in \sR(A\bullet 
B)^{r+1}
\end{equation}
pour $(f,g)\in\bar\sA(A)\times \bar\sA(B)$; 
r\'eciproquement, la donn\'ee
d'\'el\'ements $m_A$ v\'erifiant 
\eqref{cob2} implique l'existence d'une
famille $u^{(r+1)}$ comme 
dans l'\'e\-non\-c\'e du lemme \ref{l5}.

\subsubsection{Une 
relation de cocycle, II}\label{co2} Pour interpr\'eter
\eqref{coc2} 
comme une relation de
$2$-cocycle et \eqref{cob2} comme une relation 
de cobord, on proc\`ede
comme ci-dessus en introduisant un
$K\langle 
Ob(\sA)\rangle$-bimodule appropri\'e.

Notons que l'action de $\sA$ 
sur $\sR^{[r]}=\sR^r/\sR^{r+1}$ se factorise
en une action de $\bar 
\sA$. En particulier, on peut se dispenser
d'\'ecrire la section $s$ 
dans \eqref{coc2} et \eqref{cob2}.

Pour tout mot $w=A_1\dots A_m$ en 
les objets de $\sA$, notons
\[\bar\sA(w)=\bar\sA({A_1})\bullet
\dots 
\bullet \bar\sA({A_m})\subset \bar\sA(\point{w}).\]

On a 
$\bar\sA(ww')=\bar\sA(w)\bar\sA(w')$.

\begin{sloppypar}
Notons 
ensuite
$C(w)^{[r]}$ le \emph{centralisateur de $\bar\sA(w)$ 
dans
$\sR(\point{w})^{[r]}\;$}:
\[C(w)^{[r]}=\{x\in 
\sR(\point{w})^{[r]}\mid ad(a)(x)=0\,\forall a\in
\bar\sA(w)\}.\]

On 
note 
encore
\[M(w)^{[r]}=\sR(\point{w})^{[r]}/C(w)^{[r]}\]
et
\[ad^w:\sR(\point{w})^{[r]}\to M(w)^{[r]}\]
la projection canonique.

Remarquons 
que
$C(w)^{[r]}$ et
$M(w)^{[r]}$ d\'e\-pen\-dent effectivement de 
$w$, et pas seulement de
$\point{w}$. Les inclusions du type 
$\bar\sA(w_1\dots w_n)\subset
\bar\sA(\point{w_1}\dots\point{w_n})$ 
entra\^{\i}ne des inclusions du type
\[C(\point{w_1}\dots 
\point{w_n})^{[r]}\subset C(w_1\dots w_n)^{[r]}\]
et des projections 
correspondantes
\begin{equation}\label{Proj}
M(\point{w_1}\dots 
\point{w_n})^{[r]}\Surj 
M(w_1\dots
w_n)^{[r]}\end{equation}
factorisant les $ad^w$.

Les 
relations \eqref{coc2} et \eqref{cob2} se
r\'e\'ecrivent 
maintenant
\end{sloppypar}
\begin{equation}\label{coc2bis}
ad^{ABC}(r_ 
{A,B,C})=0
\end{equation}
avec
\begin{equation}\label{rcoc}
r_{A,B,C}= 
1_A\bullet
{n}_{B,C}-{n}_{A\bullet B, C}  +{n}_{A,B\bullet C} - 
{n}_{A,B}\bullet
1_C
\end{equation}
ainsi 
que
\begin{equation}\label{cob2bis}
ad^{AB}(n_{A,B}-m_A\bullet 
1_B-1_A\bullet
m_B+m_{A\bullet B})=0.
\end{equation}

L'action \`a 
gauche et \`a droite de $K\langle Ob(\sA)\rangle$ sur
$\sR^{[r]}$ en 
induit une sur son quotient
\[M^{[r]} = \prod_w 
M(w)^{[r]}.\]

\begin{sloppypar}
\subsubsection{Une relation de 
cocycle, III - fin de la preuve de a)} De
m\^eme que pour le 
lemme
\ref{l4}, \eqref{coc2} et \eqref{cob2} s'interpr\`etent 
maintenant comme
relations de
$2$-cocycle et de $2$-cobord de 
Hochschild \`a valeurs dans le $K\langle
Ob(\sA)\rangle$-bimodule 
$C^{[r]}$. La seule remarque \`a faire est
qu'\'etant donn\'e trois 
mots $w_0,w_1,w_2$, et $w=w_0w_1w_2$, la
nullit\'e de 
$ad^{\point{w_0}\point{w_1}\point{w_2}}(
r_{\point{w_0},\point{w_1},\point{w_2}})$, donn\'ee par
\eqref{coc2bis}, entra\^{\i}ne celle de 
$ad^w(r_{\point{w_0},\point{w_1},
\point{w_2}})$ via la projection 
\eqref{Proj}. On conclut 
comme
auparavant.
\end{sloppypar}

\subsubsection{Preuve de 
b)}\label{prb} Soient $s,s'$ deux sections
mo\-no\-\"{\i}\-dales strictes 
de
$\pi_\sA$. Il s'agit de montrer qu'il existe une famille 
$(u_A)_{A\in
\sA}, u_A \in 1_A +\sR(A,A)$, telle que pour 
tous
$A,B\in \sA$ on ait
$s'(f)= u_As(f)(u_A)^{-1}$ pour tout $f\in 
\bar\sA(A)$ et
\[u_{A\bullet B}=u_A\bullet u_B.\]

Nous allons 
montrer:

\begin{lemme}\label{l6} Il existe une famille 
double
$(u_A^{(r)})_{A\in\sA,r\ge 1}$, avec $u_A^{(r)}\in1+\sR(A)$, 
ayant
les propri\'et\'es suivantes:
\begin{thlist}
\item 
$u_A^{(r+1)}\equiv u_A^{(r)}\pmod{\sR(A)^r}$ pour tout
$A\in\sA$ et 
tout $r\ge 1$.
\item $s'(f)= u^{(r)}_A s(f) (u^{(r)}_A)^{-1}$ pour 
tout
$A\in \sA$, tout $f\in \bar\sA(A)$ et tout $r\ge 1$.
\item 
$u^{(r)}_{A\bullet B}\equiv 
u^{(r)}_A\bullet
u^{(r)}_B\pmod{\sR(A\bullet B)^{r}}$ pour tous 
objets $A,B$ et tout
$r\ge 1$.
\end{thlist}
\end{lemme}

Le lemme 
\ref{l6} implique b) de la m\^eme mani\`ere que les lemmes
\ref{l4} 
et \ref{l5} impliquaient c) et a). On proc\`ede toujours 
par
r\'ecurrence sur $r$. Nous aurons besoin des nouvelles 
notations
suivantes, pour un mot $w$ en les objets de $\sA$ (voir 
aussi \S
\ref{co2}):
\[
C_s(w)=\{x\in \sA(\point{w})\mid x 
s(f)=s(f)x\text{ pour tout }
f\in\bar\sA(w)\}\]
\[C_s^r(w)= 
\sR(\point{w})^r\cap 
C_s(w)\]
\[
C^{[r]}(w)=\frac{C_s^r(w)+\sR(\point{w})^{r+1}}{\sR(\point 
{w})^{r+1}}\simeq
C_s^r(w)/C_s^{r+1}(w).
\]

Le groupe $C^{[r]}(w)$ 
ne d\'epend pas du choix de $s$, puisque deux
sections sont 
conjugu\'ees par des \'el\'ements de $1+\sR$ (th\'eor\`eme
\ref{T1} 
b)) et que $1+\sR$ op\`ere trivalement sur
$\sR^{[n]}$ par 
conjugaison. On prendra garde de ne pas le
confondre avec le 
groupe
$C(w)^{[r]}$ du \S \ref{co2}: l'inclusion
\[C^{[r]}(w)\subset 
C(w)^{[r]}\]
est en g\'en\'eral stricte.

En utilisant le 
th\'eor\`eme \ref{T1} b),
choisissons une famille
$(u^{(1)}_A)_{A\in 
\sA}$, avec $u_A^{(1)}\in 1_A+\sR(A)$, v\'erifient la
propri\'et\'e 
(ii) du lemme \ref{l6}. La condition (iii) est automatique;
le lemme 
est donc vrai pour $r=1$.

Supposons maintenant $r\ge 1$ et le lemme 
connu pour $r$. Posons
\[n_{A,B}:=1_{A\bullet 
B}-(u^{(r)}_A\bullet
u^{(r)}_B)^{-1}u^{(r)}_{A\bullet 
B}\in
\sR(A\bullet B)^{r}.\]

Comme $s$ et $s'$ sont mo\-no\-\"{\i}\-dales 
strictes, on a l'identit\'e
\[ (u^{(r)}_A\bullet 
u^{(r)}_B)^{-1}u^{(r)}_{A\bullet 
B}(s(f)\bullet
s(g))(u^{(r)}_{A\bullet B})^{-1}u^{(r)}_A\bullet 
u^{(r)}_B=s(f)\bullet
s(g)\]
pour $(f,g)\in \bar\sA(A)\times\bar\sA(B)$,
ce qui \'equivaut \`a
\[n_{A,B}\in C^r(AB).\]

Par ailleurs, la d\' 
efinition de $n_{A,B}$ donne
imm\'ediatement la 
relation
\[1_A\bullet
{\bar n}_{B,C}-{\bar n}_{A\bullet B, C}  +{\bar 
n}_{A,B\bullet C} -
{\bar n}_{A,B}\bullet 1_C=0\]
o\`u $\bar n_{A,B}$ 
d\'esigne l'image de $n_{A,B}$ dans
$C^{[r]}(AB)$.

Soit maintenant 
$(u^{(r+1)}_A)$ une famille v\'erifiant les conclusions du
lemme 
\ref{l6}. Posons
\[u_A^{(r+1)}=u_A^{(r)}(1+m_A).\]

On a donc 
$m_A\in\sR(A)^r$. La condition (ii) du lemme se traduit en
\[(1+m_A) 
s(f) (1+ m_A)^{-1} =s(f)\]
pour tout $f\in \bar\sA(A)$, 
c'est-\`a-dire
\[m_A\in C^r(A).\]

La condition (iii), d'autre part, 
se traduit en
\[\bar n_{A,B}= \bar m_{A\bullet B}-1_A\bullet \bar 
m_B-\bar m_A\bullet
1_B\]
o\`u $\bar m_A$ est la projection de $m_A$ 
dans $C^{[r]}(A)$ (calcul
imm\'ediat). Inversement, la donn\'ee de 
$m_A\in \sR(A)^r$ v\'erifiant
ces deux conditions fournit une famille 
$(u^{(r+1)}_A)$ v\'erifiant les
conclusions du lemme \ref{l6}.

On 
termine comme pr\'ec\'edemment, en utilisant le 
$K\langle
Ob(\sA)\rangle$-bi\-mo\-du\-le
\[C^{[r]}=\prod_w 
C^{[r]}(w),\]
et en posant
\[\bar n_{w_0,w_1}(w) =\begin{cases}
\bar n_{\point{w_0},\point{w_1}}& \text{si $w_0
w_1=w$,}\\
0&\text{sinon.}
\end{cases}\]
Ceci conclut la preuve de b), donc celle de la proposition
\ref{p5}, et finalement celle du th\'eor\`eme \ref{T2}.\qed

\begin{rems}\label{R11/4}\
\begin{itemize}
\item[$a)$] Supposons $\sA$ ab\'elienne et engendr\'ee par un nombre fini
d'objets $A_i$, au sens o\`u tout objet est
sous-quotient d'une som\-me finie de copies de $\bullet$-mon\^omes en
ces objets. On prendra garde que $\bar\sA$
n'est pas n\'e\-ces\-sai\-re\-ment engendr\'ee par les $\pi(A_i)= A_i$. Le
point est que l'i\-ma\-ge dans
$\bar\sA$ d'un monomor\-phis\-me de
$\sA$ n'est pas n\'e\-ces\-sai\-re\-ment un monomor\-phis\-me. En fait, on
verra plus loin des exemples o\`u  $\bar \sA$
n'est engendr\'ee par aucun ensemble fini d'objets.
\item[$b)$] Dans \cite{ak(note)}, nous appliquons le th\'eor\`eme de
scindage mo\-no\-\"{\i}\-dal ci-dessus pour construire, inconditionnellement,
des r\'ea\-li\-sa\-tions de motifs nu\-m\'e\-ri\-ques, et les groupes de
Galois motiviques associ\'es.
\item[$c)$] Dans la quatri\`eme partie de ce travail, nous appliquons
le m\^eme th\'e\-o\-r\`e\-me pour d\'efinir des ``enveloppes"
pro-semi-simples ou pro-r\'e\-duc\-ti\-ves.
\end{itemize}
  \end{rems}

\subsection{Variantes} Voici un avatar mo\-no\-\"{\i}\-dal de la proposition
\ref{P3}.

\begin{sloppypar}
\begin{prop} \label{P3'} Soit $(\sB, \bullet)$ une autre ca\-t\'e\-go\-rie
mo\-no\-\"{\i}\-dale v\'erifiant les
hypoth\`eses du th\'eor\`eme \ref{T2}. Soit
$T:\sB\to
\sA$ un
$K$-foncteur mo\-no\-\"{\i}\-dal radiciel. Notons $\bar T:\bar
\sB\to \bar \sA$ le foncteur
induit. \\
Soit $s_\sA:\bar\sA\to \sA$ (\resp $s_\sB:\bar\sB\to \sB$ )
une section mo\-no\-\"{\i}\-dale de la
projection canonique. \\
Alors $s_\sA\circ \bar T$ et $T\circ s_\sB$ sont isomorphes (via un
isomor\-phis\-me mo\-no\-\"{\i}\-dal qui couvre
l'isomor\-phis\-me mo\-no\-\"{\i}\-dal identique modulo les radicaux). \\
En particulier, si $T$ et les sections
$s_A$ et $s_B$ sont strictes, il existe un
syst\`eme $(u_X)$ d'\'el\'ements de $1_{T(X)}+ \rad(\sA)(T(X),T(X))$
tel que $s_\sA\bar
T(h)=u_{X'}Ts_\sB(h)(u_X)^{-1}$ pour
tout $h\in \bar \sB(X,X')$, et $u_{X\bullet Y}=u_X\bullet u_Y$.
\end{prop}

\prf On prouve d'abord la seconde assertion, sous l'hypoth\`ese que
$\sA$ et $\sB$ sont mo\-no\-\"{\i}\-dales
strictes. La preuve est enti\`erement parall\`ele \`a celle du point $b)$
de la proposition \ref{p5}.  La
proposition \ref{P3} montre l'e\-xis\-ten\-ce d'un
syst\`eme $(u_X)$ d'\'el\'ements de $1_{T(X)}+ \rad(\sA)(T(X),T(X))$
tel que $s_\sA\bar
T(h)=u_{X'}Ts_\sB(h)(u_X)^{-1}$ pour
tout $h\in \bar \sB(X,X')$. Pour le modifier en un syst\`eme
v\'erifiant de plus $u_{X\bullet Y}=u_X\bullet
u_Y$, on construit une suite d'approximations $u^r_X$. Le point est 
de remplacer
le $K\langle Ob(\sA)\rangle$-bimodule
$C^{[r]}$ par le $K\langle Ob(\sB)\rangle$-bimodule
    $\displaystyle {C'}^{[r]} =\prod_w C'(w)^{[r]}$, o\`u $ C'(w)^{[r]}$
est d\'efini de la mani\`ere suivante:
\[
C_{Ts}(w)=\{x\in \sA(T(\point{w}))\mid x
(Ts_\sB(h))=(Ts_\sB(h))x\text{ pour tout }
h\in\bar\sB(w)\}\]
\[C_{Ts}^r(w)= \sR_\sA(T(\point{w}))^r\cap C_{Ts}(w)\]
\[
{C'}^{[r]}(w)=\frac{C_{Ts}^r(w)+\sR_\sA(T(\point{w}))^{r+1}}{\sR_\sA(T
(\point{w}))^{r+1}}\simeq
C_{Ts}^r(w)/C_{Ts}^{r+1}(w). \]
Passons maintenant au cas g\'en\'eral. D'apr\`es le th\'eor\`eme de 
conjugaison des sections mo\-no\-\"{\i}\-dales,
il est loisible de changer de sections $s_\sA, s_\sB$. Plut\^ot que 
$s_\sA, s_\sB$, c'est de
${s'}_\sA= u_\sA\circ (s_\sA)^{str} \circ \bar v_\sA,\;{s'}_\sB= 
u_\sB\circ (s_\sB)^{str} \circ \bar v_\sB$ dont
nous nous servirons (avec les notations de la construction de MacLane 
\ref{constrM}). On est alors ramen\'e \`a voir
que
$(s_\sA)^{str} \circ \bar T^{str}$ et $T^{str}\circ (s_\sB)^{str} $ 
sont isomorphes (via un
isomor\-phis\-me mo\-no\-\"{\i}\-dal qui couvre
l'isomor\-phis\-me mo\-no\-\"{\i}\-dal identique modulo les radicaux); mais on 
est alors dans le cas strict d\'ej\`a
trait\'e.\qed
\end{sloppypar}

\begin{cor} Dans la situation de la proposition \ref{P3}, supposons
$\sA$, $\sB$ strictement $K$-lin\'eaires,
\ie telles que
$\oplus$ soit strictement associative (\cf \ref{ml}). Supposons aussi
que $T$ soit strict vis-\`a-vis de
$\oplus$. Alors on peut choisir les $(u_X)$ (conjuguant $s_\sA\bar
T$ et $Ts_\sB$) de telle sorte que $u_{X\oplus Y}=u_X\oplus u_Y$.\qed

\end{cor}

Voici une variante mo\-no\-\"{\i}\-dale du compl\'ement de Malcev au
th\'eor\`eme de Wedderburn.

\begin{prop}\label{cmmalcev} Sous les
hypoth\`eses du th\'eo\-r\`eme \ref{T2},
supposons don\-n\'ee en outre une sous-ca\-t\'e\-go\-rie
mo\-no\-\"{\i}\-dale $\sB$ \emph{s\'eparable} de $\sA$ (non
n\'e\-ces\-sai\-re\-ment pleine), telle que tout objet de $\sA$ 
soit
facteur direct d'un objet de $\sB$. Alors il
existe une section 
mo\-no\-\"{\i}\-dale
$s$ de $\pi_\sA$ qui v\'erifie
$(s\circ \pi_\sA )\vert 
\sB = id_{\sB}$.
\end{prop}

\prf Quitte \`a remplacer $\sA$ et $\sB$ 
par $\sA^\natural$ et
$\sB^\natural$ respectivement, on peut
supposer 
que
$Ob(\sB)=Ob(\sA)$. Dans ce cas, le
corollaire \ref{cmalcev} 
montre l'e\-xis\-ten\-ce d'une section fonctorielle
de $\pi_\sA$ 
telle que $(s\circ \pi_\sA )_{\vert
\sB} = id_{\sB}$. Pour la 
modifier en une section mo\-no\-\"{\i}\-dale
v\'erifiant la m\^eme 
condition, on reprend la m\'ethode de preuve du
th\'eor\`eme \ref{T2} 
a), en
construisant une suite d'approximations $s_r$ (\cf le lemme 
\ref{l5} et sa
preuve). Le point est de remplacer le
$K\langle 
Ob(\sA)\rangle$-bimodule
$M^{[r]}$ par le sous-bimodule 
$\displaystyle {M'}^{[r]} =\prod_w
M'(w)^{[r]}$, o\`u
$M'(w)^{[r]}$ 
est le
sous-$K$-espace de $M(w)^{[r]}$ des \'el\'ements
dont un 
relev\'e dans $\sR(\point{w})^{[r]}$ commute \`a
$\bar\sB (w)$. On 
observe
qu'on a bien $ad^{AB}n_{AB}\in M'(AB)^{[r]}$ puisque la 
restriction
de $s_r$ \`a $\sB=\bar{\sB}$ est
mo\-no\-\"{\i}\-dale (c'est 
l'identit\'e); on
obtient ainsi des
$\nu_A$  commutant \`a 
$\bar\sB(A)$, de sorte que $s_{r+1}$ est
encore l'identit\'e sur 
$\sB=\bar{\sB}$.
\qed

\subsection{Compl\'ement} Ce compl\'ement 
servira au \S \ref{tre}.

\begin{prop}\label{t3} Sous les 
hypoth\`eses du th\'eor\`eme \ref{T2},
soit $v$ une structure 
mo\-no\-\"{\i}\-dale sur le foncteur identique de $\sA$
telle que 
$\pi_\sA(v)=1$. Alors $(Id_\sA,v)$ est mo\-no\-\"{\i}\-dalement
isomorphe 
\`a $(Id_\sA,1)$ via un isomor\-phis\-me de foncteurs mono\"{\i}daux se 
projetant sur $1$ dans $\bar\sA$.
\end{prop}

\prf On peut supposer 
$\sA$ strictement mo\-no\-\"{\i}\-dale. Il s'agit de
r\'esoudre 
l'\'equation
\[v_{A,B}=(u_A^{-1}\bullet u_B^{-1})u_{A\bullet 
B}\]
pour $A,B\in \sA$, o\`u $u=(u_A)$ est un automor\-phis\-me (non 
mo\-no\-\"{\i}\-dal)
de $Id_\sA$ se projetant sur $1$. On proc\`ede comme 
d'habitude, en supposant de plus $\sA$
strictement $K$-lin\'eaire et 
en ne traitant que d'endomor\-phis\-mes. Il
suffit de 
montrer:

\begin{lemme}\label{l18} Il existe une
famille double 
$(u^{(r)})_{r\ge 1}$ d'automor\-phis\-mes de 
$Id_\sA$,
avec
$u_A^{(r)}\allowbreak\in 1+\sR(A)$, ayant les 
propri\'et\'es
suivantes:
\begin{thlist}
\item $u_A^{(r+1)}\equiv 
u_A^{(r)}\pmod{\sR(A)^{r}}$ pour tout
$A\in\sA$ et tout $r\ge 1$.
\item $v_{A,B}\equiv ({u^{(r)}_A}^{-1}\bullet
{u^{(r)}_B}^{-1})u^{(r)}_{A 
\bullet B}\pmod{\sR(A\bullet
B)^{r}}$ pour tous objets
$A,B$ et 
tout
$r\ge 1$.
\end{thlist}
\end{lemme}

\prf R\'ecurrence sur
$r$, 
en partant de $r=1$ avec $u^{(1)}=1$. Supposons $r\ge 1$ et
trouv\'ee 
un automor\-phis\-me $u^{(r)}$. Quitte \`a remplacer 
$v_{A,B}$
par
$(u^{(r)}_A\bullet u^{(r)}_B){u^{(r)}_{A\bullet 
B}}^{-1}v_{A,B}$, on peut
supposer que $u^{(r)}=1$ et que $v\in 
\sR^{r}$. 

Pour tout mot $w$ en les objets de $\sA$, posons 

\begin{align*}
C(w)&=\{x\in \sA(\point{w})\mid x \text{ centralise 
}
\sA(w)\}\\
C^{<r>}(w)&=C(w)\cap\sR(\point{w})^r\quad (r\ge 
0)\\
\Gamma^{[r]}(w)&=C^{<r>}(w)/C^{<r+1>}(w)\\
\Gamma^{[r]}&=\prod_w 
\Gamma^{[r]}(w).
\end{align*}

L'alg\`ebre $K\langle Ob(\sA)\rangle$ 
op\`ere \`a gauche et \`a droite
sur $\Gamma^{[r]}$ de la mani\`ere 
habituelle. 

Posons $v_{A,B}=1-\nu_{A,B}$: alors $\nu_{A,B}\in 
C^{<r>}(AB)$ et on a la
congruence identique
\[\nu_{A\bullet 
B,C}+\nu_{A,B}\bullet 1_C\equiv \nu_{A,B\bullet
C}+1_A\bullet 
\nu_{B,C}\pmod{\sR(A\bullet B\bullet C)^{r+1}}.\]

Soit $u^{(r+1)}$ 
un automor\-phis\-me de $Id_\sA$ solution
de (i) et (ii). Posons 
$u_A^{(r+1)}\allowbreak=1-\mu_A$: on a donc
$\mu_A\in C^{<r>}(A)$. 
Soit $\bar\mu_A$ l'image de $\mu_A$ dans
$\Gamma^{[r]}(A)$. On a la 
congruence identique
\[\nu_{A,B}\equiv \mu_{A\bullet 
B}-1_A\bullet
\mu_B-\mu_A\bullet 1_B\pmod{\sR(A\bullet 
B)^{r+1}}.\]

\begin{sloppypar}
R\'eciproquement, si $(\mu_A)$ est 
une famille
d'\'el\'ements de $C^{<r>}(A)$ v\'erifiant ces 
congruences, alors la
famille $(u_A^{(r+1)}=1-\mu_A)$ d\'efinit un 
automor\-phis\-me de $Id_\sA$
solution de (i) et (ii). Il reste \`a 
appliquer comme pr\'ec\'edemment la
nullit\'e de $H^2(K\langle 
Ob(\sA)\rangle, 
\Gamma^{<r>})$.\qed
\end{sloppypar}

\begin{rem}\label{duaux droite} 
Si $A, B$ admettent des duaux \`a
droite (\cf \ref{s3.5}), les 
isomor\-phis\-mes canoniques 
\[\sA(\un,A^\vee\bullet
B)\cong \sA(A,B)\] 

de \eqref{eq6.2.1} passent \`a $\bar\sA$ et sont compatibles \`a 
toute
section mo\-no\-\"{\i}\-dale $s$ comme dans le th\'eor\`eme. Si tous 
les objets
de $\sA$ admettent des duaux \`a droite   (par exemple si 
$\sA$ est
sy\-m\'e\-tri\-que et rigide), on en d\'eduit que $s$ est 
d\'etermin\'ee par
ses valeurs sur les
$\bar\sA(\un , 
A)$.
\end{rem}

\section{Repr\'esentabilit\'e} \label{repr}

Dans ce 
paragraphe, nous nous proposons d'interpr\'eter les r\'esultats
des paragraphes \ref{s5} et \ref{s6} en termes plus cat\'egoriques.
Bri\`evement, on peut associer \`a une petite $K$-ca\-t\'e\-go\-rie de
Wedderburn (mo\-no\-\"{\i}\-dale) le grou\-po\-\"{\i}\-de des sections
(mo\-no\-\"{\i}\-dales) de $\pi_\sA$. Cette
construction est $2$-fonctorielle en $\sA$ pour les foncteurs pleinement
fid\`eles et les isomor\-phis\-mes naturels. Sous une hypoth\`ese de finitude
suppl\'ementaire, elle peut m\^eme s'enrichir en un
$2$-foncteur \`a valeurs dans les \emph{$K$-groupo\"{\i}des affines
scind\'es unipotents} (th\'e\-o\-r\`e\-mes \ref{gerbe} et \ref{gerbe'}).

Ces constructions ne d\'ependent pas des th\'eor\`emes \ref{T1} et
\ref{T2}. Dans le cas re\-pr\'e\-sen\-ta\-ble, elles permettent de les
r\'einterpr\'eter en des r\'esultats tr\`es forts de \emph{rationalit\'e}
pour les $K$-groupo\"{\i}des en jeu.

Nous nous pla\c cons d'abord le cadre discret, puis dans le cadre 
alg\'ebrique.

\subsection{Le paysage discret}

\subsubsection{Le cas non mo\-no\-\"{\i}\-dal}\label{i1} Soit $\sA$ une
petite $K$-ca\-t\'e\-go\-rie de Wedderburn. Consid\'erons la ca\-t\'e\-go\-rie
$\sG(\sA)$ dont les objets sont les sections de $\pi_\sA$, un mor\-phis\-me de
$s$ vers $s'$ \'etant une famille $(u_A)_{A\in\sA}$ comme dans le
th\'eor\`eme \ref{T1} b) et le compos\'e de deux mor\-phis\-mes $(u_A),(v_A)$
\'etant $(v_Au_A)$. (Comme nous l'a fait remarquer A. Bru\-gui\`e\-res,
les mor\-phis\-mes de $\sG(A)$ s'interpr\`etent aussi comme les
transformations naturelles entre sections de $\pi_\sA$ se projetant sur
l'identit\'e dans $\bar\sA$.) Alors
$\sG(\sA)$ est un petit groupo\"{\i}de: le \emph{groupo\"{\i}de
des sections de $\pi_\sA$}. Le contenu du th\'eor\`eme \ref{T1} a) (\resp
b))  est que ce groupo\"{\i}de est  \emph{non vide} (\resp \emph{connexe}).

L'interpr\'etation de la proposition \ref{P3} dans ce langage est le
peu satis\-faisant principe de fonctorialit\'e suivant: si $s_\sA\in
\sG(\sA)$ et $s_\sB\in \sG(\sB)$, d\'efinissons le \emph{transporteur de
$s_\sA$ vers $s_\sB$} comme l'ensemble $T(s_\sA,s_\sB)$ des $(u_X)$
v\'erifiant les conditions de ladite proposition: ce sont aussi les
transformations naturelles de $s_\sA\bar T$ vers $Ts_\sB$ se projetant
sur l'identit\'e. Alors, pour tout couple $(s_\sA,s_\sB)$, l'ensemble
$T(s_\sA,s_\sB)$ est \emph{non vide}.

Toutefois, si $\sA \to \sB$ est \emph{pleinement fid\`ele}, alors on a un
foncteur ca\-no\-ni\-que $\sG(\sB)\to \sG(\sA)$ dans l'autre sens. Cette
correspondance s'\'e\-tend en un $2$-foncteur: notons
\[Wed\{K\}\]
la sous-$2$-ca\-t\'e\-go\-rie non pleine de $\{K\}$ (voir \S \ref{rap}) dont
\begin{itemize}
\item les objets sont les petites $K$-ca\-t\'e\-go\-ries de Wedderburn
\item les $1$-mor\-phis\-mes sont les $K$-foncteurs pleinement fid\`eles
\item les $2$-mor\-phis\-mes sont les isomor\-phis\-mes 
naturels.
\end{itemize}

D'autre part, notons $\Gr$ la 
$2$-ca\-t\'e\-go\-rie des
petits groupo\"{\i}des (tous les foncteurs et 
toutes les transformations
naturelles sont autoris\'es). Alors on 
v\'erifie au prix d'un petit
calcul:

\begin{sorite}\label{so6} La 
construction $\sA\mapsto \sG(\sA)$ d\'efinit
un
$2$-foncteur 
$1$-con\-tra\-va\-riant et $2$-covariant
\[\sG:Wed\{K\}\to \Gr.\]
Le 
th\'eor\`eme \ref{T1} \'enonce que son image est contenue dans 
la
$2$-sous-cat\'e\-go\-rie pleine des groupo\"{\i}des non vides 
connexes.\qed
\end{sorite}

En particulier, toute $K$-\'equivalence 
de ca\-t\'e\-go\-ries $\sA\iso \sB$
d\'efinit une \'equivalence de 
ca\-t\'e\-go\-ries $\sG(\sB)\iso \sG(\sA)$, qui
n'est pas en g\'en\'eral un 
isomor\-phis\-me de ca\-t\'e\-go\-ries. Toutefois:

\begin{sorite}\label{so7} 
Pour toute $K$-ca\-t\'e\-go\-rie de Wedderburn $\sA$,
les 
foncteurs
\[
\begin{matrix}
\;&\;& \sG(\sA^\oplus) 
&\;&\;\\
\;&\nearrow\;&\;&\searrow &\;\\
 
\sG(\sA^{\oplus\natural}) &\;&\;&\;& \sG(\sA) \\
 
\;&\searrow\;&\;&\nearrow &\;\\
 
\;&\;&
\sG(\sA^\natural) &\;&\; 
\end{matrix}\]
sont des 
isomor\-phis\-mes de groupo\"{\i}des.\qed
\end{sorite}

(On peut s'en 
convaincre en remarquant que ces foncteurs d\'efinissent
des inverses 
partiels des foncteurs du sorite \ref{so01}.)

\subsubsection{Le cas 
mo\-no\-\"{\i}\-dal}\label{i2} Supposons maintenant que
$\sA$ soit 
mo\-no\-\"{\i}\-dale, et que son radical soit un id\'eal mo\-no\-\"{\i}\-dal. De 
m\^eme que dans 
\ref{i1}, on peut introduire le groupo\"{\i}de 
$\sG^{\otimes}(\sA)$
des {\it sections mo\-no\-\"{\i}\-dales de
$\pi_A$} 
(les mor\-phis\-mes \'etant les mor\-phis\-mes de foncteurs
mono\"{\i}daux se 
projetant
sur l'identit\'e dans $\bar\sA$). C'est un 
sous-groupo\"{\i}de de
$\sG(\sA)$, lui aussi connexe non 
vide\footnote{Tout foncteur d'un
groupo\"{\i}de connexe non vide vers 
un autre est essentiellement
surjectif\dots}.
      La proposition 
\ref{P3'} admet la m\^eme interpr\'etation que la
proposition 
\ref{P3} du paragraphe \ref{s5}. 
Si $\sA \to \sB$ est pleinement 
fid\`ele, on a un foncteur canonique
$\sG^\otimes(\sB)\to 
\sG^\otimes(\sA)$ dans l'autre sens. 

Notons
\[Wed^\otimes\{K\}\]
la 
$2$-ca\-t\'e\-go\-rie  dont
\begin{itemize}
\item les objets sont les 
petites $K$-ca\-t\'e\-go\-ries de Wedderburn mo\-no\-\"{\i}\-dales,
\`a 
ra\-di\-cal mo\-no\-\"{\i}\-dal
\item les $1$-mor\-phis\-mes sont les 
$K$-foncteurs mono\"{\i}daux pleinement
fi\-d\`e\-les
\item les 
$2$-mor\-phis\-mes sont les isomor\-phis\-mes naturels 
mono\"{\i}daux.
\end{itemize}

On a un $2$-foncteur \'evident 

\begin{equation}\label{eq10}
\Omega:Wed^\otimes\{K\}\to 
Wed\{K\}
\end{equation}
(oubli des structures mo\-no\-\"{\i}\-dales). Comme 
au
\S \ref{i1}, on obtient en fait:

\begin{sorite}\label{so6t} La 
construction $\sA\mapsto \sG^\otimes(\sA)$
d\'efinit un
$2$-foncteur 
$1$-con\-tra\-va\-riant et 
$2$-covariant
\[\sG^\otimes:Wed^\otimes\{K\}\to \Gr.\]
Le 
th\'eor\`eme \ref{T2} \'enonce que son image est contenue dans 
la
$2$-sous-ca\-t\'e\-go\-rie pleine des groupo\"{\i}des non vides 
connexes.\\
Il existe une $2$-transformation naturelle 
canonique
\[\sG^\otimes\Rightarrow
\sG\circ\Omega\] 
o\`u $\sG$ est 
comme dans le sorite \ref{so6} et 
$\Omega$ est comme 
en
\eqref{eq10}.\qed
\end{sorite}

\begin{sorite}\label{so7t} Pour 
toute $K$-ca\-t\'e\-go\-rie de Wedderburn $\sA$,
les 
foncteurs
\[
\begin{matrix}
\;&\;& \sG^\otimes(\sA^\oplus) 
&\;&\;\\
\;&\nearrow\;&\;&\searrow &\;\\
 
\sG^\otimes(\sA^{\oplus\natural}) &\;&\;&\;& \sG^\otimes(\sA) \\
 
\;&\searrow\;&\;&\nearrow &\;\\
 
\;&\;&
\sG^\otimes(\sA^\natural) &\;&\; 
\end{matrix}\]
sont des 
isomor\-phis\-mes de 
groupo\"{\i}des.\qed
\end{sorite}

\subsubsection{Objets compacts, 
duaux et limites inductives}\label{cdi}
Supposons
$\sA$ stable par 
limites inductives quelconques, et soit $\sA^\comp$ 
sa
sous-ca\-t\'e\-go\-rie pleine form\'ee des objets compacts. Alors on a 
un
foncteur de restriction $\sG(\sA)\to \sG(\sA^\comp)$. Si tout 
objet de
$\sA$ est isomorphe \`a une limite inductive d'objets 
compacts, on a un
foncteur en sens inverse: ces deux foncteurs sont 
des \'equivalences de
ca\-t\'e\-go\-ries, quasi-inverses l'une de 
l'autre.

M\^eme sorite dans le cas mo\-no\-\"{\i}\-dal sy\-m\'e\-tri\-que, en 
supposant $\un$
compact et en rempla\c cant ``compact" par ``ayant un 
dual".

\subsubsection{Extension des scalaires} Si $L/K$ est une 
extension, on a
des foncteurs ``extension des scalaires" (na\"{\i}ve 
et \`a la Saavedra)
$\sG(\sA)\to \sG(\sA_L)$ et $\sG(\sA)\to 
\sG(\sA_{(L)})$, et de m\^eme
dans le cas mo\-no\-\"{\i}\-dal. Si $L/K$ est 
(finie) s\'eparable, ces foncteurs
sont isomorphes en vertu du 
th\'eor\`eme \ref{t2}. Dans le cas
g\'en\'eral et en supposant $\sA$ 
stable par limites inductives
quelconques, on a un foncteur de 
restriction $\sG(\sA_{(L)})\to
\sG(\sA_L)$ commutant naturellement 
aux deux foncteurs ci-dessus. M\^eme
sorite dans le cas 
mo\-no\-\"{\i}\-dal.

\subsection{Objets en ca\-t\'e\-go\-ries, en 
groupo\"{\i}des et en actions de
grou\-pe}

\subsubsection{Objets en 
ca\-t\'e\-go\-ries et en
groupo\"{\i}des}\label{obcatgr} Soit
$\sX$ une 
ca\-t\'e\-go\-rie avec limites projectives finies. Rappelons 
qu'un
\emph{objet en ca\-t\'e\-go\-ries} dans
$\sX$ est la 
donn\'ee
$\Gamma$ de deux objets $E,S\in \sX$, d'un 
mor\-phis\-me
\[(\beta,\sigma) :  E \to S\times S\]
``but" et ``source" 
et d'un mor\-phis\-me ``loi de composition"
\[\circ : E \times_{{}^\beta 
S^\sigma} E \to E \]
telles que, pour tout objet $X\in \sX$, le 
couple $(\sX(X,E),\sX(X,S))$
et les mor\-phis\-mes correspondants 
d\'efinissent une petite ca\-t\'e\-go\-rie
$\Gamma(X)$ (objets: $S(X)$; 
mor\-phis\-mes: $E(X)$). Par le lemme de Yoneda,
cela revient
\`a exiger 
que
$\beta,\sigma,\circ$ v\'erifient les identit\'es habituelles. Un 
objet
en ca\-t\'e\-go\-ries $\Gamma$ est un \emph{objet en groupo\"{\i}des} 
si, pour
tout $X\in \sX$, $\Gamma(X)$ est un groupo\"{\i}de (m\^eme 
remarque).

Si $\Gamma=(E,S,\beta,\sigma,\circ)$ 
et
$\Gamma'=(E',S',\beta',\sigma',\circ')$ sont deux objets en 
ca\-t\'e\-go\-ries
de $\sX$, un \emph{mor\-phis\-me} de $\Gamma$ vers $\Gamma'$ 
est un couple
$T=(f,g)$, avec $f\in \sX(E,E')$, $g\in \sX(S,S')$, 
tels que les
diagrammes
\[\begin{CD}
E @>(\beta,\sigma)>> S\times 
S\\
@V{f}VV @V{g\times g}VV\\
E' @>(\beta',\sigma')>> S\times 
S
\end{CD}
\qquad
\begin{CD}
E \times_{{}^\beta S^\sigma} E @>\circ>> 
E\\
@V{f\times f}VV @V{f}VV\\
E' \times_{{}^{\beta'} {S'}^{\sigma'}} 
E' @>\circ'>> E'
\end{CD}
\]
soient commutatifs (ou, de mani\`ere 
\'equi\-va\-len\-te, que
$T(X):\Gamma(X)\to \Gamma'(X)$ soit un foncteur 
pour tout $X\in \sX$).

Si $T=(f,g)$ et $T'=(f',g')$ sont deux tels 
mor\-phis\-mes, une
\emph{homotopie} de
$T$ vers $T'$ est un mor\-phis\-me 
$u\in \sX(S,E')$ tel que les diagrammes
\[\begin{CD}
S@>u>> E'\\
\s 
(g',g)\displaystyle\searrow &&@V{(\beta',\sigma')}VV\\
&&S'\times 
S'
\end{CD}
\qquad
\begin{CD}
E@>(f,u\beta)>> E'\times_{{}^{\beta'} 
{S'}^{\sigma'}}E'\\
@V{(u\sigma,f')}VV 
@V{\circ'}VV\\
E'\times_{{}^{\beta'} {S'}^{\sigma'}}E'@>{\circ'}>> 
E'
\end{CD}
\]
soient commutatifs (ou, ce qui revient au m\^eme, que 
$u(X)$ soit une
transformation naturelle quel que soit $X\in 
\sX$).

Les objets en ca\-t\'e\-go\-ries, mor\-phis\-mes et homotopies 
d\'efinissent une
$2$-ca\-t\'e\-go\-rie $\Cat(\sX)$. Sa sous-ca\-t\'e\-go\-rie 
$1$-pleine et $2$-pleine
form\'ee desobjets en groupo\"{\i}des est 
not\'ee $\Gr(\sX)$.

\subsubsection{Objets en actions de groupe} Un 
\emph{objet en groupes}
dans $\sX$ est un couple $(G,\mu_G)$ avec 
$G\in \sX$, $\mu_G\in
\sX(G\times G,G)$ tel que $G(X)=\sX(X,G)$ 
d\'efinisse un groupe pour la
loi $\mu_G(X)=\sX(X,\mu_G)$ quel que 
soit $X\in \sX$. De m\^eme, un
\emph{objet en actions de groupe} est 
un syst\`eme
$\Delta=(G,S,\mu_G,\mu_S)$, o\`u $(G,\mu_G)$ est un 
objet en groupes,
$S\in
\sX$, $\mu_S\in
\sX(G\times S,S)$, tel que 
$(\mu_G(X),\mu_S(X))$ d\'efinisse une action
de $G(X)$ sur $S(X)$ 
pour tout $X\in \sX$.

\begin{sloppypar}\'Etant donn\'e deux objets en actions de 
groupe
$\Delta=(G,S,\mu_G,\mu_S)$,  
$\Delta'=(G',S',\mu_{G'},\mu_{S'})$, un
\emph{mor\-phis\-me} de $\Delta$ 
vers $\Delta'$ est un couple $T=(f,g)$, avec
$f\in \sX(G,G')$, $g\in 
\sX(S,S')$, tels que les diagrammes
\[\begin{CD}
G\times G@>\mu_G>> 
G\\
@V{f\times f}VV @V{f}VV\\
G'\times G'@>\mu_{G'}>> 
G'
\end{CD}
\qquad
\begin{CD}
G\times S@>\mu_S>> S\\
@V{g\times g}VV 
@V{g}VV\\
G'\times S'@>\mu_{S'}>> S'
\end{CD}
\]
soient commutatifs. \end{sloppypar}
Il y a une notion \'evidente de composition de tels 
mor\-phis\-mes.

Enfin, \'etant donn\'e deux tels mor\-phis\-mes $T=(f,g)$, 
$T'=(f',g')$, une
\emph{homotopie} de $T$ vers $T'$ est un mor\-phis\-me 
$u\in \sX(S,G')$ tel
que les diagrammes
\[\begin{CD}
S@>(u,g)>> 
G'\times S'\\
\s g'\displaystyle\searrow 
&&@V{\mu_{S'}}VV\\
&&S'
\end{CD}
\qquad
\begin{CD}
G\times 
S@>(u\circ\mu_S,f\circ p_1)>> G'\times G'\\
@V{f'\times u}VV 
@V{\mu_{G'}}VV\\
G'\times G'@>\mu_{G'}>> G'
\end{CD}\]
soient 
commutatifs (en formules: $g'(s)=u(s)g(s)$; $u(\gamma
s)=\break 
f'(\gamma)u(s)f(\gamma)^{-1}$). La composition de deux
homotopies est 
donn\'ee par 
$u'\circ u := \mu_{G'}\circ(u',u)$ ($u'\circ 
u(s)=u'(s)u(s)$).

Les objets en actions de groupe, mor\-phis\-mes et 
homotopies forment une
autre $2$-ca\-t\'e\-go\-rie not\'ee $\Acg(\sX)$.

\subsubsection{Actions de groupe et groupo\"{\i}des} \label{ex3} 
Pour
obtenir un groupo\"{\i}de $\Gamma$, il suffit de se donner un 
objet en
actions de groupe
$\Delta=(G,S,\mu_G,\mu_S)$. On d\'efinit 
$E=G\times S$, $\beta=\mu_S$,
$\sigma=$ seconde projection. On 
v\'erifie (par exemple \`a l'aide du
lemme de Yoneda) que le 
mor\-phis\-me
\[\begin{CD}
G\times G\times S@>(p_{23},1_G\times\mu_S)>> 
E\times_{{}^{\beta}
{S}^{\sigma}}E,
\end{CD}\]
o\`u $p_{23}$ est la 
projection oubliant le premier facteur, est un
isomor\-phis\-me. On 
d\'efinit alors
$\circ$ comme le 
mor\-phis\-me
\[\begin{CD}
E\times_{{}^{\beta}{S}^{\sigma}}E&(\osi)^{-1}& 
G\times G\times
S@>\mu_G\times 1_S>> G\times S=E. 
\end{CD}\] 

Un 
tel groupo\"{\i}de est dit \emph{scind\'e}.

La r\`egle ci-dessus 
s'\'etend en un $2$-foncteur
\[\Acg(\sX)\to 
\Gr(\sX).\]

\begin{sloppypar}\begin{ex}\label{e10.1} $\sX=Sch/K$ est la
ca\-t\'e\-go\-rie 
des $K$-sch\'emas (dans notre univers). Un
\emph{$K$-groupo\"{\i}de} 
est un objet en groupo\"{\i}des
$\Gamma=(S,E)$ dans $\sX$. Il est 
dit
\emph{transitif sur $S$} s'il existe
$T'$ fi\-d\`e\-le\-ment plat 
quasi-compact sur $S\times_K S$ tel que
$Hom_{S\times_K S}(T',E)\neq 
\emptyset$ (le champ associ\'e \`a $T\mapsto
(S(T), E(T),\beta, 
\sigma, \circ)$ est
alors une gerbe, \cf \cite[3.3]{de}). Il est dit 
\emph{unipotent} (\resp
\hbox{\emph{r\'eductif\dots})} si la 
restriction de $E$ \`a la diagonale
de
$S\times_K S$ est un 
$S$-groupe affine pro-unipotent 
(\resp
pro-r\'eductif\dots).

 Un \emph{$K$-espace 
homog\`ene} est un objet en actions de groupe
$(G,S)$ dans
$\sX$ 
telle que $G$ op\`ere transitivement sur $S$. On note 

$\Hmg(Sch/K)\subset  \Acg(Sch/K)$ la sous-ca\-t\'e\-go\-rie $2$-pleine 
et
$1$-pleine des $K$-espaces homog\`enes.

\'Etant 
donn\'e un objet en actions de groupe $(G,S)$ dans
$\sX$, pour que le 
mor\-phis\-me $(\beta,\gamma)$ du $K$-groupo\"{\i}de
$\Gamma$ associ\'e 
par \ref{ex3} soit fpqc, il suffit que $(G,S)$ soit un
espace 
homog\`ene. Cela implique que $\Gamma$ est transitif sur $S$
(prendre 
$T'=E$ ci-dessus). Si $G$ est pro-unipotent 
(\resp
pro-r\'eductif\dots),
$\Gamma$ est unipotent (\resp 
r\'eductif\dots): en effet, la restriction
de $E$ \`a la diagonale 
est un sous-$S$-groupe ferm\'e du $S$-groupe
constant $G\times_K 
S$.
\end{ex}\end{sloppypar}

\subsection{Repr\'esentabilit\'e alg\'ebrique} 

Soient 
$\sA$ une
$K$-ca\-t\'e\-go\-rie de Wedderburn et
$L$ une extension de
$K$. 
Par le corollaire \ref{C3.1.}, $\sA_L$ est encore de Wedderburn, 
de
radical $\sR_L=\sR\otimes L$, ce qui sugg\`ere la pr\'esence 
d'un
$K$-groupo\"{\i}de (au sens des sch\'emas) dont le 
groupo\"{\i}de
$\sG(\sA)$ du \S \ref{i1} serait le groupo\"{\i}de des 
$K$-points. De
m\^eme, si
$\sA$ est mo\-no\-\"{\i}\-dale \`a radical 
mo\-no\-\"{\i}\-dal, on peut esp\'erer que
le 
groupo\"{\i}de
$\sG^\otimes(\sA)$ du \S \ref{i2} est 
re\-pr\'e\-sen\-ta\-ble.

Nous allons montrer que c'est bien le cas, sous 
une
hypoth\`ese de finitude convenable. En fait, nous allons obtenir 
encore
mieux.
 
\subsubsection{Le cas non mo\-no\-\"{\i}\-dal}

\begin{sloppypar}
\begin{thm}\label{gerbe} Soit
$Wedf\{K\}$ la 
sous-$2$-ca\-t\'e\-go\-rie
$1$-pleine et $2$-pleine de $Wed\{K\}$ form\'ee 
des ca\-t\'e\-go\-ries $\sA$
telles que $\dim_K\sA(A,B)<\infty$ pour tout 
$(A,B)\in \sA\times\sA$.\\
a) Soit $\sA\in Wedf\{K\}$. Alors il 
existe un
$K$-groupo\"{\i}de affine scind\'e unipotent 
$\Gamma(\sA)=(E,S)$,
transitif sur
$S$, tel que pour toute 
extension
$L$ de $K$, on ait
$S(L)=Ob(\sG(\sA_L))$ 
et
$E(s_1,s_2)(L)=\sG(\sA_L)(s_1,s_2)$ pour tout couple 
de
sections
$(s_1,s_2)$ de $\pi_{\sA_L}$.\\
b) Cette construction 
provient d'un $2$-foncteur $1$-contravariant 
et
$2$-co\-va\-riant
\[\Delta:Wedf\{K\}\to \Hmg(Sch/K)\]
\`a valeurs 
dans les $K$-espaces homog\`enes unipotents, via 
le
$2$-foncteur
$\Hmg(Sch/K)\to
\Gr(Sch/K)$ de
\ref{ex3}.\\
c) Les 
$2$-foncteurs $\Delta$ et $\Gamma$ commutent \`a l'extension 
des
scalaires.
\end{thm}
\end{sloppypar}

\prf Compte tenu de 
l'hypoth\`ese de finitude, pour tout objet $A$ le
foncteur 
des
$K$-alg\`ebres commutatives vers les groupes
\[R \mapsto \;(1_A+ 
\sR(A,A)\otimes_K R, \; \circ)\]
est re\-pr\'e\-sen\-ta\-ble par un 
$K$-sch\'ema en groupes affine unipotent
$\U^A$. Soit $\displaystyle 
\U_\sA=
\prod_A \U^A$ leur produit; c'est un $K$-groupe 
pro-unipotent,
\ie une limite projective filtrante $\lim
\U_\alpha$ 
de
$K$-groupes unipotents; en tant que sch\'ema, c'est juste un 
produit
d'espaces affines.

Pour tout couple $(A,B)$ d'objets de 
$\sA$, soit $Sec(A,B)$ l'ensemble
des sections $K$-lin\'eaires de la 
projection $\sA(A,B)\to \bar
\sA(A,B)$: il est 
re\-pr\'e\-sen\-ta\-ble par un $K$-espace affine
$\Sec(A,B)$. Posons 
$\Sec_\sA=\prod_{A,B}\Sec(A,B)$: c'est aussi un
produit d'espaces 
affines.

Consid\'erons le sous-sch\'ema ferm\'e de 
$\Sec$
\begin{multline*}
S_\sA=\{(s_{A,B})\mid \forall A,B,C\in \sA, 
\forall
(f,g)\in \bar\sA(A,B)\times \bar\sA(B,C),\\
s_{A,C}(g\circ 
f)=s_{B,C}(g)\circ s_{A,B}(f)\}.
\end{multline*}

On a une action 
\'evidente de $\U_\sA$ sur $S_\sA$
\begin{align*}
\mu:\U_\sA\times_K 
S_\sA&\to S_\sA\\
(u,s)&\mapsto usu^{-1}
\end{align*}
qui est 
\emph{transitive} gr\^ace au th\'eor\`eme \ref{T1} b). 
Notons
$\Delta(\sA)$ ce $K$-espace homog\`ene unipotent: le 
groupo\"{\i}de
$\Gamma(\sA)$ cherch\'e est l'image de $\Delta(\sA)$ 
par le foncteur de
\ref{ex3}.

Il reste \`a voir que $\Delta$ 
d\'efinit un $2$-foncteur. Soit $T:\sB\to
\sA$ un foncteur pleinement 
fid\`ele: il induit un
\'epimor\-phis\-me $\U_\sA\surj \U_\sB$ de 
$K$-sch\'emas affines gr\^ace au
diagramme
\[\begin{CD}
&& 
\displaystyle\prod_{A\in\sA}\U^A\\
&&@VVV\\
\displaystyle\prod_{B\in\s 
B} \U^B@>\sim>>\displaystyle\prod_{B\in \sB}
\U^{T(B)}
\end{CD}\]
et 
de m\^eme un \'epimor\-phis\-me $\Sec_\sA\surj
\Sec_\sB$, ce dernier 
envoyant $S_\sA$ dans $S_\sB$. Ces mor\-phis\-mes
d\'efinissent 
clairement un mor\-phis\-me de $K$-espaces 
homog\`enes
$\Delta(T):\Delta(\sA)\to
\Delta(\sB)$.

Soit $T'$ un 
autre tel foncteur, et $u:T\Rightarrow T'$ un isomor\-phis\-me
naturel, 
d'o\`u un autre $\bar u:\bar T\Rightarrow \bar T'$, o\`u $\bar
T,\bar 
T':\bar \sB\to \bar \sA$ sont les foncteurs induits. Pour tout
$B\in 
\sB$ et toute section $s$ de $\pi_\sA$, l'\'el\'ement 
$s(\bar
u_B)^{-1} u_B$ appartient \`a $1+\sR(T(B),T(B))$: on en 
d\'eduit un
\'el\'ement $T^{-1}(s(\bar
u_B)^{-1} u_B)\in 1+\sR(B,B)$. 
Cette r\`egle d\'efinit un mor\-phis\-me
\[\Delta(u):S_\sA\to 
\U_\sB\]
dont on v\'erifie facilement qu'il satisfait les conditions 
d'une
homotopie $\Delta(T)\Rightarrow 
\Delta(T')$.\qed

\begin{ex}\label{typehomotop} Le th\'eor\`eme 
\ref{gerbe} montre que,
pour deux
$K$-ca\-t\'e\-go\-ries de Wedderburn 
$\sA$ et $\sB$ $K$-\'equi\-va\-len\-tes, les
groupo\"{\i}des $\Gamma(\sA)$ 
et $\Gamma(\sB)$ ont en g\'en\'eral m\^eme
type  d'homotopie (au sens 
ci-dessus), mais ne sont pas
n\'e\-ces\-sai\-re\-ment isomorphes. La 
m\^eme remarque vaut pour les
invariants plus fins
$\Delta(\sA)$ et 
$\Delta(\sB)$. La remarque \ref{r3} va nous permettre de
d\'ecrire un 
repr\'esentant explicite du type d'homotopie de
$\Delta(\sA)$ et 
$\Gamma(\sA)$ lorsque $K$ est al\-g\'e\-bri\-que\-ment
clos (ou plus 
g\'en\'eralement parfait, si tout ind\'e\-com\-po\-sa\-ble de $\sA$
le reste 
dans $\sA^\natural$ apr\`es toute extension finie $L/K$).

Il suffit 
de supposer que tout objet de $\sA$ est ind\'e\-com\-po\-sa\-ble et que
deux 
objets diff\'erents sont non isomorphes. On a alors $S=\Spec K$ 
et
$G$ est le produit des $1+\sR(A,A)$, o\`u
$A$ d\'ecrit l'ensemble 
des objets de $\sA$. En particulier, le sch\'ema
des objets de 
$\Gamma(\sA)$ est ici r\'eduit \`a un 
point.
\end{ex}

\subsubsection{Le cas 
mo\-no\-\"{\i}\-dal}\label{repmon}

\begin{sloppypar}
\begin{thm}\label{gerbe'} Soit $Wedf^\otimes\{K\}$ la
sous-$2$-ca\-t\'e\-go\-rie $1$-pleine et 
$2$-plei\-ne de $Wed^\otimes\{K\}$
form\'ee des ca\-t\'e\-go\-ries
$\sA$ 
telles que $\dim_K\sA(A,B)<\infty$ pour tout 
$(A,B)\in
\sA\times\sA$.\\
a) Soit $\sA\in Wedf^\otimes\{K\}$. Alors 
il existe un
$K$-groupo\"{\i}de affine scind\'e 
unipotent
$\Gamma^\otimes(\sA)=(E^\otimes,S^\otimes)$, transitif 
sur
$S^\otimes$, tel que pour toute extension
$L$ de $K$, on 
ait
$S^\otimes(L)=Ob(\sG^\otimes(\sA_L))$ 
et
$E^\otimes(s_1,s_2)(L)=\sG^\otimes(\sA_L)(s_1,s_2)$ pour tout 
couple de
sections mo\-no\-\"{\i}\-dales
$(s_1,s_2)$ de $\pi_{\sA_L}$.\\
b) 
Cette construction provient d'un $2$-foncteur $1$-contravariant 
et
$2$-co\-va\-riant
\[\Delta^\otimes:Wedf^\otimes\{K\}\to 
\Hmg(Sch/K)\]
\`a valeurs dans les $K$-espaces homog\`enes 
unipotents, via le
$2$-foncteur
$\Hmg(Sch/K)\to
\Gr(Sch/K)$ 
de
\ref{ex3}.\\
c) On a une $2$-transformation 
naturelle
\[\Delta^\otimes\Rightarrow\Delta\circ \Omega\] 
comme dans 
le sorite \ref{so6t}.\\
d) Les $2$-foncteurs $\Delta^\otimes$ et 
$\Gamma^\otimes$ commutent \`a
l'extension des 
scalaires.
\end{thm}
\end{sloppypar}

\prf On reprend les notations 
de la preuve du theor\`eme \ref{gerbe}.
Consid\'erons le nouvel 
espace de 
param\`etres
\[\Sec^\otimes=S\times\prod_{A,B}(1+\sR(A\bullet
B,A\bullet B)).\]

Les \'equations que doivent v\'erifier la contrainte 
mo\-no\-\"{\i}\-dale
$\tilde s$ d'une section mono\"{\i}\-dale $s$ 
d\'efinissent un sous-sch\'ema
ferm\'e $S^\otimes\subset 
\Sec^\otimes$. Soit $\U_1$ le sous-groupe
ferm\'e de $\U$ form\'e des 
$(u_A)$ tels que $u_\un=1$. Le groupe
$\U_1$ op\`ere cette fois-ci 
sur $S^\otimes$ par la 
formule
\begin{align*}
\mu^\otimes:\U_1\times_K S^\otimes&\to 
S^\otimes\\
(u,(s,\tilde s))&\mapsto (usu^{-1},\tilde 
s')
\end{align*}
o\`u $\tilde s'$ est d\'efini par le diagramme 
\eqref{eq4}. Le
th\'eor\`eme \ref{T2} b) montre que cette action est 
\emph{transitive}.
D'o\`u $\Delta^\otimes(\sA)$.

On a un mor\-phis\-me 
\'evident $\Delta^\otimes(\sA)\to \Delta(\sA)$, donn\'e
par le 
plongement de $\U_1$ dans $\U$ et par l'oubli des 
structures
mo\-no\-\"{\i}\-dales. Ceci d\'efinit la $2$-transformation 
naturelle de c)
au niveau des objets (ca\-t\'e\-go\-ries). Le reste 
se
v\'erifie facilement.\qed

\section{Sections et 
tressages}\label{tre}

\subsection{Tressages}\label{Tre} On renvoie 
\`a \cite{js} pour les 
notions de tressage et de foncteur 
mo\-no\-\"{\i}\-dal tress\'e. Contentons-nous
de les rappeler 
bri\`evement.

Un \emph{tressage} sur une ca\-t\'e\-go\-rie 
mo\-no\-\"{\i}\-dale
$(\sA,\bullet)$ est la donn\'ee, pour tout couple 
d'objets
$(A,B)\in\sA\times\sA$, d'un \'el\'ement
\[R_{A,B}\in 
\sA(A\bullet B,B\bullet A)\]
tel que le 
diagramme
\begin{equation}\label{pretress}\begin{CD}
A\bullet 
B@>f\bullet g>> C\bullet D\\
@VR_{A,B}VV @VR_{C,D}VV\\
B\bullet 
A@>g\bullet f>> D\bullet C
\end{CD}
\end{equation}
soit commutatif 
pour tous $A,B,C,D\in\sA$ et tout couple 
$(f,g)\in
\sA(A,C)\allowbreak\times\sA(B,D)$. On demande de plus que 
les $R_{A,B}$
v\'erifient les \emph{identit\'es de 
tressage}
\begin{align}\label{tress1}
R_{A,B\bullet C} 
&=a_{B,C,A}(1_B\bullet R_{A,C})a_{B,A,C}^{-1}
(R_{A,B}\bullet 
1_C)a_{A,B,C}\\ 
R_{A\bullet B, C} &=a_{C,A,B}^{-1}(R_{A,C}\bullet 
1_B)a_{A,C,B}
(1_A\bullet R_{B,C})a_{A,B,C}^{-1}\label{tress2}
\end{align}
o\`u $a$ est la contrainte d'associativit\'e.

On dit que le tressage $R$ est \emph{sy\-m\'e\-tri\-que} si on a de plus
l'identit\'e
\begin{equation}\label{tress3}
R_{B,A}=R_{A,B}^{-1}.
\end{equation}

Si $(\sB,\top)$ et $(\sA,\bullet)$ sont
tress\'ees, un foncteur mo\-no\-\"{\i}\-dal $(s,\tilde
s):(\sB,\top)\to (\sA,\bullet)$ est dit
\emph{tress\'e} s'il est compatible avec les tressages. En d'autres
termes, pour tous
$A,B\in
\sB$, on demande que le diagramme
\begin{equation}\begin{CD}\label{tress}
s(A)\bullet s(B)@>\tilde s_{A,B}>> s(A\top B)\\
@VR_{s(A),s(B)} VV @Vs(R_{A,B})VV\\
s(B)\bullet s(A)@>\tilde s_{B,A}>> s(B\top A)
\end{CD}
\end{equation}
soit commutatif.

Notons que si $s$ est strict, la condition devient
simplement $s(R_{A,B})=R_{s(A), s(B)}$.

Notons aussi qu'imposer que les tressages de $\sB$ et $\sA$
soient sy\-m\'e\-tri\-ques
n'impose aucune condition suppl\'ementaire sur $s$ (et m\^eme en retire).
Par contre, si les tressages de $\sB$ et $\sA$ sont compatibles via 
$s$, et si $s$ est essentiellement
surjective, l'un est sy\-m\'e\-tri\-que si et seulement si l'autre l'est.

\subsection{Sorites}\label{11.2} Le sorite suivant compl\`ete le sorite
\ref{so1}:

\begin{sorite}\label{so4} Soit $(s,\tilde s):(\sB,\top,R)\to
(\sA,\bullet,S)$ un foncteur mo\-no\-\"{\i}\-dal tress\'e entre deux
ca\-t\'e\-go\-ries mo\-no\-\"{\i}\-dales tress\'ees, et soit
$u:s\Rightarrow t$ un isomor\-phis\-me naturel de $s$ sur un autre foncteur
$t$. Alors l'unique structure mo\-no\-\"{\i}\-dale $\tilde t$ sur
$t$ faisant de $u$ un mor\-phis\-me de foncteurs mono\"{\i}daux (sorite
\ref{so1}) est tress\'ee.\qed
\end{sorite}

On a aussi:

\begin{sorite}\label{so5} Tout tressage $R$ satisfait les identit\'es
\begin{align*}
R_{A,B\oplus C}&=\diag(R_{A,B},R_{A,C})\\
R_{A\oplus B,C}&=\diag(R_{A,C},R_{B,C}).
\end{align*}
\end{sorite}

(Cela r\'esulte du sorite \ref{so0}, en consid\'erant $R$ comme une
transformation naturelle $\bullet \Rightarrow \bullet \circ T$, o\`u $T$
est l'\'echange des facteurs.)\qed

Rappelons \'egalement le sorite suivant \cite[rem. 2.3]{js}:

\begin{sorite}\label{so8} Notons $\sA^{sym}$ la ca\-t\'e\-go\-rie sy\-m\'e\-tri\-que
de $\sA$ (m\^emes objets, m\^emes mor\-phis\-mes mais produit d\'efini par
$A\bullet' B:=B\bullet A$, \cf \cite[I.0.1.4]{saavedra}). Alors tout
tressage $R$ sur $\sA$ d\'efinit une structure mo\-no\-\"{\i}\-dale sur le
foncteur identique $Id_\sA:\sA^{sym}\to \sA$.\qed
\end{sorite}

Enfin:

\begin{sorite}\label{so9} Soient $(\sB,\top,R)$ une ca\-t\'e\-go\-rie
mo\-no\-\"{\i}\-dale tress\'ee, $(\sA,\bullet)$ une cat\'ego\-rie mo\-no\-\"{\i}\-dale
et $(s,\tilde s):(\sB,\top)\to (\sA,\bullet)$ un foncteur mo\-no\-\"{\i}\-dal.
Pour $A,B\in\sB$, posons
\[R^s_{A,B}=\tilde s_{B,A}^{-1}s(R_{A,B})\tilde s_{A,B}:s(A)\bullet
s(B)\to s(B)\bullet s(A).\] Soit $a$ la contrainte d'associativit\'e de
$\sA$. Alors on a les identit\'es
\begin{multline*}
(\tilde s_{B,C}\bullet 1_{s(A)})^{-1}R^s_{A,B\bullet C}(1_{s(A)}\bullet
\tilde s_{B,C}) =\\
a_{s(B),s(C),s(A)}(1_{s(B)}\bullet
R^s_{A,C})a_{s(B),s(A),s(C)}^{-1} (R_{A,B}\bullet
1_{s(C)})a_{s(A),s(B),s(C)}
\end{multline*}
\begin{multline*}
(1_{s(C)}\bullet\tilde s_{A,B})^{-1}R^s_{A\bullet B, C}(\tilde
s_{A,B}\bullet 1_{s(C)}) =\\
a_{s(C),s(A),s(B)}^{-1}(R_{A,C}\bullet
1_{s(B)})a_{s(A),s(C),s(B)} (1_{s(A)}\bullet
R^s_{B,C})a_{s(A),s(B),s(C)}^{-1}.
\end{multline*}
\qed
\end{sorite}

En d'autres termes, $R^s$ v\'erifie ``presque" les relations de tressage
\eqref{tress1} et \eqref{tress2}.

\subsection{Sections}\label{sect} Reprenons maintenant les notations et
hypoth\`eses du th\'e\-o\-r\`e\-me
\ref{T2}, en pr\'esence d'un tressage\footnote{On rappelle que 
l'hypoth\`ese de commutativit\'e des
$End(\un )$-bimodules $\sA(A,B)$ est alors automatiquement 
satisfaite.}. Ce dernier reste muet sur la compatibilit\'e des 
sections
mo\-no\-\"{\i}\-da\-les
\`a un tressage \'eventuel. Nous allons voir que, dans ce cas, la
situation est tr\`es ``rigide".

Soit $R$ un tressage sur
$\sA$.
On en d\'eduit un tressage
$\bar R=\pi_\sA(R)$ sur $\bar\sA$.

Soit $(s, \tilde s)$ une section mo\-no\-\"{\i}\-dale de
$\pi_\sA$.

\begin{prop}\label{p9} a) Il existe un automor\-phis\-me (non mo\-no\-\"{\i}\-dal)
$u$ de $s$ tel que, pour tout couple d'objets $A,B\in \sA$, on ait
\begin{equation}\label{eq16}
s(\bar R_{A,B})=u_{B\bullet A}^{-1}\tilde s_{B,A} R_{A,B}
(u_A\bullet u_B) \tilde s_{A,B}^{-1}.
\end{equation}
b) On a l'identit\'e
\begin{equation}\label{eq17}
s(\bar R_{B,A})s(\bar R_{A,B})=\tilde s_{A,B}R_{B,A}R_{A,B} (u_A^2\bullet
u_B^2)\tilde s_{A,B}^{-1} u_{A\bullet B}^{-2}.
\end{equation}
\end{prop}

\prf  a) En tenant compte du sorite \ref{so8}, on applique 
la
proposition \ref{P3'} dans la situation $\sB=\sA^{sym}$, 
$T=(Id_\sA,R)$,
$s_\sA=s_\sB=s$. 

b) Remarquons les
relations de 
commutation
\[u_{B\bullet A} s(\bar R_{A,B})=s(\bar R_{A,B}) 
u_{A\bullet B},\]
due au fait que $u$ est un automor\-phis\-me du 
foncteur $s$, et
\[R_{A,B} (u_A\bullet u_B)=(u_B\bullet u_A) 
R_{A,B},\]
due au fait que $R$ est un tressage. On d\'eduit alors de 
\eqref{eq16}
\begin{equation}\label{eq17,5}
s(\bar R_{A,B})=\tilde 
s_{B,A} R_{A,B}
(u_A\bullet u_B) \tilde s_{A,B}^{-1}u_{A\bullet 
B}^{-1}
\end{equation}
d'o\`u
\begin{multline*}
s(\bar R_{B,A})s(\bar R_{A,B})\\
=u_{A\bullet B}^{-1}\tilde s_{A,B} R_{B,A}
(u_B\bullet u_A) \tilde s_{B,A}^{-1}\tilde s_{B,A} R_{A,B}
(u_A\bullet u_B) \tilde s_{A,B}^{-1}u_{A\bullet B}^{-1}\\
=u_{A\bullet B}^{-1}\tilde s_{A,B} R_{B,A}
(u_B\bullet u_A) R_{A,B}
(u_A\bullet u_B) \tilde s_{A,B}^{-1}u_{A\bullet B}^{-1}\\
=u_{A\bullet B}^{-1}\tilde s_{A,B} R_{B,A}
R_{A,B}
(u_A^2\bullet u_B^2) \tilde s_{A,B}^{-1}u_{A\bullet B}^{-1}\\
=\tilde s_{A,B} R_{B,A}
R_{A,B}
(u_A^2\bullet u_B^2) \tilde s_{A,B}^{-1}u_{A\bullet B}^{-2}
\end{multline*}
comme annonc\'e.\qed

\begin{rems}\ \label{r5}
\begin{itemize}
\item[a)] On prendra garde au fait que les
relations de commutation apparaissant dans la preuve de la proposition
\ref{p9} a) n'ont pas de raison d'\^etre vraies en
\'echangeant $R_{A,B}$ et $s(\bar R_{A,B})$. De m\^eme, $u_{A\bullet B}$
et $u_A\bullet u_B$ n'ont pas de raison de commuter.
\item[b)] Dans le sorite \ref{so8}, la condition, pour un isomor\-phis\-me
naturel $R:\bullet\Rightarrow \bullet\circ T$, de d\'efinir une structure
mo\-no\-\"{\i}\-dale sur le foncteur identique $\sA^{sym}\to \sA$, est
\emph{plus faible} que les conditions de tressage. Il en r\'esulte que
l'automor\-phis\-me $u$ n'est pas arbitraire: il v\'erifie deux conditions
suppl\'ementaires du type ``cross-effects" non commutatifs d'ordre 3, qui
reviennent \`a dire que la fonction $A\mapsto u_A$ est ``de degr\'e
$\le 2$" pour la structure mo\-no\-\"{\i}\-dale. Nous laissons au lecteur le
soin d'expliciter ces relations, que nous n'utiliserons pas.
\end{itemize}
\end{rems}

\begin{prop} \label{p10} Pour  une section mo\-no\-\"{\i}\-dale $(s,\tilde s)$,
les conditions sui\-van\-tes sont \'equi\-va\-len\-tes:
\begin{thlist}
\item $\bar R^s$ est un
tressage (notation du sorite \ref{so9}).
\item Il existe un
automor\-phis\-me (non mo\-no\-\"{\i}\-dal) $v$ de $Id_\sA$ tel qu'on ait, pour $u$
comme dans la proposition \ref{p9}, l'identit\'e
\[(u_A\bullet u_B)\tilde s_{A,B}^{-1} u_{A\bullet B}^{-1} \tilde
s_{A,B}=(v_A\bullet v_B)v_{A\bullet B}^{-1}.\]
\end{thlist}
Si ces conditions sont v\'erifi\'ees, le tressage $\bar R^s$ ne
d\'epend pas du choix de $s$: pour un choix de
$v$ comme dans {\rm (ii)}, il est donn\'e par la formule
\[\tilde 
R_{A,B}=R_{A,B} (v_A\bullet v_B)v_{A\bullet 
B}^{-1}.\]
\end{prop}

\prf En utilisant la formule \eqref{eq17,5}, 
on obtient l'identit\'e
\[R_{A,B}^{-1}\bar R^s_{A,B}
=(u_A\bullet 
u_B) \tilde s_{A,B}^{-1} u_{A\bullet B}^{-1} \tilde
s_{A,B}.\]

Si 
(i) est v\'erifi\'e, le premier membre de cette identit\'e 
d\'efinit
une structure mo\-no\-\"{\i}\-dale sur $Id_\sA:\sA^{sym}\to 
\sA^{sym}$ couvrant
la structure mo\-no\-\"{\i}\-dale triviale sur 
$Id_{\bar \sA}$. (ii) r\'esulte
alors de la proposition \ref{t3}. 
Pour la r\'eciproque, remarquons que
l'identit\'e de (ii) implique 
que $\bar R^s$ fait commuter le diagramme
\eqref{pretress} pour tous 
mor\-phis\-mes de $\sA$. En particulier, dans le
sorite \ref{so9} le 
premier membre des deux identit\'es se r\'eduit
respectivement \`a 
$\bar R^s_{A,B\bullet C}$ et \`a $\bar R^s_{A\bullet
B,C}$: $\bar 
R^s$ est donc bien un tressage.

Soit $(t,\tilde t)$ une autre 
section mo\-no\-\"{\i}\-dale. D'apr\`es
le th\'eor\`eme \ref{T2} b), $s$ et 
$t$ sont mo\-no\-\"{\i}\-dalement
conjugu\'ees. En d'autres termes, il 
existe une famille $(v_A)$
conjuguant $s$ et $t$ et v\'erifiant 
l'identit\'e
\[v_{A\bullet B}\tilde s_{A,B}=\tilde t_{A,B} v_A\bullet 
v_B.\]

On a alors
\[\bar R^t_{A,B}=(v_B^{-1}\bullet v_A^{-1})\bar 
R^s_{A,B}(v_A\bullet
v_B).\]

Ainsi, si $\bar R^s$ est un tressage, 
on a $\bar R^t=\bar R^s$.

La derni\`ere assertion est claire \`a 
partir de la formule
\eqref{eq17,5}.\qed

Nous ne connaissons pas de 
crit\`ere pour que les conditions de la
proposition \ref{p10} soient 
v\'erifi\'ees, \`a l'exception d'un cas:
celui o\`u le tressage 
r\'esiduel $\bar R$ est \emph{sy\-m\'e\-tri\-que}.
Rappelons la 
d\'efinition suivante \cite[\S 
6]{js}:

\begin{sloppypar}
\begin{defn}\label{d4} Une \emph{structure 
balanc\'ee} sur une ca\-t\'e\-go\-rie
mo\-no\-\"{\i}\-dale tress\'ee 
$(\sA,\bullet,R)$ est un automor\-phis\-me
$\theta=(\theta_A)_{A\in \sA}$ 
du foncteur identique de $\sA$ 
v\'erifiant
l'identit\'e
\[
R_{B,A}R_{A,B}=(\theta_A\bullet 
\theta_B)\theta_{A\bullet
B}^{-1}.
\]
 
\end{defn}
\end{sloppypar}

\begin{thm}\label{t4} Soit $(\sA,\bullet, 
R)$ une ca\-t\'e\-go\-rie
mo\-no\-\"{\i}\-dale tress\'ee. On suppose que le 
radical est un id\'eal
mo\-no\-\"{\i}\-dal, et on note $(\bar\sA,\bullet, 
\bar R)$ la ca\-t\'e\-go\-rie
mo\-no\-\"{\i}\-dale tress\'ee obtenue en 
quotientant par le radical.\\ a) Si
$\bar R$ est sy\-m\'e\-tri\-que, $R$ 
est balanc\'e via un automor\-phis\-me
$\theta$ v\'erifiant 
$\pi_\sA(\theta)=1$.\\ b) Si de plus $\car K\ne 2$,
les conditions de 
la proposition \ref{p10} sont v\'e\-ri\-fi\'ees. Le
tressage $\tilde 
R$ correspondant est l'unique tressage sy\-m\'e\-tri\-que de
$\sA$ tel 
que
\begin{thlist}
\item $\pi_\sA(\tilde R)=\bar 
R$;
\begin{sloppypar}
\item $(Id_\sA,\tilde R)$ est 
mo\-no\-\"{\i}\-dalement isomorphe \`a
$(Id_\sA,R)$ (\cf sorite 
\ref{so8}).
\end{sloppypar}
\end{thlist} 
c) Si de plus $R$ est 
sy\-m\'e\-tri\-que, on a
$s(\bar R)=R$ pour toute section mo\-no\-\"{\i}\-dale 
$s$.
\end{thm}

\prf \begin{sloppypar} a) On applique la proposition \ref{t3} au 
foncteur mo\-no\-\"{\i}\-dal
$(Id_\sA, R_{A,B}R_{B,A})$.

b) L'identit\'e 
\eqref{eq17} de la proposition \ref{p9} donne une
nouvelle 
identit\'e
\[(\theta_A\bullet 
\theta_B)\theta_{A\bullet
B}^{-1}=\tilde s_{A,B}^{-1} 
u_{A\bullet
B}^2\tilde s_{A,B} (u_A^{-2}\bullet u_B^{-2}).\]

Comme 
$\car K\ne 2$, l'\'el\'evation au carr\'e est bijective dans
$1+\sR$. 
En particulier,
$\theta$ a une unique racine carr\'ee $v$. Comme 
$\theta_{A\bullet B}$ et
$\theta_A\bullet\theta_B$ sont centraux dans 
$\sA(A\bullet B)$, $\tilde
s_{A,B}^{-1} u_{A\bullet B}^2\tilde 
s_{A,B}$ et $(u_A^{-2}\bullet
u_B^{-2})$ commutent. On a 
donc
\[(v_A\bullet v_B)v_{A\bullet
B}^{-1}=\tilde s_{A,B}^{-1} 
u_{A\bullet
B}\tilde s_{A,B} (u_A^{-1}\bullet u_B^{-1})\]
ce qui 
n'est autre que la condition (iii) de la proposition
\ref{p10}.

Ceci d\'emontre que $\tilde R$ v\'erifie les conditions (i) et
(ii) 
du th\'eor\`eme. Pour voir l'unicit\'e, soit $\tilde R'$ un 
autre
tressage sy\-m\'e\-tri\-que v\'erifiant ces conditions. En 
particulier,
$(Id_\sA,\tilde R')$ est mo\-no\-\"{\i}\-dalement isomorphe 
\`a $(Id_\sA,\tilde
R)$. Il existe donc $\mu\in 1+\sR$ tel 
que
\[\tilde R'_{A,B}=\tilde R_{A,B}\mu_{A\bullet 
B}(\mu_A^{-1}\bullet
\mu_B^{-1})\]
pour tout couple d'objets $(A,B)$. 
La sym\'etrie de $\tilde R$ et de
$\tilde R'$ donne maintenant 
l'identit\'e
\[\mu_{A\bullet B}^2=\mu_A^2\bullet\mu_B^2.\]

En en 
prenant la racine carr\'ee, on trouve bien que $\tilde R'=\tilde 
R$.

c) Cela r\'esulte de la derni\`ere assertion de la 
proposition
\ref{p10}. \end{sloppypar}\qed

\begin{ex}\label{ex4} \begin{sloppypar}

 Consid\'erons la quantification de 
Drinfeld-Cartier $\bar\sA[[h]]$
d'une cat\'e\-gorie mo\-no\-\"{\i}\-dale 
stricte sy\-m\'e\-tri\-que
$\bar\sA$ munie d'un ``tressage infinit\'esimal" 
$t_{A,B}$, \cf \cite[ch.
9]{krt}. Rappelons seulement que 
$\bar\sA[[h]]$ a les m\^emes objets que
$\bar\sA$, et 
que
$\bar\sA[[h]](A,B)=(\bar\sA(A,B))[[h]]$; le tressage 
de
$\bar\sA[[h]]$ est de la forme
$R_{A,B}=\bar R_{A,B}\exp(\frac{h}{ 
2} t_{A,B})$ (et la contrainte 
d'associativit\'e est construite \`a 
partir de
$t_{A,B}$ au moyen d'un associateur de Drinfeld).  En 
g\'en\'eral ce
tressage n'est pas
sy\-m\'e\-tri\-que, m\^eme modulo
$h^2$. 
Toutefois, le th\'eor\`eme \ref{t4} b) montre par un passage
\`a la 
limite que $(\bar \sA[[h]],R)$ admet une structure balanc\'ee, 
et
m\^eme que $\bar \sA[[h]]$ admet \emph{un unique tressage 
sy\-m\'e\-tri\-que
$\tilde R$} se r\'eduisant en
$\bar R$ et tel 
que
$(Id,\tilde R)$ soit mo\-no\-\"{\i}\-dalement isomorphe \`a $(Id,R)$ 
dans le
cadre du sorite \ref{so8}. 

En particulier il n'est pas 
n\'ecessaire de supposer
l'existence de duaux dans $\bar \sA$ pour 
construire une structure
balanc\'ee sur $R$, comme dans 
\cite{cartier}.
\end{sloppypar}
\end{ex}

\begin{rem} Il est probable 
que les parties b) et c) du th\'eor\`eme sont
 fausses en 
caract\'eristique 2. Il serait int\'eressant d'exhiber
un 
contre-exemple.
\end{rem}
 
 \subsection{Compl\'ements}
Pour finir, 
examinons de plus pr\`es le d\'efaut de validit\'e des
conditions de 
la proposition \ref{p10} et le d\'efaut de commutation 
du
diagramme
\eqref{tress}. On pose 
$u_{A,B}=R_{A,B}^{-1}(\tilde
s_{B,A})^{-1}s(\bar  R_{A,B})\tilde 
s_{A,B}$
et $n_{A,B}=1_{A\bullet B}-u_{A,B}$. On a $n_{A,B}\in 
C_s^1(AB)$
(\cf \S \ref{prb} pour la 
notation).

Soient
\begin{align*}
e&= e(R)=\sup\{r\mid n_{A,B}\in 
C_s^r(A B)\text{ pour tout } (A,B)\in
\sA\times \sA\}\\
\bar e&=\bar 
e(\bar R)=\sup\{r\mid [1_A\bullet m_B+m_A\bullet
1_B,n_{A,B}]\in 
\sR(A\bullet B)^r\\
&\text{ pour tout } (A,B)\in
\sA\times \sA \text{ 
et tout } (m_A,m_B)\in 
\sR(A,A)\times
\sR(B,B)\}.
\end{align*}

Notons $\bar n_{A,B}$ 
l'image de $n_{A,B}$ dans $C^{[e]}(A
B)$ (notation du \S 
\ref{prb}).

\begin{prop}\label{ab} a) L'\'el\'ement $e=e(R)\in 
\N\cup\{\infty\}$ et les \'el\'ements $\bar n_{A,B}$ 
de
$C^{[e]}(AB)$ ne d\'ependent pas du choix de $s$: ils ne 
d\'ependent
que de $R$. Si $\sA$ est stricte, on a 
les
identit\'es
\begin{align}\label{eq5}
\bar n_{A,B\bullet C} &= 
(\bar R_{A,B}\bullet 1_C)^{-1}(1_B\bullet
\bar n_{A,C})(\bar 
R_{A,B}\bullet 1_C) + \bar n_{A,B}\bullet 1_C\\
\bar n_{A\bullet B,C} 
&= (1_A\bullet \bar R_{B,C})^{-1}(\bar
n_{A,C}\bullet 1_B) 
(1_A\bullet \bar R_{B,C}) + 1_A\bullet
\bar 
n_{B,C}\label{eq6}\\
\intertext{la seconde pouvant aussi 
s'\'ecrire}
\bar n_{A\bullet B,C} &= (\bar R_{B,A}\bullet 
1_C)(1_B\bullet \bar n_{A,C})(\bar R_{B,A}\bullet
1_C)^{-1} + 
1_A\bullet
\bar n_{B,C}.\label{eq6bis}
\end{align} 
b) L'\'el\'ement 
$\bar e\in \N\cup\{\infty\}$ ne d\'epend
pas du choix de $s$: il ne 
d\'epend que de $\bar R$. On a $\bar e>e$,
et $\bar e=\infty$ si et 
seulement si les conditions de la proposition
\ref{p10} sont 
v\'erifi\'ees.\\ 
c) Soit $T: (\sA',\bullet, R')\to (\sA,\bullet, R)$ 
un foncteur 
mo\-no\-\"{\i}\-dal strict tress\'e radiciel 
(avec
$(\sA',\bullet)$ stricte). Alors $e(R')\leq 
e(R)$.
\end{prop}

\begin{sloppypar}
\prf $a)$ Soit $s'$ une autre section 
mo\-no\-\"{\i}\-dale. D'apr\`es le
th\'e\-o\-r\`e\-me
\ref{T2} b), $s$ et 
$s'$ sont conjugu\'ees par une famille $(v_A)_{A\in
\sA}$ 
d'\'e\-l\'e\-ments de $1+\sR(A)$ v\'erifiant 
l'identit\'e
$(v_{A\bullet B})\tilde s_{A,B}=\tilde 
t_{A,B}(v_A\bullet v_B)$. Si 
on
note
$u'_{A,B}=R_{A,B}^{-1}(\tilde
s'_{B,A})^{-1}s'(\bar 
R_{A,B})\tilde s'_{A,B}$, on en 
d\'eduit
l'identit\'e
\[u'_{A,B}=(v_A\bullet v_B) u_{A,B} (v_A\bullet 
v_B)^{-1}.\]
\end{sloppypar}

Si on pose $n'_{A,B}=1_{A\bullet B}-u'_{A,B}$, cette 
identit\'e implique
d'abord que $n'_{A,B}\in \sR(A\bullet B)^e$, 
ensuite que
que
$n'_{A,B}\equiv n_{A,B}\pmod{\sR(A\bullet B)^{e+1}}$, 
et enfin que
son image $\mod{\sR(A\bullet
B)^{e+1}}$ est dans 
$C^{[e]}(AB)$.

Pour \'etablir
\eqref{eq5} et
\eqref{eq6}, on peut 
supposer $s$ stricte
par le point $c)$ de la proposition
\ref{p5}.
En 
appliquant $s$ \`a \eqref{tress1} et \eqref{tress2}, on obtient les
identit\'es
\begin{align*}
s(\bar R_{A,B\bullet C}) &=(1_B\bullet s(\bar R_{A,C}))
(s(\bar R_{A,B})\bullet 1_C)\\
s(\bar R_{A\bullet B, C}) 
&=(s(\bar R_{A,C})\bullet 1_B)
(1_A\bullet s(\bar 
R_{B,C})).
\end{align*}
   
On en 
d\'eduit
\begin{align*}
u_{A,B\bullet C} &=(R_{A,B}\bullet 
1_C)^{-1}(1_B\bullet
u_{A,C})(R_{A,B}\bullet 1_C) (u_{A,B}\bullet 
1_C)\\
u_{A\bullet B, C}&=(1_A\bullet R_{B,C})^{-1}(u_{A,C}\bullet 
1_B)
(1_A\bullet R_{B,C})(1_A\bullet 
u_{B,C})
\end{align*}

\begin{sloppypar}
Les identit\'es \eqref{eq5} 
et \eqref{eq6} s'en d\'eduisent facilement.
En utilisant l'identit\'e 
$R_{B,A\bullet C} =(1_A\bullet R_{B,C})
(R_{B,A}\bullet 1_C)$ et 
l'\'egalit\'e $n_{A,C}\bullet 1_B=\bar
R_{B,A\bullet C}(1_B\bullet 
\bar n_{A,C})\bar R_{B,A\bullet C}^{-1}$, on
peut r\'e\'ecrire le 
premier terme du membre de droite de
\eqref{eq6}:
\[(\bar 
R_{B,A}\bullet 1_C)(1_B\bullet \bar n_{A,C})(\bar 
R_{B,A}\bullet
1_C)^{-1}.
\]
\end{sloppypar}

$b)$ r\'esulte 
facilement des calculs de $a)$.

$c)$ Soit $s$ (\resp $s'$) une 
section
mo\-no\-\"{\i}\-dale stricte de
$\pi_{\sA}$ (\resp $\pi_{\sA'}$) 
(il en existe d'apr\`es la 
proposition \ref{p5}). Posons
 
\[u_{T(A),T(B)}=(R_{T(A),T(B)})^{-1} s(\bar 
R_{T(A),T(B)}),\;
u'_{A,B}=(R'_{A,B})^{-1} s'(\bar R'_{A,B}) .\]

D'apr\`es la proposition \ref{P3'}, $Ts'$ et $s\bar T$ 
sont
conjugu\'ees par une famille
$(w_A)_{A\in \sA'}$
d'\'el\'ements 
de $1+\sR(T(A))$ v\'erifiant l'identit\'e
\[w_{A\bullet B}=w_A\bullet 
w_B.\]

   On d\'eduit de ces formules que
\[Tu'_{A,B}=(w_A\bullet 
w_B) u_{T(A),T(B)} (w_A\bullet w_B)^{-1}.\]
d'o\`u il est ais\'e de 
conclure.\qed

\begin{rem} \label{rdc} Les identit\'es \eqref{eq5} 
et \eqref{eq6} ne sont
autres que les relations de 
tressage
infinit\'esimales de Drinfeld-Cartier (\cf \cite[ch. 
9]{krt}),
v\'e\-ri\-fi\'ees par $ht_{A,B}$ de l'exemple
\ref{ex4}. 
Toutefois, $\bar n_{A,B}$ est ``antisy\-m\'e\-tri\-que" si $R$
est 
sy\-m\'e\-tri\-que, alors que $ht_{A,B}$ est ``sy\-m\'e\-tri\-que" dans 
la
quantification de Drinfeld-Cartier.
\end{rem}

\section{Premi\`ere 
application: ca\-t\'e\-go\-ries de Kimura}

\subsection{} Rappelons 
(d\'efinition \ref{defki}) qu'une ca\-t\'e\-go\-rie de Kimura sur un 
corps
de caract\'eristique nulle $K$ est une ca\-t\'e\-go\-rie 
$K$-lin\'eaire
mo\-no\-\"{\i}\-dale sy\-m\'e\-tri\-que rigide, v\'erifiant 
$End(\un)=K$,
pseudo-ab\'elienne, et dont tout objet est de dimension 
finie au sens de
Kimura. 

\begin{thm}\label{kimsplit} \begin{sloppypar} Soit $\sA$ une 
ca\-t\'e\-go\-rie de Kimura, de radical $\sR$. Alors\\
a) Le foncteur de 
projection $\sA\to \sA/\sR$ admet des sections
mo\-no\-\"{\i}\-dales, qui 
envoient le tressage r\'esiduel $\bar R$ de
$\sA/\sR$ sur le tressage 
$R$ de $\sA$. Deux telles sections sont
mo\-no\-\"{\i}\-dalement 
conjugu\'ees.\\ 
b) Soit $R'$ la contrainte de commutativit\'e 
de
$\sA/\sqrt[\otimes]{0}$ modifi\'ee
\`a l'aide de la 
$\Z/2$-graduation canonique de $\sA/\sqrt[\otimes]{0}$
(\cf 
th\'eor\`eme \ref{tkim.1}). Par ailleurs, soit $\pi_A^+\in
\sA/\sqrt[\otimes]{0}(A,A)$ l'idempotent (central) associ\'e
canoniquement \`a tout objet $A$ de $\sA$ (ibid.) Notons $\bar R'$ et
$\bar
\pi_A^+$ leurs images dans $\sA/\sR$. Alors, pour toute section
mo\-no\-\"{\i}\-dale $s:\sA/\sR\to \sA$ comme en a),
\begin{thlist}
\item  L'image de $s(\bar R')$ (\resp de $s(\bar\pi_A^+)$ pour tout $A$)
dans
$\sA/\sqrt[\otimes]{0}$ est \'egale \`a $R'$ (\resp \`a $\pi_A^+$).
\item La famille des $s(\bar\pi_A^+)$ d\'efinit une $\Z/2$-graduation sur
$\sA$, dont le tressage modifi\'e associ\'e est $s(\bar R')$. Pour ce
tressage sy\-m\'e\-tri\-que, la ca\-t\'e\-go\-rie $\sA$ reste rigide et la
dimension de tout objet $A$ est \'egale \`a $\kim A$; en particulier,
c'est un entier $\ge 0$.
\end{thlist}\end{sloppypar}
\end{thm}

\prf a) r\'esulte des th\'eor\`emes \ref{T2} et \ref{t4}, appliqu\'es \`a
$\sA$ en vertu du th\'eor\`eme \ref{tctki}. De m\^eme, b) r\'esulte des
m\^emes r\'ef\'erences, appliqu\'ees \`a $\sA/\sqrt[\otimes]{0}$ munie
de $R'$, en notant que la donn\'ee de $R'$ est \'equi\-va\-len\-te \`a celle de
($R$ et) des $\pi_A^+$.\qed

\newpage

\addtocontents{toc}{{\bf IV. Enveloppes}\hfill\thepage}
\
\bigskip
\begin{center} \large\bf IV. Enveloppes
\end{center}
\bigskip

Dans cette partie, on expose quelques applications ``concr\`etes" des
th\'eo\-r\`e\-mes de scindage de la partie pr\'ec\'edente \`a la
th\'eorie des re\-pr\'e\-sen\-ta\-tions ind\'e\-com\-po\-sa\-bles des 
alg\`ebres
artiniennes et des groupes alg\'ebriques 
li\-n\'e\-ai\-res. Dans ce
dernier cas, le r\'esultat-cl\'e \ref{t1} 
est une vaste
g\'en\'eralisation du th\'eor\`eme classique 
de
Jacobson-Morosov-Kostant.

\section{Le cas non mo\-no\-\"{\i}\-dal: 
enveloppes
pro-semi-simples}\label{compl}

Dans cette section, nous 
tirons les fruits du th\'eor\`eme \ref{T1}; elle
peut \^etre 
consid\'er\'ee comme un \'echauffement pr\'eparatoire aux
sections 
suivantes.

\subsection{Enveloppes pro-semi-simples} 
\label{12.1}
Soient $\sA$ une petite $K$-ca\-t\'e\-gorie et $\pi_\sA: 
\sA \to \bar\sA$
le foncteur de projection sur le quotient de
$\sA$ 
par son radical. On a alors un foncteur canonique 

\[\bar\sA\hbox{--}Modf \to \sA\hbox{--}Modf\] 
entre ca\-t\'e\-go\-ries de 
modules \`a gauche de $K$-dimension finie (nous
dirons bri\`evement: 
$K$-finis), induit par $M\mapsto M\circ \pi_\sA$,
\cf 
d\'efinition
\ref{D2 1/2}. Ce foncteur est pleinement fid\`ele (du 
fait que $\pi_\sA$
est plein et surjectif). Il en est de m\^eme du 
foncteur compos\'e avec
$\pi_{\sA{\hbox{--}}Mod}$ 

\[\bar\sA\hbox{--}Modf \to
\overline{\sA{\hbox{--}}Modf}\]
si 
$\bar\sA$-$Mod$ est ab\'elienne
semi-simple\footnote{Il y a lieu 
\'eventuellement de prendre un second
univers pour ne travailler 
qu'avec de petites ca\-t\'e\-go\-ries.}.

Prenons pour $\sA$ une alg\`ebre 
$A$ (associative unitaire) sur un corps
$K$ {\it parfait}. On a le 
foncteur d'oubli, fid\`ele et exact
\[\omega:A\hbox{--}Modf \to 
Vec_K.\]

La situation ne change pas si l'on
remplace $A$ par sa 
compl\'etion profinie (com\-pl\'e\-tion eu \'egard aux
id\'eaux 
bilat\`eres de $K$-codimension finie). On peut donc la supposer,
et 
on la supposera, profinie, ce qui permet d'\'enoncer:

\begin{lemme} 
Si $A$ est profinie, le quotient $\bar A$ de $A$ par son
radical est 
semi-simple (profinie, non n\'e\-ces\-sai\-re\-ment de
$K$-dimension finie). 
A fortiori, le foncteur
\[\bar A\hbox{--}Modf \to 
\overline{\sA{\hbox{--}}Modf}\] 
est pleinement fid\`ele. 

\end{lemme} 

\prf \'Ecrivons $A$ comme limite projective filtrante 
$\lim A_\alpha$ \`a
fl\`e\-ches de transition surjectives. Comme 
un
\'el\'ement inversible de
$A$ s'identifie \`a une collection 
d'\'el\'ements inversibles
compatibles de $A_\alpha$, on obtient, en 
revenant
\`a la d\'efinition des radicaux, que le radical $R$ de $A$ 
est la
limite $\lim R_\alpha$ des radicaux de
$A_\alpha$. On a alors 
un plongement $A/R \inj  \lim (A_\alpha/R_\alpha)$, le compos\'e avec 
la
projection sur chaque alg\`ebre semi-simple $A_\alpha/R_\alpha$ 
\'etant
surjective. Comme $A/R$ est profinie ($R$
\'etant un id\'eal 
ferm\'e), cette injection est donc un isomor\-phis\-me,
d'o\`u le 
r\'esultat.\qed

On sait que $A\hbox{--}Modf$ est une $K$-ca\-t\'e\-go\-rie 
de Wedderburn
(proposition \ref{Pwed}). On peut donc appliquer le 
th\'eor\`eme
\ref{T1}, et obtenir une section fonctorielle $s: 
\overline{A\hbox{--}Modf}\to  A\hbox{--}Modf$, deux telles 
sections
\'etant conjugu\'ees.    

D'apr\`es 
\cite[II.2.6.3]{saavedra} ou
\cite[2.2]{serre}, il existe une 
$K$-{\it alg\`ebre pro-semi-simple} $A_s$,
non n\'e\-ces\-sai\-re\-ment de 
$K$-dimension finie,  telle que le compos\'e
$\omega\circ s$ induise 
une \'equivalence de ca\-t\'e\-go\-ries  $
\overline{A{\hbox{--}}Modf}\to 
A_s\hbox{--}Modf$, et un
homomor\-phis\-me
$A\to A_s$. 
 
On rappelle 
qu'\'etant donn\'e un
homomor\-phis\-me $f: A'\to A$ de $K$-alg\`ebres 
profinies et le foncteur
associ\'e $f^\ast:A\hbox{--}Modf\to 
A'\hbox{--}Modf$, on a  (\cf
\cite[II.2.6.3]{saavedra}  pour 
l'\'enonc\'e dual):
\begin{itemize}
\item $f$ est surjectif si et 
seulement si $f^\ast$ est pleinement
fid\`ele, et pour tout 
$A$-module $K$-fini $M$, tout sous-objet de
$f^\ast(M)$ provient d'un 
sous-objet de $M$;
\item $f$ est injectif si et seulement si tout 
$A'$-module  $K$-fini
$M'$ est sous-quotient d'un module de la forme 
$f^\ast(M)$. 
\end{itemize}
 
En particulier, le mor\-phis\-me
$A\to A_s$ 
est une injection. Comme deux sections $s$ sont
conjugu\'ees, $A_s$ 
est bien d\'efinie \`a isomor\-phis\-me pr\`es. Il
d\'ecoule par ailleurs 
du lemme pr\'ec\'edent, et de {\it loc. cit.}, que
l'on a une 
surjection $A_s\to A/R$. 

\medskip On peut consid\'erer $A_s$ comme 
une ``enveloppe pro-semi-simple"
de l'al\-g\`e\-bre profinie $A$: 
d'apr\`es le lemme
\ref{l2}, \emph{les classes d'isomor\-phis\-mes de 
$A$-modules
ind\'e\-com\-po\-sa\-bles sont en bijection avec les classes 
d'isomor\-phis\-mes de
$A_s$-modules simples.}  

Cette bijection est 
donn\'ee de la mani\`ere suivante: \`a 
tout
in\-d\'e\-com\-po\-sa\-ble $M$ dans $A$-$Modf$, on peut donner 
une
structure de $A_s$-module \'etendant celle de $A$-module. On le 
note
alors $M_s$; il est simple, et $End_{A_s}M_s$ est le corps 
gauche
$\overline{End_A(M)}=End_A M/\rad(End_A M)$. 

L'alg\`ebre 
pro-semi-simple $A_s$ est semi-simple, \ie de dimension
finie sur 
$K$, si et seulement si $A$ est \emph{de type de
re\-pr\'e\-sen\-ta\-tion 
fini}. Dans ce cas, si
$M_\alpha$ parcourt un syst\`eme de 
repr\'esentants des classes
d'isomor\-phis\-mes de $A$-modules 
ind\'e\-com\-po\-sa\-bles $K$-finis,
notons $n_\alpha$ la dimension de 
$M_\alpha$ sur le corps gauche
$D_\alpha=\overline{End_A(M_\alpha)}$. 
On a $A_s\cong \prod_\alpha
M_{n_\alpha}(D_\alpha^{\rm o})$.

\begin{rems}\label{prof}\
\begin{itemize}
\item[a)] D'apr\`es 
\cite[2.17]{de}, toute $K$-ca\-t\'e\-go\-rie ab\'elienne
dont tout objet 
est de longueur finie, dont toute alg\`ebre
d'endomor\-phis\-mes est de 
$K$-dimension finie, et admettant un
g\'en\'erateur, 
est
\'equi\-va\-len\-te \`a une cat\'e\-go\-rie $A\hbox{--}Modf$ pour 
une
$K$-alg\`ebre $A$ convenable de dimension finie sur
$K$. 

\item[b)] Pour l'alg\`ebre profinie $A=K[[T_1,\dots, T_n]]$, 
la
ca\-t\'e\-go\-rie 
$A\hbox{--}Modf$ est naturellement \'equi\-va\-len\-te \`a 
la ca\-t\'e\-go\-rie
$Rep_K(\bG_a^n)$ des re\-pr\'e\-sen\-ta\-tions de dimension 
finie du groupe
vectoriel $\bG_a^n$. En effet, la cog\`ebre duale de 
$A$ s'identifie
naturellement
\`a la cog\`ebre sous-jacente \`a la 
big\`ebre affine $\sO(\bG_a^n)$, et
on conclut par \cite[3]{serre}. 

\item[c)] Soit $\sJ$ l'id\'eal de $\sA=A$-$Modf$ form\'e des 
mor\-phis\-mes
qui se factorisent \`a travers un objet projectif de 
$A$-$Modf$. Le
quotient $\sA/\sJ$ s'appelle la ca\-t\'e\-go\-rie {\it 
stable} des 
$A$-modules finis, \cf
\cite[p.23]{be}, et se note 
$A\hbox{--}\underline{Modf}$. On a un
diagramme commutatif 

\[\begin{CD} A\hbox{--}Modf@>pr>>A\hbox{--}\underline{Modf} 
\\
@VVV@VVV\\
\overline{A\hbox{--}Modf} 
@>\overline{pr}>>
\overline{A\hbox{--}\underline{Modf}}. 
\end{CD}\]

Le foncteur du bas $\overline{pr}$ est un foncteur plein 
surjectif entre
ca\-t\'e\-go\-ries semi-simples, donc admet une section. En 
fait, on peut
partitionner les classes d'isomor\-phis\-mes d'objets 
simples de
$\overline{A\hbox{--}Modf} $ en deux ensembles: les 
projectifs ($P$) et
les autres ($NP$); $\overline{pr}$ envoie les 
objets dans $P$ sur $0$, et
la restriction de $\overline{pr}$
\`a la 
sous-ca\-t\'e\-go\-rie
$K$-lin\'eaire pleine de $\overline{A\hbox{--}Modf} 
$ engendr\'ee par
les objets dans $NP$ est une
\'e\-qui\-va\-len\-ce 
de ca\-t\'e\-go\-ries. On en d\'eduit 
que
$\overline{A\hbox{--}\underline{Modf}}\cong 
\underline{A_s}\hbox{--}Modf $ pour un quotient $\underline{A_s}$ de 
$A_s$ ``correspondant" \`a $NP$. On peut
alors former le ``push-out" 
$A'=\underline{A_s}\times^{A_s}\,A/R$.

Auslander a conjectur\'e que 
si $A$ et $B$ sont deux $K$-alg\`ebres de
dimension finie telles que 
$A\hbox{--}\underline{Modf} $ et
$B\hbox{--}\underline{Modf} $ sont 
\'equi\-valentes (\'equi\-va\-lence de
Morita stable), alors il y a 
le m\^eme nombre de modules simples
non-projectifs sur $A$ et sur 
$B$, \cf
\cite[p.223]{be}. Cela \'equivaut donc \`a dire que les 
alg\`ebres
push-out $A'$ et $B'$ sont Morita-\'equi\-va\-len\-tes.

\begin{sloppypar}
Rappelons d'autre part que si $A$ est de 
dimension globale finie (au
sens de \cite{ce}), on peut lui associer 
une $K$-alg\`ebre $\hat A$ dite
r\'ep\'etitive (de dimension infinie 
sur $K$), et une \'equivalence de
ca\-t\'e\-go\-ries
\[ 
D^b(A\hbox{--}Modf)\to {\hat A\hbox{--}\underline{Modf}}\] 
\cf 
\cite[5]{ringel}, \cite{happel}; on peut d'ailleurs replacer $\hat
A$ 
par sa compl\'etion pro\-fi\-nie. Ce qui pr\'ec\`ede permet alors 
de
``d\'ecrire" le quotient de $D^b(A\hbox{--}Modf)$ par son 
radical.
\end{sloppypar}
\end{itemize}
\end{rems}

\subsection{Lien 
avec l'alg\`ebre d'Auslander.}\label{12.2} L'alg\`ebre
d'Auslander 
d'une $K$-alg\`ebre $A$ est 
\[{\bf Aus}(A):= End_A(\bigoplus_\alpha 
M_\alpha ),\]
o\`u $M_\alpha$ parcourt un syst\`eme (en g\'en\'eral 
infini) de
repr\'esentants des classes d'isomor\-phis\-mes de 
$A$-modules
ind\'e\-com\-po\-sa\-bles $K$-finis. Lorsque $A$ est de 
dimension
finie sur $K$, elle partage avec $A_s$ la propri\'et\'e que 
\emph{les
classes d'isomor\-phis\-mes de $A$-modules ind\'e\-com\-po\-sa\-bles 
sont en bijection
avec les classes d'isomor\-phis\-mes de
${\bf 
Aus}(A)$-modules simples}, \cf \cite[4.9.5]{be}.

\begin{sloppypar}
Cette bijection est la suivante: \`a tout 
ind\'e\-com\-po\-sa\-ble $M$ dans
$A\hbox{--}Modf$, on associe le ${\bf 
Aus}(A)$-module simple
$S_M:=\overline{End_A(M)}$.
 En tant que 
$K$-alg\`ebre, $\overline{End_A(M)}$ s'identifie \`a
$End_{{\bf 
Aus}(A)}S_M$ ({\it loc. cit.}).   
\end{sloppypar}

L'alg\`ebre ${\bf 
Aus}(A)$ est de dimension finie sur $K$ si et
seulement si $A$ est de 
type de re\-pr\'e\-sen\-ta\-tion fini. Dans ce cas,
\[\overline{{\bf 
Aus}(A)}:= {\bf Aus}(A)/\rad({\bf Aus}(A))=
\prod_\alpha 
D_\alpha.\]

Ainsi, $\overline{{\bf Aus}(A)}$ est 
\emph{Morita-\'equi\-va\-len\-te} \`a
$A_s$\footnote{$\overline{{\bf 
Aus}(A)}$ s'identifie d'ailleurs au
quotient par son radical de 
l'alg\`ebre des chemins du carquois
d'Auslander-Reiten de $A$ dont il 
a \'et\'e question au \S \ref{radinf},
\cf 
\eg
\cite[4.1.11]{be}.}.

\subsection{Cas des alg\`ebres 
h\'er\'editaires}\label{12.3} Rappelons
qu'une
$K$-alg\`ebre est dite 
{\it h\'er\'edi\-tai\-re} (\`a gauche) si les
conditions \'equi\-va\-len\-tes 
suivantes sont sa\-tis\-fai\-tes 
\cf
\cite[I.5]{ce}:
\begin{itemize}
\item tout id\'eal \`a gauche est 
projectif,
\item tout sous-module d'un module projectif (\`a gauche) 
est projectif.
\end{itemize}

C'est une condition Morita-invariante.

\medskip Supposons pour simplifier $K$ al\-g\'e\-bri\-que\-ment clos 
et
$A$ de dimension finie sur $K$, minimale dans sa classe 
d'\'equivalence
de Morita (ce qui revient \`a dire que $\bar A = A/R$ 
est commutative).
Alors $A$ est h\'er\'editaire si et seulement si 
$A=K\vec{\Delta}$ est
l'alg\`ebre des chemins d'un carquois fini 
$\vec{\Delta}$ (et sans boucle
orient\'ee), \cf
\cite[4.2]{be}, 
\cite[6]{cq}. 
Un $A$-module n'est donc rien d'autre qu'une 
repr\'esen\-ta\-tion
$K$-lin\'eaire de $\vec{\Delta}$, et $A$ 
s'identifie en fait \`a
l'alg\`ebre tensorielle sur un $\bar 
A$-bimodule fini, \cf \cite[6]{cq}.  

\medskip Supposons en outre 
$A$ de type de re\-pr\'e\-sen\-ta\-tion fini et
connexe. Alors par le 
th\'eor\`eme de Gabriel (\cf \cite[4.7.6]{be},
\cite[2.9]{rowen}), le 
graphe non-orient\'e ${\Delta}$ sous-jacent \`a
$\vec{\Delta}$ est un 
diagramme de Dynkin des s\'eries ADE. En outre, 
les
ind\'e\-com\-po\-sa\-bles
$M_\alpha$ sont en bijection avec les racines 
positives, et les
$n_\alpha$ ne sont autres que leurs ``longueurs" 
respectives (somme des
coefficients dans la base de racines 
standard). Cela d\'etermine
$A_s$ dans ce cas.        
 
Voici deux 
exemples. Soit d'abord $A$ la $K$-alg\`ebre des 
matrices
triangulaires sup\'erieures de taille $n$. C'est 
l'alg\`ebre
$K\vec{\Delta}$ pour le diagramme de Dynkin 
$\vec{\Delta}= A_n$
(orient\'e de gauche \`a droite). Dans la base 
standard de racines
$\alpha_1,\dots, \alpha_n$, les racines positives 
sont les
$\alpha_{ij}=\displaystyle\sum_{ i\leq k<j}\,\alpha_k, 
0<i<j\leq n+1$. On
en d\'eduit que
\[A_s \cong \displaystyle{\prod_{1\leq k\leq n}\; (M_k(K))^{n+1-k}} .\]

Prenons maintenant $A=KE_6$, l'alg\`ebre des chemins du carquois de
Dynkin (orient\'e) \`a six sommets $\vec{\Delta}=E_6$. \`A l'aide des
planches de racines \cite{tables}, on trouve
\begin{multline*} A_s \cong K^6\times M_2(K)^3 \times M_5(K) \times
M_6(K)^2\times M_7(K)^3\\
\times M_8(K)^2\times M_9(K)\times
M_{10}(K)\times M_{11}(K) .\end{multline*}

\begin{rem} Dans le cas des alg\`ebres de chemins $A=K\vec{\Delta}$, un
avatar de $A_s$ a
\'et\'e construite par voie g\'eom\'etrique dans \cite{rump}.
\end{rem}

\section{Sections mo\-no\-\"{\i}\-dales et foncteurs fibres}\label{radrep}

Dans ce paragraphe, nous appliquons le th\'eor\`eme de scindage
mo\-no\-\"{\i}\-dal sy\-m\'e\-tri\-que au cas des ca\-t\'e\-go\-ries tannakiennes.

\subsection{Un contexte g\'en\'eral}

\begin{thm}\label{c4} Soient $K$ un corps de caract\'eristique nulle
et $L$ une extension de $K$. Soit $\sA$
une (petite) ca\-t\'e\-go\-rie $K$-lin\'eaire, pseudo-ab\'elienne,
mono\-\"{\i}dale sy\-m\'e\-tri\-que rigide avec $K=End(\un )$. On suppose
qu'il existe un $K$-foncteur {\rm fid\`ele} mo\-no\-\"{\i}\-dal
sy\-m\'e\-tri\-que $\omega: \sA \to  Vec_L$.
\\ Alors $\sR=\rad(\sA)$ est mo\-no\-\"{\i}\-dal, et $\bar\sA = \sA/\sR$ est
tannakienne semi-simple sur $K$.
\\ En particulier, si $\sA$ est une ca\-t\'e\-go\-rie tannakienne sur $K$, 
de radical $\sR$, alors
$\sA/\sR$ est encore tanna\-kienne sur $K$; elle est neutre si $\sA$ est
neutre.
\end{thm}

\prf On sait par le th\'eor\`eme \ref{wedtr} que $\sA$ est de 
Wedderburn. Comme $\bar\sA$ est
$K$-lin\'eaire pseudo-ab\'elienne, elle est m\^eme
ab\'elienne (semi-simple); elle est par ailleurs rigide (\cf 
\ref{sorrig}). D'autre part, la tensorialit\'e du
radical a \'et\'e d\'emontr\'ee dans le th\'eor\`eme \ref{absJannsen}.

Les hypoth\`eses des th\'eor\`emes
\ref{T2} et \ref{t4} sont satisfaites. On peut donc trouver une section
$K$-lin\'eaire mo\-no\-\"{\i}\-dale $s$ de
$\pi_\sA$, qui est automatiquement sy\-m\'e\-tri\-que. Une telle section est
clairement fid\`ele, donc le  $K$-foncteur mo\-no\-\"{\i}\-dal
sy\-m\'e\-tri\-que $\omega_s=\omega\circ s: \bar\sA\to Vec_L$ l'est aussi.
Comme c'est un foncteur mo\-no\-\"{\i}\-dal
fid\`ele entre ca\-t\'e\-go\-ries mo\-no\-\"{\i}\-dales a\-b\'e\-lien\-nes
semi-simples, il est exact (les suites exactes se
scindent); c'est donc un foncteur fibre.

La seconde assertion d\'ecoule de l\`a.\qed

\begin{rems}\
\begin{itemize}
\item[$a)$] On pourrait aussi prouver que $\bar\sA$ est
tannakienne  si $\sA$ l'est en utilisant la
caract\'erisation interne de Deligne \cite{de}. En effet,
comme $\sR= \sN$, la dimension d'un objet de
$\sA$ ne change pas si on la calcule dans $\bar\sA$; c'est donc un
entier naturel. Mais cet argument ne montre pas que $\bar\sA$ est neutre
quand $\sA$ l'est.

\item[$b)$] L'exemple de ca\-t\'e\-go\-rie mo\-no\-\"{\i}\-dale $K$-lin\'eaire non
strictement semi-primaire vu plus haut (contre-exemple \ref{Cssp})
fournit, en caract\'e\-ris\-ti\-que nulle, un exemple d'application du
th\'eor\`eme ci-dessus (le foncteur $\omega$ \'etant induit par $x\mapsto
K$ pour tout $x\in E$).

\noindent Dans cet exem\-ple, la ca\-t\'e\-go\-rie tannakienne $\bar\sA$ 
n'est autre que
la cat\'e\-go\-rie des repr\'e\-sen\-ta\-tions du groupe diagonalisable de
groupe de carac\-t\`e\-res $E$.
\item[$c)$] L'exemple \ref{Car3} montre que le th\'eor\`eme de 
scindage mo\-no\-\"{\i}\-dal sym\'e\-tri\-que ne
s'\'etend pas, en rempla\c cant $\sR$ par $\sN$, au cas o\`u le 
radical $\sR$ n'est pas mo\-no\-\"{\i}\-dal. Dans cet
exemple,
$\sA/\sN$ est \'equi\-va\-len\-te \`a la ca\-t\'e\-go\-rie des super-espaces 
vectoriels de dimension finie, qui n'est pas
tannakienne (on a un objet de carr\'e tensoriel isomorphe \`a $\un$, 
mais de dimension $\neq 1$).
\item[$d)$] Dans \cite{ak(note)} nous d\'eduisons de ce
th\'eor\`eme, in\-d\'e\-pen\-dam\-ment de toute
conjecture, un \emph{groupe de Galois motivique pro-r\'e\-duc\-tif}
attach\'e
\`a toute cohomologie de Weil ``classique" d\'efinie sur les vari\'et\'es
projectives lisses sur un corps quelconque.
\end{itemize}
\end{rems}

\subsection{Le cas tannakien neutralis\'e}\label{neu} Ce sous-paragraphe
est pr\'eparatoire au paragraphe \ref{jm}.

On suppose encore $K$ de ca\-rac\-t\'e\-ris\-ti\-que $0$. Soient $G$ un
$K$-sch\'ema en grou\-pes affine (en bref: un $K$-groupe affine) et
$\sA=Rep_K(G)$ la ca\-t\'e\-go\-rie des
$K$-re\-pr\'e\-sen\-ta\-tions de  dimension finie de $G$. On note $\pi_G$ au lieu
de $\pi_\sA$, $\sR_G$ le radical de $\sA$, et $\omega_G$ le foncteur
fibre canonique de $\sA$ (\`a valeurs dans $Vec_K$).

Le lemme suivant est fort utile:

\begin{lemme}\label{l13.1} Soient $\sB$ une (petite) ca\-t\'e\-go\-rie 
$K$-lin\'eaire
mo\-no\-\"{\i}\-dale sym\'e\-tri\-que et $f:\sB\to Rep_K(G)$ un foncteur
mo\-no\-\"{\i}\-dal. Pour toute $K$-alg\`ebre commutative $R$, notons $f^R$ 
(\resp $\omega^R$) le foncteur
mo\-no\-\"{\i}\-dal $f$  (\resp $\omega $)
  compos\'e avec \[Vec_K\to R-Modf  \;\;\;\;(\hbox{\resp} \;\;Rep_K(G) 
\to Rep_R(G)).\] Notons
$\underline{Aut}^\otimes f$ (\resp $
\underline{Aut}^\otimes(\omega\circ f)$) le faisceau fpqc associ\'e 
au foncteur  $R\mapsto Aut^\otimes f^R
$ (\resp
$R \mapsto Aut^\otimes(\omega^R\circ f)$).
  Alors
\\ 1) $\underline{Aut}^\otimes(\omega\circ f)$ est re\-pr\'e\-sen\-ta\-ble 
par un $K$-sch\'ema en groupes affine,
  \\ 2) $\underline{Aut}^\otimes(f)$ est re\-pr\'e\-sen\-ta\-ble par le 
centralisateur de l'image de $G \allowbreak\iso
\underline{Aut}^\otimes(\omega)$ dans $\underline{Aut}^\otimes(\omega\circ 
f)$,
  \\ 3) l'homomor\-phis\-me
\begin{align*}
\underline{Aut}^\otimes(f)&\to
\underline{Aut}^\otimes(\omega\circ f)\\
(u^R)&\mapsto (\omega^R(u))
\end{align*} (vu comme homomor\-phis\-me de groupes affines) est un 
monomor\-phis\-me, \ie une immersion ferm\'ee.
\end{lemme}

\prf 1) $\underline{Aut}^\otimes(\omega\circ f)$ est re\-pr\'e\-sen\-ta\-ble 
par le sous-sch\'ema en groupes
ferm\'e de ${\bf G }:= \prod_{B\in \sB}\,GL(\omega\circ f(B))$ qui 
fixe ceux des mor\-phis\-mes de $Rep_K \bf G$ qui sont de
la forme $\omega\circ f(u)$, o\`u $u$ est un mor\-phis\-me de $\sB$.

\noindent 2) Il est clair que l'image de $\underline{Aut}^\otimes(f)$ 
dans $\underline{Aut}^\otimes(\omega\circ
f)$ centralise $\underline{Aut}^\otimes(\omega )$. Inversement, soit 
$v\in \underline{Aut}^\otimes(\omega\circ f)(R)$
centra\-lisant l'image de $\underline{Aut}^\otimes(\omega)(R')$ pour 
toute extension $R'/R$. Pour tout $B\in \sB$ et
tout $g\in G(R')$, l'\'el\'ement $v_B\in End_{R'}(\omega^{R'}(f(B)))$ 
commute \`a $g_B$. Il lui
correspond donc un unique \'el\'ement $u_B\in \sA(f^R(B),f^R(B))$ tel que
$\omega^R(u_B)=v_B$. La fid\'elit\'e de $\omega^R$ implique de plus que
$u_B$ est inversible et mo\-no\-\"{\i}\-dal. Cela montre que 
$\underline{Aut}^\otimes(\omega\circ f)$ est le centralisateur de
l'image de $
\underline{Aut}^\otimes(\omega)$ dans 
$\underline{Aut}^\otimes(\omega\circ f)$, donc re\-pr\'e\-sen\-ta\-ble.

\noindent 3) d\'ecoule de ce que les $\omega^R$ sont fid\`eles. \qed

Choisissons une section mo\-no\-\"{\i}\-dale $s$ de
$\pi_\sA$, et posons
comme ci-dessus $\omega_s=\omega_G\circ s$. Soit
$G_s=\underline{Aut}^{\otimes}(\omega_s)$: c'est un $K$-groupe 
pro-r\'eductif.
On a un homomor\-phis\-me \'evident:
\[\begin{CD}
s^\sharp:G\iso \underline{Aut}^\otimes(\omega_G) \to 
\underline{Aut}^\otimes(\omega_s)=G_s
\end{CD}\]
d'o\`u un foncteur mo\-no\-\"{\i}\-dal
\[ (s^\sharp)^*: Rep_K(G_s)\to Rep_K(G).\]

\begin{lemme} L'homomor\-phis\-me $s^\sharp$ est un monomor\-phis\-me (\ie une
immersion ferm\'ee).
\end{lemme}

\prf Soit $N=\Ker s^\sharp$. On a un diagramme commutatif de ca\-t\'e\-go\-ries
et foncteurs:
\[\begin{CD}
Rep_K(G/N)@>g^\ast>> Rep_K(G)\\
\nwarrow&&\scriptstyle (s^\sharp)^*\displaystyle\nearrow\\
&Rep_K(G_s)
\end{CD}\]

Le foncteur $g^\ast$ est pleinement fid\`ele. Comme $(s^\sharp)^*$ est
\'evidemment essentiellement surjectif, il en est de m\^eme de $g^\ast$,
qui est donc une \'equivalence de ca\-t\'e\-go\-ries, d'o\`u $N=1$.\qed

\begin{rem}\label{r2} Le foncteur
\begin{align*}
\theta_s:\bar\sA&\to Rep_K(G_s)\\
X&\mapsto (G_s\to GL(\omega_s(X)))
\end{align*}
est une \'equivalence de ca\-t\'e\-go\-ries mo\-no\-\"{\i}\-dales, mais pas en
g\'en\'eral un isomor\-phis\-me de ca\-t\'e\-go\-ries. Il poss\`ede n\'eanmoins un
quasi-inverse ca\-no\-ni\-que. On a un diagramme strictement commutatif de
ca\-t\'e\-go\-ries et foncteurs:
\[\begin{CD}
Rep_K(G)@>\omega_G>> Vec_K\\
\scriptstyle s\Big\uparrow &\omega_s\nearrow& @A\omega_{G_s}AA\\
\bar \sA@>\theta_s>> Rep_K(G_s)
\end{CD}\]

On a \'evidemment $(s^\sharp)^*\circ\theta_s=s$. En
particulier,
$\sigma_s:=\pi_G\circ (s^\sharp)^*$ est un inverse \`a
gauche de
$\theta_s$; c'en est donc aussi un quasi-inverse \`a
droite, par un
isomor\-phis\-me naturel mo\-no\-\"{\i}\-dal.
\end{rem}

Si $t$ est une autre section mo\-no\-\"{\i}\-dale de $\pi_G$,
elle d\'efinit
un autre foncteur fibre, donc un autre groupe $G_t$. D'apr\`es la
th\'eorie g\'en\'erale des ca\-t\'e\-go\-ries tannakiennes, $G_s$ et $G_t$
sont des formes int\'erieures l'un de l'autre. On a mieux, puisque $t$ est
mo\-no\-\"{\i}\-dalement conjugu\'ee \`a $s$. En effet, choisissons une telle
conjugaison $u$, de sorte que $t=usu^{-1}$. On a alors un diagramme
commutatif
\begin{equation}\label{eq13.1}
\begin{CD}
&&G_s\\
&\s s^\sharp\displaystyle \nearrow\\
G&&@V{\hat{u}}V{\wr}V\\
&\s t^\sharp\displaystyle \searrow\\
&&G_t
\end{CD}
\end{equation}
o\`u $\hat{u}:\underline{Aut}^\otimes(\omega_s)(R)\to 
\underline{Aut}^\otimes(\omega_t)(R)$ est
donn\'e par $g\mapsto ugu^{-1}$ pour toute $K$-alg\`ebre $R$. Cette
construction donne:

\begin{prop}\label{p4}
On a un foncteur canonique $T_G$ du groupo\"{\i}de connexe
$\sG_u^{\otimes}(G)$ des sections mo\-no\-\"{\i}\-dales de $\pi_G$ vers la
ca\-t\'e\-go\-rie des monomor\-phis\-mes de $G$ vers un groupe pro-r\'eductif, qui
associe \`a la section $s$ le monomor\-phis\-me $G\inj G_s$.  En particulier,
pour deux sections mo\-no\-\"{\i}\-dales $s,t$, les groupes $G_s$ et $G_t$ sont
$K$-isomorphes.\qed
  \end{prop}

\begin{para}\label{ecl} \'Eclairons et compl\'etons cette construction \`a
l'aide des r\'esultats du paragraphe \ref{repr}. Notons
$\Gamma_u(G)=(E_u,S_u)$ le groupo\"{\i}de
affine scind\'e unipotent transitif sur $S^\otimes$ associ\'e \`a
$\sA=Rep_K(G)$ par le th\'eor\`eme \ref{gerbe'}: on a
$\Gamma_u(G)(K)=\sG^\otimes(G)$. Par ailleurs, la
ca\-t\'e\-go\-rie tannakienne $\bar\sA$ \'etant neutre, sa gerbe est
repr\'esent\'ee par un $K$-groupo\"{\i}de affine que nous noterons
$\Gamma_\red(G)$: la ca\-t\'e\-go\-rie $\Gamma_\red(G)(K)$ est le
groupo\"{\i}de des $K$-foncteurs fibres sur $\bar\sA$. Consid\'erons $G$
comme groupo\"{\i}de \`a un objet. On a alors un (bi)mor\-phis\-me de
$K$-groupo\"{\i}des affines
\[\Gamma_u(G)\times_K G\to \Gamma_\red(G)\]
qui au niveau des $R$-points est d\'ecrit de la mani\`ere suivante:
\begin{align*}
\underline{Hom}^\otimes(s,t)(R)\times \underline{Aut}^\otimes(\omega)(R)&\to
\underline{Hom}^\otimes(\omega_s,\omega_t)(R)\\
(s,*)&\mapsto \omega_s=\omega\circ s.
\end{align*}

Si l'on fixe une
section ($K$-rationnelle) $s$, on obtient en particulier un
monomor\-phis\-me de
$K$-sch\'emas en groupes affines
\[U_s= U_s(G):=E_u(s,s)\inj G_s\]
avec $U_s$ pro-unipotent, o\`u $U_s$ \emph{centralise} $G$. Le diagramme
\eqref{eq13.1} se compl\`ete en un diagramme commutatif:
\[\begin{CD}
U_s\times_K G@>>> G_s\\
@V{\hat{u}\times 1}V{\wr}V @V{\hat{u}}V{\wr}V\\
U_t\times_K G@>>> G_t.
\end{CD}\]
\end{para}

\begin{prop}\label{centr}  Pour tout \'epimor\-phis\-me $\phi:G_s\surj H$, le
centralisateur de
$\phi\circ s^\sharp(G)$ dans $H$ est \'egal au produit
$\phi(U_s)\times Z(H)$.
\end{prop}

[Noter que ce produit est direct, puisque l'intersection de $\phi(U_s)$ et
de
$Z(H)$ est un groupe unipotent de type multiplicatif.]

\medskip\prf Par souci de clart\'e, notons plus pr\'e\-ci\-s\'e\-ment 
$\omega=\omega_G$.
Consid\'erons le foncteur
$f=(s^\sharp)^*\circ
\phi^*:Rep_K(H)\to Rep_K(G)$: on remarque que $\omega_G\circ f=\omega_H$.
Appliquons le lemme
\ref{l13.1}: on obtient que l'homomor\-phis\-me $\underline{Aut}^\otimes(f)\to
\underline{Aut}^\otimes(\omega_H)= H $ est injectif, d'image le 
centralisateur
de l'image de $\underline{Aut}^\otimes(\omega_G)=G $ dans $H $. Il reste \`a
calculer $\underline{Aut}^\otimes(f)$, que nous identifions ci-dessous \`a 
un
sous-groupe de $H $.

Soit $h\in \underline{Aut}^\otimes(f)(R)$, o\`u $R$ est une 
$K$-alg\`ebre commutative (on peut d'ailleurs se limiter qu
cas d'un corps, et m\^eme d'une extension alg\'ebrique de $K$, 
puisque $K$ est de caract\'eristique nulle). Alors
$\pi_G(h)\in
\underline{Aut}^\otimes(\pi_G\circ (s^\sharp)^*\circ \phi^*)(R)$. 
Comme les foncteurs $\pi_G\circ
(s^\sharp)^*$ et $\phi^*$ sont pleinement fid\`eles, leur compos\'e l'est
aussi et $\underline{Aut}^\otimes(\pi_G\circ(s^\sharp)^*\circ 
\phi^*)(R)$ s'identifie
canoniquement \`a
$\underline{Aut}^\otimes(Id_{Rep_K(H)})(R)=Z(H)(R))$. Autrement dit, 
il existe un unique
$z\in Z(H)(R)$ tel que $\pi_G(h)\allowbreak=\pi_G(z)$.
Alors $u=h\cdot z^{-1}$ v\'erifie $\pi_G(u)=1$; autrement dit, $u\in
\phi(U_s(R))$.
  D'o\`u l'assertion. \qed

\begin{prop} \label{p3} a) Pour toute section $s$ comme ci-dessus, tout
homomor\-phis\-me $\phi$ de $G$ vers un groupe pro-r\'e\-duc\-tif $H$ se
prolonge en un homomor\-phis\-me $\psi:G_s\to H$.
\end{prop}

\prf  On a un
diagramme commutatif
\[\begin{CD}
&&Rep_K(G)@>\omega_G>> Vec_K\\
&\s \phi^*\displaystyle\nearrow &@V{\pi_G}VV\\
Rep_K(H)@>\bar\phi^*>> Rep_K(G)/\sR_G
\end{CD}\]
o\`u par souci de clart\'e on note $\omega_G$ le foncteur fibre canonique
de $Rep_K(G)$. On a $\omega_H=\omega_G\circ \phi^*$.

Remarquons que $\bar\phi^*$ est exact, puisque la ca\-t\'e\-go\-rie $Rep_K(H)$
est semi-simple. Il est aussi fid\`ele: si $A\in Rep_K(H)$, le noyau de
l'homomor\-phis\-me induit par $\bar\phi^*$ de l'anneau $R$ des endomor\-phis\-mes
de $A$ vers l'anneau des endomor\-phis\-mes de $\bar \phi^*(A)$ est l'image
r\'eciproque dans $R$ du radical de l'anneau des endomor\-phis\-mes de
$\phi^*(A)$; cette image r\'eciproque est nulle, comme id\'eal nilpotent
d'une alg\`ebre semi-simple. Donc $\omega_G\circ s\circ\bar\phi^*$ est un
foncteur fibre sur $Rep_K(H)$.

Le foncteur
$\bar\phi^*$ induit un homomor\-phis\-me
\[(\bar\phi^*)^\sharp:G_s\to H'\]
avec $H':=\underline{Aut}^\otimes(\omega_G\circ s\circ\bar\phi^*)$.

Appliquons
la proposition \ref{P3'}, en remarquant que $Rep_K(H)$ est 
semi-simple (et m\^eme s\'eparable): il
existe un isomor\-phis\-me naturel mo\-no\-\"{\i}\-dal $u:\phi^*\Rightarrow s\circ
\bar\phi^*$, d'o\`u un isomor\-phis\-me
\[\begin{CD}
H\iso \underline{Aut}^\otimes(\omega_G\circ 
\phi^*)@>\widehat{\omega_G(u)}>>H'
\end{CD}\]
faisant commuter le diagramme
\[\begin{CD}
G@>\phi>> H\\
@V{s^\sharp}VV @V{\widehat{\omega_G(u)}}V{\wr}V\\
G_s@>(\bar\phi^*)^\sharp>> H'.
\end{CD}\]
D'o\`u la proposition.  \qed

\begin{lemme}\label{l3} Soient $f,g\in Hom(G,H)$, o\`u $G$ et $H$ sont
deux $K$-grou\-pes affines. Alors $f$ et $g$ sont conjugu\'es par un
\'el\'ement de $H(K)$ si et seulement s'il existe un isomor\-phis\-me naturel
mo\-no\-\"{\i}\-dal
$f^*\simtimes g^*$ au niveau des ca\-t\'e\-go\-ries de
re\-pr\'e\-sen\-ta\-tions.
\end{lemme}

\prf Soit $\theta$ un tel isomor\-phis\-me. Alors $\theta$ induit un
isomor\-phis\-me $\omega_H=\omega_G\circ f^*\simtimes\omega_G\circ
g^*=\omega_H$. Un tel isomor\-phis\-me correspond \`a un \'el\'ement $h\in
H(K)$. Cet \'el\'ement conjugue $f$ et $g$.\qed

\begin{prop}\label{p2} Soient $G$ un $K$-groupe affine $s,t$ deux
sections mo\-no\-\"{\i}\-da\-les de $\pi_G$, $G_s$ et $G_t$ les
$K$-groupes r\'eductifs attach\'es $s$ et $t$, $u$ un isomor\-phis\-me
mo\-no\-\"{\i}\-dal de $s$ sur $t$ et $\hat{u}:G_s\iso G_t$ l'isomor\-phis\-me
correspondant (\cf \eqref{eq13.1}). Soient enfin
$H$ un autre
$K$-groupe affine $H$ et $f:G_s\to H$, $g:G_t\to H$ deux homomor\-phis\-mes.\\
Supposons que $f\circ s^\sharp$ et $g\circ t^\sharp$ soient conjugu\'es
par un \'el\'ement de $H(K)$. Alors il en est de m\^eme de
$f$ et $g\circ \hat{u}$.
\end{prop}

\prf On se ram\`ene imm\'ediatement au cas o\`u $s=t$, puis au cas o\`u
$f\circ s^\sharp=g\circ s^\sharp$. On a, au niveau des ca\-t\'e\-go\-ries de
re\-pr\'e\-sen\-ta\-tions, une
\'egalit\'e de foncteurs:
\[(s^\sharp)^*\circ f^*= (s^\sharp)^*\circ g^*\]
d'o\`u, en composant \`a gauche avec $\pi_G$:
\[\sigma_s\circ f^*=\sigma_s\circ g^*\]
(\cf remarque \ref{r2}).

Comme $\sigma_s$ est un quasi-inverse mo\-no\-\"{\i}\-dal de $\theta_s$
(ibid.), on en d\'eduit l'exis\-tence d'un isomor\-phis\-me naturel
$f^*\simtimes g^*$. On conclut par le lemme \ref{l3}.\qed

\begin{sco} Avec les notations de \ref{h0}, on a \[HH_0(Rep_K(G))\cong
R_K(G_s)\otimes_\Z K.\] 
\end{sco}

\prf On a vu dans \ref{h0} que $HH_0(\sA)=
HH_0(\bar\sA)$, d'o\`u $HH_0(\sA)\allowbreak\cong HH_0(Rep_K(G_s)) $
d'apr\`es ce qui pr\'ec\`ede, et d'autre part et dans l'exemple \ref{hh0}
que
$HH_0(Rep_K(G_s)) \cong R_K(G_s)\otimes_\Z K$.
\qed

\section{Au-del\`a de Jacobson-Morozov:
enveloppes pro-r\'eductives}\label{jm}

Dans toute ce paragraphe,  \emph{$K$ est un corps de caract\'eristique
$0$} sauf mention expresse du contraire.

Soit $\frg$ une
$K$-alg\`ebre de Lie semi-simple. Le th\'eor\`eme de
Jacobson-Morozov-Kostant (\cite{jac},
\cite{mor}, \cite{kos}, \cite[p. 162, prop. 2 et 4]{lie}) \'enonce que
tout
\'el\'ement nilpotent non nul de $\frg$ est contenu dans un
$\fsl_2$-triplet de $\frg$ et que les orbites de ces deux types d'objets
sous l'action adjointe du groupe adjoint $G$ de $\frg$ sont en bijection.
En termes de groupes alg\'ebriques, cet
\'enonc\'e se traduit ainsi:
\'etant donn\'e un $K$-groupe alg\'ebrique semi-simple $G$, tout
homomor\-phis\-me non trivial $\bG_a\to G$ se prolonge en un 
homomor\-phis\-me
$SL_2\to G$ couvrant l'injection 
canonique
\begin{align}\label{eq2}\phi:\bG_a&\to 
SL_2\\
a&\mapsto\begin{pmatrix}1&a\\0&1\end{pmatrix};\notag
\end{align}
De plus, deux tels prolongements sont conjugu\'es sous l'action 
d'un
\'el\'ement de $G(K)$.

Nous allons retrouver ce r\'esultat en 
appliquant le th\'e\-o\-r\`e\-me
\ref{c4} \`a la ca\-t\'e\-go\-rie 
tannakienne neutralis\'ee des
$K$-re\-pr\'e\-sen\-ta\-tions de $\bG_a$, et 
explorer ce qui joue le r\^ole 
de $SL_2$ lorsque
$\bG_a$ est 
remplac\'e par un
$K$-groupe affine quelconque, dans le fil des 
r\'esultats du 
paragraphe pr\'ec\'edent.

\subsection{Groupes \`a 
conjugaison pr\`es}

\begin{sloppypar}
\label{i3} Soit $\Gaff_K$ la 
ca\-t\'e\-go\-rie
dont les objets sont les $K$-groupes affines et les 
mor\-phis\-mes sont les
homomor\-phis\-mes de $K$-groupes affines. Notons 
$\Gred_K$ la
sous-ca\-t\'e\-go\-rie pleine de $\Gaff_K$ form\'ee des 
groupes
pro-r\'eductifs. Ces ca\-t\'e\-go\-ries ne sont pas bien adapt\'ees 
\`a
l'interpr\'etation des r\'esultats du paragraphe pr\'ec\'edent: 
nous
devons les remplacer par des ca\-t\'e\-go\-ries plus 
grossi\`eres.
\end{sloppypar}

\begin{sloppypar}
\begin{defn}\label{d2 
} 
a) Soient $G,H$ deux $K$-groupes affines. On note
$H^1_K(G,H)$ 
l'ensemble quotient de $Hom_K(G,H)$ par la relation
d'\'equivalence 
$\sim$ telle que $f\sim g$ s'il existe $h\in H(K)$ tel
que 
$g=hfh^{-1}$ (on dit aussi que $f$ et $g$ 
sont
\emph{$H(K)$-conjugu\'es}).\\
b) On note $\bGaff_K$ la 
ca\-t\'e\-go\-rie dont les objets sont les $K$-groupes
affines et dont les 
mor\-phis\-mes sont donn\'es par les ensembles
$H^1_K(G,H)$ (ces 
mor\-phis\-mes se composant de mani\`ere \'evidente): c'est
la 
ca\-t\'e\-go\-rie des \emph{$K$-groupes affines \`a conjugaison pr\`es}. 
Sa
sous-ca\-t\'e\-go\-rie pleine form\'ee des $K$-groupes pro-r\'eductifs 
est
not\'ee $\bGred(K)$ et s'appelle ca\-t\'e\-go\-rie des 
\emph{$K$-groupes
pro-r\'eductifs
\`a conjugaison 
pr\`es}.
\end{defn}
\end{sloppypar}

\begin{rems}\
\begin{itemize}
\item[a)] Notons $\Aff_K$ la ca\-t\'e\-go\-rie des
$K$-sch\'emas affines. On a 
un diagramme commutatif
\[\begin{CD}
\Gred_K@>\iota>> \Gaff_K@>>> 
\Aff_K\\
@VVV @VVV\\
\bGred_K@>\bar\iota>> \bGaff_K
\end{CD}\]
o\`u 
les foncteurs horizontaux sont pleinement fid\`eles et les foncteurs
verticaux sont pleins et surjectifs.
\item[b)] Un certain nombre de notions ``passent" aux groupes \`a
conjugaison pr\`es (via les foncteurs de projection ci-dessus):
immersions ferm\'ees, mor\-phis\-mes fi\-d\`e\-le\-ment plats, connexit\'e,
simple connexit\'e,  (pro-)uni\-po\-ten\-ce,
(pro)-semi-simpli\-cit\'e\dots
\item[c)] La notion de sous-groupe n'a pas de sens dans $\bGaff(K)$ et
$\bGred(K)$: elle doit \^etre remplac\'ee par celle de classe de
conjugaison de sous-groupes. Par contre, la notion de sous-groupe
distingu\'e a un sens. De m\^eme, si $H$ est un sous-groupe ferm\'e (\`a
conjugaison pr\`es) de $G$, le sous-groupe (distingu\'e)
\begin{equation}\label{eq14.2}
H^\triangleleft=\{\langle gHg^{-1}\rangle\mid g\in G(K)\}
\end{equation}
a un sens.
\item[d)] On a un foncteur fid\`ele \'evident
\[\bGaff(K)\to Lien(K)\]
o\`u $Lien(K)$ d\'esigne la ca\-t\'e\-go\-rie des \emph{liens sur $K$}
\cite[II.2.1.3]{giraud}. Mais ce foncteur n'est pas plein si $K$ n'est
pas alg\'ebriquement clos.
\end{itemize}
\end{rems}

\subsection{Exactitude} Peu de limites inductives existent dans
$\Gaff_K$ et dans $\bGaff_K$. Il faut aussi prendre garde \`a ce que
les li\-mi\-tes inductives calcul\'ees dans $\Gaff_K$, $\Aff_K$ et
$\bGaff_K$, quand elles existent, ne co\"{\i}ncident pas
n\'e\-ces\-sai\-re\-ment.

Une limite inductive qui existe toujours dans $\Gaff_K$ est le
co\-\'e\-ga\-li\-sa\-teur d'un homomor\-phis\-me
$f:G\to H$ et de l'homomor\-phis\-me trivial $1:G\to H$: ce co\'egalisateur
existe aussi dans $\bGaff_K$ et co\"{\i}ncide avec le pr\'ec\'edent.
Dire que ce co\'egalisateur est trivial est plus faible
que de dire que $f$ est un \'epimor\-phis\-me.  Nous adopterons la
terminologie suivante:

\begin{defn}\label{d14.2}
a) Soient
$f:G_1\to G_2$ et
$g:G_2\to G_3$ deux mor\-phis\-mes de
$K$-groupes affines tels que $g\circ f=1$. On dit que la suite
\begin{equation}\label{eq14.3}
\begin{CD}G_1@>f>> G_2@>g>> G_3\end{CD}
\end{equation}
est \emph{exacte} (\resp \emph{faiblement exacte}) si $\Ker g=\IM f$
(\resp si $\Ker g=(\IM f)^\triangleleft$, \cf \eqref{eq14.2}).\\
b) Si $G_3=1$ dans a), on dit que $f$ est
\emph{\'epi} (\resp \emph{faiblement \'epi}) si la suite \eqref{eq14.3}
est exacte (\resp faiblement exacte).
\end{defn}

On remarquera que, 
dans b), dire que $f$ est \'epi \'equivaut
\`a dire que le mor\-phis\-me 
associ\'e de faisceaux fpqc d'ensembles est
un 
\'e\-pi\-mor\-phis\-me, ou encore que $f$ est fi\-d\`e\-le\-ment plat.

\begin{lemme}\label{l15} a) Un mor\-phis\-me $f:G_1\to G_2$ de 
$K$-groupes
affines est
faiblement \'epi si et seulement si, pour 
tout $K$-groupe affine $N$, 
l'application
d'ensembles point\'es 
$f^*:H^1_K(G_2,N)\to H^1_K(G_1,N)$ est de noyau
trivial. Si 
$G_1,G_2\in \Gred_K$, il suffit de prendre $N$
dans $\Gred_K$.\\
b) 
Une suite
\[\begin{CD}G_1@>f>> G_2@>g>> G_3@>>> 1\end{CD}\]
de 
$K$-groupes affines (avec $g\circ f=1$) est faiblement exacte si 
et 
seulement si,
pour tout $K$-groupe affine $N$, la suite d'ensembles 
point\'es
\[\begin{CD}1\to H^1_K(G_3,N)@>g^*>> 
H^1_K(G_2,N)@>f^*>>
H^1_K(G_1,N)\end{CD}\]
est exacte). Si 
$G_1,G_2,G_3\in \Gred_K$, il suffit de prendre $N$
dans 
$\Gred_K$.\\
c) L'exactitude faible d'une suite comme dans b) ne 
d\'epend que des
classes de $f$ et $g$ dans 
$\bGaff_K$.\qed
\end{lemme}

\prf $c)$ est imm\'ediat.

$a)$: 
$G_1\stackrel{f}{\to} G_2$ est faiblement \'epi $\iff \;
f(G_1)^\triangleleft = G_2 \iff \; [\forall\, G_2\stackrel{h}{\to} N$ 
(avec $N$ pro-r\'eductif si $G_2$ l'est), $h\circ f= 1\If h=1] \iff 
\;(\forall \bar h \in H^1_K(G_2,N), f^\ast 
\bar h= 1\If \bar h=1)$.

$b)$: Compte tenu de $a)$, nous n'avons \`a examiner que 
l'exactitude au milieu. Or
 $G_1\stackrel{f}{\to} 
G_2\stackrel{g}{\to}G_3$ est faiblement exacte $\iff \;
f(G_1)^\triangleleft = \ker g \iff \; [\forall\, G_2\stackrel{h}{\to} 
N $ (avec $N$
pro-r\'eductif si
$G_2$ l'est)$ , 
\;h_{\vert
f(G_1)^\triangleleft}= 1\If h_{\vert \ker g}=1] 
\iff \; 
[\forall G_2\stackrel{h}{\to} N ,\;h_{\vert
f(G_1)^\triangleleft}= 1 
\If \exists \, G_3\stackrel{i}{\to} N $ tel que $h=i\circ g)] \iff 
\;
[\forall\,
\bar h
\in H^1_K(G_2,N), f^\ast \bar h= 1\If \exists 
\,\bar i\in H^1_K(G_3,N),$ tel que $ \bar h=
g^\ast\bar i]$. 

\qed

Nous nous autoriserons de ce lemme pour parler de suites 
faiblement 
exactes dans $\bGaff_K$.

\subsection{L'enveloppe 
pro-r\'eductive}

\begin{thm}\label{t1} Le foncteur d'inclusion 
$\bar{\iota}:\bGred_K\to
\bGaff_K$ admet un adjoint \`a gauche 
$G\mapsto \Pred(G)$. Ce
foncteur est un quasi-inverse \`a gauche de 
$\bar{\iota}$.
\end{thm}

\prf Pour chaque $G\in \bGaff_K$, 
choisissons une section
mo\-no\-\"{\i}\-da\-le
$s(G)\allowbreak\in 
\sG_u^\otimes(G)$. On d\'efinit $\Pred(G)$
comme \'etant $G_{s(G)}$. 
Soit $f:G\to H$ un homomor\-phis\-me. En
appliquant la 
proposition
\ref{p3} \`a $G$ et $\phi=s(H)^\sharp\circ f:G\to 
H_{s(H)}$, on obtient
un homomor\-phis\-me $\bar\phi:G_{s(G)}\to 
H_{s(H)}$. La proposition \ref{p2}
montre que  son ima\-ge dans 
$H^1_K(\Pred(G),\Pred(H))$ ne d\'epend pas du
choix de $\bar\phi$: 
c'est $\Pred(f)$. La proposition \ref{p2} montre
aussi que
$\Pred$ 
est un foncteur. Les plongements $G\to G_{s(G)}$ d\'efinissent 
un
mor\-phis\-me de foncteurs $Id_{\bGaff_K}\to \bar \iota\circ\Pred$. Le 
fait
que ce mor\-phis\-me fasse de $\Pred$ un adjoint \`a gauche de 
$\bar\iota$
r\'esulte de la construction pr\'ec\'edente et des 
propositions \ref{p3}
et
\ref{p2}. Enfin, le fait 
que
$\Pred\circ\bar\iota\simeq Id$ est 
\'evident.\qed

\begin{defn}\label{d14.1}
Le $K$-groupe \`a 
conjugaison pr\`es $\Pred(G)$ (muni de l'injection
canonique 
$G\inj
\Pred(G)$) s'appelle \emph{l'enveloppe pro-r\'eductive 
de
$G$}.

Son sous-groupe $\PU(G)$ {(donn\'e par le sous-groupe 
unipotent $U_s$ d\'efini juste avant la prop.
\ref{centr})} s'appelle 
le
\emph{compl\'ement unipotent} de 
$G$.
\end{defn}

\begin{sloppypar}
\begin{rems} a) Si $f:G\to H$ est un 
mor\-phis\-me, nous ignorons si
$\Pred(f)(\PU(G))\subset \PU(H)$ en 
g\'en\'eral. C'est vrai si
le foncteur $f^*:Rep_K(H)\allowbreak\to 
Rep_K(G)$ est
pleinement fid\`ele, en vertu du d\'ebut du \S
\ref{i2} 
(en particulier si $f$ est un
\'epi\-mor\-phis\-me, mais aussi  dans 
le cas d'un sous-groupe
parabolique d'un groupe connexe, voir la 
d\'e\-mon\-stra\-tion de la
proposition
\ref{p7} b)). On pourrait 
aussi le d\'eduire de la proposition
\ref{centr}. Dans ce cas, nous 
ignorons si le mor\-phis\-me correspondant 
est \'epi.
\\ b) Le th\'eor\`eme montre que l'homomor\-phis\-me $\Lambda_G:G\to \Pred(G)$ est la limite projective 
des morphismes $G\to H$ (dans $\bGaff_K$), avec $H\in  \bGred_K$ (de type fini si l'on veut). En revanche, on
ne pourrait pas construire $\Pred(G)$ par une telle limite projective dans $\Gaff_K$ au lieu de $\bGaff_K$. 

\noindent Par exemple, si $G= \bG_a$, on verra plus loin que  $\Pred(G)\cong SL_2$. Consid\'erons la limite
du syst\`eme projectif des groupes $SL_2$, avec pour morphismes de transition les isomorphismes donn\'es par les
\'el\'ements de $PGL_2(K)$: cette limite vaut $SL_2$ dans $\bGred_K$, mais dans $\Gred_K$, elle est
r\'eduite au centre $\{\pm 1\}$ de $SL_2$.    
\end{rems}
\end{sloppypar}

\begin{prop}\label{p14.1} 
Supposons $K$ alg\'ebriquement clos. \`A
isomor\-phis\-me unique pr\`es 
dans
$\bGaff(K)$, l'homomor\-phis\-me $\Lambda_G:G\to \Pred(G)$ 
est
caract\'eris\'e par les deux propri\'et\'es 
suivantes:
\begin{thlist}
\item $\Pred(G)$ est pro-r\'eductif;
\item 
$\Lambda_G$ induit une bijection entre les classes d'isomor\-phis\-mes
de 
re\-pr\'e\-sen\-ta\-tions ind\'e\-com\-po\-sa\-bles de $G$ et de 
$\Pred(G)$.
\end{thlist}
\end{prop}

\prf Il est clair que 
$\Lambda_G$ a les propri\'et\'es annonc\'ees. Pour
la r\'eciproque, 
la propri\'et\'e universelle de $\Pred(G)$ nous ram\`ene
\`a 
d\'emontrer l'\'enonc\'e suivant: si $f:H\to M$ est un 
homomor\-phis\-me
de groupes pro-r\'eductifs qui induit une bijection sur 
les classes
d'isomor\-phis\-mes de re\-pr\'e\-sen\-ta\-tions irr\'eductibles, 
alors
$f$ est un isomor\-phis\-me. 

En effet, les ca\-t\'e\-go\-ries 
$Rep_K(H)$ et $Rep_K(M)$ sont
semi-simples. L'hypoth\`ese implique 
que le foncteur
$f^*:Rep_K(M)\to Rep_K(H)$ est essentiellement 
surjectif. D'autre part,
pour deux
$M$-re\-pr\'e\-sen\-ta\-tions 
ir\-r\'e\-duc\-ti\-bles $S,S'$, l'application
$Hom_M(S,S')\to 
Hom_H(S,S')$ est bijective: si $S$ et $S'$ ne sont pas
isomorphes, 
les deux membres sont nuls, et si $S=S'$ on a
$End_M(S)=End_H(S)=K$ 
(c'est ici qu'on utilise l'hypoth\`ese que $K$
est 
al\-g\'e\-bri\-que\-ment clos). Par cons\'equent, $f^*$ est 
pleinement
fid\`ele, donc c'est une \'equivalence de ca\-t\'e\-go\-ries 
et
$f$ est bien un isomor\-phis\-me.\qed

\begin{rem} Nous ignorons si 
l'hypoth\`ese que $K$ est alg\'ebriquement
clos est n\'ecessaire dans 
la proposition \ref{p14.1}.
\end{rem}

\subsection{Quelques propri\'et\'es de 
l'enveloppe pro-r\'eductive}

\begin{lemme}\label{l16} Si $G$ est 
connexe, $\Pred(G)$ est connexe.
\end{lemme}

Soit $(G_s)^0$ la 
composante neutre de
$G_s$, et soit $\Gamma=G_s/(G_s)^0$: $(G_s)^0$ 
est
pro-r\'eductif connexe et $\Gamma$ est un groupe profini
(de 
dimension $0$). Comme $G$ est connexe,
$Hom(G,\Gamma)=1$; la 
proposition \ref{p2} implique alors que
$Hom(G_s,\Gamma)=1$. Donc 
$\Gamma=1$.\qed

\begin{prop} \label{c3} Si $G$ est pro-unipotent, 
$\Pred(G)$ est
pro-semi-sim\-ple simplement connexe.
\end{prop}

\prf 
D'apr\`es le lemme \ref{l16}, $G_s$ est connexe. Pour montrer 
qu'il
est semi-simple, on raisonne de m\^eme: soit
$G_s'$ son groupe 
d\'eriv\'e. Alors $\Gamma=G_s/G_s'$ est un
pro-tore. Comme $G$ est 
pro-unipotent, $Hom(G,\Gamma)\allowbreak=1$; la
proposition \ref{p2} 
implique alors que $Hom(G_s,\Gamma)=1$. Donc $\Gamma=1$.

Enfin, soit 
$H$ un groupe pro-semi-simple connexe quelconque, et soit
$\tilde H$ 
son rev\^etement universel. On a un diagramme 
commutatif:
\[\begin{CD}
H^1_K(G_s,\tilde H)@>>> H^1_K(G,\tilde 
H)\\
@VVV @VVV\\
H^1_K(G_s,H)@>>> H^1_K(G,H).
\end{CD}\]

Les 
propositions \ref{p3} et \ref{p2} montrent que les 
fl\`eches
horizontales sont des isomor\-phis\-mes. Comme
$G$ est 
unipotent, il en est de m\^eme de la fl\`eche verticale de
droite. On 
en d\'eduit que la fl\`eche verticale de gauche est un
isomor\-phis\-me 
pour tout $H$. En appliquant ceci \`a $H=G_s$, cela
implique que la 
projection $\tilde H\to H$ admet une section \`a
conjugaison pr\`es. 
Mais ceci n'est possible que si $H$ est
pro-simplement 
connexe.
\qed

\begin{contrex} On pourrait se demander si $\Pred(U)$ 
est m\^eme
\emph{d\'eploy\'e}. Le contre-exemple suivant, pour 
$U\simeq \bG_a^4$,
nous a
\'et\'e aimablement fourni par Ulf Rehmann. 
Soit $D$ une alg\`ebre de
quaternions sur $K$. Le groupe $G=SL_{2,D}$ 
contient $U$ (avec
$U(R)=R\otimes_K D$ pour toute $K$-alg\`ebre 
commutative $R$) comme
sous-groupe triangulaire sup\'erieur strict. 
Si
$\Pred(U)$ \'etait d\'eploy\'e, son image $H$ dans $G$ le 
serait
\'egalement. Mais un tore maximal de
$G$ est donn\'e par un 
tore maximal du sous-groupe dia\-go\-nal  form\'e
des
\'el\'ements 
$(x,y)\in D^*\times
D^*$ tels que $Nrd(x)Nrd(y)=1$; de ceci on 
d\'eduit
facilement que le rang d'un tore d\'eploy\'e contenu dans 
$G$ est
$\le 1$. Ainsi $H$ serait de rang $\le 1$; mais alors il ne 
peut pas
contenir $U$.

Par contre, un groupe semi-simple anisotrope 
ne peut pas \^etre un
quotient de $\Pred(U)$, \cf 
\cite[1.5.3]{margulis}.
\end{contrex}

\begin{prop}\label{suite} a) 
Soit $G_2\to G_1$ un mor\-phis\-me faiblement
\'epi de
$K$-sch\'emas en 
groupes affines. Alors le mor\-phis\-me $\Pred(G_2)\to
\Pred(G_1)$ 
correspondant est faiblement \'epi.\\
b) Soit $G_3\to G_2\to G_1\to 
1$ une suite
faiblement exacte de
$K$-sch\'emas en groupes affines. 
Alors on a
un diagramme commutatif de suites faiblement exactes dans 
$\bGaff_K$:
\[\begin{CD}
G_3@>>> G_2@>>> G_1@>>> 1\\
@VVV @VVV 
@VVV\\
\Pred(G_3)@>>> \Pred(G_2)@>>> \Pred(G_1)@>>> 
1.
\end{CD}\]
\end{prop}

\prf Cela r\'esulte imm\'ediatement du th\'eor\`eme \ref{t1}
et du lem\-me \ref{l15} (comme un adjoint \`a gauche commute
aux limites inductives quelconques, le th\'e\-o\-r\`e\-me
\ref{t1} fournit des propri\'et\'es d'exactitude du 
foncteur $\Pred$).
  \qed

\begin{cor}\label{c7} Si $U$ est le 
radical unipotent de $G$ et si
$G^\red=G/U$, on a une suite 
faiblement exacte
\[\Pred(U)\to \Pred(G)\to G^\red\to 
1.\]\qed
\end{cor}
Voici quelques r\'esultats structurels 
suppl\'ementaires sur $\Pred(G)$:

\begin{lemme}\label{l14} Soit 
$f:G\to H$ un homomor\-phis\-me de $K$-groupes
affines. Supposons que le 
foncteur $f^*$ associ\'e soit pleinement
fid\`ele. 
Alors
\begin{thlist}
\item $\Pred(f)$ est \'epi.
\item $\Pred(f)$ est 
un isomor\-phis\-me si et seulement si $f$ est 
un
isomor\-phis\-me.
\end{thlist}
\end{lemme}

\prf (i) Si $f^*$ est 
pleinement fid\`ele, $\bar f^*$ l'est aussi; comme
c'est un foncteur 
entre ca\-t\'e\-go\-ries semi-simples, le foncteur associ\'e
sur les 
$K$-goupes affines est bien \'epi.

(ii) $f$ est un isomor\-phis\-me 
$\iff$ $f^*$ est une
\'equivalence de ca\-t\'e\-go\-ries $\iff$ $f^*$ est 
essentiellement surjectif;
de m\^eme pour $\Pred(f)$ et $\bar f^*$. 
Si $f^*$ est essentiellement
surjectif, il en est \'evidemment de 
m\^eme de $\bar f^*$. La
r\'eciproque est vraie puisque le foncteur 
de projection $\pi_G$ est
plein et conservatif.\qed

\begin{prop}\label{p7} a) Pour tout $G$ alg\'ebrique, on a une suite exacte
\[1\to \Pred(G^0)^\triangleleft\to \Pred(G)\to G/G^0\to 1\]
o\`u $G^0$ est la composante neutre de $G$.\\
b) Si $G$ est connexe, soit
$P$ un $K$-sous-groupe parabolique de $G$.
Alors
$\Pred(P)\to \Pred(G)$ est \'epi.\\
c) Si $G\to H$ est \'epi, $\Pred(G)\to \Pred(H)$ est \'epi.
\end{prop}

\prf a) \'Ecrivons $G=\lim G_i$, o\`u les $G_i=G/N_i$ sont de dimension
finie: alors, pour tout $i$, $G_i^0$ est l'image de $G^0$ dans $G_i$.
Comme le syst\`eme des $G_i^0$ est de Mittag-Leffler, on a un diagramme
commutatif de suites exactes courtes:
\[\begin{CD}
1@>>> G^0@>>> G@>>> G/G^0@>>>1\\
&&@V{A}VV @VVV @V{B}VV\\
1@>>> \lim G_i^0@>>> \lim G_i@>>> \lim G_i/G_i^0@>>> 1.
\end{CD}\]

Comme $\Ker B$ est profini et que $\Coker A$ est connexe, ce diagramme
montre que $A$ et $B$ sont des isomor\-phis\-mes. En particulier, on a
$G^0=\lim G^0/(G^0\cap N_i)$, et $N_i^0:=G^0\cap N_i$ est de codimension 
finie
dans $G^0$ et distingu\'e dans $G$.

Soit maintenant $\rho:G^0\to GL(V)$ une re\-pr\'e\-sen\-ta\-tion de $G^0$ (de
dimension finie): nous allons montrer qu'elle est isomorphe \`a un
facteur direct d'une re\-pr\'e\-sen\-ta\-tion $W$ provenant de $G$. Ceci
impliquera que
$\Pred(G^0)\allowbreak\to \Pred(G)$ est un monomor\-phis\-me, d'o\`u
l'\'enonc\'e en appliquant la proposition \ref{suite} et le lemme
\ref{l14}.

Notons $U=\Ker\rho$: alors $U$ est de codimension finie dans $G^0$, et
d'apr\`es ce qui pr\'ec\`ede il existe $V\subset G^0$, de codimension
finie et distingu\'e dans $G$. Notons $\Gamma^0=G^0/V$, de sorte qu'on a
une suite exacte
\[1\to \Gamma^0\to G/V\to G/G^0\to 1.\]

Il existe un quotient fini $\Delta$ de $G/G^0$ tel que $G/V$ soit l'image
r\'eciproque d'une extension $\Gamma$ de $\Delta$ par $\Gamma^0$. Ainsi,
on s'est ramen\'e au cas o\`u $G$ est alg\'ebrique. Dans ce cas, il suffit
de prendre pour $W$ la restriction \`a $G^0$ de $Ind_{G^0}^G V$.

\begin{sloppypar}
b) Gr\^ace au lemme \ref{l14} (i), il suffit de prouver
que le foncteur $Rep_K(G)\allowbreak\to Rep_K(P)$ est pleinement
fid\`ele. Soient
$V,W$ deux re\-pr\'e\-sen\-ta\-tions de
$G$. Le 
groupe
$G$ op\`ere sur l'espace affine $Hom(V,W)$ via un de ses 
quotients
alg\'ebriques
$\Gamma$. Soit $\bar P$ l'image de $P$ dans 
$\Gamma$: c'est un
sous-groupe parabolique de $\Gamma$. 
Comme
$\Gamma/\bar P$ est propre et  connexe, il en r\'esulte que 
$Hom_G(V,W)\to
Hom_P(V,W)$ est surjectif.\footnote{Nous remercions 
Michel Brion de nous
avoir indiqu\'e ce type d'argument; \cf 
aussi
\protect{\cite[II.4.3.3.2]{saavedra}}.}
\end{sloppypar}

c) 
Cela r\'esulte encore du lemme \ref{l14} 
(i).
\qed

\begin{sloppypar}
\begin{cor}\label{c8} Pour tout 
$K$-groupe affine $G$, on 
a
$\Pred(G^0)=(\Pred(G)^0)^\triangleleft$.
\end{cor}
\end{sloppypar}

\prf En effet, d'apr\`es le lemme \ref{l16}, $\Pred(G^0)$ 
est
connexe, donc aussi $(\Pred(G)^0)^\triangleleft$.\qed

Enfin, on 
a le r\'esultat suivant.

\begin{prop}\label{Tjm2} Pour tout 
$K$-groupe $G$ et tout $K$-groupe
pro-r\'eductif $H$, on a 
$\Pred(G\times H)=\Pred(G)\times H$.
\end{prop}

\prf On peut 
raisonner dans $\bGaff(K)$. Soit $M$ une
re\-pr\'e\-sen\-ta\-tion de 

$G\times H$. Comme $H$ est r\'eductif, la restriction de $M$ \`a $H$ 
est
semi-simple, donc l'alg\`ebre
$End_H(M)$ est semi-simple 
et
$Aut_H(M)=End_H(M)^*$ est le groupe des
$K$-points d'un groupe 
r\'eductif $H'$; on a un homomor\-phis\-me 
$G\to H'$. Ce dernier se 
prolonge en un homomor\-phis\-me $\Pred(G)\to H'$, ce
qui veut  dire que 
l'action de $G\times H$ sur
$M$ s'\'etend en une action de 
$\Pred(G)\times H$.

Consid\'erons le mor\-phis\-me canonique (dans 
$\bGred_K$)
$f:\Pred(G\times H)\to \Pred(G)\times H$. En 
consid\'erant s\'epar\'ement
les mor\-phis\-mes
$\Pred(G)\to 
\Pred(G\times H)$ et $H\to \Pred(G\times
H)$, on voit que c'est un 
\'epi. D'autre part, on vient de voir que, pour
toute 
re\-pr\'e\-sen\-ta\-tion $M$ de $\Pred(G\times H)$, 
l'application
\[H^1_K(\Pred(G)\times H,GL(M))\to H^1_K(\Pred(G\times 
H),GL(M))\]
est surjective. Par cons\'equent, la restriction de $M$ 
\`a $\Ker f$ est
triviale pour tout $M$; mais alors on a $\Ker 
f=\{1\}$.\qed

\begin{sloppypar}
\subsection{Le cas du groupe additif} \`A titre 
d'exemple embl\'ematique,
calculons $\Pred(G)$ et $\PU(G)$ pour 
$G=\bG_a$. Le th\'eor\`eme qui 
suit est une
reformulation 
pr\'ecis\'ee du th\'eor\`eme de 
Jacobson-Morozov\footnote{La preuve 
que nous
donnons de ce r\'esultat classique n'est sans doute ni la 
plus courte 
ni la plus simple!}.
\end{sloppypar}

\begin{thm}\label{Tjm} On a 
$\PU(\bG_a)=\bG_a$ et $\Pred(\bG_a)=SL_2$. 
Plus
pr\'e\-ci\-s\'e\-ment, pour toute section mo\-no\-\"{\i}\-dale $s$ 
de
$\pi_{\bG_a}$, il existe un isomor\-phis\-me de
$(\bG_a)_s$ sur
$SL_2$ 
tel que
$s^\sharp=\phi$, o\`u $\phi$ est comme en 
\eqref{eq2}.
\end{thm}

\prf Soient $V\in \sA$ la re\-pr\'e\-sen\-ta\-tion 
$\phi$ de \eqref{eq2}, et
$\bar V$ son ima\-ge dans $\bar\sA$. C'est 
une re\-pr\'e\-sen\-ta\-tion de rang $2$
de $G_s$, et elle est non triviale 
puisque sa restriction \`a $\bG_a$ via
$s^*$ est non triviale. Comme 
$G_s$ est semi-simple et connexe
(proposition \ref{c3}), son image 
dans $GL_2$ est contenue dans
$SL_2$, donc est \'egale \`a $SL_2$. 
Nous allons montrer que
l'homomor\-phis\-me $f:G_s\to SL_2$ correspondant 
est un isomor\-phis\-me.

\begin{sloppypar}
Pour cela, consid\'erons le 
foncteur $f^*:Rep_K(SL_2)\to
Rep_K(G_s)$: comme $f$ est un \'epi, il 
est pleinement fid\`ele. Son
compos\'e avec $s$ n'est autre que le 
foncteur de restriction
$\Phi:Rep_K(SL_2)\to Rep_K(\bG_a)$. Comme
$K$ 
est de caract\'eristique z\'ero, les 
re\-pr\'e\-sen\-ta\-tions
irr\'eductibles de $SL_2$ sont les puissances 
sy\-m\'e\-tri\-ques
de $V$. D'autre part, la th\'eorie des blocs de Jordan 
montre que les
re\-pr\'e\-sen\-ta\-tions ind\'e\-com\-po\-sa\-bles de $\bG_a$ sont 
\'egalement les
puissances sy\-m\'e\-tri\-ques $S^nV$ de $V$. Ainsi, $\Phi$ 
induit une
bijection entre les classes d'isomor\-phis\-mes de 
re\-pr\'e\-sen\-ta\-tions
irr\'eductibles de $SL_2$ et les classes 
d'isomor\-phis\-mes de
re\-pr\'e\-sen\-ta\-tions ind\'e\-com\-po\-sa\-bles de $\bG_a$. 
D'apr\`es le lemme
\ref{l2}, $f^*$ induit donc une bijection entre 
les classes
d'isomor\-phis\-mes d'objets irr\'eductibles de $Rep_K(SL_2)$ 
et celles de
$Rep_K(G_s)$. Comme $Rep_K(G_s)$ est semi-simple,
$f^*$ 
est essentiellement surjectif, donc est une \'equivalence 
de
ca\-t\'e\-go\-ries, d'o\`u l'assertion.
\end{sloppypar}

Enfin, 
l'\'egalit\'e $s^\sharp=\phi$ et l'\'egalit\'e $U_s=\bG_a$ 
sont
claires.  \qed

\begin{rem}\label{rga} L'espace des 
$\bG_a$-homomor\-phis\-mes $S^mV\to S^nV$
est de dimension $|P(m,n)|$, 
avec
$P(m,n)=\{j\vert \, \vert m-n\vert \leq j\leq m+n,\;j\equiv m+n 
\mod
2\}$.

En effet, comme $S^mV$ est auto-dual, ces homomor\-phis\-mes 
s'identifient
aux invariants sous $\bG_a$ dans
$S^mV\otimes S^nV$. Or 
on a la d\'ecomposition de $SL_2$-modules
$\displaystyle S^mV\otimes 
S^nV\cong
\oplus_{j\in P(m,n)} S^jV$ (Clebsch-Gordan), et chaque 
$S^jV$ n'a qu'une
droite de vecteurs invariants sous 
$\bG_a$.
\end{rem}

Pour compl\'eter l'\'etude de $Rep_K\bG_a$, 
mentionnons le r\'esultat
suivant, qui nous servira 
ult\'erieurement.

\begin{prop}\label{repG_A} $Rep_K\bG_a$ est 
strictement de Wedderburn.
\end{prop}

\prf Il s'agit de montrer que 
le radical infini est nul. Com\-men\-\c 
cons par quelques remarques 
sur les
ind\'e\-com\-po\-sa\-bles de $Rep_K \bG_a$. Ils sont de la 
forme $S^mV$,
o\`u $V$ est la re\-pr\'e\-sen\-ta\-tion fondamentale de 
dimension $2$, et
admettent une unique suite de composition (\`a 
crans de dimension $1$).
Tout homomor\-phis\-me $f$ de $S^m V$ vers 
$S^nV$ envoie le vecteur invariant
de $S^m V$ sur $0$ si
$m>n$ ou si $m=n$ et $f$ n'est pas un isomor\-phis\-me.

\noindent Soient alors $W,W'$ deux re\-pr\'e\-sen\-ta\-tions, et $f_1, \dots
f_N$ une cha\^{\i}ne d'homomor\-phis\-mes dans le radical reliant $W$ et
$W'$. Quitte \`a remplacer $f_1$ (\resp le compos\'e $f$ des $f_i$) par
le mor\-phis\-me correspondant $\tilde f_1:
\un
\to W^\vee\otimes W_1$ (\resp $\tilde f : \un \to W^\vee \otimes W'$),
et $f_i$, pour $i>1$, par $\tilde f_i= 1_{W^\vee}\otimes f_i$  (qui sont
tous dans le radical puisque le radical est un id\'eal mo\-no\-\"{\i}\-dal), on
se ram\`ene au cas $W=\un$ (noter que $\tilde f$ est le compos\'e des
$\tilde f_i$).

En outre, par additivit\'e, on peut supposer les sources et buts des
$f_i$ in\-d\'e\-com\-po\-sa\-bles, \ie de la forme $S^nV$. On a donc une
cha\^{\i}ne de non-isomor\-phis\-mes
\[\un  \stackrel{f_1}{\rightarrow} S^{n_1}V \stackrel{f_2}{\rightarrow}
S^{n_2}V \to \dots   \stackrel{f_N}{\rightarrow} S^{n_N}V \] o\`u $n_N$
est fix\'e (et inf\'erieur au produit des dimensions des $W$ et $W'$
originaux). Si $N>n_N$, on aura donc $n_{i+1}\leq n_i$ pour au moins l'un
des $i$, et alors le compos\'e des $f_j$ de $1$ \`a $i$ s'annulera
d'apr\`es l'observation pr\'ec\'edente.
   \qed

\begin{rems}\label{r1} \
\begin{itemize}
\begin{sloppypar}
\item[a)] Il suit de l'exemple \ref{typehomotop} et du th\'eor\`eme
\ref{gerbe'} que le groupo\"{\i}de
$\Gamma^\otimes(Rep_K\bG_a)$ a le m\^eme type d'homotopie que
$\bG_a\to \Spec K$ (au sens de \loccit).
\end{sloppypar}
\item[b)] Supposons $K$ de caract\'eristique $p>0$. Alors le
radical de $Rep_K(\bG_a)$ \emph{n'est pas mo\-no\-\"{\i}\-dal}.
   En effet, consid\'erons la re\-pr\'e\-sen\-ta\-tion
in\-d\'e\-com\-po\-sa\-ble standard $W$ de dimension
$p$ de  $\bG_a$ donn\'ee par $t\mapsto \exp tn_{p-1}$, o\`u $n_{p-1}$
est l'endomor\-phis\-me nilpotent d'\'echelon
$p$ (un seul bloc de Jordan). Les endomor\-phis\-mes de cette
re\-pr\'e\-sen\-ta\-tion sont donn\'es par des matrices
triangulaires dont les coefficients diagonaux sont tous \'egaux; leur
trace est donc nulle, ce qui montre
qu'avec les notations de \S4, $W$ devient nul dans
$Rep_K(\bG_a)/\sN$. Ainsi $\sR\neq \sN$, et $\sR$
n'est pas mo\-no\-\"{\i}\-dal d'apr\`es \ref{absJannsen}.
\item[c)] Revenons \`a la caract\'eristique nulle. Consid\'erons la
ca\-t\'e\-go\-rie tan\-na\-kien\-ne $\sA= Rep_K \Z$ des  re\-pr\'e\-sen\-ta\-tions de
dimension finie sur $K$ du groupe discret $\Z$. Un r\'esultat de
`folklore', bas\'e sur la  d\'ecomposition de Jordan, dit que son
enveloppe pro-alg\'ebrique est le produit de $\bG_a$ par un
$K$-groupe affine de type multiplicatif $T\times
\mu_\infty$ ($T$ est un pro-tore, $\mu_\infty$ est le groupe 
de
torsion, cyclique). Il d\'ecoule de l\`a, du
th\'eor\`eme 
pr\'ec\'edent et du corollaire \ref{suite}, que
$\bar\sA$ est 
$\otimes$-\'equi\-va\-len\-te \`a
$Rep_K(SL_2\allowbreak\times T\times 
\mu_\infty)$.
\end{itemize}
\end{rems}

\begin{sloppypar}
\subsection{Extension des scalaires} Si l'on fait 
varier le corps de
base $K$ (toujours suppos\'e de caract\'eristique 
nulle), pour un
$K$-groupe affine $G$ et une extension $L/K$, il y a 
lieu de distinguer
entre les $L$-groupes affines
$\Pred(G\times_K L)$ 
et $\Pred(G)\times_K L$. Par la propri\'et\'e
universelle de $\Pred$, 
on a toujours un mor\-phis\-me 
naturel
\begin{equation}\label{eq14.1}
\Pred(G\times_K L)\to 
\Pred(G)\times_K L.
\end{equation}
\end{sloppypar}

Nous verrons ci-dessous que ce mor\-phis\-me \emph{n'est pas un isomor\-phis\-me}
en g\'en\'eral. Par contre:

\begin{thm}\label{t14.1} Si $L/K$ est \emph{finie}, \eqref{eq14.1} est un
isomor\-phis\-me.
\end{thm}

(Il suffit en fait que $L/K$ soit alg\'ebrique, \cf th. \ref{th:basech}.)
\medskip

\prf Rappelons que dans ce cas, l'extension des scalaires \`a la Saavedra
est d\'efinie sur $\sA=Rep_K(G)$, et que l'on a un isomor\-phis\-me de
ca\-t\'e\-go\-ries $\sA_{(L)}\simeq Rep_L(G)$ (remarque \ref{r4} c)).
Consid\'erons le diagramme commutatif, o\`u une barre sup\'erieure
d\'esigne la r\'eduction modulo le radical:
\begin{equation}\begin{CD}\label{eq11}
   \sA_{(L)}\\
   \s \pi'\displaystyle\Big\downarrow
&\s\pi_{(L)}\displaystyle\searrow\\
  \overline{\sA_{(L)}}@>\alpha>> \bar\sA_{(L)}.
\end{CD}\end{equation}

Pour justifier l'existence de ce diagramme (au moyen de \ref{L1}), 
remarquons que la
ca\-t\'e\-go\-rie $\bar\sA_{(L)}$ est semi-simple en vertu de la remarque
\ref{r4} c), et que $\pi_{(L)}$ est, tout comme $\pi$, un foncteur 
plein (et d'ailleurs aussi essentiellement
surjectif).

Toute section $s$ de $\pi$ d\'efinit une section $s_{(L)}$ de
$\pi_{(L)}$, ce qui montre que $\alpha$ est plein et essentiellement
surjectif. On a une section de $\alpha$
\[\overline{s_{(L)}}=\pi'\circ s_{(L)}: \bar\sA_{(L)}\to
\overline{\sA_{(L)}}\]
qui donne un diagramme naturellement commutatif:
\[\begin{CD}
&&\sA_{(L)}@>\omega_{(L)}>>Vec_L\\
&\s s_{(L)}\displaystyle \nearrow& @A{\sigma}AA\\
\bar \sA_{(L)}@>\overline{s_{(L)}}>> \overline{\sA_{(L)}}.
\end{CD}\]

D'autre part, on a
un diagramme (naturellement) commutatif de ca\-t\'e\-go\-ries et
foncteurs:
\[\begin{CD}
Vec_K&&@>>> && Vec_L\\
@A{\omega}AA &&&& @A{\omega_{(L)}}AA\\
\sA@>>> \sA_L@>>> \sA_L^\natural @>\simeq>> \sA_{(L)}\\
@V{\pi}VV @VVV @VVV @VVV\\
\bar\sA@>>> \bar\sA_L@>>> \bar\sA_L^\natural && \bar\sA_{(L)}.
\end{CD}\]

En effet, le corollaire \ref{C3.1.} montre que le foncteur $\sA\to 
\sA_L$
est radiciel et que $(\bar \sA)_L\simeq \overline{\sA_L}$ 
(isomor\-phis\-me de
ca\-t\'e\-go\-ries). De m\^eme pour le carr\'e sui\-vant, 
par le lemme
\ref{corresp}. L'\'equivalence de ca\-t\'e\-go\-ries 
$\sA_L^\natural \cong
\sA_{(L)}$ provient du th\'eor\`eme
\ref{t2}. 
Il en r\'esulte en particulier que $\alpha$ est une
\'equivalence de 
ca\-t\'e\-go\-ries. Alors $\overline{s_{(L)}}$ est
une \'equivalence 
de ca\-t\'e\-go\-ries, ce qui termine la 
d\'emonstration.\qed

\subsection{Enveloppes pro-r\'eductives des 
groupes
pro-unipotents}\label{eprgpu}

   Si $G$ est un $K$-groupe 
pro-unipotent, on sait d\'ej\`a que 
$\Pred(G)$ est pro-semi-simple 
simplement connexe
(proposition
\ref{c3}). Nous allons voir que 
$\Pred(G)$ est en g\'en\'eral de 
dimension infinie et que son calcul 
est un
probl\`eme ``insoluble", le cas de
$G=\bG_a$ \'etant \`a cet 
\'egard exceptionnel.

Consid\'erons d'abord le cas de
$G=\bG_a\times 
\bG_a$. On peut voir que $\Pred(G)$
\emph{n'est pas de dimension 
finie} en remarquant que les homomor\-phis\-mes
$G\to
\bG_a,(a,b)\mapsto 
ax+by, (x,y\in K)$, donnent
lieu \`a une infinit\'e (param\'etr\'ee par $y/x\in \boP^1(K)$) de
re\-pr\'e\-sen\-ta\-tions ind\'e\-com\-po\-sa\-bles non
\'equi\-va\-len\-tes de dimension $2$. Elles sont du reste toutes
sous-quotients de la re\-pr\'e\-sen\-ta\-tion standard de
dimension $4$ donn\'ee par
$\begin{pmatrix}1& a\\ 0&1 \end{pmatrix}\times\begin{pmatrix}1& b\\
0&1\end{pmatrix}$.
   Cela indique aussi que le mor\-phis\-me \eqref{eq14.1} d'``extension
des scalaires" \emph{n'est pas un  isomor\-phis\-me en g\'en\'eral} (ne
serait-ce que pour raison de cardinalit\'e). On notera l'analogie
formelle entre d'une part ce fait et le th\'eor\`eme \ref{t14.1}, et
d'autre part le comportement conjectural des groupes de Galois motiviques
\cite[6.3? et remarque]{serrem}.

\medskip En fait, il s'av\`ere que \emph{la d\'etermination de
\[\bG_{a}\times \bG_{a}\inj\Pred(\bG_{a}\times \bG_{a})\] est (en un sens
convenable) un probl\`eme insoluble}. En effet,
elle inclut la classification des repr\'esen\-ta\-tions
in\-d\'e\-com\-po\-sa\-bles de
$\bG_a\times \bG_a$, ou ce qui revient au m\^eme, des $A$-modules finis
in\-d\'e\-com\-po\-sa\-bles pour
$A=K[[T_1,T_2]]$ (\cf remarque 
\ref{prof} b)). Or la classification des
$A$-modules finis 
in\-d\'e\-com\-po\-sa\-bles o\`u l'action de $A$ se
factorise par le 
cube du radical $\mathfrak{m}$ est d\'ej\`a un probl\`eme
insoluble, 
car $A/\mathfrak{m}^3$ est de type de re\-pr\'e\-sen\-ta\-tion
infini sauvage 
(\cf \eg
\cite{nathanson}\footnote{Qui donne aussi une formalisation 
de la notion
de probl\`eme de classification insoluble ou de 
difficult\'e
maximale.}); plus pr\'ecis\'ement, la th\'eorie 
des
$A/\mathfrak{m}^3$-modules finis est ind\'ecidable 
(\cf
\cite{prest}).

 Les travaux r\'ecents de 
Ginz\-burg et Panyushev sur
les paires nilpotentes \cite{gin,pan} 
permettent peut-\^etre toutefois de
d\'ecrire de nombreux quotients 
de
$\Pred(\bG_{a}\times 
\bG_{a})$.
 
\begin{thm}\label{fin} Soit $G$ un 
$K$-groupe pro-unipotent. Les
propri\'et\'es  sui\-van\-tes sont 
\'equi\-va\-len\-tes:
\\ a) $\dim\,G\leq 1$,
\\ b) $\Pred(G)$ est de 
dimension finie,
\\ c) pour toute extension $L/K$, le mor\-phis\-me 
$\Pred(G_L)\to
\Pred(G)_L$ est un isomor\-phis\-me,
\\ d) le radical 
infini $\rad^\omega(Rep_K G)$ est nul, \ie $Rep_K G$
est strictement 
de Wedderburn,
\\ e) $\rad^\omega(Rep_K G)$ est 
nilpotent.
\end{thm}

\prf D'apr\`es la discussion des cas $\bG_{a }$ 
et $\bG_{a }\times
\bG_{a }$ et le lemme \ref{l14}, il
suffit pour 
l'\'equivalence de $a),b),c)$ de d\'emontrer que $\bG_{a
}\times 
\bG_{a }$ est quotient de
$G$ d\`es que $\dim G> 1$. Passant 
aux
alg\`ebres de Lie, il s'agit de voir que pour toute alg\`ebre de 
Lie
nilpotente  $\mathfrak{n}$ de dimension $> 
1$,
$\mathfrak{n}^{ab}$ est de dimension $>1$. Mais il est bien connu 
que
$\dim\,  \mathfrak{n}^{ab}$  est une borne
inf\'erieure pour le 
nombre de g\'en\'erateurs de $\mathfrak{n}$.

L'implication $d)\If 
e)$ est triviale.

L'implication $a) \If d)$ est la proposition 
\ref{repG_A}.

Pour terminer, prouvons $e)\If a)$ par l'absurde. 
Fixons un
\'epi\-mor\-phis\-me $G\to \bG_a\times
\bG_a$. 
Identifions
$Rep_K\,
(\bG_a\times
\bG_a)$
\`a $A$-$Modf$ 
avec
$A=K[[T_1,T_2]]$ d'une part, et \`a une sous-ca\-t\'e\-go\-rie pleine 
de
$Rep_K G$ d'autre part. Soit
$\mathfrak{m}$ le radical de
$A$, et 
consid\'erons la sous-ca\-t\'e\-go\-rie pleine
$(A/\mathfrak{m}^3)$-$Modf$ de $A$-$Modf$. D'apr\`es
l'hypoth\`ese et le lemme
\ref{radinf}, le radical infini de
$(A/\mathfrak{m}^3)$-$Modf$ serait donc nilpotent, en contradiction,
d'apr\`es \cite{ks}, avec le caract\`ere sauvage de $A/\mathfrak{m}^3$.
   \qed

Au-del\`a du cas unipotent, on peut se poser la question suivante:

\begin{qn}\label{q19.1} Pour quels groupes alg\'ebriques lin\'eaires $G$
l'enveloppe pro-r\'e\-duc\-ti\-ve $\Pred(G)$ est-elle de dimension finie?
Sa formation est-elle alors compatible \`a l'extension des scalaires?
\end{qn}
 
La r\'eponse est donn\'ee dans l'appendice \ref{os}: ceci se produit si
et seulement si $G$ est de dimension finie et son radical unipotent $U$
est de dimension $\le 1$. De plus, $\Pred(G)$ est alors produit
semi-direct de $G/U$ par $SL_2$.

\medskip  Un cas \'eclairant est celui du produit semi-direct $G$ de $SL_2$
par $\bG_a\times \bG_a$ (d\'efini par la
re\-pr\'e\-sen\-ta\-tion standard de $SL_2$). Ce cas a \'et\'e examin\'e dans
\cite{piard}.
L'auteur y montre que pour tout $n\geq 1$, l'alg\`ebre de Lie
$sl_{2n+1}$ contient une copie de $Lie\,G$ mais
aucune sous-alg\`ebre de Lie semi-simple interm\'ediaire. Comme on
sait par la proposition \ref{c3} et le
corollaire \ref{c7} 
que
$\Pred(G)$ est pro-semi-simple simplement connexe, on en d\'eduit 
que 
$\Pred(G)$ admet
$SL_{2n+1}$ comme quotient. Par application du lemme de Goursat, il
admet donc aussi $\displaystyle\prod_n
SL_{2n+1}$ comme quotient.

\section{Applications aux groupes alg\'ebriques et aux
re\-pr\'e\-sen\-ta\-tions in\-d\'e\-com\-po\-sa\-bles.}\label{appl}

En d\'epit du fait que l'enveloppe pro-r\'eductive $\Pred G$ d'un
$K$-groupe alg\'ebrique lin\'eaire $G$ soit le plus souvent de dimension
infinie, son existence permet d'ob\-te\-nir des r\'esultats  concrets
sur les groupes r\'eductifs\footnote{On rappelle que dans ce texte,
r\'eductif n'implique pas connexe.} contenant
$G$ et sur les re\-pr\'e\-sen\-ta\-tions in\-d\'e\-com\-po\-sa\-bles de $G$.

Dans tout ce paragraphe, $K$ est un corps de caract\'eristique nulle.

\subsection{Applications aux groupes alg\'ebriques}

\begin{prop}[\cf \protect\cite{mor2}]\label{p13.1} Soit $H'$ un
sous-groupe r\'eductif ferm\'e \allowbreak d'un groupe r\'eductif
$H$. Alors le centralisateur de $H'$ dans $H$ est r\'e\-duc\-tif.
\end{prop}

\prf Soit $U$ le radical unipotent de ce centralisateur $C_H(H')$. Il
s'agit de montrer que tout homomor\-phis\-me
$f:\bG_a\to U$ est trivial. Un tel homomor\-phis\-me s'\'etend en un
homomor\-phis\-me $g: \bG_a\times H'\to H$,
puis, d'apr\`es la proposition \ref{Tjm2} et le th\'eor\`eme \ref{Tjm}, en
un  homomor\-phis\-me $g':SL_2\times H'\to H$. On a donc
un homomor\-phis\-me $f':SL_2\to C_H(H')$ qui prolonge $f$. La
compos\'ee de $f'$ avec la projection
$C_H(H')\to C_H(H')/U$ est triviale (puisqu'il en est de m\^eme de
$f$). On en d\'eduit que $f'$ lui-m\^eme, et
par suite $f$, est trivial.\qed

Soient $G,H$ deux $K$-groupes lin\'eaires et
$f:G\to H$ un $K$-ho\-mo\-mor\-phis\-me. (Cas le plus int\'eressant: un
monomor\-phis\-me\dots) On note $C_H(f)$ le centralisateur de $f(G)$ dans
$H$. Nous nous int\'eressons au cas o\`u $H$ est \emph{r\'eductif}.

\begin{defn}\label{d15.1} On appelle \emph{enveloppe r\'eductive de $f$} 
tout sous-groupe r\'eductif ferm\'e de $H$ contenant $f(G)$ et
\emph{minimal} pour cette propri\'et\'e. Si $f$ est un monomor\-phis\-me, on
dira aussi \emph{enveloppe r\'eductive de $G$ dans $H$}.
\end{defn}

\begin{thm}\label{T3} a) Deux enveloppes r\'eductives de $f$
sont con\-ju\-gu\'e\-es par un \'el\'ement $h \in C_H(f)(K)$.\\
b) Soit $\pi:H\to H_0$ un homomor\-phis\-me de groupes r\'eductifs. Alors
l'ima\-ge dans $H_0$ de toute enveloppe r\'eductive de $f$ est une 
enveloppe r\'eductive de $\pi\circ f$.\\ 
c) Supposons que l'enveloppe r\'eductive de $f$ soit \'egale \`a
$H$. Alors $C_H(f)$ est produit du centre de $H$ par un groupe unipotent.
\end{thm}

\prf Nous allons prouver que toute enveloppe r\'eductive $H'$ de $f$
est image dans $H$ du $K$-groupe pro-r\'eductif $G_s$ attach\'e \`a une
section mo\-no\-\"{\i}\-dale $s$ comme au \S \ref{neu}
(et r\'eciproquement). En effet, il suit de la proposition \ref{p3}
que le plongement $G\to H'$ se factorise par
$G_s$. Comme l'image de $G_s$ dans $H'$ est r\'eductive et contient
$G$, elle est \'egale \`a $H'$ par
minimalit\'e de $H'$. R\'eciproquement, le m\^eme argument montre que
pour toute factorisation de $G\to H$ \`a
travers $G_s$, l'image $H'$ de $G_s$ dans $H$ est une enveloppe
r\'eductive de $f$.

Le point b) suit imm\'ediatement de cette interpr\'etation\footnote{comme nous l'a fait observer P.
O'Sullivan, le point b) se d\'eduit aussi directement de ce que l'image inverse dans $H'$ (l'enveloppe
r\'eductive de $G$ dans $H$) de tout sous-groupe r\'eductif de $\pi(H')$ est un sous-groupe r\'eductif de
$H'$.}.

a) Ce point suit alors de la proposition \ref{p2} (avec $s=t,\;f\circ
s^\sharp=g\circ s^\sharp$).

Enfin, c) a d\'ej\`a \'et\'e d\'emontr\'e (proposition \ref{centr}).
   \qed  

Des compl\'ements \`a ce r\'esultat se trouvent dans l'appendice \ref{os}
(\ref{sec:app}).

\begin{contrex} On pourrait se demander si la propri\'et\'e c)
\emph{ca\-rac\-t\'e\-ri\-se} les enveloppes pro-r\'eductives. C'est
(totalement) faux, comme on le voit sur l'exemple de $\bG_a$ plong\'e
dans $SL_{n+1}$ par l'interm\'ediaire de la puissance sy\-m\'e\-tri\-que
$n$-i\`eme d'une re\-pr\'e\-sen\-ta\-tion fid\`ele dans $SL_2$: le
centralisateur de $\bG_a$ est alors \'egal au produit d'un groupe
isomorphe \`a
$\bG_a^n$ par le centre de $SL_{n+1}$.
\end{contrex}

\begin{rem}\label{une des dernieres rem} Pour toute enveloppe 
r\'eductive $H'$ de $G$ dans $GL(V)$, il
suit de l'interpr\'etation de $H'$ donn\'ee dans la d\'emonstration du
th\'eor\`eme \ref{T3} que la restriction
\`a $G$ des re\-pr\'e\-sen\-ta\-tions de
$H'$  induit une \emph{injection} de l'ensemble des classes
d'isomor\-phis\-mes de re\-pr\'e\-sen\-ta\-tions irr\'eductibles de
$H'$ vers l'ensemble des classes d'isomor\-phis\-mes de re\-pr\'e\-sen\-ta\-tions
ind\'e\-com\-po\-sa\-bles de $G$. Les re\-pr\'e\-sen\-ta\-tions
ind\'e\-com\-po\-sa\-bles de $G$ que l'on obtient ainsi sont, \`a
isomor\-phis\-me pr\`es, les facteurs directs ind\'e\-com\-po\-sa\-bles des sommes 
finies
$\displaystyle \oplus_i ({\check V})^{m_i}\otimes V^{n_i}$.

Plus pr\'ecis\'ement, soit $Rep(G,V)$ la sous-ca\-t\'e\-go\-rie pleine de 
$Rep_K G$ form\'ee des facteurs directs
  des sommes finies
$\displaystyle \oplus_i ({\check V})^{m_i}\otimes V^{n_i}$. C'est une 
ca\-t\'e\-go\-rie mo\-no\-\"{\i}\-dale rigide
pseudo-ab\'elienne. Son quotient par le radical n'est autre que la 
sous-ca\-t\'e\-go\-rie pleine du quotient de
$Rep_K G$ par son radi\-cal form\'ee des facteurs directs
  des sommes finies
$\displaystyle \oplus_i ({\check V})^{m_i}\otimes V^{n_i}$; c'est 
donc, \`a $\otimes$-\'equivalence pr\`es, la
sous-ca\-t\'e\-go\-rie tannakienne de
$Rep_K G_s$ engendr\'ee par $V$, qui s'identifie \`a $Rep_K \big( 
Im(G_s\to GL(V)) \big)$ d'apr\`es ce qui
pr\'ec\`ede.\end{rem}

\subsection{Applications aux re\-pr\'e\-sen\-ta\-tions
in\-d\'e\-com\-po\-sa\-bles}\label{15.2}
   Soit $V$ une re\-pr\'e\-sen\-ta\-tion fid\`ele de $G$. On note 
(abusivement) $G_V$ une enveloppe r\'educ\-tive
de $G$ dans $GL(V)$ (c'est-\`a-dire,
d'apr\`es ce qui pr\'ec\`ede, le groupe r\'educ\-tif image dans
$GL(V)$ du $K$-groupe pro-r\'eductif $G_s$ attach\'e \`a une section
mo\-no\-\"{\i}\-dale $s$ comme dans \ref{neu}).
Alors la d\'ecomposition en in\-d\'e\-com\-po\-sa\-bles de toute somme
finie
$\displaystyle \oplus_i ({\check V})^{m_i}\otimes V^{n_i}$ est
d\'etermin\'ee par la d\'ecomposi\-tion en irr\'eductibles de
$\displaystyle \oplus_i ({\check V})^{m_i}\otimes V^{n_i}$ vue comme
re\-pr\'e\-sen\-ta\-tion de $G_V$. Cela ram\`ene un certain nombre de questions
sur les re\-pr\'e\-sen\-ta\-tions \`a la d\'eter\-mi\-nation (m\^eme partielle) de
$G_V$ et \`a la th\'eorie des poids.

Voici quelques \'echantillons d'application de ce principe.

\medskip
On suppose dor\'enavant {\it $G$ connexe et $K$
al\-g\'e\-bri\-que\-ment clos} (de ca\-rac\-t\'e\-ris\-ti\-que nulle).
   Le groupe r\'eductif $G_V$ est alors connexe.

\medskip Consid\'erons le sous-groupe ab\'elien $\mathbf{a}_K(G,V)$ de
l'anneau des re\-pr\'e\-sen\-ta\-tions $\mathbf{a}_K(G)$ engendr\'e par
les classes des objets de $Rep(G,V)$, c'est-\`a-dire des facteurs 
directs des $\displaystyle
\oplus_i ({\check V})^{m_i}\otimes V^{n_i}$ (groupe de Grothendieck
vis-\`a-vis des sommes directes). C'est en fait un sous-anneau, et
m\^eme un sous-$\lambda$-anneau.

En outre, le foncteur $(s^\sharp)^\ast$ induit un isomor\-phis\-me de
$\lambda$-anneaux \[\mathbf{a}_K(G,V)\cong
\mathbf{a}_K(G_V)= R_K(G_V).\]
On rappelle que $R_K(G_V)$ s'identifie,
via l'isomor\-phis\-me ``caract\`ere" $\rm ch$, avec l'anneau des invariants
sous le groupe de Weyl de l'anneau de groupe $\Z[\Lambda]$ sur le
r\'eseau des poids $\Lambda$ (r\'eseau de rang $m=$ rang de $G_V$).

On obtient comme premi\`ere application:

\begin{thm} Si $G$ est simplement connexe, $\mathbf{a}_K(G,V)$ est un
anneau de poly\-n\^o\-mes \`a coefficients entiers en $m\leq \dim V$
ind\'etermin\'ees.
\end{thm}

\prf On montre comme dans la proposition \ref{c3} que $G_V$ est
semi-simple simplement connexe. On a alors
$\mathbf{a}_K(G_V)= R_K(G_V)= R_K(Lie\,G_V)$, et il est bien connu que
ce dernier est un anneau de poly\-n\^o\-mes en les poids fondamentaux de
$Lie\,G_V$.
\qed

\begin{rem} Si $V$ est une re\-pr\'e\-sen\-ta\-tion fid\`ele de $G$. 
D'apr\`es \ref{h0}, \ref{hh0}, et la remarque
\ref{une des dernieres rem}, on obtient \[HH_0(Rep (G,V))\cong 
HH_0(Rep_K G_V)\cong
\mathbf{a}_K(G_V)\otimes_\Z K\cong \mathbf{a}_K(G,V)\otimes_\Z K
. \]
  En particulier, l'homologie de Hochschild (en degr\'e $0$) de $Rep (G,V)$
est une alg\`ebre de poly\-n\^o\-mes sur $K$ si $G$ est simplement connexe.
\end{rem}

Pour toute re\-pr\'e\-sen\-ta\-tion $V$ (de dimension finie), on note $S^m V$ la
puissance sy\-m\'e\-tri\-que $m$-i\`eme de $V$.

\begin{thm} Les conditions
suivantes sont \'equi\-va\-len\-tes:
\\ a) $S^2V$ est ind\'e\-com\-po\-sa\-ble,
\\ b) pour tout $n\geq 0$, $S^n V$ est ind\'e\-com\-po\-sa\-ble.
\end{thm}

\prf  Supposons $a)$. Il est \'equivalent de dire que $S^nV$ est une
re\-pr\'e\-sen\-ta\-tion in\-d\'e\-com\-po\-sa\-ble de $G$ ou que vue comme
re\-pr\'e\-sen\-ta\-tion de $G_V$, $S^nV$ est irr\'eductible. On est donc
ramen\'e au cas ou $G$ est r\'eductif connexe, et le r\'esultat est alors
prouv\'e dans \cite[app.]{bbh} (en outre, \loccit montre que $G_V$ est
\'egal soit \`a $Z(G_V).SL(V)$ soit \`a $Z(G_V).Sp(V)$ pour une forme
symplectique convenable, et le centre $Z(G_V)$ de $G_V$ est r\'eduit aux
homoth\'eties de $G_V$).
\qed

\begin{prop} Soient $V,W$ deux re\-pr\'e\-sen\-ta\-tions de $G$, et $n$ un
entier $>0$. Les conditions suivantes sont
\'equi\-va\-len\-tes:
\\ i) $V\cong W$,
\\ ii) $V^{\otimes n}\cong W^{\otimes n}$.
\end{prop}

\prf Comme pr\'ec\'edemment, le probl\`eme se ram\`ene au pro\-bl\`e\-me
analogue avec $G$ remplac\'e par le groupe r\'eductif connexe $G_{V\oplus
W}$. Via l'homomor\-phis\-me ``caract\`ere", $i)$ (\resp $ii)$) se traduit
par l'\'egalit\'e ${\rm ch}(V)={\rm ch}(W)$ (\resp ${\rm ch}(V)^n={\rm
ch}(W)^n$) de poly\-n\^o\-mes de Laurent en les poids fondamentaux, \`a
coefficients entiers positifs. Or tout poly\-n\^o\-me de Laurent \`a
coefficients entiers positifs est d\'etermin\'e par sa puissance
$n$-i\`eme. \qed

\newpage

\renewcommand{\appendixname}{Appendice}

\appendix

\section{Des ca\-t\'e\-go\-ries semi-simples}\label{s10}

L'objet de cet appendice est de
clarifier la notion de semi-simplicit\'e dans les $K$-ca\-t\'e\-go\-ries (non
n\'e\-ces\-sai\-re\-ment ab\'eliennes), un anneau commutatif unitaire $K$ \'etant
fix\'e.  On renvoie au \S \ref{1.3} et au d\'ebut du \S
\ref{s2} pour les d\'efinitions de
base.

Sauf mention du contraire, \emph{$\sA$-module} signifie
\emph{$\sA$-module \`a gauche} dans tout l'appendice.

\subsection{Objets projectifs et injectifs} Par un raisonnement classique
\cite{tohoku}, la ca\-t\'e\-gorie
ab\'elienne $K$-lin\'eaire $\sA\hbox{--}Mod$ des
$\sA$-modules
(\`a gauche)
poss\`ede assez d'injectifs et de projectifs. Pour les projectifs, nous
allons retrouver ce fait de mani\`ere constructive.

\begin{lemme}\label{lA.1} Soit $f:M\to N$ un homomor\-phis\-me de
$\sA$-modules. Alors\\
a) $f$ est un \'epimor\-phis\-me si et seulement si $f(A):M(A)\to N(A)$ est
un \'epimor\-phis\-me pour tout $A\in \sA$.\\
b) $f$ est un monomor\-phis\-me si et seulement si $f(A):M(A)\to N(A)$ est
un monomor\-phis\-me pour tout $A\in \sA$.\\
c) $f$ est un isomor\-phis\-me si et seulement si $f(A):M(A)\to N(A)$ est
un isomor\-phis\-me pour tout $A\in \sA$.
\end{lemme}

\prf a) Soit $C=\Coker f$: on a $C(A)=\Coker f(A)$ pour tout $A\in \sA$,
et $C=0$ $\iff$ $C(A)=0$ pour tout $A\in \sA$.

b) R\'esulte de a) par dualit\'e (ou par le m\^eme raisonnement).

c) est clair.\qed

\begin{lemme} \label{lA.2} Pour tout objet $A\in \sA$, le $\sA$-module
$\sA_A$ est projectif.
\end{lemme}

\prf Soit $f:M\to\sA_A$ un \'epimor\-phis\-me. Par le lemme \ref{lA.1}, $f(A)$
est surjectif. Soit $m\in M(A)$ tel que $f(A)(m)=1_A$. Par le lemme de
Yoneda (proposition \ref{p1.1} a)), $m$ d\'efinit un homomor\-phis\-me
$\tilde m:\sA_A\to M$, et on voit tout de suite que $\tilde m$ est une
section de $f$.\qed

\begin{defn}\label{dA.1} Un $\sA$-module de la forme $\bigoplus
\sA_{A_i}$ est appel\'e un \emph{$\sA$-mo\-du\-le libre}.
\end{defn}

Tout $\sA$-module libre est projectif. Si $\sA$ est $K$-lin\'eaire, tout
$\sA$-module libre est de la forme $\sA_A$ pour un objet $A$ convenable.

\begin{prop}\label{pA.1} Supposons $\sA$ petite.\\
a) Tout $\sA$-module
$M$ est quotient d'un $\sA$-module libre.\\
b) Un $\sA$-module $P$ est projectif si et seulement s'il est facteur
direct d'un $\sA$-module libre.
\end{prop}

\prf a) Pour tout objet $A\in \sA$, choisissons un syst\`eme
g\'e\-n\'e\-ra\-teur
$(m_i^A)_{i\in I_A}$ de $M(A)$. Alors l'homomor\-phis\-me
\[\bigoplus_{A\in \sA}\bigoplus_{i\in I_A} \sA_A\to M\]
d\'efini par le lemme de Yoneda est un \'epimor\-phis\-me par le lemme
\ref{lA.1} a).

b) La n\'ecessit\'e r\'esulte de a). La suffisance r\'esulte du lemme
\ref{lA.2}, puis\-qu'un facteur direct d'un module projectif est
\'evidemment projectif.\qed

\begin{cor}\label{cA.1} Si $P$ est un $\sA$-module projectif,
$P_L:A\mapsto L\otimes_K P(A)$ est un $\sA_L$-module projectif pour toute
extension $L/K$.
\end{cor}

\prf D'apr\`es la proposition \ref{pA.1} b), il suffit de le voir pour
$P$ de la forme $\sA_A$, auquel cas c'est \'evident.\qed

\subsection{Ca\-t\'e\-go\-ries semi-simples}
Les d\'efinitions suivantes sont adapt\'ees de \cite{street}.

\begin{defn}\label{A.D6} Soit $\sA$ une $K$-ca\-t\'e\-go\-rie. \\
a) Un objet de $\sA$ est \emph{simple} s'il n'est pas nul et s'il n'a
pas d'autre sous-objet que lui-m\^eme et $0$.\\
b) Un objet de $\sA$ est \emph{semi-simple} s'il est somme directe
d'objets simples.\\
c) Un objet $A$ de $\sA$ est \emph{artinien} si toute cha\^{\i}ne
d\'ecroissante de sous-objets de $A$ est stationnaire.\\
d) La ca\-t\'e\-go\-rie $\sA$ est \emph{artinienne} si, pour tout $A\in
\sA$, $\sA_A$ est un objet artinien de $\sA$-$Mod$.\\
e) La ca\-t\'e\-go\-rie $\sA$ est \emph{simple} si elle l'est en tant qu'objet
de $\sA^e$-$Mod$.
\end{defn}

\begin{sloppypar}
\begin{defn}\label{A.D7} Soit $(\sA_\alpha)$ une famille de ca\-t\'e\-go\-ries
ayant les m\^emes objets.\\
a) Le \emph{produit local} des $\sA_\alpha$ est la ca\-t\'e\-go\-rie
$\prod^\lambda_\alpha \sA_\alpha$ ayant les  m\^emes objets que
les $\sA_\alpha$ et telle que, pour tout couple d'objets $(A,B)$, on ait
$(\prod^\lambda_\alpha \sA_\alpha)(A,B)=\prod_\alpha\sA_\alpha(A,B)$.\\
b) $\sA$ est \emph{somme directe locale} des $\sA_\alpha$ si elle est
produit local des $\sA_\alpha$ et que, de plus, pour tout objet $A$,
l'ensemble des $\alpha$ tels que $(\sA_\alpha)_A\ne 0$ est fini.\\
c) Supposons que les $\sA_\alpha$ soient des $K$-ca\-t\'e\-go\-ries. Le
\emph{coproduit} des $\sA_\alpha$ est la $K$-ca\-t\'e\-go\-rie
$\sA=\coprod_{\alpha}\sA_\alpha$ dont la collection d'objets est la
r\'eunion des collections d'objets des $\sA_\alpha$, avec
\[\sA(A,B)=
\begin{cases}
\sA_\alpha(A,B)&\text{si $\exists \alpha: A,B\in\sA_\alpha$}\\
0&\text{sinon}.
\end{cases}\]
\end{defn}
\end{sloppypar}

\begin{lemme}\label{A.L9} Soit $\sA$ une $K$-ca\-t\'e\-go\-rie.\\
a) Toute somme directe d'objets semi-simples est semi-simple.\\
b) Soit $A\in \sA$ un objet tel que tout sous-objet de $A$ soit facteur
direct. Alors $A$ est semi-simple dans les cas suivants:
\begin{thlist}
\item $A$ est artinien.
\item $\sA$ est de la forme $\sR$-$Mod$, o\`u $\sR$ est
une $K$-ca\-t\'e\-go\-rie.
\end{thlist}
Supposons de plus $\sA$ ab\'elienne.\\
c) Soit $A\in \sA$ un objet semi-simple et soit $B$ un
sous-objet de $\sA$. Alors il existe une famille $(S_i)_{I\in
I}$ de sous-objets simples de $A$ telle que $A=B\oplus
\bigoplus_{i\in I} S_i$.\\
d) Tout quotient, tout sous-objet d'un objet semi-simple est semi-simple.
\end{lemme}

\prf a) est \'evident. Dans b), supposons d'abord $A$ artinien. Si $A$
n'est pas semi-simple, soit $B$ un sous-objet de $A$ minimal parmi ceux
qui ne sont pas semi-simples. Alors $B$ n'est pas simple, donc poss\`ede
un sous-objet $C\ne 0,B$. Comme $C$ est facteur direct de $A$, il est
facteur direct de $B$, soit $B=C\oplus D$. Mais $C$ et $D$ sont
semi-simples, contradiction. L'autre cas se traite comme la suffisance
dans \cite[d\'em. de la prop. 4.1]{ce} (\cf \cite[bas p. 140]{street}). c)
se d\'emontre comme la n\'ecessit\'e  dans \loccit d) r\'esulte
imm\'ediatement de  c).\qed

\begin{contrex} \label{A.CeX1} Soit $\sA$ la ca\-t\'e\-go\-rie des espaces
vectoriels de dimension d\'enombrable sur un corps $k$, et soit $\bar
\sA$ la ca\-t\'e\-go\-rie quotient de $\sA$ par la sous-ca\-t\'e\-go\-rie \'epaisse
des espaces vectoriels de dimension finie. Soit $V$ un espace
vectoriel de dimension infinie, vu comme objet de $\bar \sA$. Alors
tout sous-objet de $V$ est facteur direct, mais $V$ ne contient aucun
sous-objet simple.
\end{contrex}

\begin{lemme} \label{A.L10} Soit $\sA$ une ca\-t\'e\-go\-rie ab\'elienne.
Consid\'erons les \'enonc\'es sui\-vants:
\begin{enumerate}
\item \label{abss}Tout objet de $\sA$ est semi-simple.
\item \label{abproj}Tout objet de $\sA$ est projectif.
\item \label{abinj}Tout objet de $\sA$ est injectif.
\end{enumerate}
Alors \ref{abss} $\If$ \ref{abproj} $\iff$ \ref{abinj}. De plus,
\ref{abss} $\iff$ \ref{abproj} $\iff$ \ref{abinj} dans les deux cas
suivants:
\begin{thlist}
\item Tout objet de $\sA$ est artinien.
\item $\sA$ est de la forme $\sR$-$Mod$, o\`u $\sR$ est une
$K$-ca\-t\'e\-go\-rie.
\end{thlist}
\end{lemme}

\prf Il est bien connu que \ref{abproj} $\iff$ \ref{abinj}, et \ref{abss}
$\If$ \ref{abproj} r\'esulte du lemme \ref{A.L9} a). Les implications
inverses  r\'esultent du lemme \ref{A.L9} b).
\qed

\begin{contrex} \label{A.CeX2} Dans la ca\-t\'e\-go\-rie $\bar \sA$ du
contre-exemple \ref{A.CeX1}, tout objet est projectif, mais $\bar\sA$ ne
contient aucun objet simple (\cf \cite[p. 324]{popescu}).
\end{contrex}

\begin{lemme}[\cf \protect{\cite{jannsen}}] \label{A.L11}  Soit $\sA$ une
petite ca\-t\'e\-go\-rie $K$-lin\'eaire pseudo-a\-b\'e\-lien\-ne et dont tout
objet est semi-simple.  Alors $\sA$ est ab\'elienne si et seulement si,
pour tout objet simple $S\in \sA$, l'anneau $\sA(S,S)$ est un corps.
\end{lemme}

\prf La n\'ecessit\'e r\'esulte du lemme de Schur. Suffisance: soit $T$
l'ensemble des types d'objets simples de $\sA$, et pour tout $t\in T$ soit
$\sA_t$ la sous-ca\-t\'e\-go\-rie  pleine de $\sA$ form\'ee des objets
isotypiques de type $t$. Alors $\sA_t$ v\'erifie encore l'hypoth\`ese, et
il suffit de montrer que $\sA_t$ est ab\'elienne. En d'autres termes, on
peut supposer $\sA$ isotypique. Soit $S$ un objet simple de $\sA$, et
notons $D=\sA(S,S)$ son corps d'endomor\-phis\-mes. Le foncteur $\sA^o_S$
s'enrichit en un foncteur $T$ \`a valeurs dans les $D^o$-espaces
vectoriels, o\`u $D^o$ est le corps oppos\'e \`a $D$. En utilisant
l'hypoth\`ese, on voit tout de suite que $T$ est pleinement fid\`ele
et d'image essentielle la ca\-t\'e\-go\-rie des $D^o$-espaces vectoriels de
dimension $\le c$, o\`u $c$ est le plus grand cardinal tel qu'il
existe un ensemble $I$ de cardinal $c$ tel que $S^{(I)}\in \sA$.
En particulier, $\sA$ est ab\'elienne.\qed

\begin{contrex}\label{A.CeX3} Soit $\Q[\epsilon]$ l'alg\`ebre des
nombres duaux ($\epsilon^2=0$). La ca\-t\'e\-go\-rie des $\Q[\epsilon]$-modules
libres de rang fini est $\Q$-lin\'eaire, pseudo-ab\'elienne, artinienne,
et tout objet  est semi-simple, mais elle n'est pas a\-b\'e\-lien\-ne.
\end{contrex}

\begin{lemme} \label{A.L12} Soient $\sA$ une $K$-ca\-t\'e\-go\-rie
pseudo-ab\'elienne, $A\in \sA$ et $\sA_A=\bigoplus_{\alpha\in I}
M_\alpha$ une d\'ecomposition de $\sA_A$ en somme directe de
sous-modules non nuls. Alors $I$ est fini et il existe une unique
d\'ecomposition en somme directe $A=\bigoplus A_\alpha$ telle que, pour
tout $\alpha$, $M_\alpha=\sA_{A_\alpha}$. Si $A$ est ind\'e\-com\-po\-sa\-ble,
$\sA_A$ est ind\'e\-com\-po\-sa\-ble.
\end{lemme}

\prf \'Ecrivons $Id_A=\sum e_\alpha$ avec $e_\alpha\in M_\alpha(A)$ pour
tout $\alpha$. Alors,  pour tout $B\in \sA$ et tout $f\in
\sA_A(B)=\sA(A,B)$, on a $f=\sum  fe_\alpha$. Par cons\'equent,
$M_\alpha\ne 0$ $\If$ $e_\alpha\ne  0$. Les $M_\alpha$ sont donc en nombre
fini. De plus, le lemme de  Peirce montre que les $e_\alpha$ forment un
syst\`eme  d'idempotents orthogonaux, d'o\`u la deuxi\`eme assertion en
utilisant le lemme de Yoneda. La derni\`ere assertion en r\'esulte
imm\'ediatement.\qed

Le th\'eor\`eme qui suit clarifie le lien entre diverses notions de
semi-sim\-pli\-ci\-t\'e, et montre qu'elles ne d\'ependent pas de la
$K$-structure.

\begin{thm} \label{A.P4} Soit $\sA$ une $K$-ca\-t\'e\-go\-rie.
Les conditions sui\-van\-tes sont \'equi\-va\-len\-tes:
\begin{enumerate}
\item
\label{proj} Tout objet de $\sA$-$Mod$ est projectif.
\item
\label{inj} Tout objet de $\sA$-$Mod$ est injectif.
\item
\label{mod-ss} Tout objet de $\sA$-$Mod$ est semi-simple.
\item
\label{xss} Pour tout $A\in \sA$, $\sA_A$ est semi-simple.
\item
\label{hss} $\sA$ est artinienne et le $\sA^o\boxtimes_K
\sA$-module $\sA$ est semi-simple.
\item
\label{struct} $\sA$ est \'equi\-va\-len\-te au sens de Morita \`a une
ca\-t\'e\-go\-rie de la forme
$\coprod_\alpha D_\alpha$, o\`u, pour tout $\alpha$, $D_\alpha$ est un
corps gauche.
\\
On dit alors que $\sA$ est \emph{semi-simple} (\cf d\'efinition \ref{D3}).
\\
  Si de plus $\sA$ est $K$-lin\'eaire,
  les conditions ci-dessus \'equivalent \`a:
\item
  \label{end} Pour tout objet $A\in\sA$, l'anneau $\sA(A,A)$ est 
semi-simple.
  \item
\label{rad} $\rad(\sA)=0$ et pour tout $A\in \sA$, $\sA(A,A)$ est
un anneau artinien (\`a gauche).
\\
Enfin, si $\sA$ est en outre
pseudo-ab\'elienne, les conditions ci-dessus \'equivalent \`a:
\item
\label{local}$\sA$ est somme directe locale d'une
famille de ca\-t\'e\-go\-ries artiniennes simples (ayant les m\^emes objets).
\item
  \label{ssablf} $\sA$ est ab\'elienne et tout objet de
$\sA$ est semi-simple et de longueur finie.
\item
\label{ssart} $\sA$ est ab\'elienne, artinienne, et tout
objet de $\sA$ est semi-simple.
\item
\label{sslf} $\sA$ est ab\'elienne, tout objet de
$\sA$ est de longueur finie, et $\rad(\sA)=0$.
\end{enumerate}
\end{thm}

\begin{sloppypar}
\prf Il
est clair que les
six premi\`eres conditions ne
changent pas si l'on
remplace
$\sA$
par une $K$-ca\-t\'e\-go\-rie \'equi\-va\-len\-te au sens de
Morita. De m\^eme les septi\`eme et huiti\`eme condition sont
v\'erifi\'ees pour une ca\-t\'e\-go\-rie
$K$-lin\'eaire $\sA$ si et seulement si elles le sont pour son enveloppe
pseudo-ab\'elienne $\sA^\natural$: en effet, pour tout couple d'objets
$((A,e),(A',e'))$ de $\sA^\natural$, on a
$\sA^\natural((A,e),(A',e'))=e'\sA(A,A')e$ et
$\rad(\sA^\natural)((A,e),(A',e'))=e'\rad(\sA)(A,A')e$, \cf aussi lemme
\ref{corresp} d). Ceci permet en particulier de supposer dans la suite
$\sA$ $K$-lin\'eaire
pseudo-ab\'elienne.
\end{sloppypar}

\ref{proj} $\iff$ \ref{inj} $\iff$ \ref{mod-ss} r\'esulte du
lemme \ref{A.L10}; \ref{mod-ss} $\If$ \ref{xss} est \'evident.

\ref{xss} $\If$ \ref{mod-ss}: soient $M$ un $\sA$-module \`a gauche,
$A\in \sA$ et $a\in M(A)$. Alors $Ann(a)(B)=\{f\in \sA_A(B)\mid f^*a=0\}$
d\'efinit un $A$-id\'eal \`a gauche de $\sA$ et l'application
$f\mapsto f^*a$ induit un mor\-phis\-me de $\sA$-modules $\sA_A/Ann(a)\to
M$. En faisant varier $A$ et $a$, on obtient un \'epimor\-phis\-me
$\bigoplus_{A\in \sA, a\in M(A)} \sA_A/Ann(a)\Surj M.$

Il r\'esulte alors du lemme \ref{A.L9} que $M$ est semi-simple.

\ref{xss} $\If$ \ref{ssablf}: le lemme \ref{A.L12} montre que si $S\in
\sA$ est simple, alors $\sA_S$ est simple; en r\'eappliquant ce m\^eme
lemme, on voit que tout objet de $\sA$ est semi-simple et de longueur
finie.

Enfin, pour tout objet simple $S\in\sA$,
$\sA(S,S)=\sA$-$Mod(\sA_S,\sA_S)$ est un corps,  donc $\sA$ est
ab\'elienne d'apr\`es le lemme \ref{A.L11}.

\ref{ssablf} $\If$ \ref{ssart}: si $S\in \sA$ est simple, $\sA(S,S)$
est un corps (lemme de Schur), donc ne contient aucun id\'eal $\ne 0$ et
$\sA_S$ est simple. Il en r\'esulte que $\sA_A$ est de longueur finie,
donc artinien, pour tout $A\in \sA$.

\ref{ssart} $\If$ \ref{ssablf}: le lemme de Yoneda implique que tout
objet de $\sA$ est artinien, donc de longueur finie.

\ref{ssablf} $\If$ \ref{mod-ss}: comme dans la d\'emonstration du lemme
\ref{A.L11}, on se ram\`ene au cas o\`u $\sA$ est isotypique. Soit $M$ un
$\sA$-module \`a gauche. Choisissons un objet simple $S$ de $\sA$, et soit
$D$ le  corps des endomor\-phis\-mes de $S$ (voir ci-dessus). Alors $D$
op\`ere naturellement \`a gauche sur $M(S)$. Pour tout objet $A$ de
$\sA$, on a un homomor\-phis\-me naturel
\[\sA(S,A)\otimes_D M(S)\to M(A)\]
o\`u le produit tensoriel sur $D$ est relatif \`a l'action \`a
droite de $D$ sur $\sA(S,A)$, et l'hypoth\`ese sur $\sA$ montre que
c'est un isomor\-phis\-me de groupes ab\'eliens. On d\'efinit ainsi une
\'equivalence
\[M\mapsto M(S)\]
de $\sA$-$Mod$ sur la ca\-t\'e\-go\-rie des $D$-espaces vectoriels \`a gauche,
dont tout objet est \'evidemment semi-simple.

\ref{struct} $\If$ \ref{ssablf} est \'evident. (\ref{mod-ss} et
\ref{ssart}) $\If$ \ref{struct}: la premi\`ere hypoth\`ese implique que
$\sA$ est semi-primitive au sens  de \cite[def. 2]{street} (pour toute
fl\`eche non nulle $f$ de $\sA$, il existe un module simple $S$ tel que
$S(f)\ne 0$).  La conclusion r\'esulte alors de \cite[th. 16]{street}.

\ref{struct} $\If$ \ref{hss}: pour tout $\alpha$, soit $\sI_\alpha$
l'id\'eal bilat\`ere de $\sA$ tel que $\sA/\sI_\alpha=\coprod_{\beta\ne
\alpha} D_\alpha$. Il est clair que chaque $\sI_\alpha$ est simple et
que $\sA=\bigoplus \sI_\alpha$.

\ref{hss}$\If$ \ref{local}: d\'ecomposons le bimodule $\sA$ en $\bigoplus
\sI_\alpha$, o\`u chaque id\'eal bilat\`ere $\sI_\alpha$ est simple.
Soit $\sA_\alpha=\sA/\bigoplus_{\beta\ne\alpha}\sI_\beta$. Alors
${\sA_\alpha}$ s'identifie canoniquement \`a $\sI_\alpha$, donc est
simple; en d'autres termes, $\sA_\alpha$ est simple, et elle est
\'evidemment artinienne. Le lemme \ref{A.L12} montre que $\sA$
s'identifie \`a la somme directe locale des $\sA_\alpha$.

\ref{local} $\If$ \ref{struct}: on se ram\`ene imm\'ediatement au cas
o\`u $\sA$ est simple. La conclusion r\'esulte alors de \cite[prop. 12 et
th. 16]{street}.

\ref{ssart}$\If$ \ref{end}: cela r\'esulte du lemme de Schur.

\ref{end} $\If$ \ref{ssart}: cela r\'esulte de \cite[lemma 2]{jannsen}
(o\`u l'on peut remplacer l'hypoth\`ese que $\sA(A,A)$ soit de dimension
finie sur $K$ par: $\sA(A,A)$ est artinienne.)

\ref{end} $\iff$ \ref{rad}: cela vient de ce que toute $K$-alg\`ebre
artinienne de radical nul est semi-simple.

\ref{mod-ss} $\If$ \ref{hss}: $\sA_A$ est semi-simple, donc de longueur
finie, donc artinien. Soit $I$ l'ensemble des types d'objets simples de
$\sA$-$Mod$ et, pour tout $\alpha\in I$, choisissons un module simple
$S_\alpha$ de type $\alpha$. Pour tout $\sA$-module \`a gauche
$M$, notons $Ann(M)=\{f\mid M(f)=0\}$: c'est un id\'eal bilat\`ere
de $\sA$. On va montrer l'\'egalit\'e $\sA=\bigoplus_{\alpha\in
I} Ann(\bigoplus_{\beta\ne\alpha}S_\beta)$ et que $Ann(\bigoplus_{\beta\ne
\alpha} S_\beta)$ est un id\'eal bilat\`ere simple pour tout $\alpha$.
Pour cela, \'etant donn\'e $A\in \sA$, d\'ecomposons $\sA^o_A$ en une
somme directe $\bigoplus_{\alpha\in I}M_\alpha$ de sous-modules
isotypiques. Ecrivons d'abord $Id_A=\sum e_\alpha$, avec $e_\alpha\in
M_\alpha(A)$ pour tout $\alpha$. On remarque que les $M_\alpha(A)$
sont des id\'eaux bilat\`eres de l'anneau $\sA(A,A)$; il en r\'esulte
que les $e_\alpha$ forment un syst\`eme d'idempotents orthogonaux. En
particulier, comme $S_\beta(A)$ est facteur direct de $M_\beta(A)$,
$e_\alpha$ op\`ere par $0$ sur $S_\beta(A)$ pour $\beta\ne \alpha$ et par
l'identit\'e pour $\beta=\alpha$. Soit maintenant $f\in \sA(A,B)$, et soit
$f=\sum f_\alpha$ sa d\'ecomposition sur les $M_\alpha(B)$. On a
$f_\alpha=fe_\alpha$, donc $f_\alpha$ op\`ere sur $M_\beta(B)$ par
$0$ pour $\beta\ne \alpha$ et par $f$ pour $\beta=\alpha$. On a
donc bien prouv\'e que $\sA(A,B)=\bigoplus_{\alpha\in I}
Ann(\bigoplus_{\beta\ne\alpha}S_\beta)(A,B)$. Enfin, montrons que,
pour tout $\alpha\in I$, $\sI_\alpha=Ann(\bigoplus_{\beta\ne
\alpha}S_\beta)$ est simple. Soit $\sJ$ un sous-module non nul
de $\sI_\alpha$. Choisissons $A,B\in \sA$ tels que $\sJ(A,B)\ne 0$
et soit $0\ne f\in \sJ(A,B)$. On a $S_\beta(f)=0$ pour tout
$\beta\ne\alpha$. Si $S_\alpha(f)=0$, on a donc $M(f)=0$ pour
tout $\sA$-module $M$, ce qui est absurde en choisissant
$M=\sA_A$. Par cons\'equent, $S_\alpha(f)\ne 0$.

\ref{ssablf} $\If$ \ref{sslf}: Pour tout objet semi-simple $A$,
l'alg\`ebre $\sA(A,A)$ est semi-simple donc $\rad(\sA)(A,A)=0$; on
conclut gr\^ace \`a l'additivit\'e des id\'eaux sur les objets.

\ref{sslf} $\If$ \ref{ssablf}:  Comme $\sA$ est suppos\'ee ab\'elienne et 
que
tout objet est suppos\'e de longueur finie, tout objet est somme directe 
finie
d'ind\'e\-com\-po\-sa\-bles. Tout se ram\`ene \`a montrer que tout
ind\'e\-com\-po\-sa\-ble $A$ est simple.  Cela r\'esulte du lemme
\ref{L2} appliqu\'e au monomor\-phis\-me $S\inj A$, o\`u
$S$ est un sous-objet simple de $\sA$ (dans ce cas, le lemme dit que $S$
est facteur direct de $A$, donc \'egal \`a $A$).
  \qed

\begin{cor}[\`a la d\'emonstration] Une $K$-ca\-t\'e\-go\-rie $\sA$ est simple
(d\'e\-fi\-ni\-tion \ref{A.D6}) si et seulement si elle est semi-simple et
n'a qu'un seul type de module simple.
\end{cor}

\prf Cela r\'esulte de la preuve de \ref{struct} $\If$ \ref{hss}.\qed

\begin{contrex}\label{A.CeX4} La ca\-t\'e\-go\-rie du contre-exemple
\ref{A.CeX1} est semi-simple au sens de \cite[d\'ef. 4]{street}, mais pas
au sens de la d\'efinition \ref{D3} ci-dessus.
\end{contrex}

\begin{lemme}\label{lA.3} Dans une ca\-t\'e\-go\-rie ab\'elienne semi-simple
$\sA$, tout mor\-phis\-me est somme directe d'un mor\-phis\-me nul et d'un
isomor\-phis\-me.
\end{lemme}

\prf Soit $u:A\to B$ un mor\-phis\-me de $\sA$. On peut d\'ecomposer
\begin{align*}
A&=\Ker u\oplus A'\\
B&=B'\oplus \IM u.
\end{align*}

Sur cette d\'ecomposition, $u$ a la forme d\'esir\'ee.\qed

\begin{lemme}\label{utile}
Tout foncteur additif $F: \sA\to \sB$ entre ca\-t\'e\-go\-ries
a\-b\'e\-lien\-nes semi-simples est exact. De plus,  les
conditions suivantes:
\begin{thlist}
\item $F$ est fid\`ele
\item $F$ est conservatif
\item  $F$ n'envoie aucun
objet  non nul de $\sA$ sur un objet nul de $\sB$
\end{thlist}
sont \'equi\-va\-len\-tes.
\end{lemme}

\prf Comme tout foncteur additif transforme une suite exacte courte
scind\'ee en une suite exacte courte scind\'ee, la premi\`ere affirmation
est claire. (i) $\If$ (iii) est \'evident; (iii) $\If$ (ii) et (ii)
$\If$ (i) r\'esultent imm\'ediatement du lemme \ref{lA.3}.

\subsection{Ca\-t\'e\-go\-ries s\'eparables (ou absolument semi-simples)}

\begin{thm}\label{tsep} Soit $\sA$ une $K$-ca\-t\'e\-go\-rie. Les conditions
suivantes sont \'equi\-va\-len\-tes:
\begin{enumerate}
\item \label{tsep1} Le $\sA^{\rm o}\boxtimes_K\sA$-module $\sA$ est
projectif.
\item \label{tsep1bis} Le $\sA^{\rm o}_L\boxtimes_L\sA_L$-module $\sA_L$
est projectif pour toute extension $L/K$.
\item \label{tsep2} $\dim_K \sA(A,B)<\infty$ pour tous $A,B\in \sA$ et
$\sA_L$ est semi-simple pour toute extension $L/K$.
\item \label{tsep3} $\sA_L$ est semi-simple pour toute extension $L/K$.
\item \label{tsep4} $\dim_K \sA(A,B)<\infty$ pour tous $A,B\in \sA$ et
$\rad(\sA_L)=0$ pour toute extension $L/K$.
\item \label{tsep5} La ca\-t\'e\-go\-rie $\sA^{\rm o}\boxtimes_K \sA$ est
semi-simple.\\
On dit alors que $\sA$ est \emph{s\'eparable}.\\
Si $\sA$ est $K$-lin\'eaire, ces conditions sont encore
\'equi\-va\-len\-tes \`a
\item \label{tsep6} La $K$-alg\`ebre $\sA(A,A)$ est s\'eparable pour tout
objet
$A$ de
$\sA$.
\end{enumerate}
  \end{thm}

\prf Les conditions \ref{tsep1}--\ref{tsep5} \'etant manifestement
invariantes par passage \`a l'enveloppe $K$-lin\'eaire, on peut supposer
$\sA$ $K$-lin\'eaire.

\ref{tsep2} $\If$ \ref{tsep3} et \ref{tsep2} $\If$ \ref{tsep4} sont
\'evidents. \ref{tsep3} $\iff$ \ref{tsep6} r\'esulte de la
caract\'erisation \ref{end} du th\'eor\`eme \ref{A.P4} et de la
d\'efinition d'une $K$-alg\`ebre s\'eparable (d\'efinition \ref{D4sep}
a)). \ref{tsep6} $\If$ \ref{tsep2} r\'esulte de la caract\'erisation
\ref{Psep3} des $K$-alg\`ebres s\'eparables dans la proposition
\ref{Psep}. \ref{tsep6} $\If$ \ref{tsep5} r\'esulte du lemme \ref{l2.3}.
\ref{tsep4} $\If$ \ref{tsep6} r\'esulte de la caract\'erisation
\ref{Psep5} des $K$-alg\`ebres s\'eparables dans la proposition
\ref{Psep}. \ref{tsep5} $\If$ \ref{tsep1} r\'esulte de la
caract\'erisation \ref{proj} du th\'eor\`eme \ref{A.P4}. \ref{tsep1}
$\If$ \ref{tsep1bis} r\'esulte du corollaire \ref{cA.1}.

Il reste \`a d\'emontrer que \ref{tsep1bis} $\If$ \ref{tsep3}. Il
suffit de traiter le cas $L=K$. Pour cela, nous allons reproduire la
raisonnement de \cite[ch. IX, d\'em. de la prop. 7.1]{ce} ``\`a la main".

Si $M$ et $N$ sont deux $\sA$-modules, on leur associe le $\sA^e$-module
\[\uHom_K(M,N):(A,B)\mapsto Hom_K(M(A),N(B)).\]

Le foncteur $\uHom:(\sA\hbox{--}Mod)^{\rm o}\times \sA\hbox{--}Mod\to
\sA^e\hbox{--}Mod$ est bi-exact.

Soit $0\to N\to E\to M\to 0$ une extension de $M$ par $N$. On lui
associe une extension $\tilde E$ du $\sA^e$-module $\sA$ par
$\uHom_K(M,N)$ par pull-back via le diagramme
\[\begin{CD}
0&\to& \uHom_K(M,N)@>>> \uHom_K(E,N)@>>> \uHom_K(N,N)&\to& 0\\
&&||&&@A{\rho}AA @AAA\\
0&\to& \uHom_K(M,N)@>>> \tilde E@>\pi>> \sA&\to& 0.
\end{CD}\]

Soit $s$ une section de $\pi$. Alors $\rho\circ s$ d\'efinit un
homomor\-phis\-me $\sA\to \uHom_K(E,N)$. Celui-ci correspond \`a un
homomor\-phis\-me $f:E\to N$, et la commutativit\'e du diagramme montre que
$f$ est une r\'etraction du monomor\-phis\-me $N\to E$. Ainsi, $\sA$ est
semi-simple.\qed

\section{Erratum \`a \protect\cite{ak(note)}}

O. Gabber nous a fait remarquer que la preuve de \cite[Prop. 6]{ak(note)}
\'etait (encore) incompl\`ete. Le probl\`eme est que la d\'emonstration 
donn\'ee
dans \loccit s'applique \`a un motif de la forme $h(X)$, o\`u $X$ est une
vari\'et\'e projective lisse, mais pas \emph{a priori} \`a un motif
g\'en\'eral $M$. Plus pr\'ecis\'ement, si $M$ est effectif (pour fixer
les id\'ees) et que le projecteur de K\"unneth pair de $M$ est
alg\'ebrique, il n'est pas clair qu'on puisse choisir une vari\'et\'e $X$
telle que $h(X)$ contienne $M$ en facteur direct et telle que le
projecteur de K\"unneth pair de $X$ soit alg\'ebrique. 

Voici deux mani\`eres de corriger ce point (la premi\`ere est en substance celle 
que nous a propos\'ee Gabber). Les notations sont celles de \cite{ak(note)}.

\subsection{Premi\`ere correction}Comme remarqu\'e au d\'ebut de la preuve 
de \loccit, il suffit de supposer que le corps de base $k$ est de type fini sur 
$\F_p$.

\begin{lemme} \label{lB.1} Soient $H:Mot_\rat\to Modf$-$L$
une cohomologie de Weil classique \`a coefficients dans une 
$\Q$-alg\`ebre commutative semi-simple $L$ et $Mot_H$ la ca\-t\'e\-go\-rie des
$k$-motifs $H$-homologiques \`a coefficients rationnels. Pour tout $A\in 
Mot_H$, 
notons \[H^+(A)=\bigoplus_{i\text{ pair}} H^i(A),\;\;H^-(A)=\bigoplus_{i\text{ impair}} H^i(A).\] Alors, pour
$A,B\in Mot_H$, on a
\begin{multline*}
\sR_H(A,B)=\\
\{f\in Mot_H(A,B)\mid \forall g\in
Mot_H(B,A),
\forall i\in\Z, tr H^i(g\circ f)=0\}\\
=\{f\in Mot_H(A,B)\mid \forall g\in
Mot_H(B,A),
  tr H^+(g\circ f)=tr H^-(g\circ f)=0\},
\end{multline*}
o\`u $\sR_H$ est le radical de $Mot_H$.
\end{lemme}

\prf D'apr\`es \cite[prop. 5]{ak(note)}, $Mot_H$ est semi-primaire, donc 
son radical de Kelly co\"{\i}ncide avec son radical de Gabriel. 
L'\'enonc\'e r\'esulte alors du lemme \ref{rad'}.\qed

\begin{sloppypar}
\begin{lemme}\label{lB.2} Avec les notations du lemme \ref{lB.1}, l'image 
r\'eciproque $\snum$ de $\sR_H$ dans la ca\-t\'e\-go\-rie $Mot_\rat$ des motifs de 
Chow ne d\'epend pas du choix de $H$.
\end{lemme}

\medskip
\prf \'Etant donn\'e le lemme \ref{lB.1}, il suffit de voir que, pour toute  
$k$-vari\'et\'e  projective  lisse  $X $ et toute correspondance $c\in 
CH^{\dim 
X}(X\times X)$ et tout $i\in\Z$, $tr H^i(c)$ est un nombre 
rationnel qui ne 
d\'epend pas de $H$. Par changement de base propre et lisse et par 
sp\'ecialisation des cycles alg\'ebriques, on se r\'eduit par r\'ecurrence 
sur 
le degr\'e de transcendance de $k$ sur $\F_p$ au cas o\`u $k$ est fini (\cf 
\cite[d\'em. de la prop. 5]{ak(note)}). L'assertion r\'esulte alors de 
\cite{km}.\qed
\end{sloppypar}

\begin{lemme}[\cf \protect{\cite[lemme 7]{ak(note)}}]\label{lB.3} Soient, 
pour 
fixer les id\'ees, $\ell\ne \ell'\ne p$ deux nombres premiers. Alors le 
noyau du 
foncteur $Mot_{\ell\ell'}\to 
Mot_\ell$ est un id\'eal localement nilpotent.
\end{lemme}

\prf Cela r\'esulte des lemmes \ref{lB.1} et \ref{lB.2}.\qed

\begin{lemme}\label{lB.4} Les foncteurs pleins  $Mot_\ell \to Mot_{\snum}, 
Mot_{\ell'} \to Mot_{\snum}$ induisent des bijections sur les classes 
d'isomor\-phis\-me d'objets. En particulier, pour tout objet $M$ de 
$Mot_{\ell}$,  
$b_i(M)= \dim H_{\ell}^i(M)$ ne d\'epend pas de $\ell$, et 
ne d\'epend que de la 
classe d'isomor\-phis\-me de l'image de $M$ dans $Mot_{\snum}$.
\end{lemme}

\prf Par le lemme \ref{lB.3}, ces foncteurs sont essentiellement surjectifs 
et 
conservatifs.\qed

\begin{lemme}\label{lB.5} L'image $M^\ast _{\snum}$ de $M^\ast _{\ell}$ dans 
$Mot_{\snum}$ consiste en la sous-ca\-t\'e\-go\-rie pleine form\'e des objets 
isomorphes \`a des sommes finies d'objets n'ayant qu'un seul nombre de Betti 
non 
nul. En particulier, $M^\ast _{\snum}$ est ind\'ependant de $\ell$. De 
m\^eme 
avec $M^\pm$.
\end{lemme}

\prf C'est clair par le lemme \ref{lB.4}, car on a la m\^eme propriet\'e au 
niveau de $M^\ast _{\ell}$.\qed

\begin{lemme}\label{lB.6} $M^\ast _{snum}=  M^\ast _{num}$ et $M^\pm 
_{snum}=  
M^\pm _{num}$.
\end{lemme}

\prf C'est clair.\qed

\subsection{Seconde correction} Celle-ci ne marche que pour les $M^\pm$
mais  s'applique \`a toute cohomologie de Weil, classique ou non.

D'apr\`es le th\'eor\`eme \ref{tkim.1} c), pour toute cohomologie de Weil 
$H$ on a $Mot_H^\pm=(Mot_H)_\kim$. D'autre part, la proposition
\ref{nouvellekim}  implique que le foncteur $(Mot_H)_\kim\allowbreak\to
(Mot_\num)_\kim$ est  essentiellement surjectif. Ainsi, $Mot_\num^\pm =
(Mot_\num)_\kim$ pour toute cohomologie de  Weil $H$.\qed

\subsection*{Table de concordance pour les r\'ef\'erences de 
\cite{ak(note)}}

\begin{sloppypar}
\begin{itemize}
\item Preuve du th\'eor\`eme 1: remplacer [3, 6.7.3] par [3, 
\ref{absJannsen}].
\item Preuve de la proposition 2: remplacer [3, 6.7.9 et 1.4.2b] par [3, 
\ref{new p4} et \ref{P1}b].
\item D\'ebut de {\bf 2}: remplacer [3, 9.2.1, 9.7.3, 11.3.5] par [3, 
\ref{T2}, 
\ref{P3'}, \ref{t4}].
\end{itemize}
\end{sloppypar}

\section{Finite dimensionality of reductive envelopes, by Peter
O'Sullivan}\label{os}

\renewcommand{\proofname}{Proof}

The object of this appendix is to determine when
proreductive envelopes are finite dimensional and  when
their formation is compatible with extension of scalars.
Let $k$ be a field of characteristic $0$,
and let $H$ be an affine $k$-group. Then
the proreductive envelope of $H$ is finite dimensional
if and only if $H$ is finite dimensional with prounipotent
radical of dimension $\le 1$ (Theorem~\ref{th:finite}).
This is a consequence of Theorem~\ref{th:fam}
characterising the affine $k$-groups with prounipotent radical
of dimension $\le 1$ as those whose representations are ``rigid''.
If the prounipotent radical $U$ of $H$ has dimension $1$
the proreductive
envelope of $H$ is the semidirect product of $H/U$ by $SL(2)$
(Theorem~\ref{th:SL2}).
Let $k'$ be an extension of $k$.
Then
${}^\mathrm{p}\Red(H_{k'}) \rightarrow {}^\mathrm{p}\Red(H)_{k'}$
is an isomorphism
if and only if either $k'$ is algebraic over $k$ or the
prounipotent radical of $H$ is of dimension $\le 1$
(Theorem~\ref{th:basech}).
This follows from Theorem~\ref{th:sumnd}
which gives the conditions under which
every representation of $H$ over $k'$ is a direct summand of one
defined over $k$.

Throughout this appendix $k$ denotes a field
of characteristic $0$. Unless otherwise stated,
representations are assumed to be finite dimensional.

\subsection{Families of representations}\label{sec:fam}

The object of this section is to prove
Theorems~\ref{th:fam}~and~\ref{th:sumnd}. The proofs are
based on Lemmas~\ref{le:Afam}~and~\ref{le:GL2}.
Lemma~\ref{le:GL2} is necessary also for the proof of
Theorem~\ref{th:SL2}.

\begin{lem}\label{le:te}
Suppose that $k$ is algebraically closed.
Let $M$ be a connected reductive $k$-group and let
$V$ be a representation of $M$.
Denote by $P$ the set of dominant weights of $M$ relative to some 
maximal torus $T$ and Borel subgroup
$B \supset T$ of $M$, and for each $\mu\in P$
let $V_\mu$ be a representation of $M$ with
highest weight $\mu$.
Then if $\tau \in P$ is
such that $\tau+\lambda \in P$
for each $\lambda$ in the set $\Lambda$ of weights of $V$,
there is an $M$-isomorphism
$V \otimes_k V_\tau \simeq
\bigoplus_{\lambda \in \Lambda} V_{\tau+\lambda}^{m(\lambda)}$,
where $m(\lambda)$ is the multiplicity of $\lambda$ in $V$.
\end{lem}

\begin{proof}
Denote by $\chi$ and $\chi_\mu$ for $\mu\in P$
the respective restrictions of the characters of $V$ and
$V_\mu$ to $T$. It suffices to show that
$\chi\chi_\tau =
\sum_{\lambda\in\Lambda}m(\lambda)\chi_{\tau+\lambda}$.

Let $W$ be the Weyl group of $M$ relative to $T$ and
let $\delta$ be half the sum of the positive roots of $M$
relative to $T$ and $B$. For any $\pi$ in the character
group of $T$, write $e(\pi)$ for $\pi$ regarded as an
element of $k[T]$.
Then we have
\begin{multline*}
\sum_{\lambda\in\Lambda}m(\lambda)e(\lambda)
\sum_{w\in W}(\det w)e(w(\tau+\delta))
=\\
\sum_{w\in W}(\det w)e(w(\tau+\delta))
\sum_{\lambda\in\Lambda}m(w\lambda)e(w\lambda) \\
=\sum_{\lambda\in\Lambda}m(\lambda)\sum_{w\in W}
(\det w)e(w(\tau+\lambda+\delta))
\end{multline*}
since $\Lambda$ is stable under $W$ and $m(w\lambda)=m(\lambda)$
for $w\in W$ and $\lambda\in\Lambda$. Thus if
$Q=\sum_{w\in W}(\det w)e(w\delta)$, we have by Weyl's character
formula
\[
Q\chi\chi_\tau=
Q\sum_{\lambda\in\Lambda}m(\lambda)\chi_{\tau+\lambda},
\]
whence the required result, since $Q\ne0$.
\end{proof}

\begin{lem}\label{le:A}
Suppose that $k$ is algebraically closed.
Let $M$ be a non-commutative connected reductive $k$-group
and let $V$ be a faithful irreducible
representation of $M$.
Suppose that $V$ has, relative to some
maximal torus of $M$, at most $3$ distinct non-zero weights.
Then either the derived group $M'$ of $M$ is  simple of type $A_1$
and $V$ is of dimension $2$ or $3$, or $M'$ is simple of type
$A_2$ and $V$ is of dimension $3$.
\end{lem}

\begin{proof}
We may suppose that $M=M'$ is semisimple. Let $\tilde M$ be the
universal cover of $M$, so that $V$ may also be regarded as a
representation of $\tilde M$. Let $T$ be a maximal torus and
$B \supset T$ a Borel subgroup of $M$, let $X$ be the character
group of $T$ and $W$ the Weyl group of $M$ relative to $T$,
and let $\lambda \in X$ be the highest weight of $V$ relative
to $T$ and $B$. Since $V$ has at most $3$ distinct non-zero
weights, the orbit $\Lambda$ of $\lambda$ under $W$ has at
most $3$ elements. There is a decomposition
$W = W_1 \times \cdots \times W_r$ corresponding to the simple
factors of $\tilde M$ such that the representation of $W$ on
$X \otimes \R$ is the direct sum of non-trivial irreducible
representations $X_i$ of $W_i$ for $i=1,\ldots,r$. By
faithfulness of the $M$-representation $V$, the projection
of $\lambda$ onto each $X_i$ is non-zero. Thus $r=1$ and
$\tilde M$ is simple, since otherwise $\Lambda$ would have
at least $4$ elements. The affine subspace of $X \otimes \R$
generated by $\Lambda$ is of dimension $\le2$ and is stable
under $W$, so by irreducibility the rank
$\dim_\R(X \otimes \R)$ of $\tilde M$ is $\le 2$. Since $W$
contains a rotation of order $4$ when $\tilde M$ is of type
$B_2$ and of order $6$ when $\tilde M$ is of type $G_2$,
it follows that $\tilde M$ must be of type $A_1$ or $A_2$.

If $\tilde M$ is of type $A_1$ any $\tilde M$-representation
of dimension $\ge 4$ has at least $4$ distinct non-zero weights.
If $\tilde M$ is of type $A_2$, let $\tau_1$ and $\tau_2$ be
the highest weights of the two fundamental representations of
$\tilde M$ relative to $T$ and $B$, so that the simple positive
roots of $\tilde M$ are $\mu_1 = 2\tau_1-\tau_2$ and
$\mu_2 = -\tau_1 + 2\tau_2$,
and $\lambda= m_1\tau_1 + m_2\tau_2$ for integers $m_i \ge 0$.
One of $m_1,m_2$ must be $0$, since otherwise $\Lambda$ would
contain $6$ elements. Thus either $\lambda=\tau_1$ or
$\lambda=\tau_2$, since if for example $\lambda=m_1\tau_1$
with $m_1\ge2$ then $\lambda-\mu_1$ would be a non-zero weight
of $V$ not conjugate to $\lambda$ under $W$, so $V$ would have
at least $6$ distinct non-zero weights.
\end{proof}

\begin{lem}\label{le:Vrs}
Suppose that $k$ is
algebraically closed. Let $M$ be a connected reductive $k$-group
and let $V$ be a representation of $M$ which
is either the direct sum of two non-isomorphic $1$-dimensional
representations or is irreducible of dimension $2$ or $3$.
Then there exists a family
$(V_j^i,r_j^i,s_j^i)_{i,j\in\N}$, with the $V_j^i$ irreducible
representations of $M$ and the
$r_j^i:V \otimes_k V_j^i \rightarrow V_j^{i+1}$ and
$s_j^i:V \otimes_k V_{j+1}^i \rightarrow V_j^i$ non-zero
$M$-homomorphisms, such that:

\textup{(i)} for each integer $m$ the $V_j^i$ with $i-j=m$ are pairwise
non-isomorphic;

\textup{(ii)}
the restrictions of
\begin{gather*}
r^{i+1}_j \circ (V \otimes r^i_j) :
V \otimes_k V \otimes_k V_j^i \rightarrow V_j^{i+2},  \\
s^i_{j-2} \circ (V \otimes s^i_{j-1}) :
V \otimes_k V \otimes_k V_j^i \rightarrow V_{j-2}^i,  \\
s_{j-1}^{i+1}\circ(V \otimes r_j^i)+
r_{j-1}^i\circ(V \otimes s_{j-1}^i) :
V \otimes_k V \otimes_k V_j^i \rightarrow V_{j-1}^{i+1},
\end{gather*}
respectively to
$\bigwedge^2 V \otimes_k V_j^i$,
$\bigwedge^2 V \otimes_k V_j^i$, $\bigwedge^2 V \otimes_k V_j^i$,
are $0$ for $i,j \in \N$.
\end{lem}

\begin{proof}
We may suppose that $V$ is a faithful representation of $M$.
If $V\simeq V_1\oplus V_2$ with $V_1$ and $V_2$ $1$-dimensional
and $p_1:V \rightarrow V_1$ and $p_2:V \rightarrow V_2$ are the
projections, take
$V_j^i=V_1^{\otimes i} \otimes_k (V_2^\vee)^{\otimes j}$ and
(modulo canonical isomorphisms) $r_j^i=p_1 \otimes V_j^i$ and
$s_j^i=p_2 \otimes V_{j+1}^i$. It is immediately checked that
(i) and (ii) hold.

Now suppose that $V$ is irreducible of dimension $2$ or $3$.
Let $T$ and $B\supset T$ be a maximal torus and Borel subgroup
of $M$, and write $\lambda$ and $\lambda^\vee$ for the
respective highest weights of $V$ and $V^\vee$ relative to $T$
and $B$. Let $V_j^i$ be an irreducible representation of $M$
with highest weight $i\lambda+j\lambda^\vee$. Then (i)
clearly holds because $\lambda+\lambda^\vee\ne0$.
All weights of $\bigwedge^2 V \otimes_k V_j^i$ are of the form
$(i+2)\lambda+j\lambda^\vee-\pi$ where $\pi\ne0$ is a sum of
positive roots of $M$, so
$\Hom_M(\bigwedge^2 V \otimes_k V_j^i,V_j^{i+2})=0$. Similarly
\[
\Hom_M(\bigwedge^2 V \otimes_k V_{j+2}^i,V_j^i)=
\Hom_M(V_{j+2}^i,\bigwedge^2 V^\vee \otimes_k V_j^i)=0.
\]
Thus it suffices to construct
$r^i_j \ne 0$ and $s^i_j \ne 0$ for $i,j \in \N$ such that
the restriction of the third homomorphism of (ii)
to $\bigwedge^2 V \otimes_k V_{j+1}^i$ is $0$.

Let $t_j^i\ne0$ be an
element of the highest weight space of $V_j^i$ and let $u_+\ne0$
be an element of the highest weight space of $V$ and $u_-\ne0$
an element of the lowest weight space of $V$. The weights of
$t_j^i$, $u_+$ and $u_-$ are thus respectively
$i\lambda+j\lambda^\vee$, $\lambda$ and $-\lambda^\vee$. The
$k$-linear map $V \rightarrow k$ which sends $u_-$ to $1$ and
which is $0$ on the other weight spaces of $V$ is an element
$u_+^\vee\ne0$ of the highest weight space of $V^\vee$. Now
$V \otimes_k V_j^i$ contains the weight
$(i+1)\lambda+j\lambda^\vee$ with multiplicity $1$, and all
of its other weights are of the form
$(i+1)\lambda+j\lambda^\vee-\pi$ where $\pi\ne0$ is a sum of
positive roots. Thus there is a unique $M$-homomorphism
$r_j^i:V \otimes_k V_j^i \rightarrow V_j^{i+1}$ such that
\begin{equation}\label{eq:r}
r_j^i(u_+\otimes t_j^i)=t_j^{i+1}.
\end{equation}
Similarly there is a unique $M$-homomorphism
$q_j^i:V_{j+1}^i \rightarrow V^\vee \otimes_k V_j^i$
such that $q_j^i(t_{j+1}^i)=u_+^\vee \otimes t_j^i$,
whence for any constant $C_n=C_n(M,V)$ there is a unique $M$-homomorphism
$s_j^i:V \otimes_k V_{j+1}^i \rightarrow V_j^i$
such that
\begin{equation}\label{eq:s}
s_j^i(u_-\otimes t_{j+1}^i)=C_{i+j}t_j^i.
\end{equation}
Choose $C_n$ as follows.
By Lemma~(\ref{le:A}), either the derived group $M'$ of $M$
is of type $A_1$ and $V$ is of dimension $2$ or $3$, or $M'$
is of type $A_2$ and $V$ is of dimension $3$. Set
\begin{equation}\label{eq:C}
C_n=\begin{cases}
    1/(n+2)&  \text{if $M'$ is type $A_1$ and $V$ is of
                             dimension $2$},\\
    (n+1)/(2n+3)&  \text{if $M'$ is of type $A_1$ and $V$
                              is of dimension $3$},\\
    1/(n+3)&    \text{if $M'$ is of type $A_2$}.
\end{cases}
\end{equation}

\begin{sloppypar}
The $V_j^i,r_j^i,s_j^i$ just constructed are compatible
with restriction from $M$ to $M'$, so to check
that the third homomorphism of (ii) has restriction $0$ to
$\bigwedge^2 V \otimes_k V_{j+1}^i$ we may suppose that $M'=M$.
We show first that
$\Hom_M(\bigwedge^2 V \otimes_k V_{j+1}^i,V_j^{i+1})$ is
$1$-dimensional. This is immediate when $M$ is of type $A_1$,
since then $V_j^i\simeq V_0^{i+j}$ and $\bigwedge^2 V$ is either
trivial $1$-dimensional or $M$-isomorphic to $V$ according
as $V$ is $2$-dimensional or $3$-dimensional. When $M$ is of
type $A_2$, $V$ and $V^\vee$ are the two fundamental
representations of $M$, so the dominant weights of $M$ relative
to $T$ and $B$ are the $i\lambda+j\lambda^\vee$ for $i,j\ge0$.
The weights of $V$ are $\lambda$, $-\lambda+\lambda^\vee$,
$-\lambda^\vee$, and $\bigwedge^2 V \simeq V^\vee$. Applying
Lemma~\ref{le:te} with $V$ replaced by $\bigwedge^2 V$ and
with $\tau=i\lambda+(j+1)\lambda^\vee$ shows that
$\bigwedge^2 V \otimes_k V_{j+1}^i$ contains $V_j^{i+1}$ with
multiplicity $1$ when $i\ge1,j\ge0$. Similarly, when $j\ge1$,
$\bigwedge^2 V^\vee \otimes_k V_j^1$ contains $V_{j+1}^0$ with
multiplicity $1$, whence $\bigwedge^2 V \otimes_k V_{j+1}^0$
contains $V_j^1$ with multiplicity $1$. That
$\bigwedge^2 V \otimes_k V_1^0 \simeq V^\vee \otimes_k V^\vee$
contains $V_0^1\simeq V$ with multiplicity $1$ is clear.
Thus
$\Hom_M(\bigwedge^2 V \otimes_k V_{j+1}^i,V_j^{i+1})$
is in all cases
$1$-dimensional. The restriction to
$\bigwedge^2 V \otimes_k V_{j+1}^i$ of the third homomorphism
of (iii) will thus be $0$ provided that
that the images of
$(u_+\otimes u_- -u_-\otimes u_+)\otimes t^i_{j+1}$ under
$r^i_j\circ(V \otimes s^i_j)$ and
$-s^{i+1}_j\circ(V \otimes r^i_{j+1})$
are non-zero and coincide. Now
\begin{equation}\label{eq:s0}
s^i_j(u_+ \otimes t^i_{j+1})=0
\end{equation}
because the difference between the weight
$(i+1)\lambda+(j+1)\lambda^\vee$ of $u_+ \otimes t^i_{j+1}$
and the highest weight $i\lambda+j\lambda^\vee$ of $V^i_j$
is non-zero and dominant. Thus by (\ref{eq:r}), (\ref{eq:s})
and (\ref{eq:C}) it is enough to show that
\begin{equation}\label{eq:sCt}
s^{i+1}_j(u_+ \otimes r^i_{j+1}(u_- \otimes t^i_{j+1}))
=(C_{i+j+1}-C_{i+j})t^{i+1}_j.
\end{equation}
\end{sloppypar}

Write $\mathfrak m$ for the Lie algebra of $M$.
Suppose first that $M$ is of type $A_1$ and $V$ is of dimension
$2$. Let $e$ be an element of the positive root subspace of
$\mathfrak m$, chosen so that $eu_-=u_+$, and let $f,h$ be such
that $f$ is in the negative root subspace of $\mathfrak m$ and
$[e,f]=h$, $[h,e]=2e$, and $[h,f]=-2f$. Then $hu_+=u_+$,
$hu_-=-u_-$, $fu_+=u_-$, and $ht^i_j=(i+j)t^i_j$. We have
\begin{align*}
r^i_{j+1}(u_- \otimes t^i_{j+1})&
=\frac{1}{i+j+2}ft^{i+1}_{j+1}, \\
u_+ \otimes ft^{i+1}_{j+1}&
=f(u_+ \otimes t^{i+1}_{j+1}) - u_- \otimes t^{i+1}_{j+1},
\end{align*}
where the constant $1/(i+j+2)$ in the first equality
is determined by acting with $e$ and using (\ref{eq:r}).
Using (\ref{eq:s0}), (\ref{eq:s}),  and
(\ref{eq:C}) we obtain (\ref{eq:sCt}).

When $M$ is of type $A_1$ and $V$ is of dimension
$3$, choose $e$ and $f$, $h$ such that $e^2u_-=u_+$
and $[e,f]=h$, $[h,e]=2e$, and $[h,f]=-2f$. Then
$hu_+=2u_+$, $hu_-=-2u_-$, $fu_+=2eu_-$, $f^2u_+=4u_-$,
and $ht^i_j=2(i+j)t^i_j$. We have
\begin{align*}
r^i_{j+1}(u_- \otimes t^i_{j+1})&
=\frac{1}{4(i+j+2)(2i+2j+3)}f^2t^{i+1}_{j+1}, \\
u_+ \otimes f^2t^{i+1}_{j+1}&
=f(u_+ \otimes ft^{i+1}_{j+1} - fu_+ \otimes t^{i+1}_{j+1})
+ 4u_- \otimes t^{i+1}_{j+1},
\end{align*}
where the constant in the first equality is determined
by acting with $e^2$. The images under $s^{i+1}_j$ of
$u_+ \otimes ft^{i+1}_{j+1}$ and
$fu_+ \otimes t^{i+1}_{j+1}$  are $0$
by weights, so using (\ref{eq:s}) and
(\ref{eq:C}) we obtain (\ref{eq:sCt}).

Finally suppose that $M$ is of type $A_2$. Then the simple
positive roots of $M$ are $\mu=2\lambda-\lambda^\vee$ and
$\mu^\vee=-\lambda+2\lambda^\vee$. Let $e$ and $e^\vee$ be
elements of the root subspaces of $\mathfrak m$ corresponding
respectively to $\mu$ and $\mu^\vee$. We may suppose $e$ and
$e^\vee$ chosen so that $ee^\vee u_-=u_+$. Let $f$, $f^\vee$,
$h$, $h^\vee$ be the unique elements of $\mathfrak m$ such that
$f$ and $f^\vee$ lie in the root subspaces corresponding
respectively to $-\mu$ and $-\mu^\vee$ and $[e,f]=h$, $[h,e]=2e$,
$[h,f]=-2f$, $[e^\vee,f^\vee]=h^\vee$, $[h^\vee,e^\vee]=2e^\vee$,
$[h^\vee,f^\vee]=-2f^\vee$. Then $[h,h^\vee]=0$,
$[h,e^\vee]=-e^\vee$, $[h^\vee,e]=-e$, $[h,f^\vee]=f^\vee$,
$[h^\vee,f]=f$, and $[e,f^\vee]=[e^\vee,f]=0$. Also
$eu_-=f^\vee u_+=0$, $e^\vee u_-=fu_+$, $hu_+=u_+$, $h^\vee u_+=0$,
$hu_-=0$, $h^\vee u_-=-u_-$, $f^\vee fu_+=u_-$, and
$ht^i_j=it^i_j$, $h^\vee t^i_j=jt^i_j$, $et^i_j=e^\vee t^i_j=0$.
We have
\begin{align*}
r^i_{j+1}(u_- \otimes t^i_{j+1})&
=\frac{1}{(i+1)(i+j+3)}
((i+2)f^\vee f - (i+1)ff^\vee)t^{i+1}_{j+1}, \\
u_+ \otimes f^\vee f t^{i+1}_{j+1}&
=f^\vee(u_+ \otimes ft^{i+1}_{j+1}), \\
u_+ \otimes ff^\vee t^{i+1}_{j+1}&
=f(u_+ \otimes f^\vee t^{i+1}_{j+1})-
f^\vee(fu_+ \otimes t^{i+1}_{j+1})+
u_- \otimes t^{i+1}_{j+1},
\end{align*}
where the coefficients in the first equality
are determined by acting with $ee^\vee$ and $e$.
The respective weights of $u_+ \otimes ft^{i+1}_{j+1}$,
$u_+ \otimes f^\vee t^{i+1}_{j+1}$, and
$fu_+ \otimes t^{i+1}_{j+1}$ differ from the
highest weight of $V^{i+1}_j$ by
the positive roots $\mu^\vee$, $\mu$ and $\mu^\vee$,
so the image under $s^{i+1}_j$ of all three elements
is $0$. Using (\ref{eq:s}) and
(\ref{eq:C}) we obtain (\ref{eq:sCt}).
\end{proof}

Let $H$ be an affine $k$-group and let $S$ be a $k$-scheme.
A family of representations of $H$ over $S$, or simply
a representation of $H$ over $S$,
is a locally free $\sO_S$-module $\sV$ of finite type
together with
a homomorphism $H_S \rightarrow GL(\sV)$ over $S$.
Suppose that $H$ is
a semidirect product $U \rtimes M$. Then $M$ acts on $U$ by group automorphisms. A representation of $H$ over
$S$ is the same as a representation $\sV$ of $M$ together with an
$M$-homomorphism
$\rho:U_S \rightarrow GL(\sV)$
over $S$. If $U$ is of finite type and $V$
is its Lie algebra, the
$M$-structure on $U$ defines a structure of $M$-Lie algebra on $V$.
>From $\rho$ we then obtain a homomorphism
$\sigma:V \otimes_k \sO_S \rightarrow \underline\End_{\sO_S}(\sV)$
of $(M,\sO_S)$-Lie algebras.

Suppose that $U$ is a vector group and
$S$ is quasi-compact. Then the $\rho$ correspond
uniquely to the $\sigma$ for which
$\sigma^n=0$ for some $n$. It follows that
a representation of $H$ over $S$ is the same
as a pair $(\sV,\mu)$, where $\sV$ is a representation of $M$
over $S$ and where $\mu$ is an $M$-homomorphism
\begin{equation}\label{eq:lrep}
\mu:V \otimes_k \sV \rightarrow \sV
\end{equation}
of $\sO_S$-modules such that the
\begin{equation}\label{eq:iter}
\mu^{(n)}: V^{\otimes n} \otimes_k \sV \rightarrow \sV,
\end{equation}
defined inductively by $\mu^{(0)} = 1_{\sV}$ and
$\mu^{(n)} = \mu \circ (V \otimes_k \mu^{(n-1)})$,
are $0$ for $n$ large,
and such that the restriction of $\mu^{(2)}$ to
$\bigwedge^2 V \otimes_k \sV$ is $0$
(corresponding to $[\sigma,\sigma]=0$).

When $S=\Spec(k)$ the images of the $\mu^{(n)}$ define an
$H$-invariant filtration on $\sV$.
Further $H$ acts on the associated graded $\GR_{\mu}\sV$ through $M$,
and $\mu$ induces an $M$-homomorphism
$V \otimes_k \GR_\mu^n\sV \rightarrow \GR_\mu^{n+1}\sV$
for each $n$.

\begin{lem}\label{le:Afam}
Suppose that $k$ is algebraically closed, and
let $H$ be a connected affine $k$-group with prounipotent
radical of dimension $>1$. Then there exists a
representation of $H$ over
the affine line $\bA^1$ with fibres at $k$-points
of $\bA^1$ pairwise non-isomorphic.
\end{lem}

\enlargethispage*{20pt}

\begin{proof}
We may suppose that $H$ is of finite type over $k$. If $U$
is the prounipotent radical of $H$ and $Z_1\subset Z_2\subset \dots
\subset Z_{m-1}\subset Z_m=U$ is the upper central series of $U$,
then the $Z_i$ are normal subgroups of $H$, and $Z_m/Z_{m-1}$ is
of dimension $>1$ since $U$ is of dimension $>1$. Factoring out
$Z_{m-1}$, we may suppose that $U$ is commutative. We have
$H=U\rtimes M$ with $M$ a connected reductive subgroup of $H$.
Write $V$ for the Lie algebra of $U$.
The action of $M$ on $U$ defines a representation of $M$ on
$V$. If $\bar U$ is an $M$-quotient
of $U$ of dimension $>1$ we may replace $H=U\rtimes M$ by
its quotient $\bar U \rtimes M$, so we may suppose that
$V$ is either irreducible of dimension $>1$ or is the
direct sum of two $1$-dimensional representations. It will be convenient
to distinguish the following three cases:
\begin{itemize}
\item[I] $V$ is the direct sum of two isomorphic
$1$-dimensional representations;
\item[II] $V$ is either the direct sum of two
non-isomorphic $1$-dimensional representations, or is
irreducible of dimension $2$ or $3$;
\item[III] $V$ is irreducible of dimension $>3$.
\end{itemize}

Suppose I holds. Then we may identify $U$ with $\bG_a^2$,
where $M$ acts on both factors as the same character $\chi$.
With $\bA^1=\Spec(k[t])$, we may identify the endomorphism
ring of $(\bG_a)_{\bA^1}$ over $\bA^1$ with $k[t]$.
If $M$ acts on $(\bG_a)_{\bA^1}$ as $\chi$,
then the semidirect product over $\bA^1$ of
$(1,-t):U_{\bA^1}=(\bG_a)_{\bA^1}^2
\rightarrow (\bG_a)_{\bA^1}$ with $M_{\bA^1}$
is a homomorphism $p:H_{\bA^1} \rightarrow
(\bG_a \rtimes M)_{\bA^1}$ over $\bA^1$. Given a faithful
representation $\tau:\bG_a \rtimes M \rightarrow GL(n)$,
the kernel of the fibre above $x \in \bA^1(k)=k$ of
$\tau_{\bA^1} \circ p: H_{\bA^1} \rightarrow GL(n)_{\bA^1}$
is the
subgroup $\binom{x}{1}:\bG_a \rightarrow \bG_a^2 = U$ of
$U \subset H$. The fibres of the representation
$\tau_{\bA^1} \circ p$ of $H$ over $\bA^1$
are thus pairwise non-isomorphic.

\begin{sloppypar}
Now suppose II holds. Apply Lemma~\ref{le:Vrs}
to obtain a family $(V^i_j,r^i_j,s^i_j)_{i,j\in\N}$,
with $r^i_j:V \otimes_k V^i_j \rightarrow V^{i+1}_j$
and $s^i_j:V \otimes_k V^i_{j+1} \rightarrow V^i_j$,
such that (i) and (ii) of Lemma~\ref{le:Vrs} hold.
Write $\bA^1=\Spec(k[t])$, and define
a family $(\sE^i_j,f^i_j,g^i_j)_{i,j\in\N}$
of free $\sO_{\bA^1}$-modules of finite type $\sE^i_j$
and $\sO_{\bA^1}$-homomorphisms
$f^i_j:\sE^i_j \rightarrow \sE^{i+1}_j$,
$g^i_j:\sE^i_{j+1} \rightarrow \sE^i_j$ by
$\sE^0_0=\sE^1_0=\sE^0_1=\sE^2_1=\sE^1_2=\sE^2_2=\sO_{\bA^1}$,
$\sE^1_1=\sO_{\bA^1}^2$, $\sE^i_j=0$ for all other $i,j$, and
$f^0_0=g^0_0=f^1_2=g^2_1=1$ and $f^0_1=\binom{1}{0}$,
$g^1_0=(1,-t)$, $g^1_1=\binom{0}{1}$, $f^1_1=(-1,1)$. The non-zero
$\sE$, $f$, $g$ thus form a commutative diagram
\begin{equation}\label{di:OA}
\begin{CD}
{} @.  \sO_{\bA^1}  @>1>>  \sO_{\bA^1} \\
@.     @VV{\binom{0}{1}}V     @VV1V \\
\sO_{\bA^1}  @>{\binom{1}{0}}>> \sO_{\bA^1}^2
@>{(-1,1)}>>  \sO_{\bA^1}  \\
@VV1V    @VV{(1,-t)}V   @. \\
\sO_{\bA^1} @>1>>  \sO_{\bA^1} @. {}
\end{CD}
\end{equation}
\end{sloppypar}

Now set $\sV=\bigoplus_{i,j\in\N}V^i_j \otimes \sE^i_j$ and let
$\mu:V \otimes_k \sV \rightarrow \sV$ have component
\[
V \otimes_k V^i_j \otimes_k \sE^i_j
\rightarrow V^{i'}_{j'} \otimes \sE^{i'}_{j'}
\]
given by
$r^i_j \otimes f^i_j$ if $(i',j')=(i+1,j)$,
$s^i_j \otimes g^i_j$ if $(i',j'+1)=(i,j)$, and $0$ in all
other cases.  Clearly $\mu^{(3)}=0$, and
the restriction of $\mu^{(2)}$ to
$\bigwedge^2 V \otimes_k \sV$ is $0$ by (ii) of
Lemma~\ref{le:Vrs} together with the fact that (\ref{di:OA})
commutes. Thus $(\sV,\mu)$ is a representation
of $H$ over $\bA^1$.

We now show that an isomorphism between the fibres $(\sV_x,\mu_x)$
and $(\sV_y,\mu_y)$ of $(\sV,\mu)$ at the $k$-points $x$ and $y$
of $\bA^1$ induces an isomorphism between the fibres at $x$ and
$y$ of the diagram (\ref{di:OA}). This will give an isomorphism
$k^2=(\sE^1_1)_x \simeq (\sE^1_1)_y=k^2$ which sends $\Img(g^1_1)_x$,
$\Ker(f^1_1)_x$, $\Img(f^0_1)_x$, $\Ker(g^1_0)_x$ respectively to
$\Img(g^1_1)_y$, $\Ker(f^1_1)_y$, $\Img(f^0_1)_y$, $\Ker(g^1_0)_y$,
hence an isomorphism $\bP^1 \simeq \bP^1$ which sends
$0,1,\infty,x$ respectively to $0,1,\infty,y$, and so will imply
$x=y$ as required.

For any $x$ we have
$\sV_x=\bigoplus_{i,j\in\N}V^i_j \otimes_k E^i_j$ with
$E^i_j$ the fibre of $\sE^i_j$ at $x$. Write
$\sV_x^n=\bigoplus_{i-j=n}V^i_j \otimes_k E^i_j$,
so that $\sV_x = \bigoplus_{n\in\Z}\sV_x^n$. Then
$\mu_x:V \otimes_k \sV_x \rightarrow \sV_x$
restricts to an $M$-homomorphism
\[
\mu_x^n: V \otimes_k \sV_x^n \rightarrow
\sV_x^{n+1}
\]
for each $n$. Further $\mu_x^n$
is an epimorphism for $n\ge-1$ because its composite with each
projection $\sV_x^{n+1} \rightarrow V^i_j \otimes_k E^i_j$
is an epimorphism by (\ref{di:OA}), and the $V^i_j$ with $i-j=n+1$
are pairwise non-$M$-isomorphic
by (i) of Lemma~\ref{le:Vrs}. Since $\mu_x^n=0$ for
$n<-1$, it follows that the image of
$\mu_x^{(n)}: V^{\otimes n} \otimes_k \sV_x \rightarrow \sV_x$
is $\bigoplus_{m\ge n-1} \sV_x^m$. Hence
$\GR_{\mu_x}^n \sV_x = \sV_x^{n-1}$, and
$V \otimes_k \GR_{\mu_x}^n \sV_x \rightarrow \GR_{\mu_x}^{n+1} \sV_x$
induced by $\mu_x$ coincides with $\mu_x^{n-1}$. Thus any
isomorphism $\theta:(\sV_x,\mu_x) \rightarrow (\sV_y,\mu_y)$
induces on passage to the associated graded
an isomorphism
\[
\theta^n:\sV_x^n \rightarrow \sV_y^n
\]
for each $n$ such that the diagram
\begin{equation}\label{di:theta}
\begin{CD}
V \otimes_k \sV_x^n @>{\mu_x^n}>> \sV_x^{n+1}\\
@VV{V \otimes \theta^n}V   @VV{\theta^{n+1}}V\\
V \otimes_k \sV_y^n @>{\mu_y^n}>> \sV_y^{n+1}
\end{CD}
\end{equation}
commutes. Since for each $n$ the $V^i_j$ with $i-j=n$ are
pairwise non-$M$-isomorphic by (i) of Lemma~\ref{le:Vrs},
there are $k$-isomorphisms $\theta^i_j:E^i_j \rightarrow E^i_j$
such that
\[
\theta^n=V^i_j \otimes_k \theta^i_j.
\]
 From (\ref{di:theta}) it then follows that
\begin{align*}
r^i_j \otimes (\theta^{i+1}_j(f^i_j)_x)&
=r^i_j \otimes ((f^i_j)_y\theta^i_j)\\
s^i_j \otimes (\theta^i_j(g^i_j)_x)&
=s^i_j \otimes ((g^i_j)_y\theta^i_{j+1})
\end{align*}
for $i,j\ge0$. Since the $r^i_j$ and $s^i_j$ are non-zero,
\begin{align*}
\theta^{i+1}_j(f^i_j)_x&
=(f^i_j)_y\theta^i_j\\
\theta^i_j(g^i_j)_x&
=(g^i_j)_y\theta^i_{j+1},
\end{align*}
so the $\theta^i_j$ give an isomorphism between the fibres of
(\ref{di:OA}) at $x$ and $y$.

Finally suppose III holds. In the notation of Lemma~\ref{le:te},
choose $\tau\in P$ such
that $\tau+\lambda\in P$ and $\tau + \lambda \notin \Lambda$ for
$\lambda\in\Lambda$. Then
$V \otimes_k V_\tau \simeq
\bigoplus_{\lambda \in \Lambda} V_{\tau+\lambda}^{m(\lambda)}$.
Choose by Lemma~\ref{le:A} distinct non-zero elements $\lambda_1$,
$\lambda_2$, $\lambda_3$, $\lambda_4$ of $\Lambda$, and write
$V_0=V_\tau$ and $V_j=V_{\tau+\lambda_j}$, $j=1,2,3,4$. The $V_i$,
$i=0,1,2,3,4$ are then pairwise non-isomorphic irreducible
representations of $M$, and $V_j$ is a
direct summand of $V \otimes_k V_0$ for $j=1,2,3,4$.
Let $r_j:V \otimes_k V_0 \rightarrow V_j$ be a non-zero
$M$-homomorphism for $j=1,2,3,4$. With $\bA^1 = \Spec(k[t])$,
set $\sE_0=\sO_{\bA^1}^2$,
$\sE_j=\sO_{\bA^1}$ for $j=1,2,3,4$, and define
$f_j:\sE_j \rightarrow \sE_0$ for $j=1,2,3,4$ by
$f_1 = (1,0)$, $f_2 = (-1,1)$, $f_3 = (0,1)$,
$f_4 = (1,-t)$. Now set
$\sV  =\bigoplus_{i=0}^4 V_i \otimes_k \sE_i$ and let
$\mu:V \otimes_k \sV \rightarrow \sV$ have component
\[
V \otimes_k V_i \otimes \sE_i
\rightarrow V_{i'} \otimes \sE_{i'}
\]
given by $r_j \otimes f_j$ when $i=0,i'=j$ and $j=1,2,3,4$, with
all other components $0$.
Then $\mu^{(2)}=0$ so $(\sV,\mu)$ is a representation
of $H$ over $\bA^1$.

Suppose that $\theta:\sV_x \rightarrow \sV_y$
is an isomorphism between the fibres $(\sV_x,\mu_x)$ and
$(\sV_y,\mu_y)$ at $x$ and $y$. Since
$\sV_x=\bigoplus_{i=0}^4 V_i \otimes_k E_i$
with $E_i$ the fibre of $\sE_i$ at $x$, and since
the $V_i$ are pairwise non-$M$-isomorphic, we have
$\theta=\bigoplus_{i=0}^4 V_i \otimes \theta_i$
for isomorphisms $\theta_i:E_i \rightarrow E_i$. Compatibility
of $\theta$ with $\mu_x$ and $\mu_y$ then shows that for
$j=1,2,3,4$
\[
r_j \otimes (\theta_j(f_j)_x) = r_j \otimes ((f_j)_y\theta_0)
\]
whence
\[
\theta_j(f_j)_x  = (f_j)_y\theta_0
\]
since $r_j\ne0$. Thus the $\theta_i$ give an isomorphism
$k^2=E_0\simeq E_0=k^2$ which sends $\Ker (f_j)_x$ to
$\Ker (f_j)_y$ for $j=1,2,3,4$, and hence an isomorphism
$\bP^1\simeq\bP^1$ which sends $0,1,\infty,x$ respectively
to $0,1,\infty,y$. Thus $x=y$ as required.
\end{proof}

\begin{lem}\label{le:GL2}
Let $H$ be an affine $k$-group with prounipotent
radical $U$ of dimension $\le 1$. Then
there exists a $k$-homomorphism $H \rightarrow GL(2)$
with restriction to $U$ an embedding. If $g$ is such a
$k$-homomorphism and $p:H \rightarrow H/U$ is the projection
then every indecomposable representation of $H$ is isomorphic
to $g^*V \otimes_k p^*W$ for some representation $V$ of $GL(2)$
and $W$ of $H/U$.
\end{lem}
\begin{proof}
We may suppose that $U$ is of dimension $1$.
Then $M=H/U$ acts on $U$ as a character $\chi$. Taking
a Levi decomposition of $H$ and identifying $U$ with $\bG_a$
we may suppose $H=\bG_a \rtimes_\chi M$.
The restriction to $U$ of the $k$-homomorphism
$H \rightarrow GL(2)$ defined by
\[
am \mapsto
{\begin{pmatrix}
\chi(m)&a\\ 0&1
\end{pmatrix}}
\]
for points $a$ of $\bG_a$ and $m$ of $M$ is then an embedding.

A representation of $H$ on a $k$-vector space $E$ is a
representation of $M$ on $E$ together with an $M$-homomorphism
$E \otimes_k \chi \rightarrow E$ such that
$E \otimes_k \chi^r \rightarrow E$ induced by iteration is $0$
for $r$ large. Equivalently, a representation of $H$ on $E$ is
a representation of $M$ on $E$ together with a structure of
module on $E$ over the polynomial ring $k[n]$ such that $n$ acts
nilpotently and $mne=\chi(m)nme$
for each $k$-algebra $k'$, $m \in M(k')$ and
$e \in E \otimes_k k'$.

If $D$ is the standard $2$-dimensional
representation of $GL(2)$, then the kernel of
$g^*D \rightarrow g^*D$ defined by the
action of $n$ is $1$-dimensional and stable under $M$.
Let $d_0\ne0$ be an element of a $1$-dimensional
$M$-stable subspace of $g^*D$ complementary to this kernel. Then
$d_1=nd_0\ne0$, $nd_1=0$, and $md_0=\chi_0(m)d_1$,
$md_1=\chi(m)\chi_0(m)d_1$ for points $m$ of $M$ and some character
$\chi_0$ of $M$.
Thus if $V$ is the $r$th symmetric power of $D$, there is
a basis $v_0,\dots,v_r$ of $g^*V$ for which
$mv_i=\chi^i(m)\chi_0^r(m)v_i$ and $nv_i=v_{i+1}$
if $i<r$ and $nv_r=0$.

Now let $E$ be a representation of $H$. Let $E_0$ be an
$M$-subrepresentation of $E$, and suppose that $n^{r+1}E=0$.
Then if $W = E_0 \otimes \chi_0^{-r}$, the assignment
$v_i \otimes w \mapsto n^iw$
defines an $H$-homomorphism
$\phi:g^*V \otimes_k p^*W \rightarrow E$.
For the second statement, it will thus suffice to
show that if $E$ is $H$-indecomposable, and if $r$ and
$E_0$ are chosen so that $n^rE\ne0$, $n^{r+1}E=0$,
and $E_0 \oplus nE = E$,
then $\phi$ is an isomorphism,
i.e.\ each $e\in E$ has a unique decomposition
\begin{equation}\label{eq:dec}
e=e_0+ne_1+\dots+n^re_r
\end{equation}
with $e_i\in E_0$ for $i=0,\dots,r$. In fact from
$E_0 \oplus nE = E$
it follows inductively that $E=E_0+\dots+n^iE_0+n^{i+1}E$ so taking
$i=r$ shows that each $e\in E$ has a decomposition (\ref{eq:dec}).
To prove the uniqueness, write $E''$ for the kernel of $n^r$
restricted
to $E_0$ and let $E'$ be an $M$-subspace of $E_0$ such that
$E_0=E' \oplus E''$. Since $n^rE\ne0$ we have $E'\ne0$. If
$e' = \sum_{i=0}^r n^i e'_i$ and $e'' = \sum_{i=0}^r n^i e''_i$
with $e'_i \in E'$ and $e''_i \in E''$, then $e'=e''$ implies
$e'_i=0$ for $i=0,1,\dots,r$, as is seen by induction on $i$:
given that $e'_0=e'_1=\cdots=e'_{i-1}=0$, from
$n^{r-i}e'=n^{r-i}e''$ we have $n^re'_i=n^re''_i=0$ whence
$e'_i=0$. Thus $k[n]E' \cap k[n]E''=0$ and
$E=k[n]E_0=k[n]E' \oplus k[n]E''$. Since $E$ is $H$-indecomposable
and $k[n]E'$ and $k[n]E''$ are $H$-subspaces of $E$, we have $E''=0$
and $E'=E_0$, whence the required uniqueness since $e'=0$ implies
$e'_i=0$ for $i=0,1,\dots,r$.
\end{proof}

\begin{lem}\label{le:sumnd}
Let $H$ be an affine $k$-group, let $k'$ be an extension of $k$,
and let $V'$ be a representation of $H$ over $k'$.

\textup{(i)} If $H$ is proreductive then $V'$
is a direct summand of some $V \otimes_k k'$.

\textup{(ii)} If $k'$ is algebraic over $k$ then $V'$
is a direct summand of some $V \otimes_k k'$.

\textup{(iii)} If $k$ is algebraically closed and
$V'$ is a direct summand of $V \otimes_k k'$, then
$V' \simeq V_0 \otimes_k k'$ for some direct summand $V_0$ of $V$.
\end{lem}

\begin{proof}
There is a canonical epimorphism
\begin{equation}\label{eq:epi}
V' \otimes_k k' \rightarrow V' \rightarrow 0
\end{equation}
of (not necessarily finite dimensional) $H$-representations over $k'$.
Since $V'$ is
finite-dimensional over $k'$ there is a finite dimensional
$H$-sub\-re\-pre\-sen\-ta\-tion $V$ of $V'$ over $k$ such that the
restriction of (\ref{eq:epi}) to $V \otimes_k k'$ is
an epimorphism. If $H$ is reductive, $V'$ is a direct summand of
$V \otimes_k k'$, whence (i).
To prove (ii) we may suppose $k'$ finite over $k$. Then $V'$ is
finite dimensional over $k$ and there is
a canonical splitting of (\ref{eq:epi}), defined using the trace
from $k'$ to $k$. For (iii), note that if $A=\End_H(V)$,
any idempotent $e$ in
$A \otimes_k k'$ is conjugate to one in $A$. Indeed $e$ generates a
semisimple $k'$-subalgebra of $A \otimes_k k'$ and so some conjugate
lies in the extension to $k'$ of a maximal semisimple subalgebra
of $A$. With $k$ algebraically closed, it is then obvious that
some conjugate lies in $A$.
\end{proof}

\begin{Th}\label{th:fam}
If $H$ is an affine $k$-group the following conditions are
equivalent:

\textup{(a)} the prounipotent radical of $H$ is of dimension $\le 1$;

\textup{(b)} for every $k$-scheme $S$ of finite type and family
$\sV$ of representations of $H$ over $S$, there is
a finite extension $k'$ of $k$ and a stratification of $S_{k'}$ by
locally closed subschemes such that the family $\sV_{k'}$ of
representations of $H_{k'}$ over $S_{k'}$
is constant along each stratum;

\textup{(c)} for every $k$-scheme $S$ of finite type and family
$\sV$ of representations
of $H$ over $S$, the set of
$H$-isomorphism classes of the fibres of $\sV$ at $k$-points of
$S$ is finite;

\textup{(d)} for every family $\sV$ representations of $H$
over $\bA^1$
there exists a representation $V$ of $H$ over $k$
such that $\sV_x \simeq V$
for an infinite set of $x\in\bA^1(k)$.
\end{Th}

\begin{proof}
(a) $\Longrightarrow$ (b): We may suppose that $k$ is
algebraically closed. It is enough to show then that, assuming
the prounipotent  radical $U$ of $H$ is of dimension $\le 1$,
each representation $\sV$
of $H$ over a non-empty reduced and irreducible
$k$-scheme $S$ of finite type is constant along some non-empty
open subscheme of $S$. Let $V_1$ be the generic fibre of $\sV$.
If $K=GL(2) \times M/U$ and $k'$ is the function field of $S$
there is by Lemma~\ref{le:GL2} a
$k$-homomorphism $f: H \rightarrow K$ such that
$V_1 \simeq f_{k'}^*V'$ for some representation $V'$ of $K_{k'}$.
Further $V' \simeq V \otimes_k k'$ for some representation
$V$ of $K$ by (i) and (iii) of Lemma~\ref{le:sumnd}, so there
is an $H$-isomorphism $V_1 \simeq f^*V \otimes_k k'$. This
extends to an $H$-isomorphism of $\sV$ with the constant
family $f^*V$ along some non-empty open subscheme of $S$.

(b) $\Longrightarrow$ (c): It is enough to note that
representations $V$ and $W$ of $H$ which become isomorphic
over $k'$ are themselves isomorphic, because the
determinant is then
not identically zero on $\Hom_H(V,W)$.

(c) $\Longrightarrow$ (d) is immediate.

(d) $\Longrightarrow$ (a): Call the fibres of a representation
$\sV$ of $H$ over $\bA^1$ almost distinct if
for each representation $V$ of $H$ over $k$ we have
$\sV_x \simeq V$ for only a finite set of $x \in \bA^1(k)$.
The fibres of $\sV$ are
almost distinct provided that those of the representation obtained
from $\sV$ by either extension of scalars, pullback along some
$H' \rightarrow H$, or passage to a direct summand,
are almost distinct. This is immediate in the first two cases
and in the case of direct summands follows from the fact that
by Krull--Schmidt any $\sV_x$ has only a finite set of pairwise
non-$H$-isomorphic direct summands.
Now suppose that (a) does not hold. Then we show that there
exists a $\sV$ whose fibres are almost distinct, so that (d)
does not hold. To do this
we may suppose by passing to a quotient that $H$ is of
finite type.

Consider first the case where $k$ is algebraically closed.
If $H^0$ is the connected component of $H$ there is then by
Lemma~\ref{le:Afam} a representation $\sV^0$ of $H^0$ over
$\bA^1$ with almost distinct fibres. Since $\sV^0$ is a direct
$H^0$-summand of the induced $H$-representation
$\Ind_{H^0}^H\sV^0$, the fibres of
$\Ind_{H^0}^H\sV^0$ are almost distinct.

For arbitrary $k$, with algebraic closure $\bar k$, there is thus
a representation $\bar\sV$ of $H$ over $\bA_{\bar k}^1$ with
almost distinct fibres. We have $\bar\sV=\sV' \otimes_{k'} \bar k$
for some finite extension $k'$ of $k$ and $H$-representation $\sV'$
over $\bA_{k'}$. If $p:\bA_{k'}^1 \rightarrow \bA^1$ is the
projection and we write $\sV=p_*\sV'$ then $\sV'$ is a direct
$H$-summand of $\sV \otimes_k k' = p^*p_*\sV'$. The fibres of
$\sV'$ and hence $\sV$ are thus almost distinct.
\end{proof}

\begin{Th}\label{th:sumnd}
Let $H$ be an affine $k$-group and let $k'$ be an extension of $k$.
Then every representation
of $H$ over $k'$ is a direct summand of a
representation defined over $k$
if and only if
either $k'$ is algebraic over $k$ or
the prounipotent radical of $H$ has dimension $\le 1$.
\end{Th}

\begin{proof}
If $k'$ is algebraic over $k$ then
each representation $V'$ of $H$ over $k'$ is a
direct summand of some $V \otimes_k k'$ by Lemma~\ref{le:sumnd}(ii).
Suppose that the prounipotent
radical $U$ of $H$ has dimension $\le 1$.
If $K=GL(2) \times M/U$
there is then by Lemma~\ref{le:GL2} a
$k$-homomorphism $f: H \rightarrow K$ such that
any representation $V'$ of $H_{k'}$ is isomorphic to
$f_{k'}^*W'$ for some representation $W'$ of $K_{k'}$.
Then $W'$ is a direct summand of $W \otimes_k k'$ for some
representation $W$ of $K$ by Lemma~\ref{le:sumnd}(i), so
$V'$ is a direct summand of $f^*W \otimes_k k'$.

For the converse, suppose that
$k'$ is transcendental over $k$ and that the prounipotent
radical of $H$ has dimension $>1$.
Then by the equivalence
of (a) and (d) of Theorem~\ref{th:fam}, there is a representation
$\sV$ of $H$ over $\bA^1$ with an infinite set of pairwise
non-isomorphic fibres at $k$-points of $\bA^1$.
Choose an embedding of $k(t)$ in $k'$, and let $V_1$ be
the generic fibre of $\sV$. We show that
$V_1 \otimes_{k(t)} k'$ is
not a direct summand of any $V \otimes_k k'$.
In fact suppose the contrary.
Then for some $V$ the identity of $V_1$ lies
in the image of the composition homomorphism
\[
\Hom_{H,k(t)}(V \otimes_k k(t),V_1) \otimes_{k(t)}
\Hom_{H,k(t)}(V_1,V \otimes_k k(t)) \rightarrow
\End_{H,k(t)}V_1,
\]
since this is so after extending scalars from $k(t)$ to $k'$.
Thus for some $n$ there are $H$-homomorphisms
$r_i:V \otimes_k k(t) \rightarrow V_1$,
$s_i:V_1 \rightarrow V \otimes_k k(t)$, $i=1,\dots,n$
over $k(t)$ such that $r_1s_1+\dots+r_ns_n=1_{V_1}$.
This implies that $V_1$ is a direct summand of
$V^n \otimes_k k(t)$. Thus $\sV$ is a direct summand of
the constant family $V^n$ along some non-empty open subscheme
of $\bA^1$. By Krull--Schmidt this contradicts the fact that
$\sV$ has an infinite set of pairwise non-isomorphic fibres.
\end{proof}

\subsection{Application to proreductive envelopes}\label{sec:app}

In connection with proreductive envelopes it will be convenient
to use the following terminology. Let $f:H \rightarrow K$
be a $k$-homomorphism from an affine $k$-group to a proreductive
$k$-group. Call $f$ \emph{universal} if each
$k$-homomorphism $H \rightarrow K'$ to a proreductive $K'$
factors, uniquely up to conjugation by a $k$-point of $K'$,
through $f$. Equivalently, $f$ is universal if its conjugacy class
is the embedding of $H$ into its proreductive
envelope. Th\'eor\`eme \ref{t1} states that for each $H$
there exists a universal $f$. Call $f$
\emph{minimal} if $f$ factors through no proper proreductive
subgroup of $K$. If $L$ is proreductive, $h:H \rightarrow L$ is
universal, and $f=lh$, then $f$ is universal if and only if
$l$ is an isomorphism and $f$ is minimal if and only if $l$
is faithfully flat.

By considering $k$-homomorphisms to general linear groups
it can be seen that if $f$ is universal then $f^*$ induces
a bijection on isomorphism classes of objects (Proposition
\ref{p14.1}). We note also the following fact, which generalises
a result due to Kostant \cite[Theorem 3.6]{kos}. If $f$ is universal,
then any $k$-homomorphism $h:H \rightarrow K'$ to a
proreductive $K'$ factors through $f$ uniquely up to conjugation
by a \emph{unique} $k$-point of the prounipotent radical
$R_uC_{K'}(h)$ of the
centraliser $C_{K'}(h)$ of $h$. Equivalently, if $h=lf$ then
$C_{K'}(l)$ is a Levi subgroup of $C_{K'}(h)$, i.e.\ $C_{K'}(h)$
is the semidirect product of $R_uC_{K'}(h)$ and $C_{K'}(l)$. To verify
this, let $L$ be a Levi subgroup of $C_{K'}(h)$. Then
$h$ factors through the proreductive subgroup $C_{K'}(L)$ of
$K'$, so by universality of $f$ the conjugate of
$l$ by some $k$-point $z$ of $C_{K'}(h)$
factors through $C_{K'}(L)$. Thus
$L \subset zC_{K'}(l)z^{-1} \subset C_{K'}(h)$ and
$zC_{K'}(l)z^{-1} = L$ since $C_{K'}(l)$ is proreductive.

The main facts required about universal and minimal
$k$-homomorphisms are contained in the following lemma.

\begin{lem}\label{le:univ}
Let $f:H \rightarrow K$ be a $k$-homomorphism from an affine
$k$-group to a proreductive $k$-group. Then

\textup{(i)} $f$ is universal if and only if $f$ is minimal and
each representation of $H$ is a direct summand of $f^*V$
for some representation $V$ of $K$.

\textup{(ii)}
If $f$ is
minimal then $f_{k'}:H_{k'} \rightarrow K_{k'}$ is minimal
for any extension $k'$ of $k$.

\textup{(iii)} If $f$ is
universal then $f_{k'}:H_{k'} \rightarrow K_{k'}$ is universal
for any algebraic extension $k'$ of $k$.
\end{lem}

\begin{proof}
Factor $f$ as
$H \stackrel{g}{\rightarrow} L \stackrel{l}{\rightarrow} K$
with $L$ proreductive and $g$ universal. Then $g^*$ induces
a bijection on isomorphism classes of objects, and $f$ is
universal (resp.\ minimal) if and only if $l$ is an isomorphism
(resp.\ faithfully flat).

Thus $f$ is universal $\Longleftrightarrow$ $f$
is minimal and $l$ is a closed immersion. Further
$l$ is a closed immersion $\Longleftrightarrow$ each representation
of $L$ is a direct summand of some $l^*V$ (because $L$ is
proreductive) $\Longleftrightarrow$ each representation
of $H$ is a direct summand of some $f^*V$. This gives (i).

\begin{sloppypar}
If $W$ is representation  of a
$k$-group $G$, write $m_G(W)$ for the multiplicity of the trivial
$G$-representation $k$ in a decomposition of $W$ into
$G$-in\-de\-com\-pos\-ables. Clearly $m_G(W)$ is the rank of the
canonical $k$-homomorphism $W^G \rightarrow W_G$, so it is
preserved by extension of the scalars. Then
$l$ is surjective
$\Longleftrightarrow$ $l^*$ is fully faithful (because $L$ is proreductive)
$\Longleftrightarrow$ $\dim_k\Hom_K(k,V)=\dim_k\Hom_L(k,l^*V)$ for each
representation $V$ of $K$
$\Longleftrightarrow$ $m_K(V)=m_L(l^*V)$ for each $V$ $\Longleftrightarrow$
$m_K(V)=m_H(f^*V)$ for each $V$
$\Longleftrightarrow$ each $V$ is a direct summand of some $V_1$
for which $m_K(V_1)=m_H(f^*V_1)$ (because $m_K(V) \le m_H(f^*V)$).
By Lemma~\ref{le:sumnd}(i), if the last condition is satisfied
it is also satisfied with $f$, $H$ and $K$ replaced by $f_{k'}$,
$H_{k'}$ and $K_{k'}$. This gives (ii).
\end{sloppypar}

For (iii) it is enough by (i) and (ii) to show if each
representation of $H$ is a direct summand of some $f^*V$
then each representation of $H_{k'}$ is a direct summand
of some $f_{k'}^*V'$. This is clear from Lemma~\ref{le:sumnd}(ii).
\end{proof}

\begin{Th}\label{th:basech}
Let $H$ be an affine $k$-group and  let $k'$ be an extension of $k$.
Then the canonical morphism
${}^\mathrm{p}\Red(H_{k'}) \rightarrow {}^\mathrm{p}\Red(H)_{k'}$
is an isomorphism if and only if either $k'$ is algebraic over $k$
or the prounipotent radical of $H$ is of dimension $\le 1$.
\end{Th}

\begin{sloppypar}
\begin{proof}
Let $f:H \rightarrow K$, with $K$ proreductive, be universal.
Then
${}^\mathrm{p}\Red(H_{k'}) \rightarrow {}^\mathrm{p}\Red(H)_{k'}$
is an isomorphism $\Longleftrightarrow$
$f_{k'}$ is universal $\Longleftrightarrow$ each representation of
$H_{k'}$ is a direct summand of some $f_{k'}^*W'$
(Lemma~\ref{le:univ})
$\Longleftrightarrow$ each representation of $H_{k'}$ is a direct
summand of some $f^*W \otimes_k k'$ (Lemma~\ref{le:sumnd}(i))
$\Longleftrightarrow$
each representation of $H_{k'}$ is a direct summand of some
$V \otimes_k k'$. The result thus follows from
Theorem~\ref{th:sumnd}.
\end{proof}
\end{sloppypar}

\begin{Th}\label{th:SL2}
Let $U$ be a unipotent $k$-group of dimension $1$, equipped
with an action of the proreductive $k$-group $M$. Then

\textup{(i)} There exists an $M$-embedding $U \rightarrow SL(2)$
for some action of $M$ on $SL(2)$.

\textup{(ii)} If $f:U \rightarrow SL(2)$ is an $M$-embedding and
$h:U \rightarrow K$ is an $M$-homomorphism with $K$ proreductive,
then there is an $M$-homomorphism $l:SL(2) \rightarrow K$,
unique up to conjugation by a $k$-point of $K$ fixed by $M$,
such that $h=lf$.

\textup{(iii)} If $f:U \rightarrow SL(2)$ is an $M$-embedding
then (the conjugacy class of)
$f \rtimes M:U \rtimes M \rightarrow SL(2) \rtimes M$
is the embedding of $U \rtimes M$ into its proreductive envelope.
\end{Th}

\begin{proof}
By Lemma~\ref{le:GL2} there is a $k$-homomorphism
$g:U \rtimes M \rightarrow GL(2)$ non-trivial on $U$.
The restriction of $g$ to $U$ factors through an $M$-embedding
$f:U \rightarrow SL(2)$, where $M$ acts on $SL(2)$ through
$g$ and the inner action of $GL(2)$. This gives (i).

It will suffice to prove (ii) and (iii) for the $f$
just constructed. Consider first (iii).
Let $L$ be a proreductive subgroup of $SL(2) \rtimes M$ through
which $f \rtimes M$ factors. Then $L \cap SL(2) \supset U$ is
reductive since $SL(2)$ is normal in $SL(2) \rtimes M$, so
$L \cap SL(2) = L$ and $L \supset SL(2)$. Also $L \supset M$,
so $L = SL(2) \rtimes M$. Thus $f \rtimes M$ is minimal.
The projection $U \rtimes M \rightarrow M$ factors through
$f \rtimes M:U \rtimes M \rightarrow SL(2) \rtimes M$, and $g$ factors
through $f \rtimes M$ by construction. Thus
$(f \rtimes M)^*$ is essentially surjective by Lemma~\ref{le:GL2},
and so universal by Lemma~\ref{le:univ}(i).

Given $h$ as in (ii),
$h \rtimes M$ is by universality of $f \rtimes M$ the
composite of $f \rtimes M$
with a $k$-homomorphism $SL(2) \rtimes M \rightarrow K \rtimes M$,
necessarily of the form $l \rtimes M$ for some
$M$-homomorphism $l$. Thus $h=lf$. If also $h=l'f$ for an
$M$-homomorphism $l'$, then  $l$ and $l'$ are by universality of
$f$ conjugate by a unique $k$-point $\alpha$ of the prounipotent
radical of the centraliser of $h$. This implies that
$l_{k'}$ and $l'_{k'}$
are conjugate by $m\alpha$ for any extension $k'$ of $k$
and $m\in M(k')$. Since $M$ stabilises the centraliser of
$h$ and so the prounipotent radical of this centraliser,
$M$ fixes $\alpha$ by uniqueness.
\end{proof}

An affine $k$-group $G$ will
be called finite dimensional if its quotients of finite type are
of bounded dimension. It is equivalent to require that the connected
component of $G$ be an extension of a group of finite type by a
pro\'etale group.

\begin{lem}\label{le:fin}
If $k$ is algebraically closed then
a connected finite dimensional reductive $k$-group
has only a finite number of non-isomorphic
irreducible representations of given dimension and determinant.
\end{lem}

\begin{proof}
The connected semisimple case follows from
Weyl's dimension formula.
In general, a connected finite dimensional reductive group
is a quotient of $T \times G$ for some protorus $T$
and connected semisimple $G$. An irreducible representation
of $T \times G$ of dimension $n$ and determinant $\chi$ is
the tensor product of an $n$-dimensional irreducible
representation of $G$ with a $1$-dimensional representation
of $T$ defined by a character $\psi$ with $\psi^n = \chi$.
Since there is at most one such $\psi$, the result follows
from the connected semisimple case.
\end{proof}

\begin{Th}\label{th:finite}
The proreductive envelope of an affine $k$-group
$H$ is of finite type over $k$ (resp.\ finite dimensional)
if and only if $H$ is of finite type over $k$
(resp.\ finite dimensional) and the prounipotent radical of $H$ is of
dimension $\le1$.
\end{Th}

\begin{proof}
The ``if'' is clear by Theorem~\ref{th:SL2}(iii).
Let $f:H \rightarrow K$, with $K$ proreductive, be universal.
For the ``only if'', it suffices to show that if the
prounipotent radical of $H$ has dimension $>1$ then
$K$ is not finite dimensional.
By Lemma~\ref{le:univ}(iii) we may suppose that $k$ is algebraically
closed. If $\bar H$ is a quotient of $H$ and
$\bar H \rightarrow \bar K$, with $\bar K$ proreductive, is
universal, then $H \rightarrow \bar H \rightarrow \bar K$
is minimal, and so factors through a faithfully flat
$K \rightarrow \bar K$.
Replacing $H$ by an appropriate $\bar H$,
we may thus further suppose $H$ to
be of finite type (see also Proposition 19.4.7 c)).
Let $j:H^0 \rightarrow H$ be the embedding of the
connected component of $H$ and
$g:H^0 \rightarrow L$, with $L$ proreductive, be universal,
so that $fj=lg$ for some $l:L \rightarrow K$. Then $l$ is a
closed immersion. Indeed otherwise some representation of $L$
would not be a direct summand of any $l^*V$, whence the same
would hold with $L$ and $l$ replaced by $H^0$ and $j$,
since $f^*$ and $g^*$ induce bijections on
isomorphism classes of objects. But each representation $W$ of
$H^0$ is a direct summand of $j^*\Ind_{H^0}^H W$.
Replacing $H$ by $H^0$ we may thus suppose $H$ connected
(see also Proposition 19.4.7 a)).
Then $K$ also is connected, as is seen by considering
$k$-homomorphisms to finite groups.

By Theorem~\ref{th:fam} there is a
family $\sV$ of representations of $H$, pa\-ra\-me\-tri\-sed by $\bA^1$,
such that the set of $H$-isomorphism classes of the fibres $\sV_x$
for $x \in \bA^1(k)$ is infinite.
There thus exists an infinite set $\sS$ of pairwise
non-isomorphic irreducible representations $V$ of $H$ such that each
$V \in \sS$ is a direct summand of some $\sV_x$ for $x\in\bA^1(k)$.
Let $T$ be the radical of a a Levi subgroup of $H$.
The restrictions to $T$ of the $\sV_x$ are isomorphic, so among
direct summands of the $\sV_x$ there are only a finite number
of possible $T$-isomorphism classes. Since the determinant of
a representation of $H$ is defined uniquely by its restriction
to $T$, it follows that there is an infinite subset $\sS'$ of $\sS$
such that each $V \in \sS'$ has the same determinant $\chi$. Then
$\chi=\psi \circ f$ for a unique character $\psi$ of $K$.
Each $V \in \sS'$ is isomorphic to $f^*W$ for some representation
$W$ of $K$, and each such $W$ has determinant $\psi$ and
dimension at most that of the fibres of $\sV$.
Lemma~\ref{le:fin} then shows that $K$ is not finite dimensional.
\end{proof}

By a similar argument it can be shown
for example that if $k$ is uncountable the proreductive envelope
of an affine $k$-group with prounipotent radical of dimension $>1$
is not a countable limit of $k$-groups of finite type.

\newpage

\section*{Index terminologique}

{\Small

$A$-id\'eal \`a gauche \dotfill \ref{d1}, 
\ref{corresp}, \ref{L4.2}

$\sA$-module (\`a gauche) \dotfill \ref{D2 
1/2}, \ref{corresp}, \ref{D3},
\ref{prodtens}, appendice 
\ref{s10}

$\sA$-bimodule \dotfill \ref{prodtens}

anneau de 
re\-pr\'e\-sen\-ta\-tions \dotfill \ref{hh0}, \ref{15.2}

artinien \dotfill 
\ref{P4/3}, \ref{R1}, \ref{P3/2}, \ref{A.D6}, \ref{A.L9}

auto-dual 
\dotfill \ref{P1}, \ref{P4/3}, \ref{D4sep1}, \ref{R1}, 
\ref{rga}

Auslander (alg\`ebre d')  \dotfill 
\ref{12.2}\bigskip

biproduit \dotfill \ref{biprod}\bigskip

compact 
\dotfill \ref{d1.1}, \ref{r4}, \ref{l11}, \ref{l12}, 
\ref{p6},
\ref{t2}, \ref{cdi}

compl\'ement unipotent  \dotfill 
\ref{d14.1}

conservatif  \dotfill \ref{P1}, \ref{L1}, \ref{meg1}, 
\ref{P3/2},
\ref{s.saa}, \ref{new p3/2}, \ref{l14}, 
\ref{utile}

coproduit (de ca\-t\'e\-go\-ries)  \dotfill 
\ref{A.D7}

cor\'etraction  \dotfill \ref{P1}, \ref{L2}, 
\ref{Ldim}\bigskip

dimension (rigide) \dotfill \ref{traces}, \ref{a 
la Deligne},
\ref{nouvelle}, \ref{nil, tr}, \ref{dim 0 ou 1}, 
\ref{rmot}, \ref{new
p5/2}, \ref{nouveau}, \ref{Z/2}, \ref{pkim.2}, 
\ref{tkim.1}

dimension finie au sens de Kimura \dotfill \ref{dkim.1}

dimension de Kimura (ou kimension) \dotfill \ref{dkim.2} 

dual 
\dotfill \ref{l7}, \ref{6.2}, \ref{hh0}, \ref{duaux 
droite},
\ref{cdi}, \ref{prof} \bigskip

enveloppe pro-semi-simple 
\dotfill \ref{12.1}

enveloppe r\'eductive  \dotfill 
\ref{d15.1}

enveloppe pro-r\'eductive  \dotfill 
\ref{d14.1}

\'epaisse \dotfill\ref{ep-loc}\bigskip

faiblement 
(exact, \'epi)  \dotfill \ref{d14.2}, \ref{suite}, \ref{c7}
\bigskip

gerbe \dotfill \ref{e10.1}, \ref{ecl},

groupe \`a conjugaison pr\`es  \dotfill \ref{i3}, 
\ref{d14.1}\bigskip

h\'er\'editaire  \dotfill \ref{12.3}

Hochschild 
(homologie) \dotfill \ref{h0}, \ref{rhoch}

Hochschild (cohomologie) 
\dotfill \ref{7.2}
\bigskip

id\'eal (bilat\`ere) d'une ca\-t\'e\-go\-rie 
additive  \dotfill \ref{1.3}

indice \dotfill 
\ref{nouvelle}
\bigskip

$K$-ca\-t\'e\-go\-rie, $K$-foncteur, 
$K$-lin\'eaire \dotfill \ref{biprod}

Kimura (ca\-t\'e\-go\-rie de) 
\dotfill \ref{defki} \bigskip   

mo\-no\-\"{\i}\-dal (id\'eal)  \dotfill 
\ref{d6.1}

mo\-no\-\"{\i}\-dal (foncteur)   \dotfill \ref{9.1}

Morita 
(\'equivalence de)  \dotfill \ref{morita}

\bigskip

pr\'esentation 
finie  \dotfill \ref{r1.1}

produit local (de ca\-t\'e\-go\-ries) \dotfill 
\ref{A.D7}

profinie (alg\`ebre)  \dotfill 
\ref{12.1}

pseudo-ab\'elien   \dotfill \ref{biprod}

puissance 
ext\'erieure, sy\-m\'e\-tri\-que  \dotfill \ref{alt, sym}\bigskip

radical 
\dotfill \ref{D1}

$\otimes$-radical  \dotfill \ref{d6.2}

radical 
infini \dotfill \ref{d3.1}

radiciel (mon mari) \dotfill \ref{D2}, 
\ref{NewL1}, \ref{P3/2},
\ref{radinf}, \ref{C3.1.}, \ref{P3}, 
\ref{P3'}, \ref{ab} 

r\'etraction \dotfill \ref{P1}, 
\ref{L2}

rigide  \dotfill \ref{6.2}, \ref{new p1}, \ref{rtss}, 
\ref{rmot},
\ref{new p2}, \ref{c2}, \ref{absJannsen}, \ref{new p4}, 
\ref{3contrex},
\ref{c4}\bigskip

semi-primaire  \dotfill 
\ref{D4}

semi-simple  \dotfill \ref{D3}, \ref{A.D6}

s\'eparable 
(ca\-t\'e\-go\-rie)  \dotfill \ref{D4sep1}

simple (objet, ca\-t\'e\-go\-rie) 
\dotfill \ref{A.D6}

sorite \dotfill \ref{sorital}

stricte (ca\-t\'e\-go\-rie 
mo\-no\-\"{\i}\-dale, foncteur mo\-no\-\"{\i}\-dal) 
\dotfill 
\ref{ml}

strictement (semi-primaire, de Wedderburn, $K$-lin\'eaire) 

\dotfill \ref{D4swed}

structure balanc\'ee  \dotfill \ref{d4}, 
\ref{ex4}\bigskip

trace\dotfill \ref{traces},  \ref{tr et pext}, 

\ref{p6.4.1}, \ref{a la Kimura}, \ref{nil, tr}, \ref{a la 
bru},
\ref{rtss},
\ref{3contrex},\ref{h0}

tressage \dotfill 
\ref{6.2}, \ref{traces}, \ref{new p1}, \ref{Tre},
\ref{11.2}, 
\ref{p10}, \ref{t4}, \ref{ex4}

tressage infinit\'esimal \dotfill 
\ref{ex4}, \ref{rdc}

type de re\-pr\'e\-sen\-ta\-tion  \dotfill \ref{d3.1}, 
\ref{ccm}, \ref{t3.1},
\ref{12.2}, \ref{12.3}, 
\ref{eprgpu}\bigskip

Wedderburn (ca\-t\'e\-go\-rie de) \dotfill 
\ref{D4wed}  } 
\bigskip\bigskip\bigskip\bigskip

\newpage

\begin{thebibliography}{I}

\bibitem{ak(note)} Y. Andr\'e, B. Kahn, Construction inconditionnelle 
de groupes de Galois motiviques, C.R. Acad. Sci. Paris, S\'er I.\, {\bf 331} 
(2002), 989--994.

\bibitem{beilinson} A. Beilinson, Height pairing between algebraic
cycles, in: $K$-theory, Arithmetic and
Geometry, Lect.notes in Math. {\bf 1289}, Springer (1987), 27--41.

\bibitem{be} D. Benson, {\it Representations and cohomology I}, Cambridge
studies {\bf 30}, Cambridge Univ. Press,1995.

\bibitem{bbh} F. Beukers, D. Brownawell, G. Heckman, {  Siegel
normality}, Annals of Math. {\bf 127} (1988), 279--308.

\bibitem{alg} N. Bourbaki, {\it Alg\`ebre}, chapitre VIII, Hermann, 1958.

\bibitem{tables} N. Bourbaki, {\it Groupes et alg\`ebres de Lie}, chapitre 
VI, Hermann/CCLS, 1975.

\bibitem{lie} N. Bourbaki, {\it Groupes et alg\`ebres de Lie}, chapitre 
VIII, Hermann/CCLS, 1975.

\bibitem{breen} L. Breen, { Tannakian categories}, {\it in} Motives,
Proc. Symposia pure Math. {\bf 55} (I), AMS, 1994, 337--376.

\bibitem{bru1} A. Brugui\`eres, Th\'eorie tannakienne non
commutative, Comm. in  Algebra, {\bf 22}(14) (1994), 5817--5860.

\bibitem{bru2} A. Brugui\`eres, Tresses et structure enti\`ere sur la
ca\-t\'e\-go\-rie des re\-pr\'e\-sen\-ta\-tions de $SL_N$ quantique, Comm. in Algebra
{\bf 28} (2000), 1989--2028.

\bibitem{ce} H. Cartan, S. Eilenberg, {\it Homological algebra}, Princeton
Univ. Press, 1956.

\bibitem{cartier} P. Cartier, { Construction combinatoire des
invariants de Vassiliev-Kontsevich des n\oe uds}, C. R. Acad. Sci. Paris
{\bf 316} (1993), 1205--1210.

\bibitem{cq} J. Cuntz, D. Quillen, { Algebra extensions and
nonsingularity}, Journal A.M.S. {\bf 8} 2 (1995), 251--289.

\bibitem{de} P. Deligne, {Ca\-t\'e\-go\-ries tannakiennes}, in:
The Grothendieck Festschrift, vol. 2, Birkh\"auser P.M. {\bf 87} (1990),
111--198.

\bibitem{qfs} P. Deligne et {\it al.}, {\it Quantum fields and strings: a
Course for Mathematicians}, AMS, 1999.

\bibitem{gabriel} P. Gabriel, { Des ca\-t\'e\-go\-ries ab\'eliennes}, Bull.
Soc. Math. France {\bf 90} (1962), 323--448.

\bibitem{gabrielrev} P. Gabriel, { Probl\`emes actuels de th\'eorie
des re\-pr\'e\-sen\-ta\-tions}, L'Ens. Math.  {\bf 20} (1974), 323--332.

\bibitem{giraud} J. Giraud, {\it Cohomologie non ab\'elienne}, Springer, 
1971.

\bibitem{tohoku} A. Grothendieck, { Sur quelques points d'alg\`ebre 
homologique}, Tohoku Math. J. (2), {\bf 9} (1957), 119--221.

\bibitem{gin} V. Ginzburg, { Principal nilpotent pairs in a semisimple Lie 
algebra. I}, Invent. Math. {\bf 140} (2000), 511--561.

\bibitem{gu-pe} V. Guletskii, C. Pedrini, { The Chow motive of the
Godeaux surface}, pr\'epublication (2001).

\bibitem{happel} D. Happel, {\it Triangulated categories in the
representation theory of finite dimensional algebras}, London Math. Soc.
Lect. notes {\bf 119} (1988), Cambridge Univ. Press.

\bibitem{higman}  G. Higman, { On a conjecture of Nagata},
Proc. Camb. Philos. Soc. {\bf 52} (1956), 1--4.

\bibitem{jac} N. Jacobson, Ann. of Math. {\bf 36} (1935), 875--881.

\bibitem{jannsen} U. Jannsen, { Motives, numerical equivalence and
semi-simplicity}, Invent. Math. {\bf 107} (1992), 447--452.

\bibitem{jannsen2} U. Jannsen, { Equivalence relations on algebraic
cycles}, {\it in} The arithmetic and geometry of algebraic cycles (Banff,
AB, 1998), 225--260, NATO Sci. Ser. C Math. Phys. Sci., {\bf 548}, Kluwer
Acad. Publ., Dordrecht, 2000.



\bibitem{js} A. Joyal, R. Street, { Braided monoidal categories}, Adv.
in Math. {\bf 102} (1993), 20--78.

\bibitem{krt} C. Kassel, M. Rosso, V. Turaev, {\it Quantum groups and knot
invariants}, S.M.F. Panoramas et synth\`eses {\bf 5} (1997).

\bibitem{km} N. Katz, W. Messing, { Some consequences of the Riemann 
hypothesis for varieties over finite fields}, Invent. Math. {\bf 23} (1974), 
73--77.

\bibitem{kelly} G. M. Kelly, { On the radical of a category}, J.
Australian Math. Soc. {\bf 4} (1964), 299--307.

\bibitem{ks} O. Kerner, A. Skowro\'nski, { On module categories with
nilpotent infinite radical}, Compos. Math. {\bf 77} 3 (1991), 313--333.

\bibitem{ki} S.I. Kimura,  { Chow motives can be finite-dimensional, in
some 
sense}, \`a para\^{\i}tre au J. of Alg. Geom.

\bibitem{kos} B. Kostant, { The principal three-dimensional subgroup
and the Betti numbers of a complex simple Lie group},  Amer. J. Math.
{\bf 81} (1959), 973--1032.

\bibitem{ku} K. K\"unneman, { On the Chow motive of an abelian 
scheme}, 
in: Motives,
Proc. Symposia pure Math. {\bf 55} (I), AMS, 1994, 189--205.



\bibitem{leduc} P.Y. Leduc, { ca\-t\'e\-go\-ries semi-simples et
ca\-t\'e\-go\-ries primitives}, Canad. J. Math. {\bf 20} (1968), 612--628.

\bibitem{maclane} S. Mac Lane, {\it Categories for the working
mathematician}, 2\`eme \'ed., Springer GTM {\bf 5} (1998).

\bibitem{margulis} G. A. Margulis, {\it Discrete subgroups of semisimple Lie
groups}, Springer, Berlin, 1991.

\bibitem{mitchell} B. Mitchell, { Rings with several objects}, Adv.
Math. {\bf 8} (1972), 1--161.

\bibitem{mor} V.V. Morozov, { Sur un \'el\'ement nilpotent dans une
alg\`ebre de Lie semi-simple} (en russe), Dokl. Akad. Nauk SSSR {\bf 36}
(1942), 83--86.

\bibitem{mor2} V.V. Morozov, { Sur le centralisateur d'une
sous-alg\`ebre semi-simple d'une alg\`ebre de Lie semi-simple} (en
russe), Dokl. Akad. Nauk SSSR {\bf 36} (1942), 259--261.

\bibitem{murre} J.P. Murre, { On a conjectural filtration on the Chow
groups 
of an algebraic variety}, parts I and II, Indag. Math. {\bf 4}
(1993), 
177--201.

\bibitem{nagata} M. Nagata, { On the nilpotency of nil-algebras}, 
J. Math. Soc. Japan {\bf 4} (1952), 296--301.

\bibitem{nathanson} M. Nathanson, { Classification problems in
$K$-categories}, Fund. Math. {\bf 105} 3 (1979/80), 187--197.

\bibitem{os} P. O'Sullivan, lettres aux auteurs, 29 avril et 12 mai 2002.

\bibitem{pan} D. I. Panyushev, { Nilpotent pairs
in semisimple Lie algebras and their characteristics},  Internat.
Math. Res. Notices {\bf 2000} 1--21.

\bibitem{piard} A. Piard, { Indecomposable representations of a 
 semi-direct product $sl(2)\ltimes\sA$ and semi-simple 
groups containing $sl(2)\ltimes\sA$}, in: Sympos. Math. XXXI 
(Roma, 1988), 185--195, Acad.  Press, 1990.  

\bibitem{platon} Platon, {\it Ph\'edon}, \S $\;$ LIII.

\bibitem{popescu} N. Popescu, {\it Abelian categories with applications
to rings and modules}, Acad. Press, 1973.

\bibitem{prest} M. Prest, {\it Model theory and modules}, L.M.S. Lecture note
series {\bf 130}, Cambridge Univ. Press 1988.

\bibitem{ringel} C. M. Ringel, {  Recent advances in the
representation theory of finite dimensional algebras}, in Representation
theory of finite groups and finite-dimensional algebras, Birkh\"auser
Progress in Math. {\bf 95} (1991).

\bibitem{rowen} L. Rowen, {\it Ring theory}, vol. 1, Acad. Press, 1988.

\bibitem{rump} W. Rump, { Doubling a path algebra, or: how to extend
indecomposable modules to simple modules}, in: Representation theory of
groups, algebras and orders (Costan\c{t}a, 1995), An. \c{S}tiin\c{t}.
Univ. Ovidius Constan\c{t}a Ser Mat. {\bf 4} 2 (1996), 174--185. 

\bibitem{saavedra} N. Saavedra Rivano, {\it ca\-t\'e\-go\-ries tannakiennes},
Lect. Notes in Math. {\bf 265}, Springer, 1972.

\bibitem{serre} J.-P. Serre, { G\`ebres}, L'Ens. Math. {\bf 39}
(1993), 33--85.

\bibitem{serrem} J.-P. Serre, { Propri\'et\'es conjecturales des
groupes de Galois motiviques et des re\-pr\'e\-sen\-ta\-tions
$l$-adiques}, {\it in} Motives, Proc. Symposia pure Math. {\bf 55} (I),
AMS, 1994, 377--400.

\bibitem{sisk} D. Simson, A. Skowro\'nski, { The Jacobson radical
power series of module categories and the representation type}, Bol. Soc.
Mat. Mexicana {\bf 5} 2 (1999), 223--236.

\bibitem{street} R. Street, { Ideals, radicals and structure of
additive categories},  Appl. Cat. Structures {\bf 3} (1995), 139--149.

\bibitem{thomason} R. Thomason, { The classification of triangulated
categories}, Compos. Math. {\bf 105} (1997), 
1--27.

\bibitem{voevodsky} V. Voevodsky, { A nilpotence theorem for cycles
algebraically equivalent to zero}, International Mathematics Research
Notices {\bf 4} (1995), 1--12.

\bibitem{vogel} P. Vogel, { Invariants de Witten-Reshetikin-Turaev
  et th\'eories quantiques des champs}, in: Panoramas et synth\`eses {\bf 7}
(1999), 117-143.
\end{thebibliography}
 \end{document}